# Symmetry, topology and the maximum number of mutually pairwise touching infinite cylinders: complete configuration classification

*Peter V. Pikhitsa[a)] and Stanislaw Pikhitsa*

[a)] *Seoul National University, Seoul, 151-742, Korea*

peter@snu.ac.kr

**Abstract**

We provide a complete classification of possible configurations of mutually pairwise touching infinite cylinders in Euclidian 3D space. It turns out that there is a maximum number of such cylinders possible in 3D independently on the shape of the cylinder cross-sections. We give the explanation of the uniqueness of the non-trivial configuration of seven equal mutually touching round infinite cylinders found earlier. Some results obtained for the chirality matrix which is equivalent to the Seidel adjacency matrix may be found useful for the theory of graphs.

## 1. Introduction.

The problem of mutually touching infinite cylinders begets a non-trivial geometry leading to expansive or auxetic behavior of a regular network of cylinders that may be interesting for Physics [1-3]. A configuration of mutually touching $n$ infinite cylinders (named an $n$-knot in [1-3] and, alternatively, $n$-cross in [3]) lies in the core of the auxetic properties because each mutual touching removes a degree of freedom from the configuration until for the configurations in the network only expansion or contraction as a whole is allowed of all possible motions [1].

Since the discovery of the mutually touching 7,8, and 9 round infinite cylinders [1-4] many questions have remained open. First of all there is the question of configuration classification, that is the development of a method that would give a possibility to distinguish configurations to find new ones. Another question is of how many cylinders of arbitrary adjustable cross-sections could mutually contact. Before only round cylinders have been under consideration [1-4]. And finally, is the configuration with seven all equal round cylinders (named here **a89** for the reasons explained below) first found in [2] and then rediscovered in [4] unique or there are other non-trivial configurations of the seven cylinders? We will answer these questions by considering the symmetry and topology of configurations of oriented lines in 3D as well as by utilizing the restrictions imposed by the conditions of mutual contacts.

We use the normalized *chirality matrix* that we introduced in [3], which describes the topology of a configuration of infinite oriented lines in 3D space and reads

$$P_{i,k} \equiv \begin{pmatrix} 0 & \frac{R_{ik}[n_i \times n_k]}{|R_{ik}[n_i \times n_k]|} & \cdots & \cdots \\ \frac{R_{ik}[n_i \times n_k]}{|R_{ik}[n_i \times n_k]|} & 0 & \cdots & \vdots \\ \vdots & \vdots & \ddots & \vdots \\ \vdots & \cdots & \cdots & 0 \end{pmatrix}, \quad (1)$$

where $\boldsymbol{n}_i$ is a unit vector along the $i^{\text{th}}$ line (Fig. 1) and $\boldsymbol{R}_{ik}$ is the vector of the shortest distance between the $i^{\text{th}}$ and the $k^{\text{th}}$ lines, collinear with $[\boldsymbol{n}_i \times \boldsymbol{n}_k]$. Note that chirality matrix $P$ as a symmetrical matrix of the entries +1 and -1 and the zeroes on its diagonal is identical to the Seidel adjacency matrix used in the theory of graphs. Therefore some results that we obtain below for the chirality matrix may be useful for the graph theory. The chirality matrix has its inner symmetries: changing directions of lines leads to simultaneous changing the sign of $i^{th}$ column and $i^{th}$ row and the line number permutations lead to $i^{th}$ and $j^{th}$ permutations of corresponding columns and rows. These operations are the similarity transformations that do not change the determinant of the chirality matrix which remains the topological invariant [3].

As far as the chirality matrix $P$ is not enough to tell all the nuances of the topology of configurations of mutually pairwise touching cylinders, we introduce a novel matrix $\mathcal{R}$ that we name the *ring matrix*. Together with the chirality matrix $P$ it can unambiguously distinguish the topology of non-trivial configurations of lines in 3D in general and thus can be applied for the classification of the configurations of 7, 8, and 9 mutually pairwise touching infinite cylinders (called 7, 8, and 9-knots in [3]). Along with $\mathcal{R}$ we will generate a matrix $Q$ from the chirality matrix $P$ that will ignore the direction of oriented lines because the orientations are irrelevant for the configurations of the cylinders. Such situation is known for the nematic liquid crystals where the director rather than the vector determines their configurations.

## 2. Ring matrix and complete classification of the configurations.

One can write

$$\boldsymbol{R}_{ik} = |\boldsymbol{R}_{ik}|\boldsymbol{r}_{ik},$$

$$\boldsymbol{r}_{ik} = P_{ik}[\boldsymbol{n}_i \times \boldsymbol{n}_k]/|P_{ik}[\boldsymbol{n}_i \times \boldsymbol{n}_k]| \quad . \tag{2}$$

The unit vectors $\boldsymbol{r}_{ik}$ are shown in Fig. 1 where the projection along the direction of the $i^{th}$ line of the configuration of four lines is given. To get the criterion of whether the $i^{th}$ line is encaged by the other three lines we project

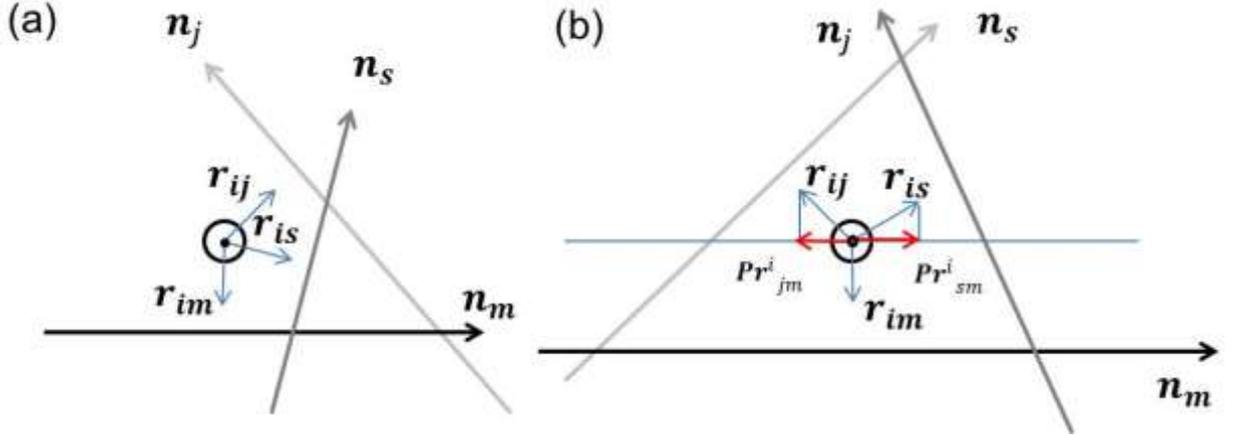

Fig. 1. Top view of the $i^{th}$ line encaged by three lines that we call a ring. The line direction towards the reader is denoted by a circle with a dot. (a), no encaging; (b) the line encaged by the other three lines. The red arrows show the projections of the unit vectors in the direction of $r_{im}$.

the unit vectors to define

$$Pr^i_{jm} = r_{im} - (r_{im} r_{ij}) r_{ij} \quad . \qquad (3)$$

When the projections point in different directions for each of the three lines, then the $i^{th}$ line is surrounded by the other lines as shown in Fig. 1b and thus is encaged (entangled). The criterion would be

$$S(j,s,m;i) + S(m,s,j;i) + S(j,m,s;i) \equiv -3, \qquad (4)$$

where

$$S(j,s,m;i) = sign(Pr^i_{jm} Pr^i_{sm}). \qquad (4)$$

Calculated for each line it can form a vector

$$\mathfrak{R}_i = \sum_{s>m>j} 1, if\ [S(j,s,m;i) + S(m,s,j;i) + S(j,m,s;i) = -3], 0, otherwise\ . \quad (5)$$

Such a vector tells how many rings of three lines encage each of the line in a given configuration. However, we found that for the unambiguous discrimination of the topology of line configurations the vector, being a topological invariant, is still insufficient. As we mentioned earlier, the complete description stems from the *ring matrix* $\mathcal{R}$. To define $\mathcal{R}$ we give the vector $\mathfrak{R}$ the structure where in Eq. (5) instead of 1s (the scalars) we sum up vectors

$$\mathcal{R}_{i,k} = \sum_{s>m>j}\{1,\ if[(s=k)+(j=k)+(m=k)], 0\ otherwise\}, if\ [S(j,s,m;i) + S(m,s,j;i) + S(j,m,s;i) = -3], 0\ otherwise\ . \qquad (6)$$

that indicate how many times each line enters into each of rings encaging a given line. For illustration we give the vector $\mathfrak{R}$ and the matrix $\mathcal{R}$ calculated for the configuration **a89** for $7^*$-knot [2, 3] with the help of Eq. (5) and Eq. (6), respectively:

$$\mathfrak{R}(a89) = \begin{pmatrix} 4 \\ 8 \\ 0 \\ 0 \\ 4 \\ 4 \\ 0 \end{pmatrix}, \quad (7)$$

$$\mathcal{R}(a89) = \begin{pmatrix} 0 & 1 & 1 & 4 & 1 & 1 & 4 \\ 4 & 0 & 4 & 4 & 4 & 4 & 4 \\ 0 & 0 & 0 & 0 & 0 & 0 & 0 \\ 0 & 0 & 0 & 0 & 0 & 0 & 0 \\ 1 & 1 & 4 & 4 & 0 & 1 & 1 \\ 1 & 1 & 1 & 1 & 4 & 0 & 4 \\ 0 & 0 & 0 & 0 & 0 & 0 & 0 \end{pmatrix}. \quad (8)$$

One can notice that the sum of the numbers in the first row in Eq. (8) is three times the first number in the vector in Eq. (7). The same is true for all the rows in Eq. (8). It is so because each ring encaging a line consists of three other lines. The matrix $\mathcal{R}$ is a topological invariant which together with the chirality matrix [3]

$$P(a89) = \begin{pmatrix} 0 & +1 & +1 & +1 & +1 & +1 & +1 \\ +1 & 0 & +1 & +1 & +1 & -1 & +1 \\ +1 & +1 & 0 & -1 & -1 & -1 & +1 \\ +1 & +1 & -1 & 0 & -1 & +1 & +1 \\ +1 & +1 & -1 & -1 & 0 & +1 & -1 \\ +1 & -1 & -1 & +1 & +1 & 0 & +1 \\ +1 & +1 & +1 & +1 & -1 & +1 & 0 \end{pmatrix} \quad (9)$$

completely and unambiguously describes the topology of the line configuration. It is remarkable that one can reproduce the chirality matrix and the ring matrix just by inspecting any given configuration while establishing directions along the cylinders and marking the chirality at each contact, along with the counting the number of rings around each cylinder that contain this cylinder. The reverse problem can be solved as well: being given both matrices one can reproduce the topologically equivalent line configuration. However, for practical purposes of finding new configurations instead of inspecting the two matrices it would be easier to have a numerical invariant that could distinguish different topologies.

## 3. A numerical invariant for the complete classification of configurations.

Unfortunately, the chirality matrix is sensitive to the orientations of the lines. To avoid this dependence we derive a new matrix $Q$ from $P$ which is invariant with respect to the line orientation as follows (we are using Eq. (9) as an example):

> One has to transform the first row in $P(a89)$ into all +1s (it is already so) and to sum up numbers in the corresponding columns and put the sums into the first row.
> One has to transform the second row in $P(a89)$ into all +1s (here by reverting sign in the sixth column and then in the sixth row) then to sum up numbers in the corresponding columns and put the sums into the second row.

After proceeding in the same way through the whole matrix of Eq. (9) one gets the symmetric matrix:

$$Q(P(a89)) = \begin{pmatrix} 6 & 4 & 0 & 2 & 0 & 2 & 4 \\ 4 & 6 & 2 & 0 & -2 & -2 & 2 \\ 0 & 2 & 6 & -2 & 0 & 2 & 2 \\ 2 & 0 & -2 & 6 & -2 & 2 & 4 \\ 0 & -2 & 0 & -2 & 6 & 0 & 0 \\ 2 & -2 & 2 & 2 & 0 & 6 & 0 \\ 4 & 2 & 2 & 4 & 0 & 0 & 6 \end{pmatrix} . \quad (10)$$

An additional advantage of $Q$ is that it is different for the mirror configuration

$$Q(-P(a89)) = \begin{pmatrix} 6 & -2 & 2 & 0 & 2 & 0 & -2 \\ -2 & 6 & 0 & 2 & 4 & 4 & 0 \\ 2 & 0 & 6 & 4 & 2 & 0 & 0 \\ 0 & 2 & 4 & 6 & 4 & 0 & -2 \\ 2 & 4 & 2 & 4 & 6 & 2 & 2 \\ 0 & 4 & 0 & 0 & 2 & 6 & 2 \\ -2 & 0 & 0 & -2 & 2 & 2 & 6 \end{pmatrix} . \quad (11)$$

Eq. (10) and Eq. (11) satisfy the equation that holds for all 7-knots:

$$Q(P) + Q(-P) = \begin{pmatrix} 12 & 2 & 2 & 2 & 2 & 2 & 2 \\ 2 & 12 & 2 & 2 & 2 & 2 & 2 \\ 2 & 2 & 12 & 2 & 2 & 2 & 2 \\ 2 & 2 & 2 & 12 & 2 & 2 & 2 \\ 2 & 2 & 2 & 2 & 12 & 2 & 2 \\ 2 & 2 & 2 & 2 & 2 & 12 & 2 \\ 2 & 2 & 2 & 2 & 2 & 2 & 12 \end{pmatrix} . \quad (12)$$

The form of Eq. (12) also holds for all 8- and 9- knots where the matrices are 8x8 and 9x9 and the numbers on the diagonal of the corresponding matrix are 14 and 16, respectively.

All possible numerical invariants could be formed from the traces of corresponding matrices such as: $tr(\mathcal{R}^2), tr(\mathcal{R}^3), \ldots tr(\mathcal{R}^7)$. Yet, in order to distinguish mirror configurations we have to include $Q$ as well, for example, by adding $tr(Q\mathcal{R}), tr(Q\mathcal{R}^2), tr(Q^2\mathcal{R})$, etc. Such numbers could be combined in a single numerical for each configuration. Here we suggest a simple formal summation over $n = 0,1,2 \ldots$ that proved to be a good practical numerical invariant that could tell the topologically different configurations:

$$\wp(P, R) = \sum_n tr(Q\mathcal{R}^n) = tr[Q(1 - \mathcal{R})^{-1}] . \quad (13)$$

One can introduce an additional ring invariant that would have the same value for the configurations with the same ring matrix which is insensitive to mirroring:

$$\wp_\mathcal{R} = \wp(P, R) + \wp(-P, R). \quad (14)$$

It directly follows from Eq. (13) along with the Eq. (12) that Eq. (14) is determined only by $\mathcal{R}$. Below we will give the results of the calculations of invariants for a number of configurations.

### 4. The chirality matrices forbidden for mutual pairwise touching.

We will show that much information about $n$-knot configurations can be extracted solely from the analysis of the chirality matrix. The most important is that there is a rigorous sufficiency criterion as to chirality matrices that lead to impossible configurations, that is to those configurations that never allow mutual pairwise contacts of all cylinders. Namely,

- *if P contains a subset of 5 cylinders that have a $5 \times 5$ submatrix (which we conventionally call $K5$ from the theory of graphs)*

$$K5 = \begin{pmatrix} 0 & 1 & 1 & 1 & 1 \\ 1 & 0 & 1 & 1 & 1 \\ 1 & 1 & 0 & 1 & 1 \\ 1 & 1 & 1 & 0 & 1 \\ 1 & 1 & 1 & 1 & 0 \end{pmatrix} \quad (15)$$

or any of its transformation with changing directions and cylinder permutations (simultaneous changing the sign of $i^{th}$ column and $i^{th}$ row or/and $i^{th}$ and $j^{th}$ permutations), as well as mirror transformation that changes the signs of all entries of $K5$, then this configuration is impossible for mutual pairwise contacts.

To prove the impossibility of $K5$ to have all the cylinders in the mutual contacts we will use a projection diagram to form the matrix $K5$. First, one can draw a diagram element to illustrate the assignment of signs +1 and -1 for the projection of two cylinders (Fig. 2): nearest rotation from red to blue going counter-clockwise gives +1 for the corresponding element of the chirality matrix, otherwise -1. Here in the diagram the red line should run over the blue line.

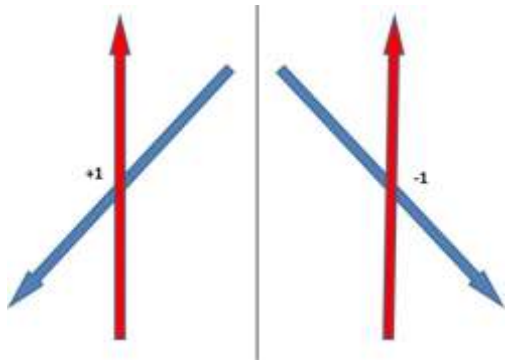

Fig. 2. The elementary diagrams illustrating the chirality +1 (left) and -1 (right) of the configuration of two oriented lines in 3D.

It is not difficult to construct a diagram for $K5$ and prove the theorem by direct inspection of the configuration which we give below in Fig. 3. The determinant of $K5$ is either 4 or -4 (for the mirror configuration) and only $K5$ can have such determinants. This property can be used to

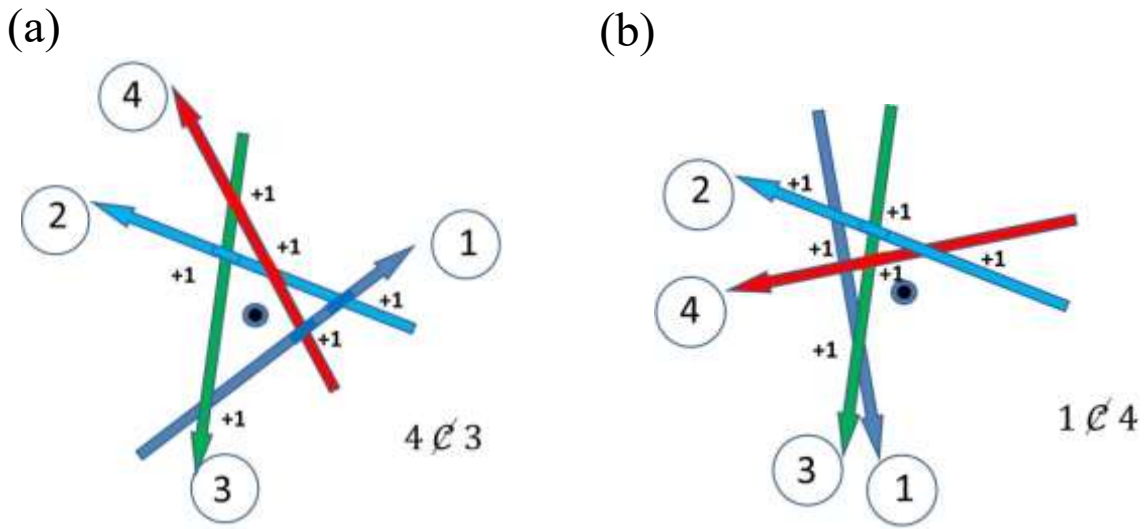

Fig. 3. (a) Illustration of $K5$ where cylinder 4 does not contact cylinder 3 marked with the crossed letter C. $0^{th}$ cylinder points out of the picture. (b) The other case where $1^{st}$ cylinder does not contact with the $4^{th}$ cylinder.

find the presence of $K5$ in any chirality matrix. The requirement that $K5$ should be absent selects out of all possible chirality matrices $P$ for 7x7 matrices with all possible determinants which absolute values are

2, 6, 10, 14, 18, 22, 26, 30, 34, 42, 46, 50, 54, 66, 70, 78, 90, 102, 150, 162, 250 (16)

only the matrices of determinants with absolute values

2, 10, 18, 42, 54, 66, 102, 150, 162, **250**                (17)

which do not contain $K5$. We found existing configurations for all determinants of Eq. (17) (see Fig. 7 and Appendices 1,3,4) but for 250.

Of course there may exist configurations which do not contain $K5$ but still cannot have all the cylinders in the mutual contacts. Namely, of all the matrices of Eq. (17) only for $|P| = 250$ all cylinders cannot be in mutual contacts. Below we show that $P250$ in spite of not having $K5$ still cannot be realized with all mutual contacts for positive radii of the cylinders: at least two of the cylinders cannot contact.

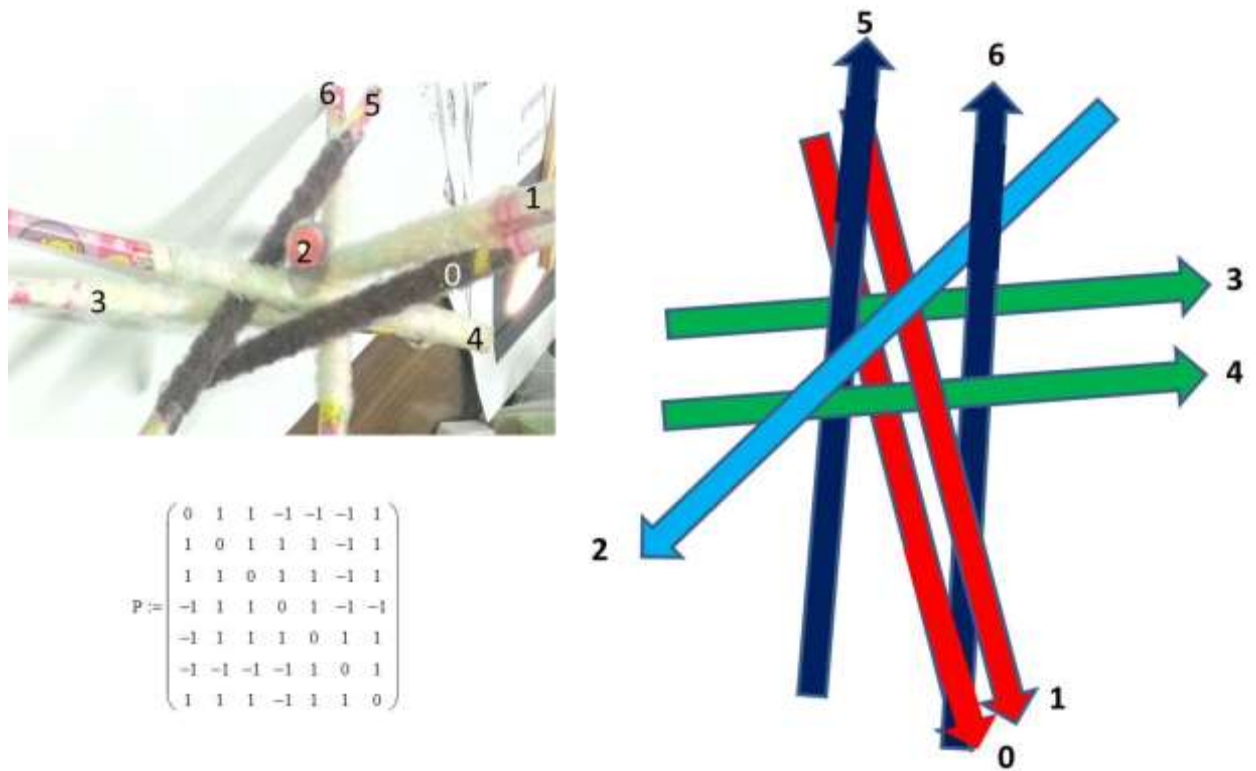

Fig. 4. A pencil model of the configuration of cylinders $P250$ (left) and its schematic presentation with directed lines. The numbers of lines correspond to the numbers on the pencils.

Let us introduce a useful definition. Two cylinders are said to be in *equal environment* (EE) if two rows/columns that correspond to these two cylinders in the chirality matrix are identical or can be made identical by multiplying by -1. Of course one should ignore the diagonal zeroes. The necessary condition for a matrix to have such a property is to have 1 or -1 among its eigenvalues. Yet, some chirality matrices (for example, the one with $|P| = 42$) while having 1 and -1 among its eigenvalues do not have EE. The EE property is important for understanding why $P250$ configuration is impossible for mutual contacts. Later we will show that EE also helps to select out those configurations that cannot have all equal cylinders.

From the matrix $P250$ shown in Fig. 4 one can see that the first and the second cylinder are in EE. All other chiral matrices out of the list of Eq. (17), for example, for the configuration **b14** with the chiral matrix determinant -102 containing cylinders in EE, possess a property:

*If one switches the chirality between the very two EE cylinders in the matrix, then the matrix turns into the matrix which contains $K5$.*

The only exception is the chiral matrix $P250$: if one switches the chirality matrix entry between the first and the second cylinders in EE (which is $P_{1,2} = 1$ in Fig. 4) to $P_{1,2} = -1$, then the determinant will switch to -250 and will not contain $K5$, unlike other chiral matrices of Eq. (17) that have EE. This specific property greatly restricts the possible configurations of the cylinders corresponding to $P250$. Indeed, the two cylinders in EE cannot go through any ring similar to the one in Fig. 3a which is created by the lines 1,2,3. We call such a ring *right wreath*. In this picture the lines 0 and 4 are in EE and both punch the ring. The configuration is $K5$, yet if we switch the chirality between lines 0 and 4 then there is no $K5$. That means that if the two EE cylinders in a configuration with $|P| = 250$ punched any of such rings, then after switching their chirality, $K5$ would appear as it happens for other chiral matrices. On the other hand, if the two EE cylinders are located outside the ring made by lines 1, 2, 3 in Fig. 3b (we call such a ring *right scissors*) then again they would be in $K5$ configuration. Switching the chirality between the two removes $K5$ which never happens for $P250$. The restrictions of not having such rings to punch and to pass by, reduce the configuration for $P250$ down to the one shown in Fig. 4. The schematic in Fig.4 with lines demonstrate that the 2$^{nd}$ cylinder is in the see-saw configuration balancing on the 1$^{st}$ cylinder and is being able to touch either the 3$^{rd}$ or the 4$^{th}$ cylinder but never both.

Let us note that while investigating the matrix $P250$ we serendipitously came across a remarkable property of the matrix: we found that the determinant $|P_{i,k}(r_i + r_k)| \geq 0$ for all $r_i \geq 0$. We proved it analytically by direct calculation finding the series of all positive terms for the determinant. A similar determinant with all $P_{i,k} = 1$ was considered long ago in 1880 on page 221 in item 29 of the section "Examples on the methods of the text" in the book "A treatise on the theory of determinants", Cambridge, Univ. Press by Robert Forsyth Scott.

## 5. The maximum number of mutually pairwise touching straight cylinders of arbitrary cross-section.

The question of the number of cylinders of arbitrary cross-section arises when one tries to use the arguments of the degrees of freedom that we used for 7,8 and 9 round cylinders [1-3]. Solely on the arguments it would be reasonable to expect that while arbitrary cross-sections have infinitely large number of shapes and therefore the degrees of freedom, then the number of mutually pairwise touching straight cylinders of arbitrary cross-section could be infinite. However, the theorems given below claim that this is not the case: the number of mutually pairwise touching straight cylinders is restricted from above at least by the number 14.

The notion of "pairwise" means here that the point of a contact is shared by two and only two cylinders. The degree of freedom arguments used in our previous works [1-3] give a definite prediction as to the maximum possible number of cylinders of arbitrary radii (which is 7) provided that the radii are not adjusted. If the radii are involved into calculations then the maximum number increases up to 9 because each adjustable radius adds a degree of freedom. However, there is a question: what if the cross-section of cylinders is not a circle but is possible to be varied and adjusted? It is easy to imagine an elliptical cross-section which immediately produces up to two additional degrees of freedom for each cylinder (rotation around the cylinder axis and the aspect ratio). It can be generalized for infinitely many degrees of freedom for arbitrary shaped cross-sections. As we said earlier, one would expect an arbitrary number of mutually touching cylinders.   With the help of the chirality matrix we will show that it is not the case for pairwise touching, namely:

*the maximum number of mutually pairwise touching infinitely long straight cylinders with arbitrary adjustable cross-sections is restricted from above.*

The degeneracy of a straight line when several lines pass through one single point in space is not considered because in such a degenerate case the notion of "mutual contact" does not have a meaning of a pairwise contact which we explore. Still it is possible to consider infinite straight stripes (blades) with straight line edges if only the pairwise contact is allowed which we will publish elsewhere.

$P250$ plays the role of $K5$ for larger than 7 number of cylinders. The presence of $K5$ or $P250$ is a sufficient condition for any configuration of cylinders of arbitrary cross-section not to have all the cylinders in the mutual pairwise contacts. We use the properties of $K5$ and $P250$ to prove important theorems:

1. The number of arbitrary cylinders that can be in mutual pairwise contacts cannot exceed 18. That means that starting from $n > 18$ any chirality matrix contains $K5$.
2. Moreover, the number of arbitrary cylinders in mutual pairwise contacts cannot exceed 14. This happens because starting from $n > 14$ any chirality matrix contains $P250$ as a submatrix.

We prove these theorems by direct calculations showing that:

*All chiral matrices larger than 18x18 contain $K$*5*.  There is only one prototype matrix 18x18 (up to the similarity transformations) that does not contain* 5 *:*

$P = $

|    | 0  | 1  | 2  | 3  | 4  | 5  | 6  | 7  | 8  | 9  | 10 | 11 | 12 | 13 | 14 | 15 | 16 | 17 |
|----|----|----|----|----|----|----|----|----|----|----|----|----|----|----|----|----|----|----|
| 0  | 0  | 1  | 1  | 1  | -1 | -1 | -1 | 1  | 1  | 1  | -1 | -1 | -1 | 1  | 1  | 1  | -1 | -1 |
| 1  | 1  | 0  | 1  | 1  | -1 | -1 | 1  | -1 | 1  | 1  | -1 | 1  | 1  | -1 | -1 | -1 | 1  | -1 |
| 2  | 1  | 1  | 0  | 1  | 1  | -1 | -1 | 1  | -1 | -1 | 1  | -1 | 1  | 1  | -1 | -1 | 1  | -1 |
| 3  | 1  | 1  | 1  | 0  | 1  | 1  | 1  | -1 | 1  | -1 | -1 | -1 | 1  | 1  | -1 | 1  | -1 | 1  |
| 4  | -1 | -1 | 1  | 1  | 0  | 1  | 1  | 1  | -1 | 1  | 1  | -1 | 1  | -1 | -1 | 1  | -1 | -1 |
| 5  | -1 | -1 | -1 | 1  | 1  | 0  | 1  | 1  | 1  | -1 | -1 | 1  | -1 | 1  | -1 | 1  | 1  | -1 |
| 6  | -1 | 1  | -1 | 1  | 1  | 1  | 0  | 1  | 1  | 1  | -1 | -1 | 1  | -1 | 1  | -1 | 1  | 1  |
| 7  | 1  | -1 | 1  | -1 | 1  | 1  | 1  | 0  | 1  | 1  | 1  | -1 | -1 | 1  | 1  | -1 | 1  | -1 |
| 8  | 1  | 1  | -1 | 1  | -1 | 1  | 1  | 1  | 0  | 1  | 1  | 1  | -1 | 1  | -1 | -1 | -1 | 1  |
| 9  | 1  | 1  | -1 | -1 | 1  | -1 | 1  | 1  | 1  | 0  | 1  | 1  | 1  | -1 | 1  | 1  | -1 | -1 |
| 10 | -1 | -1 | 1  | -1 | 1  | -1 | -1 | 1  | 1  | 1  | 0  | 1  | 1  | 1  | -1 | -1 | -1 | 1  |
| 11 | -1 | 1  | -1 | -1 | -1 | 1  | -1 | -1 | 1  | 1  | 1  | 0  | 1  | 1  | -1 | 1  | 1  | -1 |
| 12 | -1 | 1  | 1  | 1  | 1  | -1 | 1  | -1 | -1 | 1  | 1  | 1  | 0  | 1  | 1  | 1  | 1  | 1  |
| 13 | 1  | -1 | 1  | 1  | -1 | 1  | -1 | 1  | 1  | -1 | 1  | 1  | 1  | 0  | 1  | 1  | 1  | 1  |
| 14 | 1  | -1 | -1 | -1 | -1 | -1 | 1  | 1  | -1 | 1  | -1 | -1 | 1  | 1  | 0  | 1  | 1  | 1  |
| 15 | 1  | -1 | -1 | 1  | 1  | 1  | -1 | -1 | -1 | 1  | -1 | 1  | 1  | 1  | 1  | 0  | -1 | -1 |
| 16 | -1 | 1  | 1  | -1 | -1 | 1  | 1  | 1  | -1 | -1 | -1 | 1  | 1  | 1  | 1  | -1 | 0  | -1 |
| 17 | -1 | -1 | -1 | 1  | -1 | -1 | 1  | -1 | 1  | -1 | 1  | -1 | 1  | 1  | 1  | -1 | -1 | 0  |

This matrix has a simple characteristic polynomial $(x^2 - 17)^9$.

*All chirality matrices larger than 14x14 contain P250. There is only one prototype matrix 14x14 (up to the similarity transformations) that does not contain both K5 and 250 :*

$P = $

|    | 0  | 1  | 2  | 3  | 4  | 5  | 6  | 7  | 8  | 9  | 10 | 11 | 12 | 13 |
|----|----|----|----|----|----|----|----|----|----|----|----|----|----|----|
| 0  | 0  | 1  | -1 | 1  | 1  | -1 | -1 | 1  | 1  | -1 | -1 | 1  | 1  | -1 |
| 1  | 1  | 0  | 1  | -1 | 1  | -1 | 1  | -1 | 1  | -1 | -1 | -1 | 1  | 1  |
| 2  | -1 | 1  | 0  | 1  | 1  | -1 | 1  | -1 | -1 | 1  | 1  | -1 | 1  | -1 |
| 3  | 1  | -1 | 1  | 0  | 1  | -1 | 1  | 1  | -1 | 1  | -1 | 1  | -1 | -1 |
| 4  | 1  | 1  | 1  | 1  | 0  | 1  | 1  | -1 | 1  | -1 | 1  | 1  | -1 | -1 |
| 5  | -1 | -1 | -1 | -1 | 1  | 0  | 1  | 1  | 1  | -1 | 1  | -1 | -1 | -1 |
| 6  | -1 | 1  | 1  | 1  | 1  | 1  | 0  | 1  | 1  | 1  | -1 | -1 | -1 | 1  |
| 7  | 1  | -1 | -1 | 1  | -1 | 1  | 1  | 0  | 1  | 1  | -1 | -1 | 1  | -1 |
| 8  | 1  | 1  | -1 | -1 | 1  | 1  | 1  | 1  | 0  | 1  | 1  | 1  | 1  | 1  |
| 9  | -1 | -1 | 1  | 1  | -1 | -1 | 1  | 1  | 1  | 0  | 1  | 1  | 1  | 1  |
| 10 | -1 | -1 | 1  | -1 | 1  | 1  | -1 | -1 | 1  | 1  | 0  | 1  | 1  | -1 |
| 11 | 1  | -1 | -1 | 1  | 1  | -1 | -1 | -1 | 1  | 1  | 1  | 0  | -1 | 1  |
| 12 | 1  | 1  | 1  | -1 | -1 | -1 | -1 | 1  | 1  | 1  | 1  | -1 | 0  | -1 |
| 13 | -1 | 1  | -1 | -1 | -1 | -1 | 1  | -1 | 1  | 1  | -1 | 1  | -1 | 0  |

This matrix has a simple characteristic polynomial $(x^2 - 13)^7$. The direct calculation was performed by using a MathCad11 solver programmed to find *K*5 and *P*250 matrices in all possible chirality matrices of successively growing size. The text of the program for finding *K*5 is given in Appendix 5. A complete calculation up to the matrix 18x18 may take up to 4 hours depending on the CPU.

## 6. Collection of configurations of 7,8, and 9 mutually touching round infinite cylinders.

In order to analyze and classify different configurations we produced a collection of 7-knots (see Fig. 8 below and Appendices 1,3,4). For this purpose we first calculated several 9-knots (in addition to the one first published in [3] with the determinant of the chirality matrix (see Eq. (1) for the definition) $|P| = 0$, which we designated first as **a** and then as **a0a** (when the number of configurations found grew big) with determinants -52, 16, 4, 0 which we first designated **b**, **c**, **d**, and **e** when they were few and then changed to **am52a**, **a16**, **a4a**, and **a0b**, correspondingly, and others (see Fig. 5 and Appendix 2 for the complete list). Here **a=a0a** is different from **e= a0b** in spite of having the same zero determinant of $P$. Then we extracted a number of 7-knot configurations by removing any two cylinders, designating the resulting configuration as **a89**, for example, when it was extracted from **a** by removing $8^{th}$ and $9^{th}$ cylinders or say **e69** if it were $6^{th}$ and $9^{th}$ cylinders from **e**. At the same time we produced 8-knots in the same manner and designated them accordingly (see Fig. 8), for example **a8** when the $8^{th}$ cylinder was removed from **a=a0a**. Except those configurations that originate from different 9-knots we produced a number of configurations of 7-knots and 8-knots independently. They have their own labels usually containing the information of the determinant of the chirality matrix of a certain configuration. The letter **m** before the name of a configuration means a **mirror** configuration of the one without the letter.

Many topologically different 7-knot configurations have the same determinant $|P|$, positive or negative, which had the form $4n + 2$ where $n$ is an integer. The chirality matrix turned out to be different for many configurations with the same determinant but for all cases except $n = 0, -1$, the matrices could be transformed into each other by permutations and/or by changing the sign of the row/column which corresponds to permutations between the cylinders or to a change of the direction of the orientation of a given cylinder which transformation does not change the topology of the configuration. The case $n = 0, -1$ was slightly different because the chirality matrices could have a different characteristic polynomial while having the same determinant. As we said before, the chirality matrix and its determinant could distinguish only configurations that belong to different classes with different determinants. In our paper [3] we distinguished by this means two configurations (now we identify one of them as **e69**) with $|P| = 10$, first of any 7-knot configurations discovered in 2004 [1] and **a89** (it was denoted as $7^*$-knot in [3]) with equal radii with $|P| = -18$, discovered in 2009 [2, 3] and rediscovered in [4]. Note that for the 8-knot and 9-knot the determinants take values $4n + 1$ and $4n$, respectively.

For 9-cross we obtained several cases with determinants with absolute values 0, 4, 16, 52, 80, 84, 104 (encrypted in the names of configurations) which invariants are given in Fig. 5:

$$
\begin{aligned}
a0a &:= 20.837393218 \\
a0b &:= 31.278680391 \\
a0c &:= 10.329027872 \\
a0d &:= 34.588654781 \\
a0e &:= 40.209551479 \\
a0f &:= -19.696390243 \\
a0g &:= -17.864406779 \\
a0i &:= -16.683884892 \\
a0j &:= 24.775095554 \\
a4a &:= 25.352726356 \\
a4b &:= 27.017463933 \\
a16 &:= 27.627505709 \\
am52a &:= 24.019357049 \\
a16b &:= 26.265792952 \\
a4c &:= 25.748401074 \\
a4d &:= 17.370024424 \\
a4e &:= 22.662334217 \\
a4f &:= 21.401499765 \\
am80a &:= 28.273868340 \\
am80b &:= 24.614050270
\end{aligned}
\qquad
\begin{aligned}
ma0a &:= 23.843846756 \\
ma0b &:= 25.423787672 \\
ma0c &:= 46.413324269 \\
ma0d &:= 22.057568881 \\
ma0e &:= 15.967677501 \\
ma0f &:= 65.506146341 \\
ma0g &:= 63.086956521 \\
ma0i &:= 63.345227817 \\
ma0j &:= 20.801649567
\end{aligned}
\qquad
\begin{aligned}
am52b &:= 32.409748667 \\
am84a &:= 30.082683463 \\
a16c &:= 17.573207219 \\
a4g &:= -6.239582169 \\
am80c &:= 23.919756097 \\
a4h &:= 29.864447357 \\
a4i &:= 24.230326797 \\
a4j &:= 34.110579410 \\
a4k &:= 28.729273902 \\
a0h &:= 43.885428506 \\
a16d &:= 29.764023934 \\
a16e &:= -229.801282051 \\
a16f &:= 33.223866917 \\
a16g &:= 21.961822660 \\
a104a &:= 23.231469769
\end{aligned}
$$

Fig. 5. The invariants of all 9-knot configurations found.

The first configuration ever found **a=a0a** was given in [3] (see Appendix 2).

For the 7[*]-knot or **a89** we obtain $\wp = 14.(07317)$ and $\wp = 7.(80487)$ for the mirror configuration **ma89**. In Fig. 8 we give $\wp$ for all 7-knots with $|P| = -18,$ extracted from the known 9-knots (for other determinants see Appendices 1 and 3,4 for the parameters of each configuration). We give in Fig. 6 as an example of one of the tables from Appendix 4 the unique configuration **a89.** The same configuration is shown in Fig. 7 in two different projections to emphasize a beautiful property of any of 7-knots to have projections where cylinders look parallel in pairs.

Configuration label **a89**

$$\begin{pmatrix} \text{orange} \\ \text{red} \\ \text{blue} \\ \text{green} \\ \text{cyan} \\ \text{magenta} \\ \text{gray} \end{pmatrix} \begin{pmatrix} 0 & 0 & 0 & 1 \\ t1 & p1 & z1 & r1 \\ t2 & p2 & z2 & r2 \\ t3 & p3 & z3 & r3 \\ t4 & p4 & z4 & r4 \\ t5 & p5 & z5 & r5 \\ t6 & p6 & z6 & r6 \end{pmatrix} = \begin{pmatrix} 0 & 0 & 0 & 1 \\ 0.3255445253 & 3.1415926536 & 0 & 1 \\ 0.3800224809 & 9.3515809945 & 13.4293807485 & 1 \\ 2.4244469645 & 6.4336249988 & -23.5948605815 & 1 \\ 0.3316040913 & 9.1928411285 & 4.8892717374 & 1 \\ 0.2843728713 & 9.2680129687 & -9.9091055237 & 1 \\ 2.8406082265 & 6.0133479798 & -10.6464648447 & 1 \end{pmatrix}$$

$$P = \begin{pmatrix} 0 & 1 & 1 & 1 & 1 & 1 & 1 \\ 1 & 0 & 1 & 1 & 1 & -1 & 1 \\ 1 & 1 & 0 & -1 & -1 & -1 & 1 \\ 1 & 1 & -1 & 0 & -1 & 1 & 1 \\ 1 & 1 & -1 & -1 & 0 & 1 & -1 \\ 1 & -1 & -1 & 1 & 1 & 0 & 1 \\ 1 & 1 & 1 & 1 & -1 & 1 & 0 \end{pmatrix} \qquad R = \begin{pmatrix} 0 & 1 & 1 & 4 & 1 & 1 & 4 \\ 4 & 0 & 4 & 4 & 4 & 4 & 4 \\ 0 & 0 & 0 & 0 & 0 & 0 & 0 \\ 0 & 0 & 0 & 0 & 0 & 0 & 0 \\ 1 & 1 & 4 & 4 & 0 & 1 & 1 \\ 1 & 1 & 1 & 1 & 4 & 0 & 4 \\ 0 & 0 & 0 & 0 & 0 & 0 & 0 \end{pmatrix}$$

I3(P, R) = 14.0731707317

I3(−P, R) = 7.8048780488

$|P| = -18$

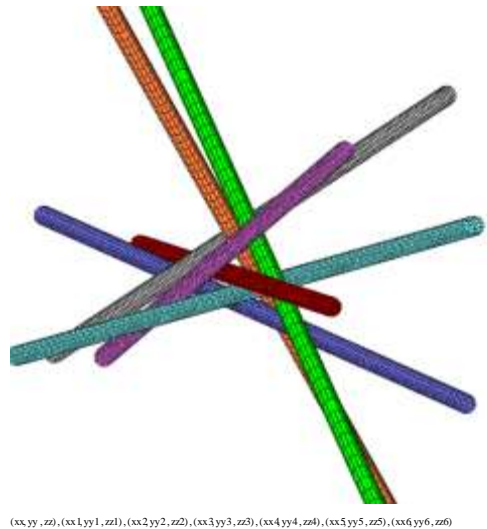

(xx,yy,zz), (xx1,yy1,zz1), (xx2,yy2,zz2), (xx3,yy3,zz3), (xx4,yy4,zz4), (xx5,yy5,zz5), (xx6,yy6,zz6)

Fig. 6. The characteristics of **a89** configuration. The numerical invariant $\wp$ is designated here as I3(P,R).

Here the matrix has the entries $t1, t2, \ldots, t6$ and $p1, p2, \ldots, p6$ which are latitude and longitude angles of the spherical coordinates of the oriented axes of the cylinders with corresponding numbers, the zero cylinder is the vertical pivot cylinder in this coordinate system that lies along the axis $z$. This coordinate system is identical to the one used in [3]. The heights of the touching points of the cylinders with the pivot $0^{th}$ cylinder are $z1, z2, \ldots, z6$. The radii of the cylinders are $r1, r2, \ldots r6$ correspondingly, and the radius of the $0^{th}$ cylinder is always 1.

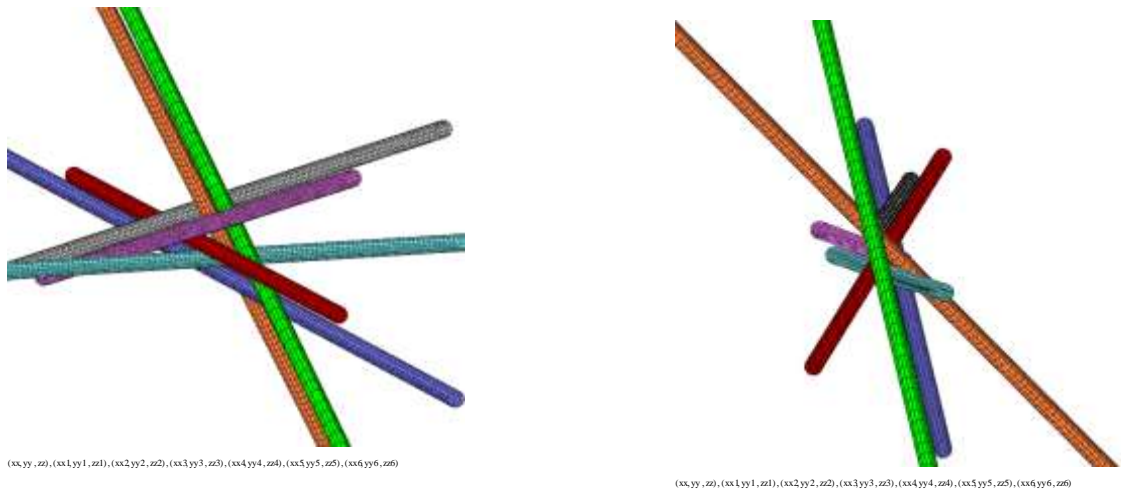

(xx,yy,zz),(xx1,yy1,zz1),(xx2,yy2,zz2),(xx3,yy3,zz3),(xx4,yy4,zz4),(xx5,yy5,zz5),(xx6,yy6,zz6)

(xx,yy,zz),(xx1,yy1,zz1),(xx2,yy2,zz2),(xx3,yy3,zz3),(xx4,yy4,zz4),(xx5,yy5,zz5),(xx6,yy6,zz6)

Fig. 7. Projections display a beautiful property of any 7-knot configuration that three unit vectors orthogonal to the directions of corresponding pairs of the cylinders be nearly coplanar. Quantitatively, the enclosed volume as the coplanar check for the three unit vectors is less than circa 0.01 for the arbitrary scissors angle $t1$. At a certain scissors angle this value can be exactly zero. Both projections are shown on the example of one and the same **a89** configuration.

| | | | |
|---|---|---|---|
| a24 := 16.736842105 | m18a := 16.563909774 | m18v := 15.377777777 | m18at := 19.0312 |
| a48 := 15.059210526 | m18b := 14.075471698 | m18w := 10.573208722 | m18au := 12.622691292 |
| a89 := 14.073170731 | m18c := 18.217105263 | m18x := 20.550724637 | m18av := 22.739018087 |
| ma28 := 20.973821989 | m18d := 11.843137254 | m18y := 13.264705882 | m18aw := 8.273972602 |
| ma29 := 16.378048780 | m18e := 4.823529411 | m18z := 12.906862745 | m18ax := 16.150819672 |
| ma49 := −1 | m18f := 20.790960453 | m18aa := 11.5 | m18ay := 19.064935064 |
| b27 := 17.503546099 | m18g := 20.307692307 | m18ab := 14.563176895 | m18az := 17.795454545 |
| b45 := 12.344827586 | m18h := 16.943396226 | m18ac := 17.296296296 | m18ba := 9.238095238 |
| mb56 := 18.164502164 | m18i := 12.451127819 | m18ad := 14.046511627 | m18bb := 16.146853146 |
| mb57 := 7.401459854 | m18j := 19.629629629 | m18ae := 19.158974359 | m18bc := 18.615384615 |
| mb68 := 18.368794326 | m18k := 15.14102564 | m18af := 13.72413793 | m18bd := 17.974025974 |
| c14 := 9.8 | m18l := 17.350993377 | m18ag := 18.8 | m18be := 16.536082474 |
| c24 := 16.526315789 | m18m := 19.523341523 | m18ah := 5.928571428 | m18bf := 13.573770491 |
| c27 := 16.543689320 | m18n := 11.287356321 | m18ai := 20.678899082 | m18bg := 7.379746835 |
| c46 := 0.526315789 | m18o := 19.200980392 | m18aj := 16.156626506 | m18bh := 19.242424242 |
| c48 := 16.353383458 | m18p := 4.466019417 | m18ak := 19.829787234 | m18bi := 16.679687 |
| mc13 := 14.269230769 | f56 := 11.875776397 | m18al := 15.4383202 | m18bj := 19.409836065 |
| d18 := 17.8437 | f57 := 7.307692307 | m18am := 12.508250825 | m18bk := 20.397849462 |
| d28 := 22.2 | f67 := 21.714285714 | m18an := −0.285714285 | m18bl := 11.512820512 |
| md16 := 20.120689655 | m18q := 13.945205479 | m18ao := 18.250626566 | m18bm := 4.746268656 |
| md24 := 17.375939849 | m18r := 13.467048710 | m18ap := 18.090909090 | m18bn := 23.3203125 |
| e48 := 17.165413533 | m18s := 11.609756097 | m18aq := 13.398412698 | m18bo := 20.382978723 |
| me28 := 15.283950617 | m18t := 19.674740484 | m18ar := 20.481818181 | m18bp := 11.003115264 |
| me29 := 13.310344827 | m18u := 13.753694581 | m18as := 18.194174757 | m18bq := 8.798301486 |
| | | | m18br := 17.355140186 |

Fig. 8. The list of different invariants of 7-knot configurations with the chirality matrix determinant -18.

We put in Fig. 8 only non-equivalent configurations with different invariants. Equivalent configurations like **a89** and **b47** should be the same as far as they have the identical $\wp =$ 14.(07317) so we omitted **b47** and the like from the list. A direct inspection of the image of the configurations confirms that they are topologically identical. The same is true for other coincidences.

Using the number invariant of Eq. (13) gives a powerful tool of complete classification and control over topologically different configurations. We can illustrate the selectivity properties of the number invariant of Eq. (14) on the example of **e9** and **d9** configurations of 8-knots that have different invariants and even have different determinants $|P| = 9$ and $|P| = 1$, respectively, as one can see from Fig. 9. Still $\wp_{\mathcal{R}}(e9) = \wp_{\mathcal{R}}(d9) = 40.4$ and if the numerical invariant works well it should indicate that both configurations have the equal ring matrices. We make sure that they indeed have one and the same ring matrix

$$\mathcal{R}(e9) = \mathcal{R}(d9) = \begin{pmatrix} 0 & 3 & 3 & 7 & 5 & 3 & 7 & 5 \\ 5 & 0 & 7 & 5 & 7 & 5 & 5 & 5 \\ 0 & 0 & 0 & 0 & 0 & 0 & 0 & 0 \\ 0 & 0 & 0 & 0 & 0 & 0 & 0 & 0 \\ 1 & 1 & 5 & 5 & 0 & 1 & 1 & 1 \\ 0 & 0 & 0 & 0 & 0 & 0 & 0 & 0 \\ 0 & 0 & 0 & 0 & 0 & 0 & 0 & 0 \\ 2 & 2 & 4 & 2 & 2 & 8 & 4 & 0 \end{pmatrix}. \quad (18)$$

|     |                        |                      |     |                        |                      |
|-----|------------------------|----------------------|-----|------------------------|----------------------|
| 1   | a1 := 13.727687742     | ma1 := 22.921022597  | 41  | c1 := 28.844580777     | mc1 := 17.0817995    |
| 9   | a2 := 14.329527991     | ma2 := 23.672228320  | 9   | c2 := 21.201474201     | mc2 := 23.7076167    |
| 1   | a3 := 14.527522935     | ma3 := 23.325688073  | 9   | c3 := 21.201474201     | mc3 := 23.7076167    |
| 9   | a4 := 20.428571428     | ma4 := 21            | −7  | c4 := 23.696428571     | mc4 := 20.4821428    |
| 1   | a5 := 24.442477876     | ma5 := 20.707964601  | 17  | c5 := 22.115207373     | mc5 := 22.88940092   |
| 25  | a6 := 13.126179563     | ma6 := 24.778916149  | 17  | c6 := 14.117647058     | mc6 := 28.6256684    |
| 1   | a7 := 26.811158798     | ma7 := 5.450643776   | 17  | c7 := 31.98            | mc7 := 3.835         |
| 9   | a8 := 15.754385964     | ma8 := 21.182456140  | 41  | c8 := 21.928746928     | mc8 := 22.9803439    |
| 9   | a9 := 15.991416309     | ma9 := 16.270386266  | 1   | c9 := 10.7777777777    | mc9 := 29.6222222    |
| 105 | b1 := 20.817931034     | mb1 := 16.575172413  | 17  | d1 := 20.395744680     | md1 := 24.8127659    |
| 185 | b2 := 13.812532637     | mb2 := 25.683028720  | 41  | d2 := 23.200772200     | md2 := 21.7104247    |
| 1   | b3 := 14.527522935     | mb3 := 23.325688073  | 9   | d3 := 21.201474201     | md3 := 23.7076167    |
| 65  | b4 := 26.866666666     | mb4 := 14.561904761  | 17  | d4 := 18               | md4 := 26.1785714    |
| 17  | b5 := 22.115207373     | mb5 := 22.889400921  | 1   | d5 := 24.442477876     | md5 := 20.7079646    |
| 209 | b6 := 16.466009064     | mb6 := 22.842975206  | 41  | d6 := 11.262086514     | md6 := 28.1959287    |
| 17  | b7 := 31.697399527     | mb7 := 2.869976359   | 1   | d7 := 29.115839243     | md7 := 5.45153664    |
| 185 | b8 := 16.499202551     | mb8 := 22.582137161  | 9   | d8 := 22.251518833     | md8 := 22.7424058    |
| 9   | b9 := 15.991416309     | mb9 := 16.270386266  | 1   | d9 := 10.7777777777    | md9 := 29.6222222    |
| 1   | e1 := 24.341523341     | me1 := 20.449631449  | 65  | x65a := 22.662579695   | m65a := 16.3751839   |
| 9   | e2 := 25.531001589     | me2 := 20.640699523  | 17  | x17a := 23.232558139   | m17a := 15.0418604   |
| 1   | e3 := 20.584633853     | me3 := 24.623049219  | 185 | x185a := 8.635220125   | m185a := 30.874213   |
| 9   | e4 := 24.196428571     | me4 := 19.982142857  | 41  | x41a := 22.005815305   | m41a := 16.680856    |
| 1   | e5 := 24.442477876     | me5 := 20.707964601  | 185 | x185b := 31.257861635  | m185b := 8.157232    |
| 25  | e6 := 12.213592233     | me6 := 26.825242718  | 17  | x17b := 14.523809523   | m17b := 26.904761    |
| 1   | e7 := 26.811158798     | me7 := 5.450643776   |     |                        |                      |
| 9   | e8 := 21.751243781     | me8 := 22.935323383  |     |                        |                      |
| 9   | e9 := 23.044444444     | me9 := 17.355555555  |     |                        |                      |

Fig. 9. The invariants $\wp$ for 8-knots. Before the name of each configuration we give the determinant of the corresponding chirality matrix. Note that in case of the even chirality matrix its determinant is the same for the mirror configuration. Stray 8-knot configurations obtained in various ways not from any of 9-knots are marked with the letter x and also given with their mirror ones. For reader's attention we put some of the names of topologically equal configurations: **a3**=**b3**; **a9**=**b9**; **c2**=**c3**=**d3**; **b5**=**c5**; **a5**=**d5**=**e5**; **c9**=**d9**; **a7**=**e7**.

The ring matrix gives the information as to how many of cylinders are not knotted, which is the number of rows with all zeroes. There are configurations of 7-knots that have only two unknotted cylinders (see Appendix 4) with configurations: **m18r**, **m18ab**, **m18aq**, **m18cc**, **x42aa**, **x42aj**, **m102k**, **m150n**, **m150o**, **m150q**, and **m102g**.

7. **Configurations of equal radii for the 7-knot.**

Surprisingly, the chiral matrix alone can tell whether a configuration can have equal radii for all cylinders or not. It is easy to understand that the configuration with EE cannot have all cylinders being of equal radius because even the two cylinders being in EE cannot have equal radii. Suppose they have. The scissors angle between them can be made arbitrary small since no other cylinder hinders their mutual rotation at least in one direction as far as they are in EE. But then all other cylinders are forced to lie on these two cylinders in two planes from both sides. In one of the planes there are at least three cylinders. One of them lying in the plane and being between the other two inevitably blocks the other two from the mutual contact. Not having EE is the necessary condition for the configuration to have all equal cylinders.

The configurations with the chirality matrices that do not have EE can have equal radii of the cylinders but only asymptotically, except the **a89**. As we discuss below this very configuration does not have a possibility for *parallel* cylinders of equal radii to exist because of the topological impossibility of the Necker cube configuration. A collection of equal radii configurations is given in Appendix 3. They are not many and some of them are non-trivial though contain a pair of parallel cylinders: those are **ma37** and **ma49**. Interestingly, when the radii of the cylinders in both configurations approach 1 the large triangle extends self-similarly with the sharp angle of circa 15 degrees for **ma37** and 22 degrees for **ma49**.

8. **The uniqueness of the non-trivial configuration of mutually touching 7 round infinite cylinders of equal radii.**

The uniqueness of configuration **a89** of mutually touching 7 round infinite cylinders of equal radii has been suspected after the publication of [4] that independently rediscovered in 2013 the same configuration that we reported in 2009 [2]. Now we will produce the arguments that support its uniqueness. Any other configuration with equal radii should approach a degenerate one with some cylinders nearly parallel that we call "cube" (see Fig. 10a and Appendix 3). Many configurations, such as **c24**, **md24**, **m18c**, etc. with $|P| = -18$ are such.

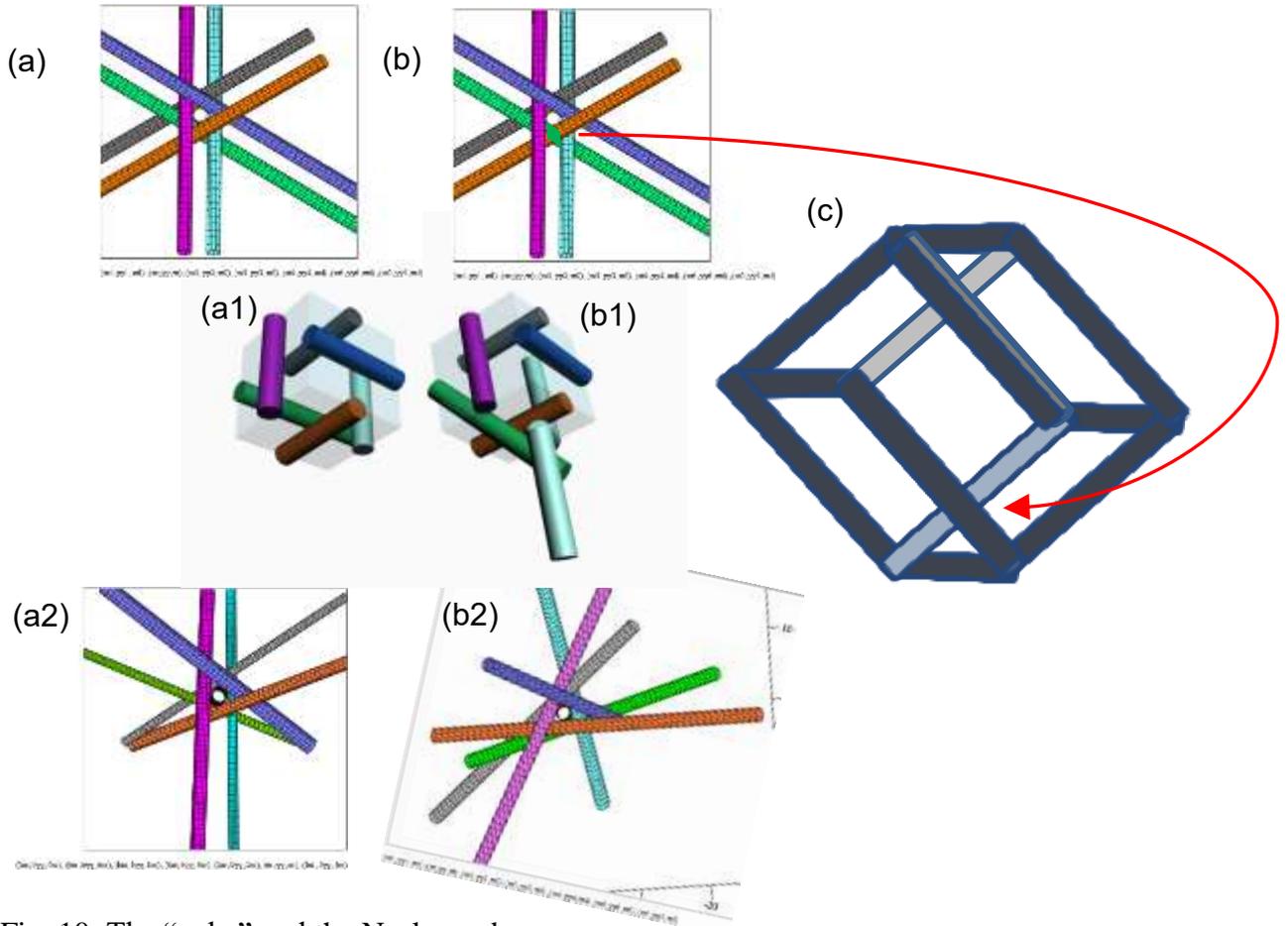

Fig. 10. The "cube" and the Necker cube.

The specific topology of configuration **a89** makes it similar to an Impossible cube (the Necker cube [5]) given in Fig. 10c. Let us show this with the help of the plates in Fig. 10. If one compares Fig. 10a and Fig. 10b one finds that the only difference is that the green cylinder goes over the brown one in Fig. 10b. Fig. 10b1 shows that it is impossible for the parallel cylinders unless there is a break in the cylinder. The red arrow connector shows the similarity to the Necker cube in Fig. 10c. The plates Fig. 10a2 and Fig. 10b2 show the real configurations of mutually touching cylinders topologically equivalent to Fig. 10a and Fig. 10b, respectively. Again, the only difference in both configurations is the switching between the green and the brown cylinders. The unique topology of **a89** is revealed even more if one notices that it fits the impossible Penrose triangle [6,7] (Fig. 11).

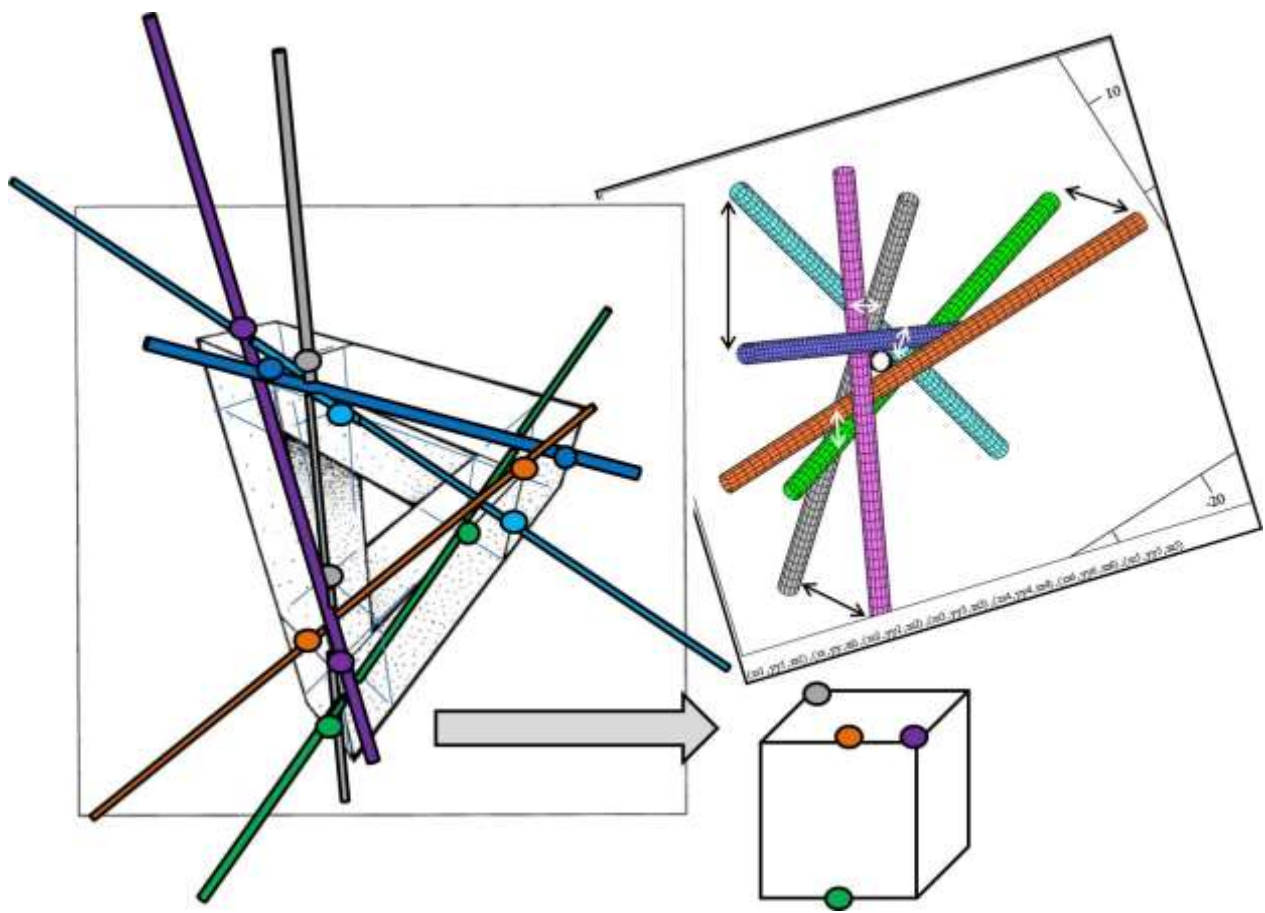

Fig. 11. Schematic of **a89** (given in the right plate) on the background of the impossible Penrose triangle (left).

Decoration of corner cubes emphasizes the position of the lines on the triangle. On the right plate the arrows show the pairs of the cylinders in a top view of the pivot cylinder of **a89** so that the white arrows mark the parts of one pair that is squeezed between the "legs" of the other pair marked with the black arrows. Configuration **a89** has its mirror **ma89** and no other configurations with non-parallel all equal cylinders are possible. The configuration **m18c** that can have asymptotically equal cylinders given in Appendix 3 is the closest to **a89** geometrically because it has a similar chirality matrix and a similar ring matrix. However as one can see its invariant is different and its image is identical to the "normal cube" of Fig. 10a with parallel cylinders.

9. Conclusion

We suggested a classification scheme for mutually pairwise touching infinite cylinders. Our approach relies on the properties of the chirality matrix and the newly introduced ring matrix. The chirality matrix contains enough information to restrict possible cylinder configurations and establish the upper limit for the number of mutually pairwise touching infinite cylinders of adjustable cross-sections. The chirality matrix helps to point out to those configurations that may have equal radii in the case of round cylinders. The Necker cube topology of the non-trivial equal radii configuration of the seven cylinders found earlier stipulates its uniqueness. Finally, revealing the topology characteristics of cylinder configurations presented here may help to understand and quantify the difference and similarities in auxetic behaviour between our approach of the touching cylinder regular network [1-3] and the auxetic geometry of expanding periodic bar-and-joint frameworks recently published in [8].

# Supplementary Material

Appendix 1
7-cross invariants

−2
| | | |
|---|---|---|
| a12 := 3.1278688525 | m2a := −6 | m2ab := 0.7904191617 |
| a15 := 7.4242424242 | m2b := −4.8235294118 | m2ac := 19.7142857143 |
| a59 := 6.8780487805 | m2c := −10.380952381 | m2ad := 14.2330097087 |
| ma34 := −0.4854368932 | m2d := 31.8979591837 | m2ae := 32.3235294118 |
| ma37 := 5.5714285714 | m2e := 21.1666666667 | m2af := 26.9693877551 |
| a13 := 13.4923076923 | m2f := 22.7636363636 | m2ag := 11.1827956989 |
| a19 := 4.4324324324 | m2g := 18.3333333333 | m2ah := 29.0882352941 |
| a23 := 13.7272727273 | m2h := 17.3333333333 | m2ai := 13.5342465753 |
| a25 := 18.5744680851 | m2i := 16.8965517241 | m2aj := 2.3636363636 |
| a57 := 26.0810810811 | m2j := 9.6904761905 | m2ak := 7.6229508197 |
| ma17 := 1.4705882353 | m2k := 8.4523809524 | m2al := 33.3768115942 |
| ma18 := 9.3179190751 | m2l := 12.2272727273 | m2am := 11.2057416268 |
| ma35 := 14.7142857143 | m2m := 21.6354166667 | m2an := 1.6271186441 |
| ma38 := 26.2340425532 | m2n := 21.7604166667 | m2ao := 24.2023809524 |
| ma45 := 18.8965517241 | m2o := −8.5531914894 | m2ap := 17.6666666667 |
| ma47 := 16.5365853659 | f25 := 16.1578947368 | m2aq := 19.68 |
| a17m2 := 19.2156862745 | f26 := 26.5 | m2ar := 1.9642857143 |
| b17 := 30.4385026738 | f28 := −9.2380952381 | m2as := 1.2307692308 |
| b18 := 16.2160278746 | m2p := 10.2307692308 | m2at := 23.9191919192 |
| c17 := 32.1005025126 | m2q := 22.1090909091 | m2au := 5.0277777778 |
| mc19 := 14.2142857143 | m2r := 16.7341772152 | m2av := 22.7628205128 |
| mc26 := 6.3448275862 | m2s := 32.0869565217 | m2aw := 26.2804878049 |
| c68 := 16.0194174757 | m2t := 19.5769230769 | m2ax := 23.1052631579 |
| c69 := 7.8 | m2u := 16.793814433 | m2ay := 19.8461538462 |
| md68 := 6.0491803279 | m2v := 11.841509434 | m2az := 6.8472222222 |
| e13 := 15.9166666667 | m2w := 12 | |
| me17 := 4.3037974684 | m2x := 15.4330708661 | |
| me18 := 7.5483870968 | m2y := 28.9438202247 | |
| e19 := 15.7142857143 | m2z := −2.393442623 | |
| me38 := 29 | m2aa := −1.2173913043 | |

−10
| | | |
|---|---|---|
| a14 := 7.9270072993 | m10a := 0.75 | m10w := 0.4 |
| a16 := 8.180878553 | m10b := 6.4230769231 | m10x := 6.3770491803 |
| a46 := −0.8823529412 | m10c := 8.9818181818 | m10y := −0.8072289157 |
| a56 := 8.4242424242 | m10d := 8.4727272727 | m10z := 4.3577981651 |
| a58 := 1.7142857143 | m10e := −0.0618556701 | m10aa := 6.9827586207 |
| ma26 := 12.6612903226 | m10f := 12.0260416667 | m10ab := −1.3617021277 |
| ma27 := −0.7073170732 | m10g := 9.0175438596 | m10ac := 6.2597402597 |
| ma36 := 10.7094017094 | m10h := −0.2745098039 | m10ad := 7.2307692308 |
| ma39 := 2.5573770492 | m10i := 14.4934383202 | m10ae := 7.5333333333 |
| ma67 := 8.3770491803 | f29 := 8.3571428571 | m10af := 15.9887955182 |
| ma69 := 11.2151898734 | m10j := 7.4835164835 | m10ag := 2.8645833333 |
| c12 := 7.5172413793 | m10k := 11.4090909091 | m10ah := 1.1523178808 |
| c29 := −0.2068965517 | m10l := 4.3939393939 | m10ai := 11.0309278351 |
| mc18 := 6.2571428571 | m10m := 3.417721519 | m10aj := 9.2234042553 |
| md27 := −1.0294117647 | m10n := −0.3311258278 | m10ak := 5.93 |

| | |
|---|---|
| md67 := 10.1384615385 | m10o := −5.2941176471 |
| e14 := 1.2 | m10p := 7.4805194805 |
| e16 := 7.7313432836 | m10q := 13.676 |
| me26 := 11.4652777778 | m10r := 11.7487437186 |
| me36 := 10.8601398601 | m10s := 7.3440860215 |
| me39 := −0.3448275862 | m10t := 1.0782608696 |
| e46 := 3.8947368421 | m10u := 11.7666666667 |
| me69 := 7.9 | m10v := 8.71875 |

−18

| | | |
|---|---|---|
| a24 := 16.7368421053 | m18n := 11.2873563218 | m18av := 22.7390180879 |
| a48 := 15.0592105263 | m18o := 19.2009803922 | m18aw := 8.2739726027 |
| a89 := 14.0731707317 | m18p := 4.4660194175 | m18ax := 16.1508196721 |
| ma28 := 20.9738219895 | f56 := 11.8757763975 | m18ay := 19.0649350649 |
| ma29 := 16.3780487805 | f57 := 7.3076923077 | m18az := 17.7954545455 |
| ma49 := −1 | f67 := 21.7142857143 | m18ba := 9.2380952381 |
| b27 := 17.5035460993 | m18q := 13.9452054795 | m18bb := 16.1468531469 |
| b45 := 12.3448275862 | m18r := 13.4670487106 | m18bc := 18.6153846154 |
| mb56 := 18.1645021645 | m18s := 11.6097560976 | m18bd := 17.974025974 |
| mb57 := 7.401459854 | m18t := 19.6747404844 | m18be := 16.5360824742 |
| mb68 := 18.3687943262 | m18u := 13.7536945813 | m18bf := 13.5737704918 |
| c14 := 9.8 | m18v := 15.3777777778 | m18bg := 7.3797468354 |
| c24 := 16.5263157895 | m18w := 10.5732087227 | m18bh := 19.2424242424 |
| c27 := 16.5436893204 | m18x := 20.5507246377 | m18bi := 16.6796875 |
| c46 := 0.5263157895 | m18y := 13.2647058824 | m18bj := 19.4098360656 |
| c48 := 16.3533834586 | m18z := 12.9068627451 | m18bk := 20.3978494624 |
| mc13 := 14.2692307692 | m18aa := 11.5 | m18bl := 11.5128205128 |
| d18 := 17.84375 | m18ab := 14.5631768953 | m18bm := 4.7462686567 |
| d28 := 22.2 | m18ac := 17.2962962963 | m18bn := 23.3203125 |
| md16 := 20.1206896552 | m18ad := 14.0465116279 | m18bo := 20.3829787234 |
| md24 := 17.3759398496 | m18ae := 19.158974359 | m18bp := 11.0031152648 |
| e48 := 17.1654135338 | m18af := 13.724137931 | m18bq := 8.7983014862 |
| me28 := 15.2839506173 | m18ag := 18.8 | m18br := 17.3551401869 |
| me29 := 13.3103448276 | m18ah := 5.9285714286 | m18bs := 22.695890411 |
| m18a := 16.5639097744 | m18ai := 20.6788990826 | m18bt := 11.6097560976 |
| m18b := 14.0754716981 | m18aj := 16.156626506 | m18bu := 4.6376811594 |
| m18c := 18.2171052632 | m18ak := 19.829787234 | m18bv := 4.102739726 |
| m18d := 11.8431372549 | m18al := 15.43832021 | m18bw := 12.9342105263 |
| m18e := 4.8235294118 | m18am := 12.5082508251 | m18bx := 20.0970873786 |
| m18f := 20.790960452 | m18an := −0.2857142857 | m18by := 14.8235294118 |
| m18g := 20.3076923077 | m18ao := 18.2506265664 | m18bz := 5.1475409836 |
| m18h := 16.9433962264 | m18ap := 18.0909090909 | m18ca := 8.1428571429 |
| m18i := 12.4511278195 | m18aq := 13.3984126984 | m18cb := 20.5714285714 |
| m18j := 19.6296296296 | m18ar := 20.4818181818 | m18cc := 16.1330108827 |
| m18k := 15.141025641 | m18as := 18.1941747573 | m18cd := 12.9407894737 |
| m18l := 17.3509933775 | m18at := 19.03125 | m18ce := 9.6437659033 |
| m18m := 19.5233415233 | m18au := 12.6226912929 | |
| d12 := 5.6981132075 | x42h := 8.3862660944 | x42aa := 15.0678127985 |
| d14 := −6.6 | x42i := 5.2948717949 | x42ab := −0.198757764 |
| d46 := −5.0526315789 | x42j := 12.2 | x42ac := −0.2857142857 |
| b12 := 7.6444444444 | x42k := 4.7483443709 | x42ad := 9.6504065041 |
| b15 := 11.625 | x42l := 10.3921568627 | x42ae := 11.375 |
| b58 := −2.6206896552 | x42m := −0.0952380952 | x42af := 10.6616541353 |
| md26 := 0.598540146 | x42n := 11.0406091371 | x42ag := 7.8156028369 |
| mc16 := 14.5105263158 | x42o := 1.6603773585 | x42ah := 11.2146596859 |
| mc67 := 10.188976378 | x42p := 4.5853658537 | x42ai := 5.5451713396 |
| mc78 := 9.2857142857 | x42q := 8.2727272727 | x42aj := 6.5258358663 |
| mb26 := 13.3502538071 | x42r := 5.1976047904 | x42ak := 2.1134020619 |
| mb78 := 7.4237288136 | x42s := 5.4545454545 | |
| x42a := −10.3529411765 | x42t := 4.5373134328 | |
| x42b := 12.3786407767 | x42u := 16.2103990326 | |
| x42c := 12.3333333333 | x42v := 13.3300492611 | |
| x42d := 2.4607329843 | x42w := 4.5490196078 | |
| x42e := 3.198019802 | x42x := −0.3404255319 | |
| x42f := 16.6864864865 | x42y := 4.4754098361 | |
| x42g := 12.3349753695 | x42z := 12.2178217822 | |

−54

    m54a := 3.3150684932
    m54b := 1.3675213675
    m54c := −3.86
    m54d := −5.4143646409
    m54e := 1.9404761905
    m54f := 3.5344827586
    m54g := −4.9827586207

    a2466 := 19.05907173
    x66a := 11.2810810811
    x66b := 15.0476190476
    x66c := 16.8034188034
    x66d := 22.9705882353
    x66e := 15.3210526316
    x66f := 21.5294117647
    x66g := 20.1696428571
    x66h := 17.3717948718
    x66i := 8.8195488722

−102

    a24m102 := 14.5185185185
    ma24102a := 17.6
    ma24102b := 14.2941176471
    b14 := 14.7299270073
    b16 := 17.4590909091
    b46 := 14.2941176471
    m102a := 13.7023809524
    m102b := 15.8468899522
    m102c := 15.2205882353
    m102d := 15.3805970149
    m102e := 18.775
    m102f := 15.9375
    m102g := 15.1657848325
    m102h := 16.0277777778
    m102i := 18.8242424242
    m102j := 16.8080808081
    m102k := 15.8082706767
    m102l := 9.4736842105
    m102m := 16.19
    m102n := 16.2112676056
    m102o := 14.455026455
    m102p := 14.7878787879
    m102q := 15.8770949721
    m102r := 17.3617886179
    m102s := 17.0418250951
    m102t := 14.3214285714
    m102u := 17.391025641

−150

    b24 := 19.4962406015
    mb28 := 25.9052132701
    b48 := 18.2763157895
    a48m150 := 19.2857142857
    m150a := 19.7402597403
    m150b := 15.8482142857
    m150c := 13.1189189189
    m150d := 15.223628692
    m150e := 11.6714975845
    m150f := 12.9111111111
    m150g := 14.9160305344
    m150h := 23.8510638298
    m150i := 18.0212765957
    m150j := 12.2274881517
    m150k := 13.2955223881
    m150l := 16.4615384615
    m150m := 21.4342105263
    m150n := 11.1624548736
    m150o := 19.0739614995
    m150p := 14.525198939
    m150q := 11.5339412361

−162

    m162a := 15.8409090909
    m162b := 12.557319224
    m162c := 14.2978723404
    m162d := 13.3173652695
    m162e := 11.1226765799
    m162f := 13.4059405941
    m162g := 14.7235772358
    m162h := 12.0238095238
    m162i := 14.380952381
    m162j := 13.4487179487

# Appendix 2

## 9-cross invariants

a0a := 20.8373932182
a0b := 31.2786803919
a0c := 10.3290278722
a0d := 34.5886547812
a0e := 40.2095514799
a0f := −19.6963902439
a0g := −17.8644067797
a0i := −16.6838848921
a0j := 24.7750955542
a4a := 25.3527263569
a4b := 27.0174639332
a16 := 27.6275057092
am52a := 24.0193570498
a16b := 26.2657929528
a4c := 25.7484010746
a4d := 17.3700244249
a4e := 22.6623342175
a4f := 21.4014997657
am80a := 28.2738683402
am80b := 24.6140502709

ma0a := 23.8438467564
ma0b := 25.4237876721
ma0c := 46.4133242692
ma0d := 22.0575688817
ma0e := 15.9676775014
ma0f := 65.5061463415
ma0g := 63.0869565217
ma0i := 63.3452278177
ma0j := 20.8016495675

am52b := 32.4097486672
am84a := 30.0826834633
a16c := 17.5732072194
a4g := −6.2395821696
am80c := 23.9197560976
a4h := 29.8644473579
a4i := 24.2303267974
a4j := 34.1105794102
a4k := 28.7292739021
a0h := 43.8854285065
a16d := 29.7640239342
a16e := −229.8012820513
a16f := 33.2238669172
a16g := 21.9618226601
a104a := 23.2314697698

Appendix 2
9-cross
Configuration label a104a

$$\begin{pmatrix} \text{orange} \\ \text{red} \\ \text{blue} \\ \text{green} \\ \text{cyan} \\ \text{magenta} \\ \text{gray} \\ \text{olive} \\ \text{pink} \end{pmatrix} \begin{pmatrix} 0 & 0 & 0 & 1 \\ t1 & p1 & z1 & r1 \\ t2 & p2 & z2 & r2 \\ t3 & p3 & z3 & r3 \\ t4 & p4 & z4 & r4 \\ t5 & p5 & z5 & r5 \\ t6 & p6 & z6 & r6 \\ t7 & p7 & z7 & r7 \\ t8 & p8 & z8 & r8 \end{pmatrix} = \begin{pmatrix} 0 & 0 & 0 & 1 \\ 3.6159045172 & 3.1415926536 & 0 & 12.7464268152 \\ 3.6471234187 & 15.7076980451 & -28.1964771012 & 0.0264410812 \\ 4.4610324534 & 14.2687792623 & -29.9442332468 & 0.0215838157 \\ 3.5024946052 & 10.8291087486 & -11.7316631231 & 2.7631602455 \\ 4.3896094306 & 4.5325796511 & -26.1949005192 & 0.0414398674 \\ 5.2966856921 & 9.1488334071 & -27.7808437608 & 1.2792883208 \\ 0.4753260593 & 3.2079035145 & -27.9222944335 & 0.0160532311 \\ 1.7154517039 & 1.66791117 & -29.4018268684 & 0.9115071088 \end{pmatrix}$$

$$P = \begin{pmatrix} 0 & -1 & -1 & -1 & -1 & -1 & -1 & 1 & 1 \\ -1 & 0 & -1 & -1 & 1 & 1 & -1 & -1 & 1 \\ -1 & -1 & 0 & 1 & -1 & -1 & 1 & 1 & -1 \\ -1 & -1 & 1 & 0 & 1 & 1 & -1 & -1 & -1 \\ -1 & 1 & -1 & 1 & 0 & -1 & -1 & 1 & -1 \\ -1 & 1 & -1 & 1 & -1 & 0 & -1 & -1 & 1 \\ -1 & -1 & 1 & -1 & -1 & -1 & 0 & -1 & 1 \\ 1 & -1 & 1 & -1 & 1 & -1 & -1 & 0 & 1 \\ 1 & 1 & -1 & -1 & -1 & 1 & 1 & 1 & 0 \end{pmatrix} \qquad R = \begin{pmatrix} 0 & 0 & 0 & 0 & 0 & 0 & 0 & 0 & 0 \\ 0 & 0 & 0 & 0 & 0 & 0 & 0 & 0 & 0 \\ 0 & 0 & 0 & 0 & 0 & 0 & 0 & 0 & 0 \\ 7 & 10 & 5 & 0 & 6 & 5 & 7 & 5 & 9 \\ 0 & 0 & 0 & 0 & 0 & 0 & 0 & 0 & 0 \\ 3 & 3 & 4 & 3 & 12 & 0 & 4 & 3 & 4 \\ 8 & 8 & 7 & 7 & 8 & 8 & 0 & 7 & 7 \\ 7 & 4 & 11 & 4 & 5 & 8 & 5 & 0 & 4 \\ 2 & 10 & 5 & 2 & 2 & 2 & 2 & 5 & 0 \end{pmatrix}$$

$|P| = 104$

$I3(P,R) = 23.2314697698$

$I3(-P,R) = 33.6931444473$

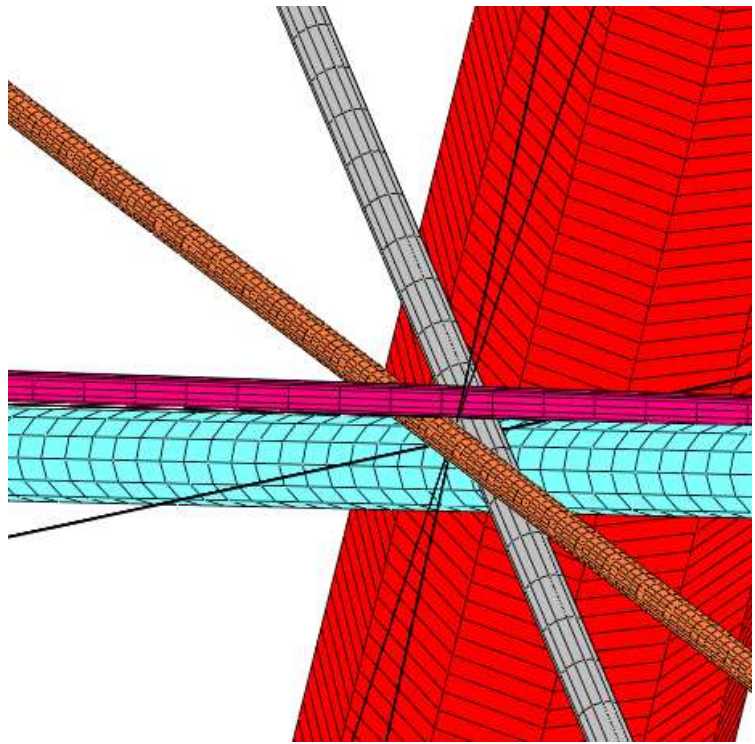

Appendix 2
9-cross
Configuration label a0j

$$\begin{pmatrix} \text{orange} \\ \text{red} \\ \text{blue} \\ \text{green} \\ \text{cyan} \\ \text{magenta} \\ \text{gray} \\ \text{olive} \\ \text{pink} \end{pmatrix} \begin{pmatrix} 0 & 0 & 0 & 1 \\ t_1 & p_1 & z_1 & r_1 \\ t_2 & p_2 & z_2 & r_2 \\ t_3 & p_3 & z_3 & r_3 \\ t_4 & p_4 & z_4 & r_4 \\ t_5 & p_5 & z_5 & r_5 \\ t_6 & p_6 & z_6 & r_6 \\ t_7 & p_7 & z_7 & r_7 \\ t_8 & p_8 & z_8 & r_8 \end{pmatrix} = \begin{pmatrix} 0 & 0 & 0 & 1 \\ 1.9586383944 & 3.1415926536 & 0 & 1.3075218005 \\ 2.4312039098 & -0.4960015764 & -3.9413620072 & 5.51215675 \\ 1.5257358336 & -1.2834566523 & 0.8861892933 & 0.0794574576 \\ 2.5339757481 & -1.595664017 & -8.0991142589 & 0.9653327888 \\ 0.5904057601 & 8.5334827203 & -6.7999227653 & 0.1061532441 \\ 0.3475231492 & 4.1301093819 & -1.6395400906 & 0.4466675943 \\ 1.5937336823 & 2.1083459158 & -4.4226586756 & 4.6006495871 \\ 0.5148958868 & 5.0283842979 & 10.5241855909 & 1.7084311671 \end{pmatrix}$$

$$P = \begin{pmatrix} 0 & 1 & 1 & 1 & 1 & 1 & 1 & 1 & 1 \\ 1 & 0 & -1 & -1 & 1 & -1 & 1 & -1 & -1 \\ 1 & -1 & 0 & 1 & -1 & 1 & 1 & -1 & 1 \\ 1 & -1 & 1 & 0 & -1 & 1 & -1 & -1 & -1 \\ 1 & 1 & -1 & -1 & 0 & 1 & 1 & -1 & -1 \\ 1 & -1 & 1 & 1 & 1 & 0 & 1 & 1 & 1 \\ 1 & 1 & 1 & -1 & 1 & 1 & 0 & 1 & -1 \\ 1 & -1 & -1 & -1 & -1 & 1 & 1 & 0 & 1 \\ 1 & -1 & 1 & -1 & -1 & 1 & -1 & 1 & 0 \end{pmatrix} \quad R = \begin{pmatrix} 0 & 2 & 10 & 2 & 2 & 5 & 2 & 5 & 2 \\ 6 & 0 & 6 & 9 & 6 & 6 & 6 & 6 & 9 \\ 0 & 0 & 0 & 0 & 0 & 0 & 0 & 0 & 0 \\ 5 & 2 & 2 & 0 & 2 & 5 & 2 & 2 & 10 \\ 2 & 2 & 2 & 5 & 0 & 2 & 2 & 10 & 5 \\ 4 & 4 & 4 & 3 & 3 & 0 & 3 & 12 & 3 \\ 7 & 5 & 11 & 4 & 8 & 4 & 0 & 4 & 5 \\ 0 & 0 & 0 & 0 & 0 & 0 & 0 & 0 & 0 \\ 0 & 0 & 0 & 0 & 0 & 0 & 0 & 0 & 0 \end{pmatrix}$$

$|P| = -3.7470027081 \times 10^{-15}$

$I3(P, R) = 24.7750955542$

$I3(-P, R) = 20.8016495675$

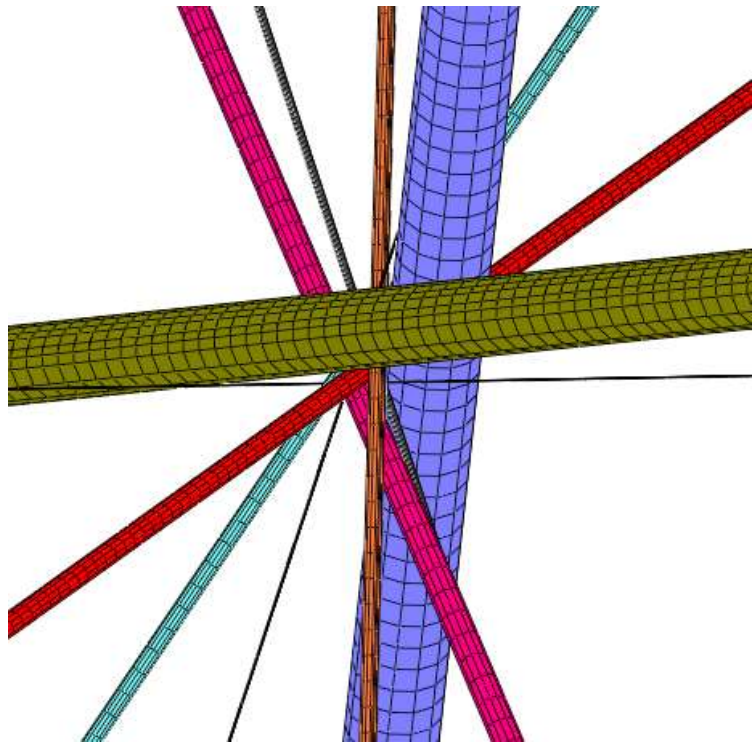

Appendix 2
9-cross
Configuration label a16g

$$\begin{pmatrix} \text{orange} \\ \text{red} \\ \text{blue} \\ \text{green} \\ \text{cyan} \\ \text{magenta} \\ \text{gray} \\ \text{olive} \\ \text{pink} \end{pmatrix} \begin{pmatrix} 0 & 0 & 0 & 1 \\ t_1 & p_1 & z_1 & r_1 \\ t_2 & p_2 & z_2 & r_2 \\ t_3 & p_3 & z_3 & r_3 \\ t_4 & p_4 & z_4 & r_4 \\ t_5 & p_5 & z_5 & r_5 \\ t_6 & p_6 & z_6 & r_6 \\ t_7 & p_7 & z_7 & r_7 \\ t_8 & p_8 & z_8 & r_8 \end{pmatrix} = \begin{pmatrix} 0 & 0 & 0 & 1 \\ 0.8041996236 & 3.1415926536 & 0 & 0.0434500445 \\ 0.5817171932 & 3.16792351 & 1.9381760752 & 0.7953349257 \\ 0.9250165926 & 2.4665439012 & -0.4254425452 & 0.1481285692 \\ 1.8969547221 & 4.8662589349 & 1.595850801 & 0.3503249865 \\ 2.3798663618 & 4.6665177898 & -0.4257947425 & 0.2911489214 \\ 0.7376527843 & 2.2297534273 & -1.2121984795 & 0.055731755 \\ 0.5955478753 & 2.5319425226 & -0.9353491824 & 0.0258429553 \\ 0.6430160735 & 2.5035894535 & -0.708333231 & 0.0100000003 \end{pmatrix}$$

$$P = \begin{pmatrix} 0 & 1 & 1 & 1 & 1 & 1 & 1 & 1 & 1 \\ 1 & 0 & -1 & 1 & -1 & 1 & -1 & -1 & 1 \\ 1 & -1 & 0 & -1 & 1 & 1 & -1 & -1 & -1 \\ 1 & 1 & -1 & 0 & -1 & -1 & -1 & 1 & 1 \\ 1 & -1 & 1 & -1 & 0 & -1 & 1 & 1 & 1 \\ 1 & 1 & 1 & -1 & -1 & 0 & -1 & 1 & 1 \\ 1 & -1 & -1 & -1 & 1 & -1 & 0 & 1 & -1 \\ 1 & -1 & -1 & 1 & 1 & 1 & 1 & 0 & 1 \\ 1 & 1 & -1 & 1 & 1 & 1 & -1 & 1 & 0 \end{pmatrix} \qquad R = \begin{pmatrix} 0 & 0 & 0 & 0 & 0 & 0 & 0 & 0 & 0 \\ 11 & 0 & 4 & 4 & 4 & 5 & 8 & 5 & 7 \\ 0 & 0 & 0 & 0 & 0 & 0 & 0 & 0 & 0 \\ 4 & 7 & 8 & 0 & 5 & 5 & 11 & 4 & 4 \\ 0 & 0 & 0 & 0 & 0 & 0 & 0 & 0 & 0 \\ 2 & 2 & 2 & 5 & 10 & 0 & 5 & 2 & 2 \\ 0 & 0 & 0 & 0 & 0 & 0 & 0 & 0 & 0 \\ 7 & 7 & 6 & 5 & 10 & 9 & 5 & 0 & 5 \\ 8 & 4 & 11 & 7 & 5 & 4 & 4 & 5 & 0 \end{pmatrix}$$

$|P| = 16$

$I3(P, R) = 21.9618226601$

$I3(-P, R) = 35.1369634061$

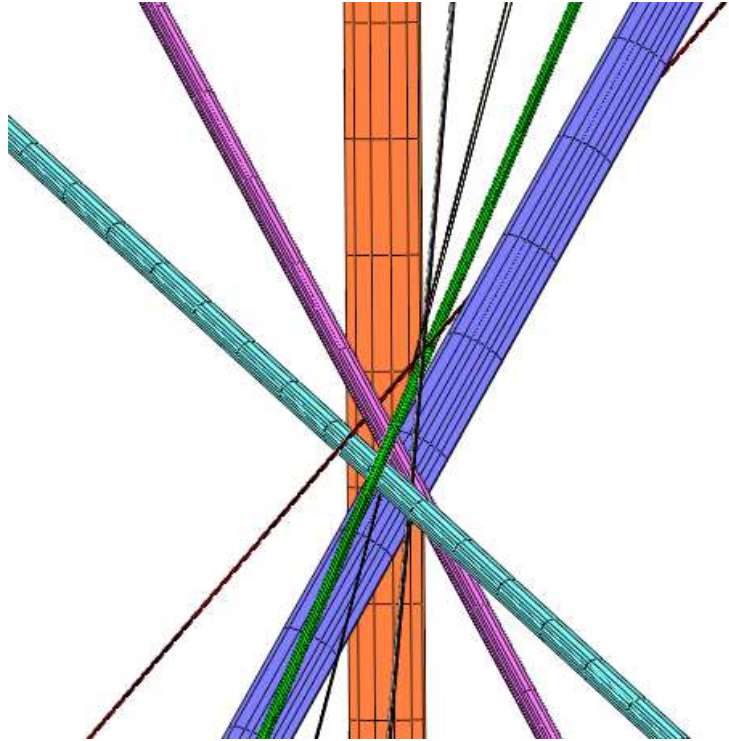

Appendix 2
9-cross
Configuration label a16f

$$\begin{pmatrix} \text{orange} \\ \text{red} \\ \text{blue} \\ \text{green} \\ \text{cyan} \\ \text{magenta} \\ \text{gray} \\ \text{olive} \\ \text{pink} \end{pmatrix} \begin{pmatrix} 0 & 0 & 0 & 1 \\ t_1 & p_1 & z_1 & r_1 \\ t_2 & p_2 & z_2 & r_2 \\ t_3 & p_3 & z_3 & r_3 \\ t_4 & p_4 & z_4 & r_4 \\ t_5 & p_5 & z_5 & r_5 \\ t_6 & p_6 & z_6 & r_6 \\ t_7 & p_7 & z_7 & r_7 \\ t_8 & p_8 & z_8 & r_8 \end{pmatrix} = \begin{pmatrix} 0 & 0 & 0 & 1 \\ 1.6955521322 & 3.1415926536 & 0 & 0.0260027908 \\ 2.1256365591 & 3.1196854471 & 1.1866196019 & 0.1338203207 \\ 0.6165463826 & 3.3306952633 & 0.5547247366 & 0.036683492 \\ 2.9472293108 & 1.7832438493 & 4.9828529395 & 0.1184071733 \\ 0.1584179202 & 2.0278148404 & -4.5943575093 & 0.1046700367 \\ 2.0934372625 & 3.4527932085 & -0.6307324638 & 1.0393544902 \\ 2.5797877454 & 3.2753516767 & -0.6889038469 & 0.0115084535 \\ 0.108897054 & 1.5664445131 & -11.1282948661 & 0.1744861218 \end{pmatrix}$$

$$P = \begin{pmatrix} 0 & 1 & 1 & 1 & 1 & 1 & 1 & 1 & 1 \\ 1 & 0 & 1 & -1 & 1 & -1 & 1 & 1 & -1 \\ 1 & 1 & 0 & 1 & -1 & -1 & 1 & 1 & -1 \\ 1 & -1 & 1 & 0 & 1 & 1 & 1 & -1 & 1 \\ 1 & 1 & -1 & 1 & 0 & 1 & -1 & -1 & -1 \\ 1 & -1 & -1 & 1 & 1 & 0 & 1 & -1 & -1 \\ 1 & 1 & 1 & 1 & -1 & 1 & 0 & -1 & 1 \\ 1 & 1 & 1 & -1 & -1 & -1 & -1 & 0 & -1 \\ 1 & -1 & -1 & 1 & -1 & -1 & 1 & -1 & 0 \end{pmatrix} \qquad R = \begin{pmatrix} 0 & 0 & 0 & 0 & 0 & 0 & 0 & 0 & 0 \\ 12 & 0 & 4 & 4 & 4 & 3 & 3 & 3 & 3 \\ 0 & 0 & 0 & 0 & 0 & 0 & 0 & 0 & 0 \\ 5 & 8 & 7 & 0 & 4 & 4 & 11 & 5 & 4 \\ 3 & 3 & 4 & 3 & 0 & 3 & 4 & 4 & 12 \\ 6 & 5 & 5 & 7 & 9 & 0 & 7 & 5 & 10 \\ 0 & 0 & 0 & 0 & 0 & 0 & 0 & 0 & 0 \\ 8 & 8 & 8 & 7 & 8 & 7 & 7 & 0 & 7 \\ 0 & 0 & 0 & 0 & 0 & 0 & 0 & 0 & 0 \end{pmatrix}$$

$|P| = 16$

$I3(P,R) = 33.2238669172$

$I3(-P,R) = 24.7434123092$

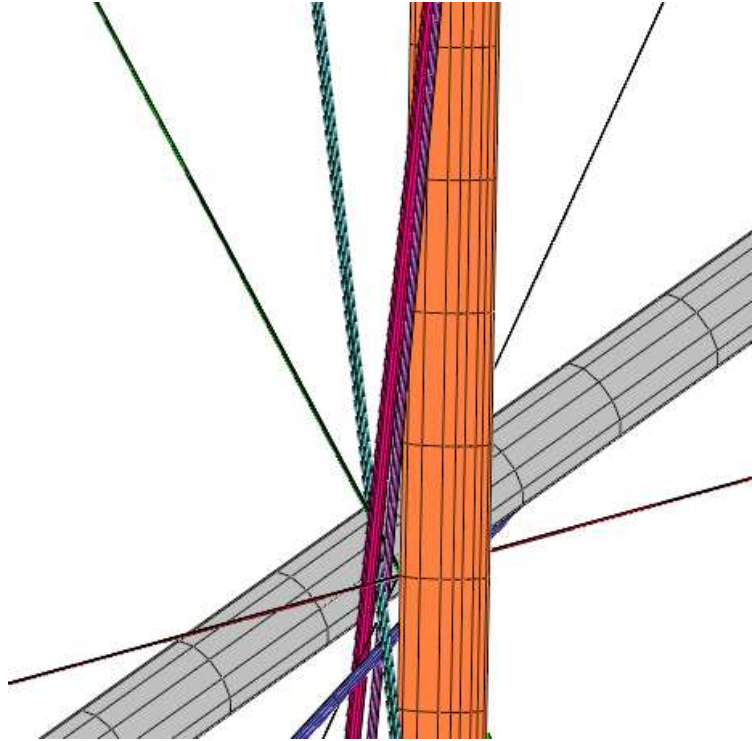

Appendix 2
9-cross
Configuration label a0i

$$\begin{pmatrix} \text{orange} \\ \text{red} \\ \text{blue} \\ \text{green} \\ \text{cyan} \\ \text{magenta} \\ \text{gray} \\ \text{olive} \\ \text{pink} \end{pmatrix} \begin{pmatrix} 0 & 0 & 0 & 1 \\ t_1 & p_1 & z_1 & r_1 \\ t_2 & p_2 & z_2 & r_2 \\ t_3 & p_3 & z_3 & r_3 \\ t_4 & p_4 & z_4 & r_4 \\ t_5 & p_5 & z_5 & r_5 \\ t_6 & p_6 & z_6 & r_6 \\ t_7 & p_7 & z_7 & r_7 \\ t_8 & p_8 & z_8 & r_8 \end{pmatrix} = \begin{pmatrix} 0 & 0 & 0 & 1 \\ 1.5989896879 & 3.1415926536 & 0 & 0.3900345267 \\ 1.7178216079 & 3.1371629465 & 21.1483270324 & 1.4573608218 \\ 1.1128163648 & 0.038926593 & -23.8393094004 & 0.362599636 \\ 1.768236628 & 2.8954323991 & 6.691610052 & 10.8313139628 \\ -1.705699145 & 6.1348109238 & -14.7960910148 & 30.8054598997 \\ 2.8624333709 & 5.8694832165 & -36.7494146048 & 0.7819383766 \\ 1.9632443888 & 3.1350724308 & -18.3987280763 & 0.1554806057 \\ 1.9844738257 & 3.1581236904 & -19.4304301885 & 1.3468389228 \end{pmatrix}$$

$$P = \begin{pmatrix} 0 & 1 & 1 & 1 & 1 & -1 & 1 & 1 & 1 \\ 1 & 0 & 1 & 1 & 1 & -1 & 1 & -1 & 1 \\ 1 & 1 & 0 & -1 & -1 & -1 & 1 & -1 & 1 \\ 1 & 1 & -1 & 0 & -1 & -1 & -1 & -1 & 1 \\ 1 & 1 & -1 & -1 & 0 & 1 & -1 & 1 & 1 \\ -1 & -1 & -1 & -1 & 1 & 0 & 1 & 1 & 1 \\ 1 & 1 & 1 & -1 & -1 & 1 & 0 & -1 & -1 \\ 1 & -1 & -1 & -1 & 1 & 1 & -1 & 0 & 1 \\ 1 & 1 & 1 & 1 & 1 & 1 & -1 & 1 & 0 \end{pmatrix} \quad R = \begin{pmatrix} 0 & 3 & 3 & 9 & 6 & 9 & 6 & 3 & 3 \\ 6 & 0 & 9 & 7 & 8 & 7 & 8 & 6 & 9 \\ 0 & 0 & 0 & 0 & 0 & 0 & 0 & 0 & 0 \\ 1 & 1 & 6 & 0 & 1 & 1 & 1 & 1 & 6 \\ 3 & 3 & 12 & 4 & 0 & 3 & 3 & 4 & 4 \\ 0 & 0 & 0 & 0 & 0 & 0 & 0 & 0 & 0 \\ 2 & 2 & 2 & 2 & 2 & 10 & 0 & 5 & 5 \\ 6 & 6 & 6 & 9 & 6 & 6 & 6 & 0 & 9 \\ 0 & 0 & 0 & 0 & 0 & 0 & 0 & 0 & 0 \end{pmatrix}$$

$|P| = 0$

$I3(P, R) = -16.6838848921$

$I3(-P, R) = 63.3452278177$

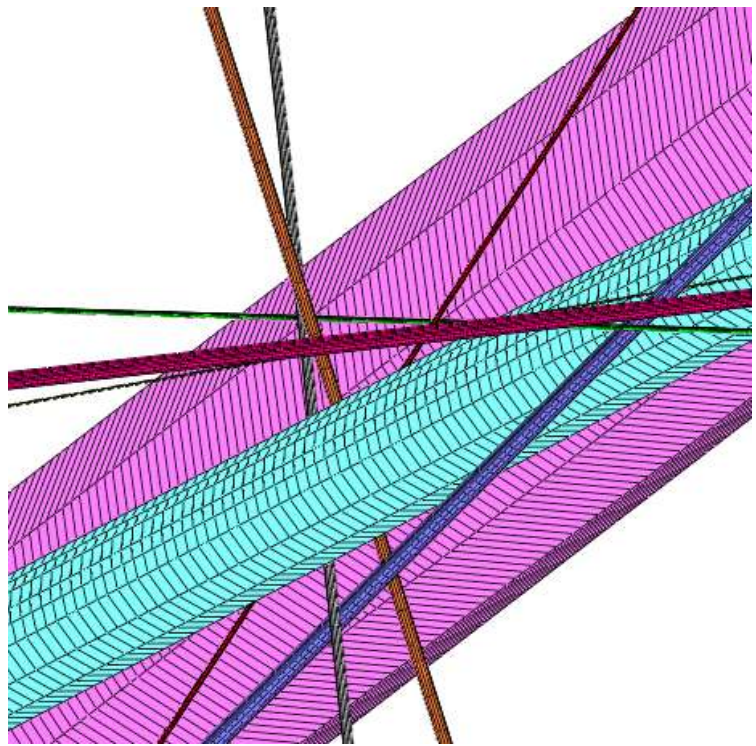

Appendix 2
9-cross
Configuration label a16e

$$\begin{pmatrix} \text{orange} \\ \text{red} \\ \text{blue} \\ \text{green} \\ \text{cyan} \\ \text{magenta} \\ \text{gray} \\ \text{olive} \\ \text{pink} \end{pmatrix} \begin{pmatrix} 0 & 0 & 0 & 1 \\ t_1 & p_1 & z_1 & r_1 \\ t_2 & p_2 & z_2 & r_2 \\ t_3 & p_3 & z_3 & r_3 \\ t_4 & p_4 & z_4 & r_4 \\ t_5 & p_5 & z_5 & r_5 \\ t_6 & p_6 & z_6 & r_6 \\ t_7 & p_7 & z_7 & r_7 \\ t_8 & p_8 & z_8 & r_8 \end{pmatrix} = \begin{pmatrix} 0 & 0 & 0 & 1 \\ 1.7846537154 & 3.1415926536 & 0 & 0.2976733205 \\ 2.149017431 & 4.8763675791 & -0.8507159308 & 0.0151779391 \\ 2.0546859163 & 5.8961117366 & -3.3645659432 & 1.1067467124 \\ 1.955117469 & 5.37560695 & -0.4596866099 & 1.0206370654 \\ 2.4029758679 & 4.8791741969 & -5.6847061704 & 1.7116774838 \\ 0.942557675 & 2.599286595 & -1.4413285361 & 0.3973400647 \\ 1.9633295594 & 4.8598834268 & -1.2843978616 & 0.0610508236 \\ 2.1015764012 & 5.3141935544 & -1.3106762568 & 0.112325474 \end{pmatrix}$$

$$P = \begin{pmatrix} 0 & 1 & 1 & 1 & 1 & 1 & 1 & 1 & 1 \\ 1 & 0 & -1 & -1 & -1 & 1 & -1 & 1 & -1 \\ 1 & -1 & 0 & 1 & -1 & 1 & -1 & -1 & 1 \\ 1 & -1 & 1 & 0 & 1 & -1 & 1 & 1 & 1 \\ 1 & -1 & -1 & 1 & 0 & -1 & 1 & -1 & -1 \\ 1 & 1 & 1 & -1 & -1 & 0 & 1 & 1 & -1 \\ 1 & -1 & -1 & 1 & 1 & 1 & 0 & 1 & 1 \\ 1 & 1 & -1 & 1 & -1 & 1 & 1 & 0 & -1 \\ 1 & -1 & 1 & 1 & -1 & -1 & 1 & -1 & 0 \end{pmatrix} \quad R = \begin{pmatrix} 0 & 1 & 1 & 6 & 1 & 1 & 6 & 1 & 1 \\ 6 & 0 & 6 & 6 & 6 & 9 & 9 & 6 & 6 \\ 6 & 1 & 0 & 1 & 6 & 1 & 1 & 1 & 1 \\ 0 & 0 & 0 & 0 & 0 & 0 & 0 & 0 & 0 \\ 0 & 0 & 0 & 0 & 0 & 0 & 0 & 0 & 0 \\ 0 & 0 & 0 & 0 & 0 & 0 & 0 & 0 & 0 \\ 1 & 1 & 6 & 1 & 1 & 6 & 0 & 1 & 1 \\ 7 & 5 & 4 & 11 & 5 & 8 & 4 & 0 & 4 \\ 6 & 5 & 10 & 6 & 10 & 5 & 6 & 6 & 0 \end{pmatrix}$$

$|P| = 16$

$I3(P, R) = -229.8012820513$

$I3(-P, R) = -42.8429487179$

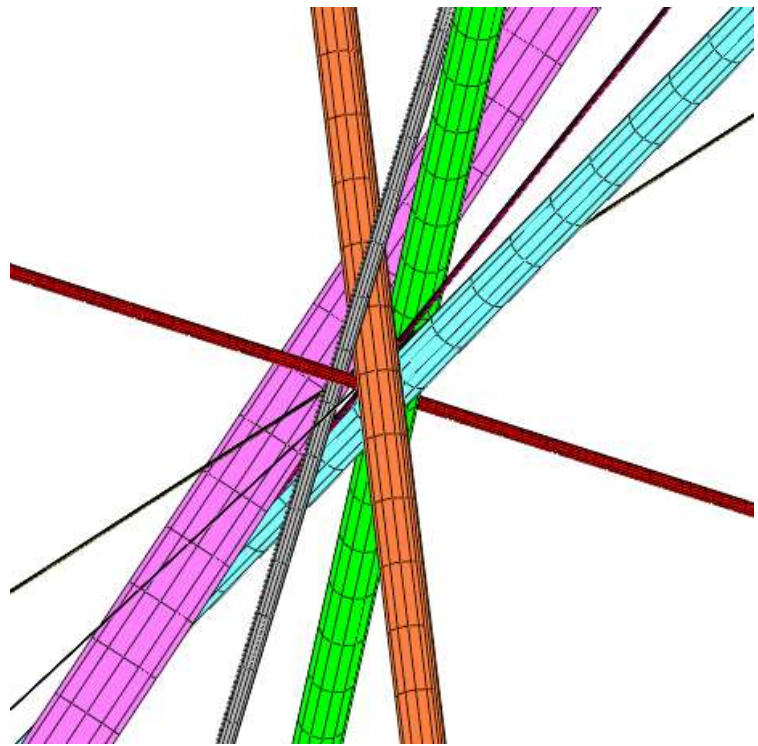

Appendix 2
9-cross
Configuration label a16d

$$\begin{pmatrix} \text{orange} \\ \text{red} \\ \text{blue} \\ \text{green} \\ \text{cyan} \\ \text{magenta} \\ \text{gray} \\ \text{olive} \\ \text{pink} \end{pmatrix} \begin{pmatrix} 0 & 0 & 0 & 1 \\ t_1 & p_1 & z_1 & r_1 \\ t_2 & p_2 & z_2 & r_2 \\ t_3 & p_3 & z_3 & r_3 \\ t_4 & p_4 & z_4 & r_4 \\ t_5 & p_5 & z_5 & r_5 \\ t_6 & p_6 & z_6 & r_6 \\ t_7 & p_7 & z_7 & r_7 \\ t_8 & p_8 & z_8 & r_8 \end{pmatrix} = \begin{pmatrix} 0 & 0 & 0 & 1 \\ 2.5622734476 & 3.1415926536 & 0 & 0.6826712717 \\ 1.0727666193 & 3.4281886493 & 9.82454416 & 1.1384585694 \\ 2.2771822054 & 4.1202230939 & 2.765542016 & 4.9913950652 \\ 0.2500820235 & 1.4593669402 & -8.8046490914 & 0.2283086078 \\ 0.4755337624 & 1.4958704478 & 0.4026814723 & 0.0837579078 \\ 2.7634546352 & 4.1126329069 & -1.625988671 & 0.0100000003 \\ 0.2005623996 & 1.5386179518 & -10.6887623332 & 0.0806573671 \\ 2.9090525788 & 4.0823844398 & -9.0221361324 & 0.0664625337 \end{pmatrix}$$

$$P = \begin{pmatrix} 0 & 1 & 1 & 1 & 1 & 1 & 1 & 1 & 1 \\ 1 & 0 & -1 & -1 & -1 & 1 & -1 & -1 & 1 \\ 1 & -1 & 0 & 1 & -1 & -1 & 1 & -1 & 1 \\ 1 & -1 & 1 & 0 & 1 & 1 & -1 & 1 & -1 \\ 1 & -1 & -1 & 1 & 0 & -1 & -1 & 1 & -1 \\ 1 & 1 & -1 & 1 & -1 & 0 & -1 & 1 & 1 \\ 1 & -1 & 1 & -1 & -1 & -1 & 0 & -1 & -1 \\ 1 & -1 & -1 & 1 & 1 & 1 & -1 & 0 & 1 \\ 1 & 1 & 1 & -1 & -1 & 1 & -1 & 1 & 0 \end{pmatrix} \qquad R = \begin{pmatrix} 0 & 0 & 0 & 0 & 0 & 0 & 0 & 0 & 0 \\ 12 & 0 & 3 & 3 & 4 & 3 & 3 & 4 & 4 \\ 0 & 0 & 0 & 0 & 0 & 0 & 0 & 0 & 0 \\ 0 & 0 & 0 & 0 & 0 & 0 & 0 & 0 & 0 \\ 0 & 0 & 0 & 0 & 0 & 0 & 0 & 0 & 0 \\ 5 & 6 & 10 & 5 & 7 & 0 & 9 & 5 & 7 \\ 6 & 3 & 9 & 9 & 3 & 3 & 0 & 6 & 3 \\ 6 & 6 & 5 & 10 & 6 & 6 & 5 & 0 & 10 \\ 3 & 3 & 3 & 4 & 12 & 4 & 4 & 3 & 0 \end{pmatrix}$$

$|P| = 16$

$I3(P, R) = 29.7640239342$

$I3(-P, R) = 27.7499065071$

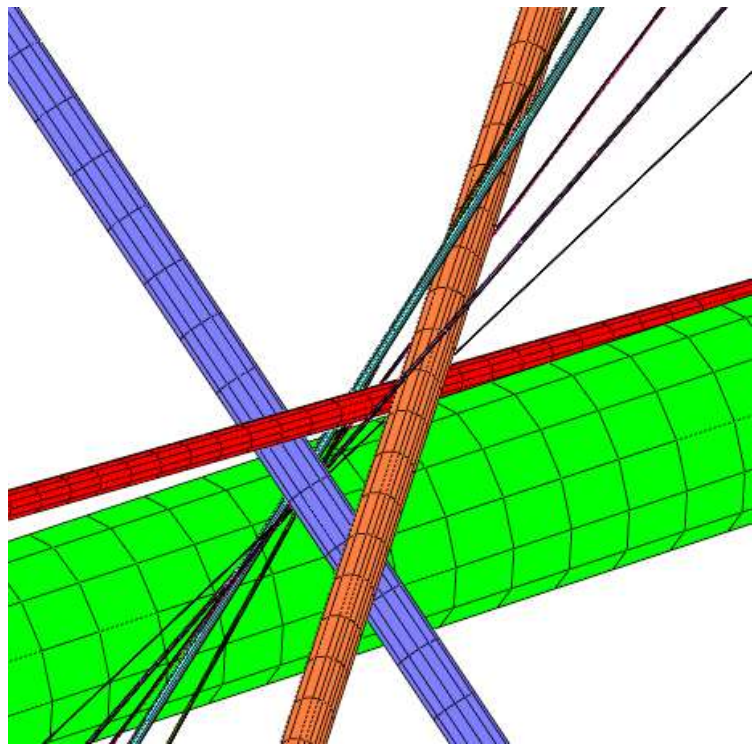

Appendix 2
9-cross
Configuration label a0h

$$\begin{pmatrix} \text{orange} \\ \text{red} \\ \text{blue} \\ \text{green} \\ \text{cyan} \\ \text{magenta} \\ \text{gray} \\ \text{olive} \\ \text{pink} \end{pmatrix} \begin{pmatrix} 0 & 0 & 0 & 1 \\ t1 & p1 & z1 & r1 \\ t2 & p2 & z2 & r2 \\ t3 & p3 & z3 & r3 \\ t4 & p4 & z4 & r4 \\ t5 & p5 & z5 & r5 \\ t6 & p6 & z6 & r6 \\ t7 & p7 & z7 & r7 \\ t8 & p8 & z8 & r8 \end{pmatrix} = \begin{pmatrix} 0 & 0 & 0 & 1 \\ 0.5457246655 & 3.1415926536 & 0 & 0.1825734315 \\ 0.7327597662 & 9.4757393687 & -4.4339207685 & 6.1475445899 \\ 1.1043126937 & 7.228500969 & -3.9906007363 & 0.2523992152 \\ 0.5665860022 & 8.8817078492 & -0.3370613761 & 0.1898117496 \\ 0.5022607733 & 8.9036595843 & -1.5690765114 & 0.0543667867 \\ 2.6146475015 & 5.5471029923 & -1.1712232517 & 0.3129178368 \\ 2.7103309344 & -0.9910201149 & -2.300399585 & 0.3149959406 \\ 2.4684794754 & 4.6516026407 & -0.3474402013 & 0.0618381786 \end{pmatrix}$$

$$P = \begin{pmatrix} 0 & 1 & 1 & 1 & 1 & 1 & 1 & 1 & 1 \\ 1 & 0 & 1 & 1 & 1 & -1 & 1 & 1 & 1 \\ 1 & 1 & 0 & 1 & 1 & 1 & -1 & -1 & -1 \\ 1 & 1 & 1 & 0 & -1 & 1 & 1 & -1 & 1 \\ 1 & 1 & 1 & -1 & 0 & 1 & -1 & -1 & 1 \\ 1 & -1 & 1 & 1 & 1 & 0 & 1 & -1 & -1 \\ 1 & 1 & -1 & 1 & -1 & 1 & 0 & -1 & -1 \\ 1 & 1 & -1 & -1 & -1 & -1 & -1 & 0 & 1 \\ 1 & 1 & -1 & 1 & 1 & -1 & -1 & 1 & 0 \end{pmatrix} \quad R = \begin{pmatrix} 0 & 3 & 3 & 12 & 3 & 3 & 4 & 4 & 4 \\ 6 & 0 & 6 & 9 & 9 & 3 & 3 & 3 & 3 \\ 0 & 0 & 0 & 0 & 0 & 0 & 0 & 0 & 0 \\ 0 & 0 & 0 & 0 & 0 & 0 & 0 & 0 & 0 \\ 0 & 0 & 0 & 0 & 0 & 0 & 0 & 0 & 0 \\ 5 & 4 & 11 & 4 & 7 & 0 & 8 & 5 & 4 \\ 0 & 0 & 0 & 0 & 0 & 0 & 0 & 0 & 0 \\ 3 & 3 & 4 & 3 & 4 & 4 & 12 & 0 & 3 \\ 7 & 8 & 7 & 7 & 8 & 7 & 8 & 8 & 0 \end{pmatrix}$$

$|P| = -4.1633363423 \times 10^{-15}$

$I3(P,R) = 43.8854285065$

$I3(-P,R) = 12.6262029249$

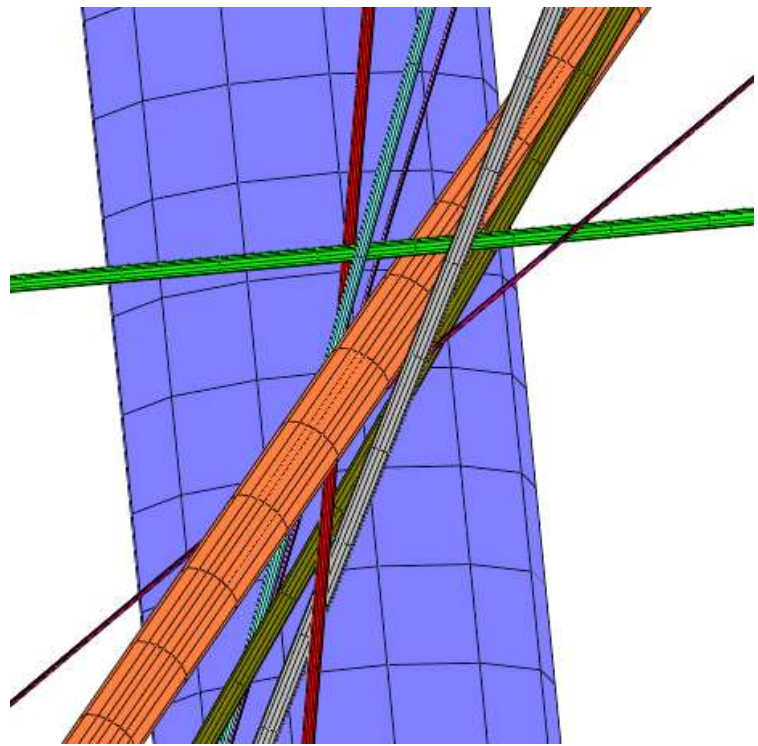

Appendix 2
9-cross
Configuration label a4k

$$\begin{pmatrix} \text{orange} \\ \text{red} \\ \text{blue} \\ \text{green} \\ \text{cyan} \\ \text{magenta} \\ \text{gray} \\ \text{olive} \\ \text{pink} \end{pmatrix} \begin{pmatrix} 0 & 0 & 0 & 1 \\ t_1 & p_1 & z_1 & r_1 \\ t_2 & p_2 & z_2 & r_2 \\ t_3 & p_3 & z_3 & r_3 \\ t_4 & p_4 & z_4 & r_4 \\ t_5 & p_5 & z_5 & r_5 \\ t_6 & p_6 & z_6 & r_6 \\ t_7 & p_7 & z_7 & r_7 \\ t_8 & p_8 & z_8 & r_8 \end{pmatrix} = \begin{pmatrix} 0 & 0 & 0 & 1 \\ 2.994723984 & 3.1415926536 & 0 & 0.1849609631 \\ 2.8598837075 & 2.6900657369 & 4.2576672553 & 0.0767172779 \\ 2.9473798622 & 9.6482790877 & 1.825562307 & 0.0750189535 \\ 2.929440114 & 2.9000897964 & 6.0726784104 & 0.2116212282 \\ 2.9742886539 & -2.9921450169 & 1.6324493756 & 0.0218147246 \\ 0.1815482625 & 1.7900077977 & 1.7419886291 & 0.010000179 \\ 2.8530561706 & 2.9410130425 & 2.9099958667 & 0.0466844304 \\ 2.7136560471 & 3.0288609384 & 3.9340229253 & 0.0407533452 \end{pmatrix}$$

$$P = \begin{pmatrix} 0 & 1 & 1 & 1 & 1 & 1 & 1 & 1 & 1 \\ 1 & 0 & 1 & -1 & 1 & -1 & 1 & 1 & 1 \\ 1 & 1 & 0 & -1 & -1 & -1 & -1 & 1 & -1 \\ 1 & -1 & -1 & 0 & 1 & -1 & 1 & -1 & 1 \\ 1 & 1 & -1 & 1 & 0 & 1 & -1 & 1 & 1 \\ 1 & -1 & -1 & -1 & 1 & 0 & -1 & -1 & 1 \\ 1 & 1 & -1 & 1 & -1 & -1 & 0 & 1 & 1 \\ 1 & 1 & 1 & -1 & 1 & -1 & 1 & 0 & -1 \\ 1 & 1 & -1 & 1 & 1 & 1 & 1 & -1 & 0 \end{pmatrix} \quad R = \begin{pmatrix} 0 & 0 & 0 & 0 & 0 & 0 & 0 & 0 & 0 \\ 0 & 0 & 0 & 0 & 0 & 0 & 0 & 0 & 0 \\ 9 & 7 & 0 & 5 & 5 & 5 & 7 & 10 & 6 \\ 2 & 2 & 5 & 0 & 10 & 2 & 2 & 5 & 2 \\ 0 & 0 & 0 & 0 & 0 & 0 & 0 & 0 & 0 \\ 7 & 6 & 5 & 9 & 10 & 0 & 5 & 5 & 7 \\ 2 & 10 & 2 & 2 & 5 & 5 & 0 & 2 & 2 \\ 6 & 6 & 1 & 1 & 1 & 1 & 1 & 0 & 1 \\ 8 & 4 & 5 & 4 & 7 & 4 & 11 & 5 & 0 \end{pmatrix}$$

$|P| = 4$

$I3(P, R) = 28.7292739021$

$I3(-P, R) = 17.7995838353$

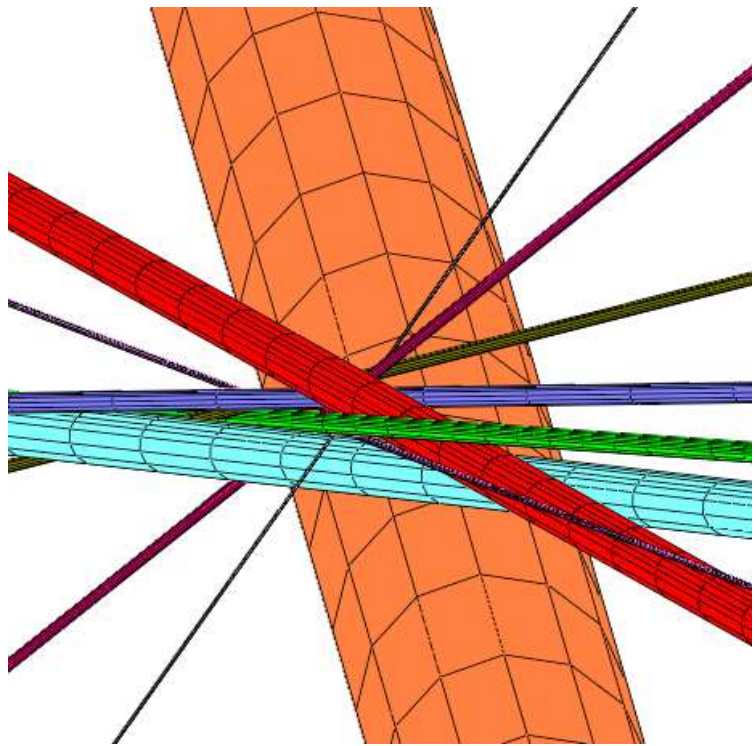

Appendix 2
9-cross
Configuration label a0g

$$\begin{pmatrix} \text{orange} \\ \text{red} \\ \text{blue} \\ \text{green} \\ \text{cyan} \\ \text{magenta} \\ \text{gray} \\ \text{olive} \\ \text{pink} \end{pmatrix} \begin{pmatrix} 0 & 0 & 0 & 1 \\ t_1 & p_1 & z_1 & r_1 \\ t_2 & p_2 & z_2 & r_2 \\ t_3 & p_3 & z_3 & r_3 \\ t_4 & p_4 & z_4 & r_4 \\ t_5 & p_5 & z_5 & r_5 \\ t_6 & p_6 & z_6 & r_6 \\ t_7 & p_7 & z_7 & r_7 \\ t_8 & p_8 & z_8 & r_8 \end{pmatrix} = \begin{pmatrix} 0 & 0 & 0 & 1 \\ 1.747329354 & 3.1415926536 & 0 & 1.0115499624 \\ 2.1247109916 & 3.3938393128 & 0.4269362765 & 0.0685147831 \\ 1.5903285771 & -3.5467816347 & 0.8843541831 & 0.0820356347 \\ 1.9618957607 & 3.2451355147 & 1.2620077595 & 0.2268695551 \\ 1.6497964417 & 9.2286660755 & 0.7346764212 & 0.0395819055 \\ 2.6902276163 & 4.4119404014 & -1.3891274311 & 0.0100000003 \\ 2.1730765077 & 3.2244995418 & 0.4326610702 & 0.0311033041 \\ 2.2887055058 & 3.1904326972 & 0.8080783793 & 0.0279803562 \end{pmatrix}$$

$$P = \begin{pmatrix} 0 & 1 & 1 & 1 & 1 & 1 & 1 & 1 & 1 \\ 1 & 0 & -1 & 1 & -1 & 1 & -1 & -1 & -1 \\ 1 & -1 & 0 & 1 & 1 & 1 & 1 & -1 & 1 \\ 1 & 1 & 1 & 0 & -1 & 1 & -1 & 1 & -1 \\ 1 & -1 & 1 & -1 & 0 & -1 & 1 & -1 & -1 \\ 1 & 1 & 1 & 1 & -1 & 0 & 1 & 1 & -1 \\ 1 & -1 & 1 & -1 & 1 & 1 & 0 & -1 & 1 \\ 1 & -1 & -1 & 1 & -1 & 1 & -1 & 0 & 1 \\ 1 & -1 & 1 & -1 & -1 & -1 & 1 & 1 & 0 \end{pmatrix} \quad R = \begin{pmatrix} 0 & 0 & 0 & 0 & 0 & 0 & 0 & 0 & 0 \\ 0 & 0 & 0 & 0 & 0 & 0 & 0 & 0 & 0 \\ 9 & 7 & 0 & 5 & 5 & 5 & 7 & 10 & 6 \\ 2 & 2 & 5 & 0 & 10 & 2 & 2 & 5 & 2 \\ 0 & 0 & 0 & 0 & 0 & 0 & 0 & 0 & 0 \\ 7 & 6 & 5 & 9 & 10 & 0 & 5 & 5 & 7 \\ 4 & 11 & 4 & 4 & 5 & 5 & 0 & 7 & 8 \\ 6 & 6 & 1 & 1 & 1 & 1 & 1 & 0 & 1 \\ 6 & 1 & 1 & 1 & 6 & 1 & 1 & 1 & 0 \end{pmatrix}$$

$|P| = -4.0800696155 \times 10^{-15}$

$I3(P, R) = -17.8644067797$

$I3(-P, R) = 63.0869565217$

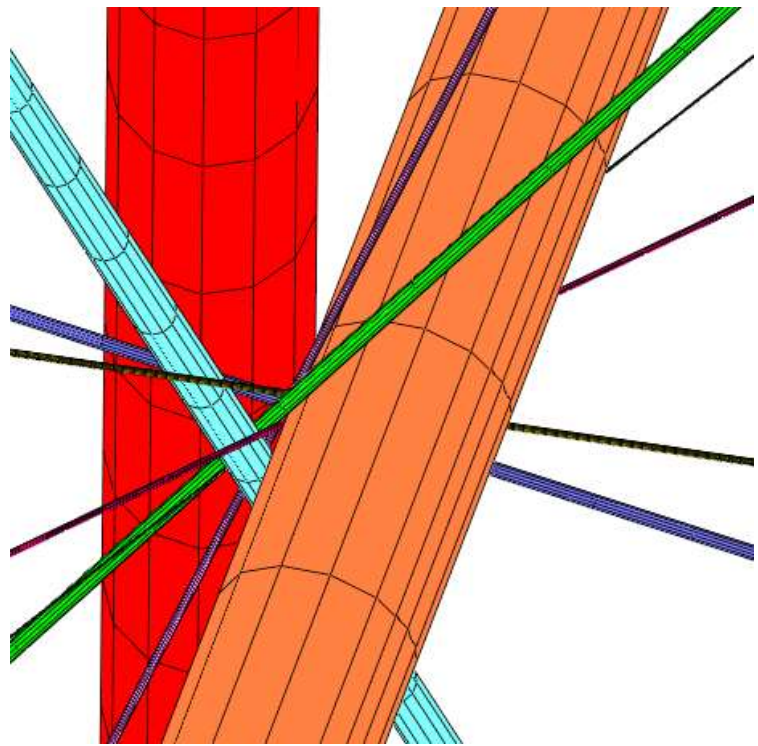

Appendix 2
9-cross
Configuration label a0f

$$\begin{pmatrix} \text{orange} \\ \text{red} \\ \text{blue} \\ \text{green} \\ \text{cyan} \\ \text{magenta} \\ \text{gray} \\ \text{olive} \\ \text{pink} \end{pmatrix} \begin{pmatrix} 0 & 0 & 0 & 1 \\ t1 & p1 & z1 & r1 \\ t2 & p2 & z2 & r2 \\ t3 & p3 & z3 & r3 \\ t4 & p4 & z4 & r4 \\ t5 & p5 & z5 & r5 \\ t6 & p6 & z6 & r6 \\ t7 & p7 & z7 & r7 \\ t8 & p8 & z8 & r8 \end{pmatrix} = \begin{pmatrix} 0 & 0 & 0 & 1 \\ 0.1954140012 & 3.1415926536 & 0 & 0.4465748441 \\ 0.5033958776 & 4.2221494641 & 8.864027872 & 0.2726595369 \\ 0.3619147647 & -3.475242815 & 5.3396098657 & 0.0420129909 \\ 0.3502652565 & 3.8104489481 & 11.955977224 & 0.7517904552 \\ 0.3327828688 & 9.1122487433 & 5.1424823572 & 0.0100000004 \\ 2.7908772789 & 5.0376834993 & 5.3605280806 & 0.0308820687 \\ 0.7942601842 & 3.6206712338 & 6.6392338438 & 0.2444111998 \\ 0.5610232046 & 3.4248999887 & 8.1875189853 & 0.2664259746 \end{pmatrix}$$

$$P = \begin{pmatrix} 0 & 1 & 1 & 1 & 1 & 1 & 1 & 1 & 1 \\ 1 & 0 & -1 & 1 & -1 & 1 & -1 & -1 & -1 \\ 1 & -1 & 0 & 1 & 1 & 1 & 1 & -1 & 1 \\ 1 & 1 & 1 & 0 & -1 & 1 & -1 & 1 & -1 \\ 1 & -1 & 1 & -1 & 0 & -1 & 1 & -1 & -1 \\ 1 & 1 & 1 & 1 & -1 & 0 & 1 & 1 & -1 \\ 1 & -1 & 1 & -1 & 1 & 1 & 0 & -1 & 1 \\ 1 & -1 & -1 & 1 & -1 & 1 & -1 & 0 & 1 \\ 1 & -1 & 1 & -1 & -1 & -1 & 1 & 1 & 0 \end{pmatrix} \qquad R = \begin{pmatrix} 0 & 0 & 0 & 0 & 0 & 0 & 0 & 0 & 0 \\ 0 & 0 & 0 & 0 & 0 & 0 & 0 & 0 & 0 \\ 8 & 8 & 0 & 4 & 4 & 4 & 8 & 8 & 4 \\ 3 & 3 & 6 & 0 & 9 & 3 & 3 & 9 & 6 \\ 0 & 0 & 0 & 0 & 0 & 0 & 0 & 0 & 0 \\ 8 & 7 & 6 & 8 & 9 & 0 & 6 & 7 & 9 \\ 3 & 12 & 3 & 3 & 4 & 4 & 0 & 3 & 4 \\ 6 & 6 & 1 & 1 & 1 & 1 & 1 & 0 & 1 \\ 6 & 1 & 1 & 1 & 6 & 1 & 1 & 1 & 0 \end{pmatrix}$$

$|P| = -4.0800696155 \times 10^{-15}$

$I3(P, R) = -19.6963902439$

$I3(-P, R) = 65.5061463415$

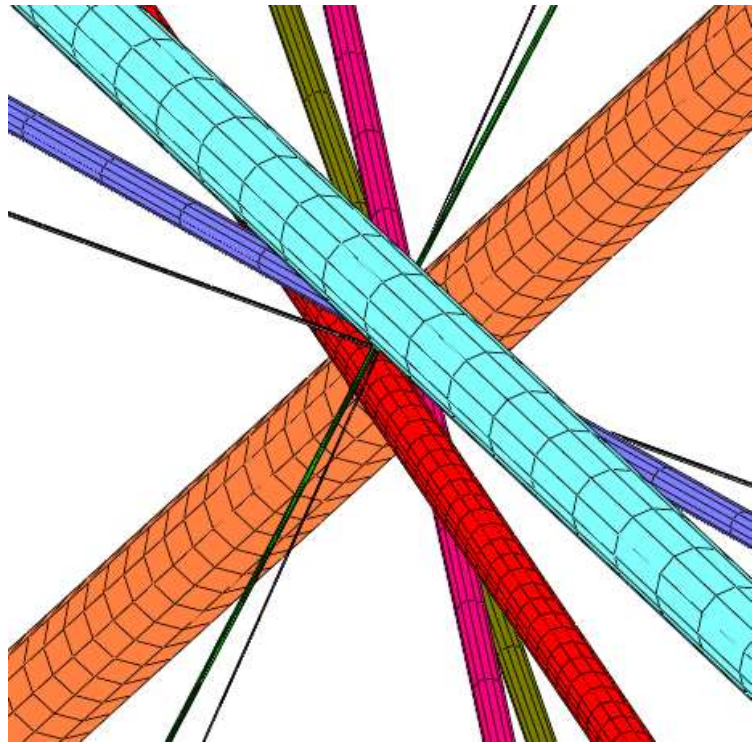

Appendix 2
9-cross
Configuration label a4j

$$\begin{pmatrix} \text{orange} \\ \text{red} \\ \text{blue} \\ \text{green} \\ \text{cyan} \\ \text{magenta} \\ \text{gray} \\ \text{olive} \\ \text{pink} \end{pmatrix} \begin{pmatrix} 0 & 0 & 0 & 1 \\ t_1 & p_1 & z_1 & r_1 \\ t_2 & p_2 & z_2 & r_2 \\ t_3 & p_3 & z_3 & r_3 \\ t_4 & p_4 & z_4 & r_4 \\ t_5 & p_5 & z_5 & r_5 \\ t_6 & p_6 & z_6 & r_6 \\ t_7 & p_7 & z_7 & r_7 \\ t_8 & p_8 & z_8 & r_8 \end{pmatrix} = \begin{pmatrix} 0 & 0 & 0 & 1 \\ 1.6189785457 & 3.1415926536 & 0 & 0.0284748449 \\ 1.8987257976 & 9.4163562807 & 1.9916307431 & 0.9330070813 \\ 1.0008040669 & 6.8905335063 & -1.5405549848 & 0.551500874 \\ 0.890925264 & 8.3897776378 & -1.9300903512 & 1.6684866678 \\ 1.986084849 & 4.9128169609 & -0.9147276674 & 0.2190708722 \\ 2.2368181696 & -2.336920652 & -0.4530903535 & 0.035986277 \\ 1.9330836533 & 3.5278662821 & 0.4420586225 & 0.2394180799 \\ 2.0978544521 & 3.454748937 & -0.2495902219 & 0.1182289542 \end{pmatrix}$$

$$P = \begin{pmatrix} 0 & 1 & 1 & 1 & 1 & 1 & 1 & 1 & 1 \\ 1 & 0 & 1 & 1 & -1 & 1 & 1 & -1 & 1 \\ 1 & 1 & 0 & -1 & -1 & 1 & 1 & 1 & 1 \\ 1 & 1 & -1 & 0 & 1 & 1 & -1 & 1 & -1 \\ 1 & -1 & -1 & 1 & 0 & 1 & -1 & -1 & -1 \\ 1 & 1 & 1 & 1 & 1 & 0 & -1 & 1 & -1 \\ 1 & 1 & 1 & -1 & -1 & -1 & 0 & 1 & -1 \\ 1 & -1 & 1 & 1 & -1 & 1 & 1 & 0 & -1 \\ 1 & 1 & 1 & -1 & -1 & -1 & -1 & -1 & 0 \end{pmatrix} \qquad R = \begin{pmatrix} 0 & 2 & 2 & 10 & 2 & 5 & 5 & 2 & 2 \\ 8 & 0 & 8 & 9 & 6 & 7 & 6 & 9 & 7 \\ 0 & 0 & 0 & 0 & 0 & 0 & 0 & 0 & 0 \\ 0 & 0 & 0 & 0 & 0 & 0 & 0 & 0 & 0 \\ 0 & 0 & 0 & 0 & 0 & 0 & 0 & 0 & 0 \\ 0 & 0 & 0 & 0 & 0 & 0 & 0 & 0 & 0 \\ 3 & 3 & 6 & 3 & 9 & 9 & 0 & 3 & 6 \\ 4 & 3 & 12 & 4 & 3 & 3 & 3 & 0 & 4 \\ 5 & 5 & 4 & 8 & 11 & 4 & 4 & 7 & 0 \end{pmatrix}$$

$|P| = 4$

$I3(P,R) = 34.1105794102$

$I3(-P,R) = 21.9848680807$

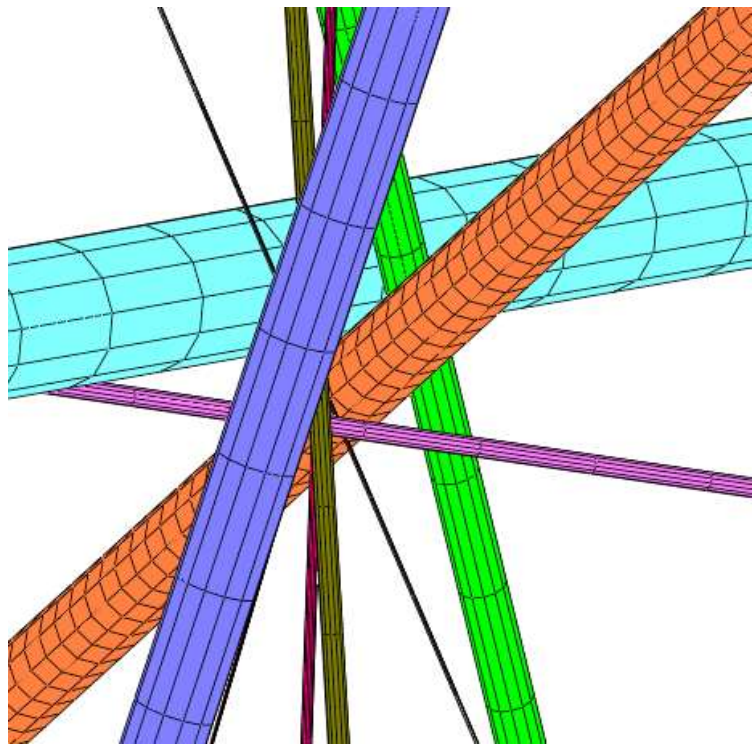

Appendix 2
9-cross
Configuration label a0e

$$\begin{pmatrix} \text{orange} \\ \text{red} \\ \text{blue} \\ \text{green} \\ \text{cyan} \\ \text{magenta} \\ \text{gray} \\ \text{olive} \\ \text{pink} \end{pmatrix} \begin{pmatrix} 0 & 0 & 0 & 1 \\ t_1 & p_1 & z_1 & r_1 \\ t_2 & p_2 & z_2 & r_2 \\ t_3 & p_3 & z_3 & r_3 \\ t_4 & p_4 & z_4 & r_4 \\ t_5 & p_5 & z_5 & r_5 \\ t_6 & p_6 & z_6 & r_6 \\ t_7 & p_7 & z_7 & r_7 \\ t_8 & p_8 & z_8 & r_8 \end{pmatrix} = \begin{pmatrix} 0 & 0 & 0 & 1 \\ 1.2216350113 & 3.1415926536 & 0 & 0.0494811013 \\ 1.8827858834 & 9.4028039268 & 3.1104069948 & 1.0157161473 \\ 0.8993685414 & 7.5501888487 & -0.7522687077 & 1.0098115716 \\ 0.9178786422 & 8.7477954774 & -0.9813708975 & 0.5107153102 \\ 1.9855914701 & 5.0514936152 & -0.5630849428 & 0.432615911 \\ 2.3265090141 & -1.5387116626 & -0.8120487686 & 0.0915751043 \\ 1.8173595491 & 3.9915916461 & 1.0739890829 & 0.5082164055 \\ 2.1717045399 & 3.0133667691 & 0.8614571593 & 0.3069422699 \end{pmatrix}$$

$$P = \begin{pmatrix} 0 & 1 & 1 & 1 & 1 & 1 & 1 & 1 & 1 \\ 1 & 0 & 1 & 1 & -1 & 1 & 1 & -1 & 1 \\ 1 & 1 & 0 & -1 & -1 & 1 & 1 & 1 & -1 \\ 1 & 1 & -1 & 0 & 1 & 1 & -1 & 1 & -1 \\ 1 & -1 & -1 & 1 & 0 & 1 & -1 & -1 & -1 \\ 1 & 1 & 1 & 1 & 1 & 0 & -1 & 1 & -1 \\ 1 & 1 & 1 & -1 & -1 & -1 & 0 & 1 & -1 \\ 1 & -1 & 1 & 1 & -1 & 1 & 1 & 0 & -1 \\ 1 & 1 & -1 & -1 & -1 & -1 & -1 & -1 & 0 \end{pmatrix} \qquad R = \begin{pmatrix} 0 & 2 & 2 & 10 & 2 & 5 & 5 & 2 & 2 \\ 8 & 0 & 7 & 8 & 7 & 8 & 7 & 8 & 7 \\ 0 & 0 & 0 & 0 & 0 & 0 & 0 & 0 & 0 \\ 0 & 0 & 0 & 0 & 0 & 0 & 0 & 0 & 0 \\ 0 & 0 & 0 & 0 & 0 & 0 & 0 & 0 & 0 \\ 0 & 0 & 0 & 0 & 0 & 0 & 0 & 0 & 0 \\ 3 & 3 & 6 & 3 & 9 & 9 & 0 & 3 & 6 \\ 4 & 3 & 12 & 4 & 3 & 3 & 3 & 0 & 4 \\ 6 & 6 & 6 & 10 & 10 & 5 & 5 & 6 & 0 \end{pmatrix}$$

$|P| = 1.082467449 \times 10^{-14}$

$I3(P, R) = 40.2095514799$

$I3(-P, R) = 15.9676775014$

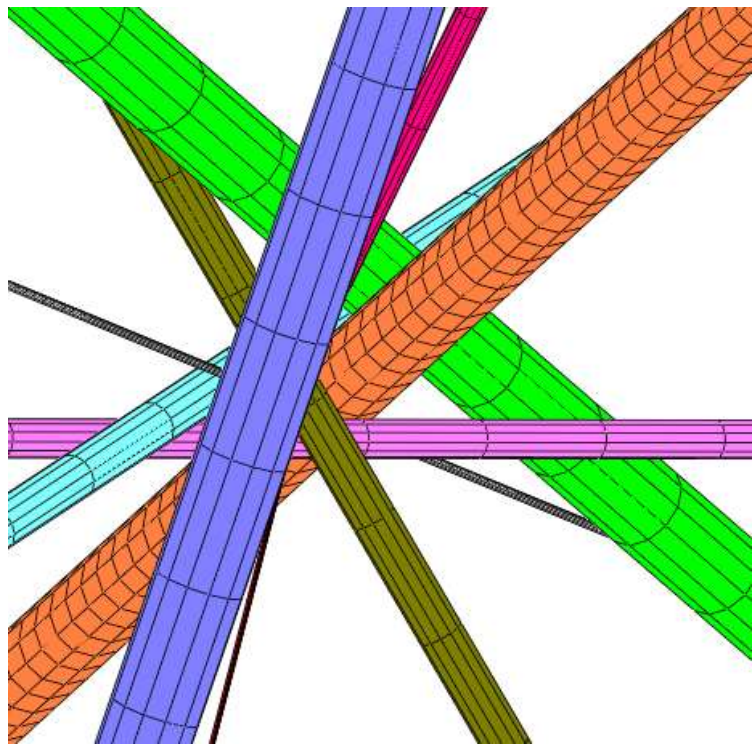

Appendix 2
9-cross
Configuration label a4h

$$\begin{pmatrix} \text{orange} \\ \text{red} \\ \text{blue} \\ \text{green} \\ \text{cyan} \\ \text{magenta} \\ \text{gray} \\ \text{olive} \\ \text{pink} \end{pmatrix} \begin{pmatrix} 0 & 0 & 0 & 1 \\ t_1 & p_1 & z_1 & r_1 \\ t_2 & p_2 & z_2 & r_2 \\ t_3 & p_3 & z_3 & r_3 \\ t_4 & p_4 & z_4 & r_4 \\ t_5 & p_5 & z_5 & r_5 \\ t_6 & p_6 & z_6 & r_6 \\ t_7 & p_7 & z_7 & r_7 \\ t_8 & p_8 & z_8 & r_8 \end{pmatrix} = \begin{pmatrix} 0 & 0 & 0 & 1 \\ 1.1613746904 & 3.1415926536 & 0 & 0.1541911841 \\ 2.0068122055 & 2.4299464456 & 0.5054587138 & 0.1261791964 \\ 1.345772384 & 0.7431013078 & -0.7054617068 & 1.9472431103 \\ 2.0377895112 & 2.0109528049 & 0.4481952418 & 0.1324239577 \\ 1.7151309138 & 2.4137649889 & -4.5452243299 & 2.0549051762 \\ -2.1236726842 & 5.4101747992 & 1.249799213 & 1.8199608034 \\ 2.1472665252 & 2.4177497498 & -0.0309324097 & 0.0177909955 \\ 2.3766917975 & 0.8198428215 & 1.9338716184 & 0.0787348425 \end{pmatrix}$$

$$P = \begin{pmatrix} 0 & 1 & 1 & 1 & 1 & 1 & -1 & 1 & 1 \\ 1 & 0 & 1 & 1 & 1 & -1 & -1 & -1 & 1 \\ 1 & 1 & 0 & -1 & -1 & -1 & -1 & -1 & -1 \\ 1 & 1 & -1 & 0 & -1 & 1 & -1 & -1 & -1 \\ 1 & 1 & -1 & -1 & 0 & 1 & -1 & 1 & 1 \\ 1 & -1 & -1 & 1 & 1 & 0 & -1 & -1 & 1 \\ -1 & -1 & -1 & -1 & -1 & -1 & 0 & 1 & 1 \\ 1 & -1 & -1 & -1 & 1 & -1 & 1 & 0 & 1 \\ 1 & 1 & -1 & -1 & 1 & 1 & 1 & 1 & 0 \end{pmatrix} \qquad R = \begin{pmatrix} 0 & 3 & 3 & 4 & 4 & 3 & 12 & 3 & 4 \\ 4 & 0 & 8 & 4 & 5 & 11 & 7 & 4 & 5 \\ 0 & 0 & 0 & 0 & 0 & 0 & 0 & 0 & 0 \\ 0 & 0 & 0 & 0 & 0 & 0 & 0 & 0 & 0 \\ 6 & 5 & 10 & 6 & 0 & 5 & 6 & 6 & 10 \\ 0 & 0 & 0 & 0 & 0 & 0 & 0 & 0 & 0 \\ 0 & 0 & 0 & 0 & 0 & 0 & 0 & 0 & 0 \\ 8 & 4 & 4 & 11 & 5 & 7 & 5 & 0 & 4 \\ 2 & 2 & 10 & 5 & 2 & 2 & 2 & 5 & 0 \end{pmatrix}$$

$|P| = 4$

$I3(P,R) = 29.8644473579$

$I3(-P,R) = 26.9479564499$

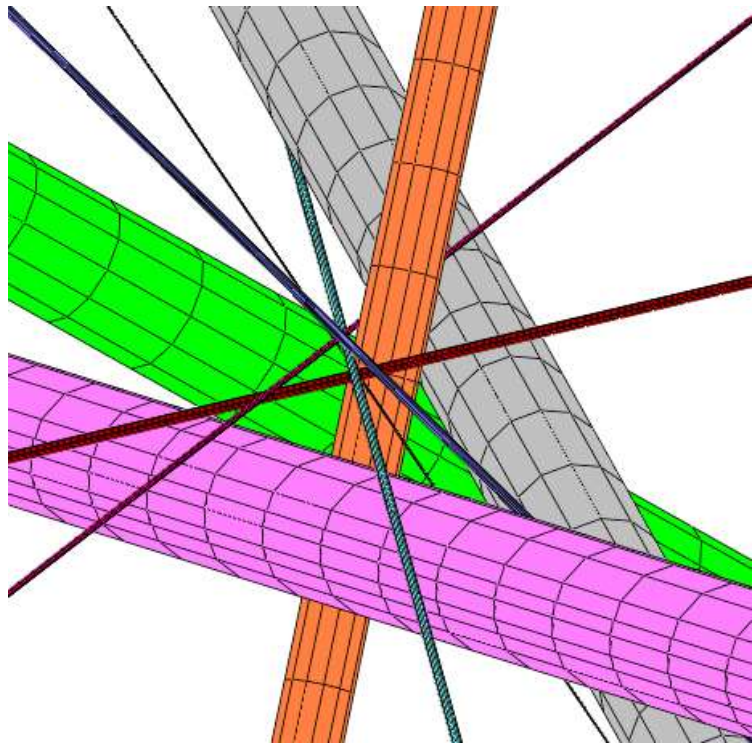

Appendix 2
9-cross
Configuration label am80c

$$\begin{pmatrix} \text{orange} \\ \text{red} \\ \text{blue} \\ \text{green} \\ \text{cyan} \\ \text{magenta} \\ \text{gray} \\ \text{olive} \\ \text{pink} \end{pmatrix} \begin{pmatrix} 0 & 0 & 0 & 1 \\ t1 & p1 & z1 & r1 \\ t2 & p2 & z2 & r2 \\ t3 & p3 & z3 & r3 \\ t4 & p4 & z4 & r4 \\ t5 & p5 & z5 & r5 \\ t6 & p6 & z6 & r6 \\ t7 & p7 & z7 & r7 \\ t8 & p8 & z8 & r8 \end{pmatrix} = \begin{pmatrix} 0 & 0 & 0 & 1 \\ 1.7351306481 & 3.1415926536 & 0 & 0.2502162721 \\ 1.4221116719 & 3.4745244228 & 0.5819273119 & 0.0100000003 \\ 2.321756933 & 6.1240167852 & -13.3368127606 & 3.1664268278 \\ 1.5159695906 & 3.5610345529 & 0.4952474563 & 0.1800093032 \\ 1.5652724673 & 3.4747367374 & -13.0136265056 & 7.7672086413 \\ 4.8807520806 & 0.3589629579 & -1.45422382 & 7.13141202 \\ 1.3586558994 & 3.4729167157 & -0.2351605208 & 0.055265482 \\ 1.2838667417 & 6.1006138557 & 2.1904770191 & 0.1472747638 \end{pmatrix}$$

$$P = \begin{pmatrix} 0 & 1 & 1 & 1 & 1 & 1 & -1 & 1 & 1 \\ 1 & 0 & -1 & -1 & -1 & 1 & 1 & 1 & -1 \\ 1 & -1 & 0 & 1 & 1 & 1 & 1 & -1 & 1 \\ 1 & -1 & 1 & 0 & 1 & -1 & 1 & 1 & 1 \\ 1 & -1 & 1 & 1 & 0 & -1 & 1 & -1 & -1 \\ 1 & 1 & 1 & -1 & -1 & 0 & 1 & 1 & -1 \\ -1 & 1 & 1 & 1 & 1 & 1 & 0 & -1 & -1 \\ 1 & 1 & -1 & 1 & -1 & 1 & -1 & 0 & -1 \\ 1 & -1 & 1 & 1 & -1 & -1 & -1 & -1 & 0 \end{pmatrix} \qquad R = \begin{pmatrix} 0 & 3 & 3 & 4 & 4 & 3 & 12 & 3 & 4 \\ 4 & 0 & 8 & 4 & 5 & 11 & 7 & 4 & 5 \\ 0 & 0 & 0 & 0 & 0 & 0 & 0 & 0 & 0 \\ 0 & 0 & 0 & 0 & 0 & 0 & 0 & 0 & 0 \\ 6 & 5 & 10 & 6 & 0 & 5 & 6 & 6 & 10 \\ 0 & 0 & 0 & 0 & 0 & 0 & 0 & 0 & 0 \\ 0 & 0 & 0 & 0 & 0 & 0 & 0 & 0 & 0 \\ 7 & 5 & 4 & 11 & 4 & 8 & 4 & 0 & 5 \\ 2 & 2 & 10 & 5 & 2 & 2 & 2 & 5 & 0 \end{pmatrix}$$

$|P| = -80$

$I3(P, R) = 23.9197560976$

$I3(-P, R) = 32.651097561$

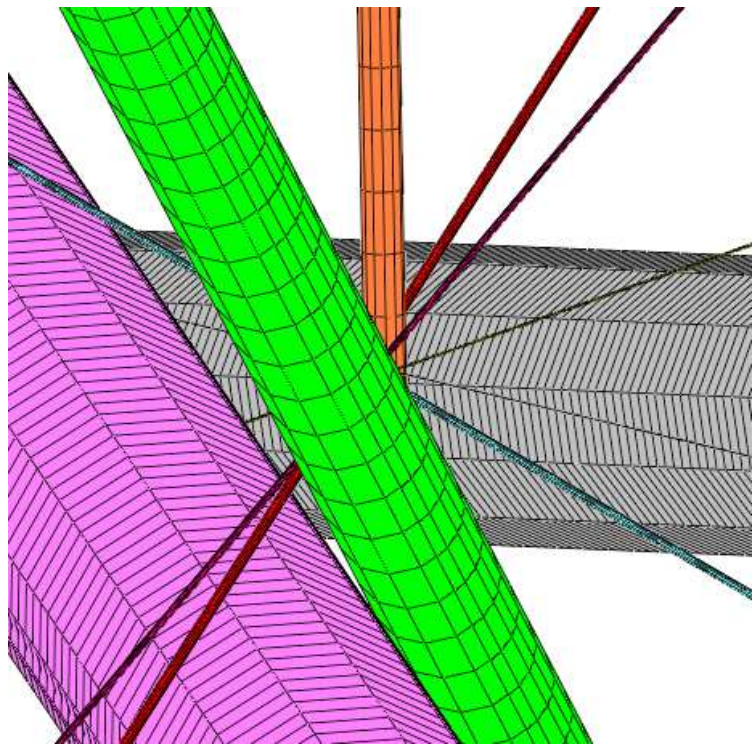

Appendix 2
9-cross
Configuration label a4g

$$\begin{pmatrix} \text{orange} \\ \text{red} \\ \text{blue} \\ \text{green} \\ \text{cyan} \\ \text{magenta} \\ \text{gray} \\ \text{olive} \\ \text{pink} \end{pmatrix} \begin{pmatrix} 0 & 0 & 0 & 1 \\ t_1 & p_1 & z_1 & r_1 \\ t_2 & p_2 & z_2 & r_2 \\ t_3 & p_3 & z_3 & r_3 \\ t_4 & p_4 & z_4 & r_4 \\ t_5 & p_5 & z_5 & r_5 \\ t_6 & p_6 & z_6 & r_6 \\ t_7 & p_7 & z_7 & r_7 \\ t_8 & p_8 & z_8 & r_8 \end{pmatrix} = \begin{pmatrix} 0 & 0 & 0 & 1 \\ 2.7011866007 & 3.1415926536 & 0 & 0.9431928298 \\ 0.3746072364 & -3.002666633 & 63.7214304416 & 3.0737840201 \\ 2.952238415 & -0.1405453591 & -103.0406267752 & 4.3770550371 \\ 0.1697330374 & -0.8437273541 & 66.4477951877 & 5.3414577362 \\ 3.0939857773 & -1.0515125919 & -116.3200364676 & 0.7133861389 \\ 0.2295694918 & 2.4794290356 & -39.5364947458 & 19.8376582975 \\ 0.0707623408 & 8.8846694326 & -62.2870255057 & 0.1336254991 \\ 2.992247585 & 4.5501503515 & 4.5621191927 & 2.0233801993 \end{pmatrix}$$

$$P = \begin{pmatrix} 0 & 1 & 1 & 1 & 1 & 1 & 1 & 1 & 1 \\ 1 & 0 & -1 & -1 & -1 & 1 & -1 & -1 & -1 \\ 1 & -1 & 0 & 1 & 1 & 1 & -1 & -1 & 1 \\ 1 & -1 & 1 & 0 & 1 & -1 & -1 & 1 & 1 \\ 1 & -1 & 1 & 1 & 0 & -1 & 1 & 1 & -1 \\ 1 & 1 & 1 & -1 & -1 & 0 & -1 & 1 & 1 \\ 1 & -1 & -1 & -1 & 1 & -1 & 0 & 1 & -1 \\ 1 & -1 & -1 & 1 & 1 & 1 & 1 & 0 & -1 \\ 1 & -1 & 1 & 1 & -1 & 1 & -1 & -1 & 0 \end{pmatrix} \qquad R = \begin{pmatrix} 0 & 2 & 2 & 10 & 2 & 2 & 5 & 5 & 2 \\ 9 & 0 & 6 & 9 & 6 & 6 & 6 & 6 & 6 \\ 0 & 0 & 0 & 0 & 0 & 0 & 0 & 0 & 0 \\ 0 & 0 & 0 & 0 & 0 & 0 & 0 & 0 & 0 \\ 1 & 1 & 6 & 6 & 0 & 1 & 1 & 1 & 1 \\ 1 & 1 & 1 & 1 & 6 & 0 & 6 & 1 & 1 \\ 0 & 0 & 0 & 0 & 0 & 0 & 0 & 0 & 0 \\ 4 & 4 & 5 & 4 & 7 & 11 & 8 & 0 & 5 \\ 8 & 5 & 11 & 4 & 7 & 4 & 5 & 4 & 0 \end{pmatrix}$$

$|P| = 4$

$I3(P, R) = -6.2395821696$

$I3(-P, R) = 45.6298178974$

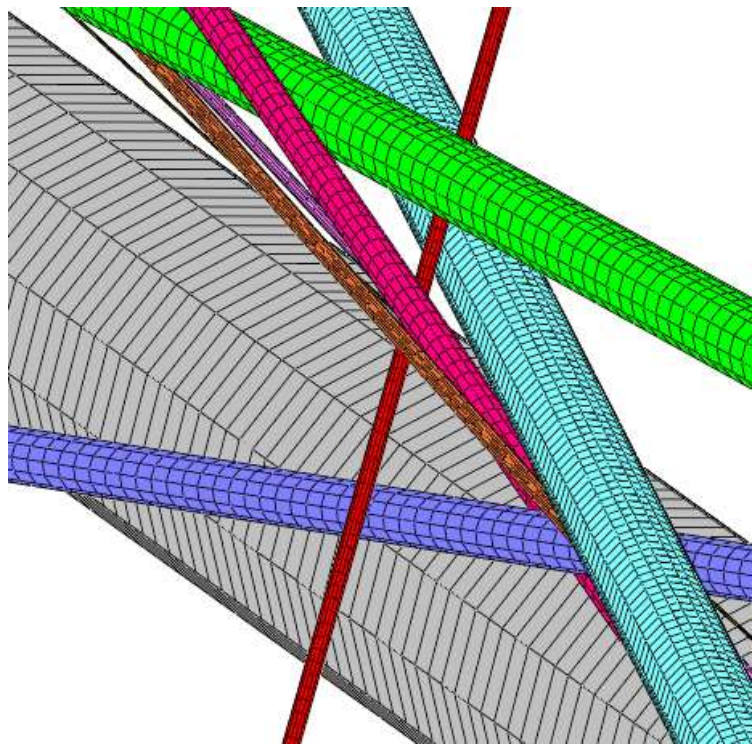

Appendix 2
9-cross
Configuration label a16c

$$\begin{pmatrix} \text{orange} \\ \text{red} \\ \text{blue} \\ \text{green} \\ \text{cyan} \\ \text{magenta} \\ \text{gray} \\ \text{olive} \\ \text{pink} \end{pmatrix} \begin{pmatrix} 0 & 0 & 0 & 1 \\ t_1 & p_1 & z_1 & r_1 \\ t_2 & p_2 & z_2 & r_2 \\ t_3 & p_3 & z_3 & r_3 \\ t_4 & p_4 & z_4 & r_4 \\ t_5 & p_5 & z_5 & r_5 \\ t_6 & p_6 & z_6 & r_6 \\ t_7 & p_7 & z_7 & r_7 \\ t_8 & p_8 & z_8 & r_8 \end{pmatrix} = \begin{pmatrix} 0 & 0 & 0 & 1 \\ 0.4549880834 & 3.1415926536 & 0 & 0.1872004942 \\ 0.2488415013 & -3.0319086824 & 7.0946189827 & 0.478239907 \\ 3.0158204404 & -0.5437851537 & -19.4026004488 & 0.5328271151 \\ 0.4254734395 & -2.1071078624 & 7.236906291 & 1.096092478 \\ 2.9905471528 & -2.1734541452 & -6.0778082599 & 0.0507942918 \\ 0.1208584003 & 2.0354765338 & -10.8364712401 & 0.2056825376 \\ 0.0997358752 & 8.5134521892 & -10.2507441737 & 0.0273328981 \\ 2.8806805443 & 4.1611992352 & -0.7607081587 & 0.0100000003 \end{pmatrix}$$

$$P = \begin{pmatrix} 0 & 1 & 1 & 1 & 1 & 1 & 1 & 1 & 1 \\ 1 & 0 & -1 & -1 & -1 & 1 & -1 & -1 & -1 \\ 1 & -1 & 0 & 1 & 1 & 1 & -1 & -1 & 1 \\ 1 & -1 & 1 & 0 & 1 & -1 & -1 & 1 & 1 \\ 1 & -1 & 1 & 1 & 0 & -1 & 1 & 1 & -1 \\ 1 & 1 & 1 & -1 & -1 & 0 & -1 & 1 & -1 \\ 1 & -1 & -1 & -1 & 1 & -1 & 0 & 1 & -1 \\ 1 & -1 & -1 & 1 & 1 & 1 & 1 & 0 & -1 \\ 1 & -1 & 1 & 1 & -1 & -1 & -1 & -1 & 0 \end{pmatrix} \qquad R = \begin{pmatrix} 0 & 2 & 2 & 10 & 2 & 2 & 5 & 5 & 2 \\ 9 & 0 & 6 & 9 & 6 & 6 & 6 & 6 & 6 \\ 0 & 0 & 0 & 0 & 0 & 0 & 0 & 0 & 0 \\ 0 & 0 & 0 & 0 & 0 & 0 & 0 & 0 & 0 \\ 1 & 1 & 6 & 6 & 0 & 1 & 1 & 1 & 1 \\ 2 & 2 & 2 & 2 & 5 & 0 & 10 & 2 & 5 \\ 0 & 0 & 0 & 0 & 0 & 0 & 0 & 0 & 0 \\ 4 & 4 & 5 & 4 & 7 & 11 & 8 & 0 & 5 \\ 7 & 4 & 11 & 5 & 8 & 4 & 4 & 5 & 0 \end{pmatrix}$$

$|P| = 16$

$I3(P, R) = 17.5732072194$

$I3(-P, R) = 25.832159058$

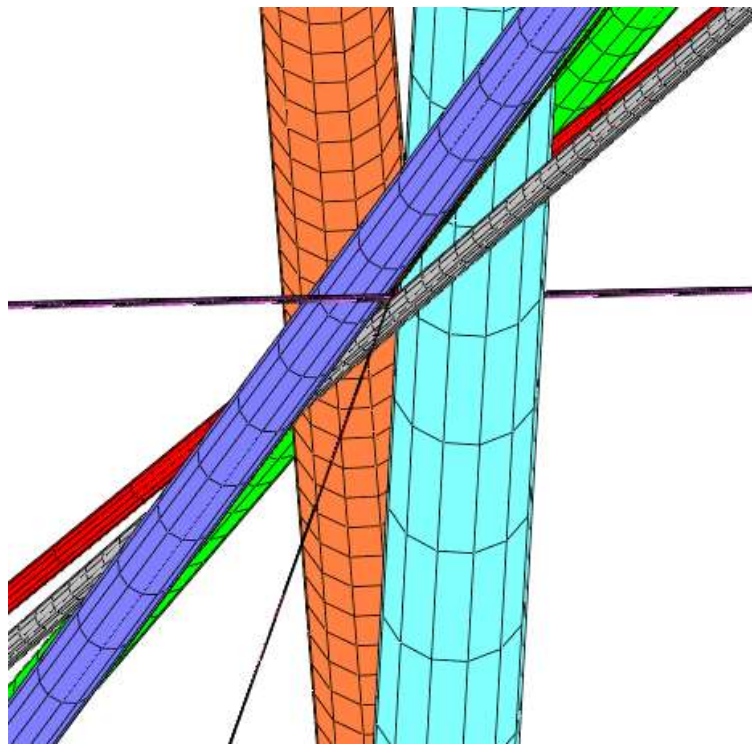

Appendix 2
9-cross
Configuration label am84a

$$\begin{pmatrix} \text{orange} \\ \text{red} \\ \text{blue} \\ \text{green} \\ \text{cyan} \\ \text{magenta} \\ \text{gray} \\ \text{olive} \\ \text{pink} \end{pmatrix} \begin{pmatrix} 0 & 0 & 0 & 1 \\ t1 & p1 & z1 & r1 \\ t2 & p2 & z2 & r2 \\ t3 & p3 & z3 & r3 \\ t4 & p4 & z4 & r4 \\ t5 & p5 & z5 & r5 \\ t6 & p6 & z6 & r6 \\ t7 & p7 & z7 & r7 \\ t8 & p8 & z8 & r8 \end{pmatrix} = \begin{pmatrix} 0 & 0 & 0 & 1 \\ 1.0354026338 & 3.1415926536 & 0 & 0.2862693429 \\ 2.4888353334 & 2.6701354048 & 2.8703977215 & 1.0039734989 \\ 1.5531584976 & 1.1073194686 & 1.2249102553 & 3.6749966949 \\ 0.1454554833 & 1.1149639368 & -16.5644796297 & 0.01 \\ 2.2884486933 & 3.5534936595 & 0.0935640111 & 0.0137031255 \\ 0.3151900167 & 1.1416512709 & -5.579992459 & 0.0384922992 \\ 2.8837157252 & 3.4944832539 & -7.8395625615 & 0.0815984379 \\ 2.7688665815 & 3.562604393 & -4.998310941 & 0.416740625 \end{pmatrix}$$

$$P = \begin{pmatrix} 0 & 1 & 1 & 1 & 1 & 1 & 1 & 1 & 1 \\ 1 & 0 & 1 & 1 & -1 & -1 & 1 & 1 & 1 \\ 1 & 1 & 0 & -1 & -1 & 1 & -1 & 1 & 1 \\ 1 & 1 & -1 & 0 & 1 & -1 & 1 & -1 & -1 \\ 1 & -1 & -1 & 1 & 0 & -1 & -1 & -1 & 1 \\ 1 & -1 & 1 & -1 & -1 & 0 & 1 & -1 & 1 \\ 1 & 1 & -1 & 1 & -1 & 1 & 0 & 1 & 1 \\ 1 & 1 & 1 & -1 & -1 & -1 & 1 & 0 & -1 \\ 1 & 1 & 1 & -1 & 1 & 1 & 1 & -1 & 0 \end{pmatrix} \qquad R = \begin{pmatrix} 0 & 0 & 0 & 0 & 0 & 0 & 0 & 0 & 0 \\ 12 & 0 & 3 & 3 & 4 & 3 & 3 & 4 & 4 \\ 0 & 0 & 0 & 0 & 0 & 0 & 0 & 0 & 0 \\ 0 & 0 & 0 & 0 & 0 & 0 & 0 & 0 & 0 \\ 2 & 2 & 2 & 5 & 0 & 5 & 2 & 2 & 10 \\ 5 & 6 & 10 & 5 & 7 & 0 & 9 & 5 & 7 \\ 7 & 4 & 11 & 8 & 5 & 4 & 0 & 5 & 4 \\ 6 & 6 & 5 & 10 & 10 & 6 & 5 & 0 & 6 \\ 0 & 0 & 0 & 0 & 0 & 0 & 0 & 0 & 0 \end{pmatrix}$$

$|P| = -84$

$I3(P,R) = 30.0826834633$

$I3(-P,R) = 27.4861027794$

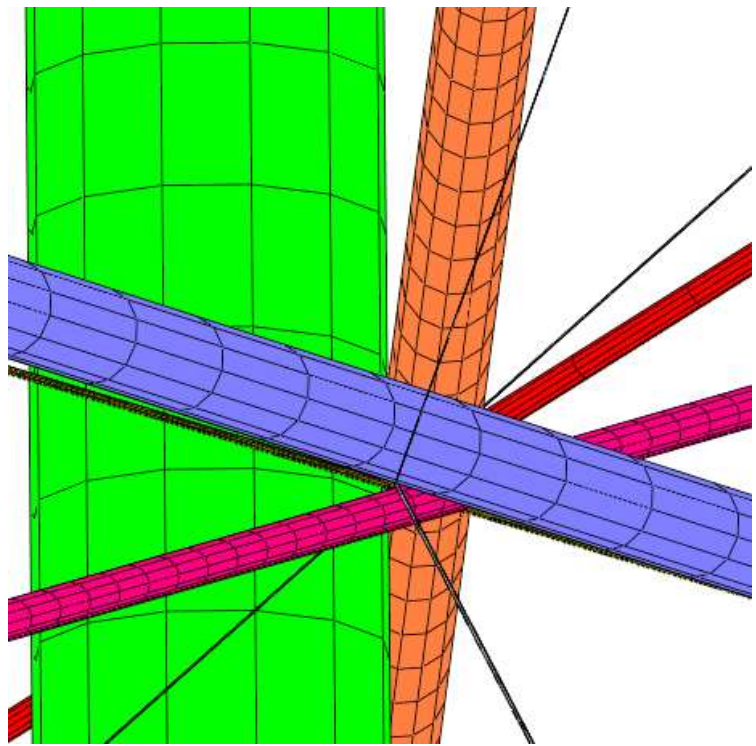

Appendix 2
9-cross
Configuration label am52b

$$\begin{pmatrix} \text{orange} \\ \text{red} \\ \text{blue} \\ \text{green} \\ \text{cyan} \\ \text{magenta} \\ \text{gray} \\ \text{olive} \\ \text{pink} \end{pmatrix} \begin{pmatrix} 0 & 0 & 0 & 1 \\ t_1 & p_1 & z_1 & r_1 \\ t_2 & p_2 & z_2 & r_2 \\ t_3 & p_3 & z_3 & r_3 \\ t_4 & p_4 & z_4 & r_4 \\ t_5 & p_5 & z_5 & r_5 \\ t_6 & p_6 & z_6 & r_6 \\ t_7 & p_7 & z_7 & r_7 \\ t_8 & p_8 & z_8 & r_8 \end{pmatrix} = \begin{pmatrix} 0 & 0 & 0 & 1 \\ 0.3721316833 & 3.1415926536 & 0 & 0.0160364801 \\ 1.0534193296 & 4.3151608465 & 2.0242474847 & 0.0225177048 \\ 0.0707211132 & 2.9242823815 & -4.2431513008 & 0.3677476218 \\ -2.507140559 & 8.2555404112 & -3.0798038877 & 0.0871378284 \\ 0.3210012399 & 3.0894073368 & -0.4360887622 & 0.0135682913 \\ 0.3414550799 & 3.1569554893 & -0.3485614095 & 0.0116944258 \\ 2.8099830966 & 4.3505661598 & -2.5317646169 & 1.9961939809 \\ 0.3937567996 & 3.1982623248 & 0.023537858 & 0.0128594411 \end{pmatrix}$$

$$P = \begin{pmatrix} 0 & 1 & 1 & 1 & -1 & 1 & 1 & 1 & 1 \\ 1 & 0 & 1 & -1 & 1 & -1 & 1 & 1 & 1 \\ 1 & 1 & 0 & 1 & 1 & 1 & -1 & 1 & -1 \\ 1 & -1 & 1 & 0 & -1 & -1 & -1 & 1 & -1 \\ -1 & 1 & 1 & -1 & 0 & -1 & -1 & 1 & 1 \\ 1 & -1 & 1 & -1 & -1 & 0 & 1 & 1 & -1 \\ 1 & 1 & -1 & -1 & -1 & 1 & 0 & 1 & -1 \\ 1 & 1 & 1 & 1 & 1 & 1 & 1 & 0 & 1 \\ 1 & 1 & -1 & -1 & 1 & -1 & -1 & 1 & 0 \end{pmatrix} \quad R = \begin{pmatrix} 0 & 0 & 0 & 0 & 0 & 0 & 0 & 0 & 0 \\ 0 & 0 & 0 & 0 & 0 & 0 & 0 & 0 & 0 \\ 4 & 8 & 0 & 5 & 4 & 5 & 7 & 11 & 4 \\ 2 & 2 & 5 & 0 & 10 & 2 & 2 & 5 & 2 \\ 0 & 0 & 0 & 0 & 0 & 0 & 0 & 0 & 0 \\ 6 & 7 & 5 & 9 & 10 & 0 & 5 & 5 & 7 \\ 9 & 3 & 3 & 6 & 3 & 3 & 0 & 9 & 6 \\ 0 & 0 & 0 & 0 & 0 & 0 & 0 & 0 & 0 \\ 9 & 9 & 6 & 3 & 6 & 3 & 3 & 3 & 0 \end{pmatrix}$$

$|P| = -52$

$I3(P, R) = 32.4097486672$

$I3(-P, R) = 24.8189895913$

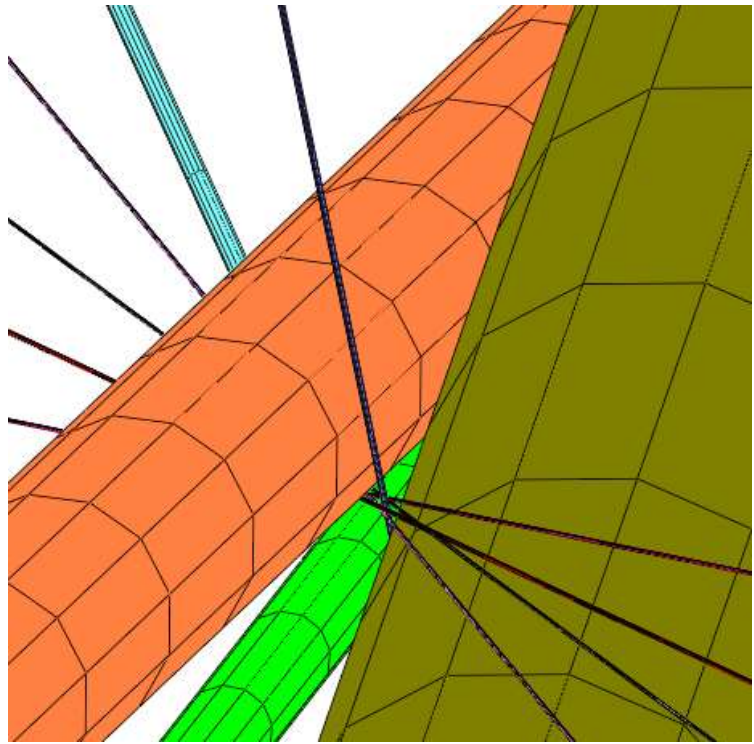

Appendix 2
9-cross
Configuration label am80b

$$\begin{pmatrix} \text{orange} \\ \text{red} \\ \text{blue} \\ \text{green} \\ \text{cyan} \\ \text{magenta} \\ \text{gray} \\ \text{olive} \\ \text{pink} \end{pmatrix} \begin{pmatrix} 0 & 0 & 0 & 1 \\ t_1 & p_1 & z_1 & r_1 \\ t_2 & p_2 & z_2 & r_2 \\ t_3 & p_3 & z_3 & r_3 \\ t_4 & p_4 & z_4 & r_4 \\ t_5 & p_5 & z_5 & r_5 \\ t_6 & p_6 & z_6 & r_6 \\ t_7 & p_7 & z_7 & r_7 \\ t_8 & p_8 & z_8 & r_8 \end{pmatrix} = \begin{pmatrix} 0 & 0 & 0 & 1 \\ 0.1515509911 & 3.1415926536 & 0 & 0.6258654451 \\ 2.3347020777 & 1.3469560216 & -3.0484200651 & 0.8046677748 \\ 2.0821640109 & 0.6028817738 & -16.5248443871 & 21.4232953336 \\ 1.1588602766 & 4.971489346 & 5.9361823071 & 1.6821351477 \\ 2.54285859 & 4.8995084721 & 95.3498362282 & 6.5480428404 \\ -2.3188905994 & 5.6594329156 & 21.2005537557 & 10.877881584 \\ 2.2551350516 & 1.7098038956 & 0.8016071944 & 3.3599245632 \\ 1.0094898969 & 4.6053889986 & 2.2078388278 & 0.4615618363 \end{pmatrix}$$

$$P = \begin{pmatrix} 0 & 1 & 1 & 1 & 1 & 1 & -1 & 1 & 1 \\ 1 & 0 & 1 & 1 & 1 & -1 & -1 & 1 & 1 \\ 1 & 1 & 0 & -1 & 1 & 1 & -1 & -1 & -1 \\ 1 & 1 & -1 & 0 & 1 & 1 & -1 & -1 & 1 \\ 1 & 1 & 1 & 1 & 0 & 1 & 1 & -1 & -1 \\ 1 & -1 & 1 & 1 & 1 & 0 & -1 & -1 & -1 \\ -1 & -1 & -1 & -1 & 1 & -1 & 0 & -1 & 1 \\ 1 & 1 & -1 & -1 & -1 & -1 & -1 & 0 & -1 \\ 1 & 1 & -1 & 1 & -1 & -1 & 1 & -1 & 0 \end{pmatrix} \quad R = \begin{pmatrix} 0 & 11 & 5 & 5 & 4 & 4 & 4 & 8 & 7 \\ 1 & 0 & 1 & 1 & 1 & 1 & 1 & 6 & 6 \\ 4 & 4 & 0 & 11 & 5 & 8 & 4 & 5 & 7 \\ 0 & 0 & 0 & 0 & 0 & 0 & 0 & 0 & 0 \\ 6 & 6 & 6 & 6 & 0 & 9 & 9 & 6 & 6 \\ 0 & 0 & 0 & 0 & 0 & 0 & 0 & 0 & 0 \\ 5 & 5 & 2 & 2 & 2 & 10 & 0 & 2 & 2 \\ 0 & 0 & 0 & 0 & 0 & 0 & 0 & 0 & 0 \\ 1 & 1 & 1 & 6 & 1 & 1 & 6 & 1 & 0 \end{pmatrix}$$

$|P| = -80$

$I3(P, R) = 24.6140502709$

$I3(-P, R) = 25.7315623119$

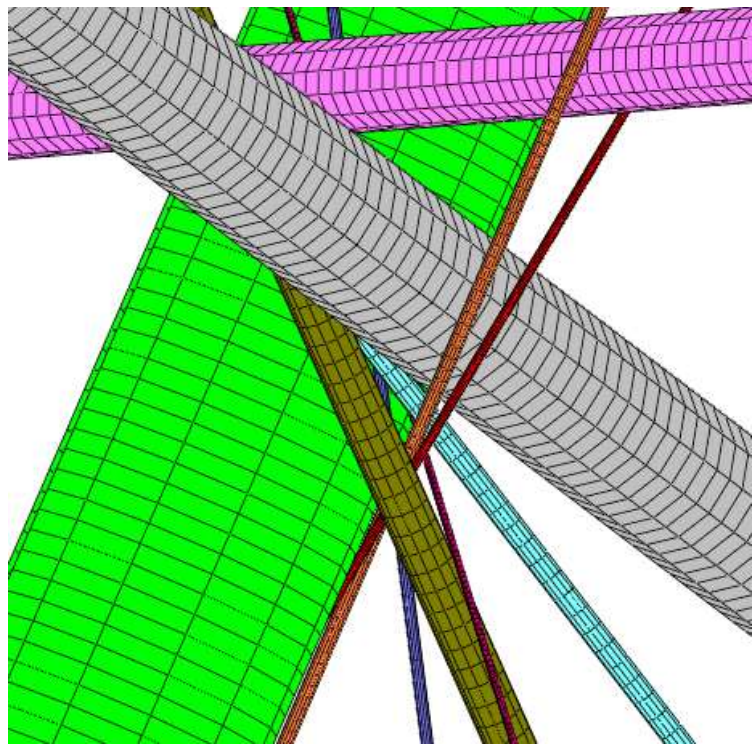

Appendix 2
9-cross
Configuration label am80a

$$\begin{pmatrix} \text{orange} \\ \text{red} \\ \text{blue} \\ \text{green} \\ \text{cyan} \\ \text{magenta} \\ \text{gray} \\ \text{olive} \\ \text{pink} \end{pmatrix} \begin{pmatrix} 0 & 0 & 0 & 1 \\ t_1 & p_1 & z_1 & r_1 \\ t_2 & p_2 & z_2 & r_2 \\ t_3 & p_3 & z_3 & r_3 \\ t_4 & p_4 & z_4 & r_4 \\ t_5 & p_5 & z_5 & r_5 \\ t_6 & p_6 & z_6 & r_6 \\ t_7 & p_7 & z_7 & r_7 \\ t_8 & p_8 & z_8 & r_8 \end{pmatrix} = \begin{pmatrix} 0 & 0 & 0 & 1 \\ 0.1116411462 & 3.1415926536 & 0 & 0.3717177653 \\ 2.6240996785 & 0.9049604018 & -12.6416451699 & 0.4400785253 \\ 2.5567598882 & 0.6870265463 & -18.5747439625 & 5.2728149752 \\ 0.8660129465 & 4.407131826 & 2.9045922335 & 1.4362482396 \\ 2.252258735 & 4.3322993921 & 66.8188169264 & 2.9712550596 \\ -2.5848106684 & 5.6917708188 & 32.7331455038 & 9.8526448918 \\ 2.6145671983 & 1.1236234164 & -8.382567729 & 2.2690360078 \\ 2.7223774326 & 1.116697513 & -15.6400637065 & 0.0475113874 \end{pmatrix}$$

$$P = \begin{pmatrix} 0 & 1 & 1 & 1 & 1 & 1 & -1 & 1 & 1 \\ 1 & 0 & 1 & 1 & 1 & -1 & -1 & 1 & -1 \\ 1 & 1 & 0 & -1 & 1 & 1 & -1 & -1 & 1 \\ 1 & 1 & -1 & 0 & 1 & 1 & -1 & -1 & -1 \\ 1 & 1 & 1 & 1 & 0 & 1 & 1 & -1 & -1 \\ 1 & -1 & 1 & 1 & 1 & 0 & -1 & -1 & -1 \\ -1 & -1 & -1 & -1 & 1 & -1 & 0 & -1 & -1 \\ 1 & 1 & -1 & -1 & -1 & -1 & -1 & 0 & -1 \\ 1 & -1 & 1 & -1 & -1 & -1 & -1 & -1 & 0 \end{pmatrix} \qquad R = \begin{pmatrix} 0 & 12 & 3 & 3 & 4 & 4 & 4 & 3 & 3 \\ 0 & 0 & 0 & 0 & 0 & 0 & 0 & 0 & 0 \\ 7 & 7 & 0 & 10 & 5 & 5 & 6 & 5 & 9 \\ 0 & 0 & 0 & 0 & 0 & 0 & 0 & 0 & 0 \\ 4 & 4 & 8 & 8 & 0 & 8 & 8 & 4 & 4 \\ 0 & 0 & 0 & 0 & 0 & 0 & 0 & 0 & 0 \\ 5 & 5 & 2 & 2 & 2 & 10 & 0 & 2 & 2 \\ 0 & 0 & 0 & 0 & 0 & 0 & 0 & 0 & 0 \\ 4 & 7 & 4 & 8 & 5 & 5 & 4 & 11 & 0 \end{pmatrix}$$

$|P| = -80$

$I3(P, R) = 28.2738683402$

$I3(-P, R) = 29.5840566164$

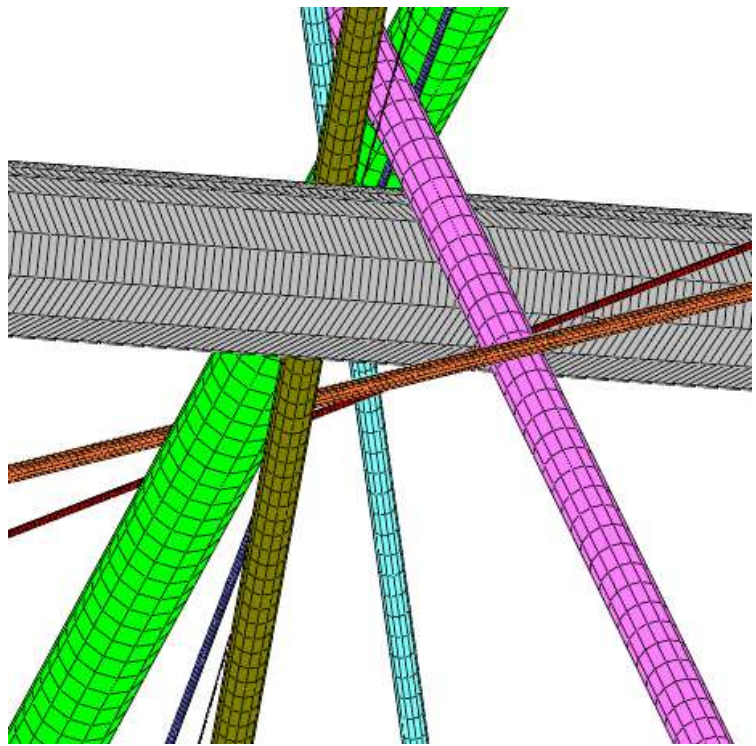

Appendix 2
9-cross
Configuration label a4f

$$\begin{pmatrix} \text{orange} \\ \text{red} \\ \text{blue} \\ \text{green} \\ \text{cyan} \\ \text{magenta} \\ \text{gray} \\ \text{olive} \\ \text{pink} \end{pmatrix} \begin{pmatrix} 0 & 0 & 0 & 1 \\ t_1 & p_1 & z_1 & r_1 \\ t_2 & p_2 & z_2 & r_2 \\ t_3 & p_3 & z_3 & r_3 \\ t_4 & p_4 & z_4 & r_4 \\ t_5 & p_5 & z_5 & r_5 \\ t_6 & p_6 & z_6 & r_6 \\ t_7 & p_7 & z_7 & r_7 \\ t_8 & p_8 & z_8 & r_8 \end{pmatrix} = \begin{pmatrix} 0 & 0 & 0 & 1 \\ 2.7040942541 & 3.1415926536 & 0 & 1.0162125668 \\ 3.0701273951 & 4.3086550849 & -9.7144982343 & 0.0399060706 \\ 3.1136349429 & -0.9797115344 & -59.151897546 & 0.6592097137 \\ 3.0506589457 & -1.9031047538 & -6.9600090151 & 0.1679623742 \\ 3.0815089674 & -2.0502025301 & -54.6730999731 & 0.6322296118 \\ 0.0627484769 & 1.3596736564 & -82.416584213 & 1.919125505 \\ 0.030125705 & 7.9677247019 & -102.2844075415 & 0.2202476032 \\ 2.9749445689 & 4.2273296132 & -16.7837578073 & 0.0100000001 \end{pmatrix}$$

$$P = \begin{pmatrix} 0 & 1 & 1 & 1 & 1 & 1 & 1 & 1 & 1 \\ 1 & 0 & -1 & -1 & -1 & 1 & -1 & -1 & 1 \\ 1 & -1 & 0 & 1 & -1 & -1 & -1 & -1 & -1 \\ 1 & -1 & 1 & 0 & 1 & -1 & -1 & 1 & 1 \\ 1 & -1 & -1 & 1 & 0 & -1 & 1 & 1 & -1 \\ 1 & 1 & -1 & -1 & -1 & 0 & -1 & 1 & 1 \\ 1 & -1 & -1 & -1 & 1 & -1 & 0 & 1 & 1 \\ 1 & -1 & -1 & 1 & 1 & 1 & 1 & 0 & 1 \\ 1 & 1 & -1 & 1 & -1 & 1 & 1 & 1 & 0 \end{pmatrix} \quad R = \begin{pmatrix} 0 & 2 & 2 & 10 & 2 & 2 & 5 & 5 & 2 \\ 6 & 0 & 9 & 6 & 9 & 6 & 6 & 6 & 6 \\ 6 & 1 & 0 & 1 & 6 & 1 & 1 & 1 & 1 \\ 0 & 0 & 0 & 0 & 0 & 0 & 0 & 0 & 0 \\ 0 & 0 & 0 & 0 & 0 & 0 & 0 & 0 & 0 \\ 2 & 2 & 5 & 2 & 5 & 0 & 10 & 2 & 2 \\ 0 & 0 & 0 & 0 & 0 & 0 & 0 & 0 & 0 \\ 4 & 4 & 7 & 4 & 5 & 8 & 11 & 0 & 5 \\ 7 & 5 & 4 & 11 & 5 & 8 & 4 & 4 & 0 \end{pmatrix}$$

$|P| = 4$

$I3(P, R) = 21.4014997657$

$I3(-P, R) = 22.9702129876$

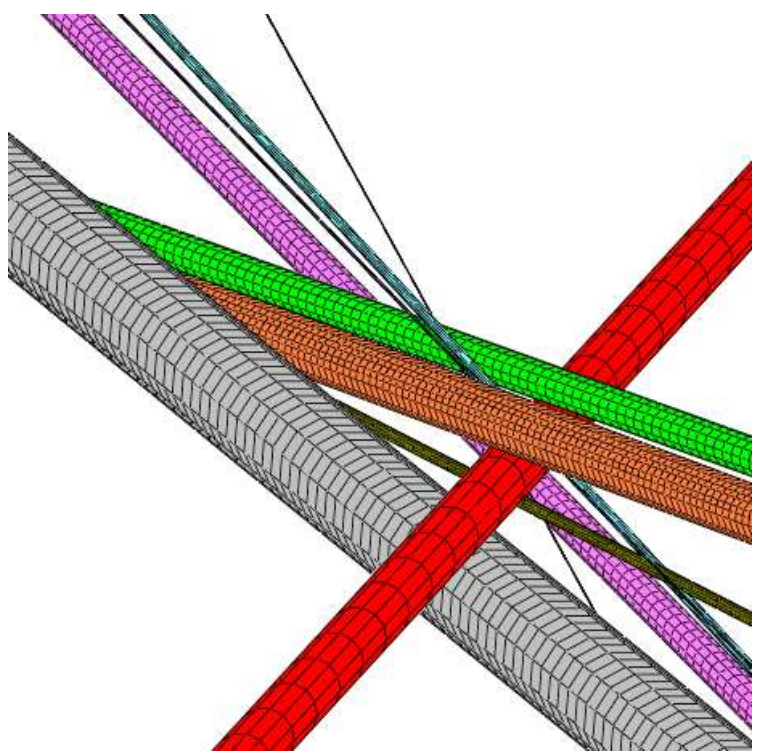

Appendix 2
Appendix 2
9-cross
Configuration label a0d

$$\begin{pmatrix} \text{orange} \\ \text{red} \\ \text{blue} \\ \text{green} \\ \text{cyan} \\ \text{magenta} \\ \text{gray} \\ \text{olive} \\ \text{pink} \end{pmatrix} \begin{pmatrix} 0 & 0 & 0 & 1 \\ t_1 & p_1 & z_1 & r_1 \\ t_2 & p_2 & z_2 & r_2 \\ t_3 & p_3 & z_3 & r_3 \\ t_4 & p_4 & z_4 & r_4 \\ t_5 & p_5 & z_5 & r_5 \\ t_6 & p_6 & z_6 & r_6 \\ t_7 & p_7 & z_7 & r_7 \\ t_8 & p_8 & z_8 & r_8 \end{pmatrix} = \begin{pmatrix} 0 & 0 & 0 & 1 \\ -0.298234814 & 3.1415926536 & 0 & 0.6709065337 \\ -2.0832362636 & 1.0001005257 & 0.9795605159 & 13.2387485882 \\ -0.8501493831 & 3.7749654038 & -44.4175175804 & 4.7740708205 \\ 0.2994069419 & 6.196847276 & 3.508701753 & 0.4987627935 \\ -0.2439755461 & 2.7513674982 & -0.4566915874 & 0.1406982758 \\ 0.4273052862 & 0.591934356 & 18.8582581225 & 2.0417844538 \\ -2.3936221243 & 0.5864538552 & 6.1539559136 & 0.0310374388 \\ 0.490045626 & 0.5915282126 & 3.0206837734 & 0.01 \end{pmatrix}$$

$$P = \begin{pmatrix} 0 & -1 & -1 & -1 & 1 & -1 & 1 & -1 & 1 \\ -1 & 0 & -1 & 1 & 1 & -1 & -1 & -1 & 1 \\ -1 & -1 & 0 & 1 & -1 & -1 & -1 & 1 & -1 \\ -1 & 1 & 1 & 0 & 1 & 1 & -1 & -1 & -1 \\ 1 & 1 & -1 & 1 & 0 & 1 & -1 & 1 & 1 \\ -1 & -1 & -1 & 1 & 1 & 0 & -1 & 1 & 1 \\ 1 & -1 & -1 & -1 & -1 & -1 & 0 & 1 & -1 \\ -1 & -1 & 1 & -1 & 1 & 1 & 1 & 0 & -1 \\ 1 & 1 & -1 & -1 & 1 & 1 & -1 & -1 & 0 \end{pmatrix} \quad R = \begin{pmatrix} 0 & 8 & 5 & 11 & 7 & 5 & 4 & 4 & 4 \\ 0 & 0 & 0 & 0 & 0 & 0 & 0 & 0 & 0 \\ 0 & 0 & 0 & 0 & 0 & 0 & 0 & 0 & 0 \\ 0 & 0 & 0 & 0 & 0 & 0 & 0 & 0 & 0 \\ 4 & 4 & 5 & 8 & 0 & 4 & 11 & 7 & 5 \\ 4 & 8 & 8 & 4 & 4 & 0 & 4 & 8 & 8 \\ 0 & 0 & 0 & 0 & 0 & 0 & 0 & 0 & 0 \\ 5 & 5 & 2 & 2 & 2 & 2 & 10 & 0 & 2 \\ 3 & 3 & 12 & 3 & 4 & 3 & 4 & 4 & 0 \end{pmatrix}$$

$|P| = -6.9944050551 \times 10^{-15}$

$I3(P, R) = 34.5886547812$

$I3(-P, R) = 22.0575688817$

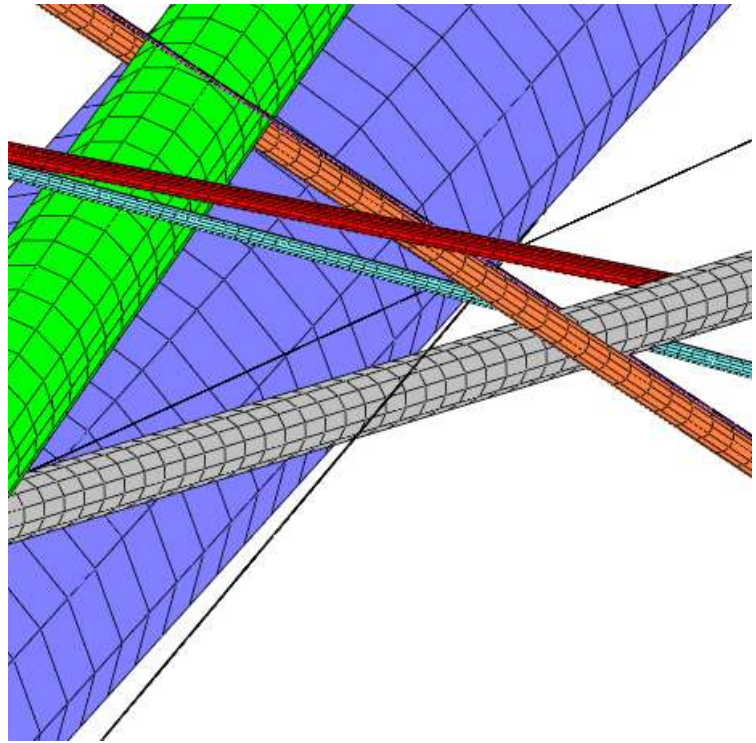

Appendix 2
9-cross
Configuration label am52a

$$\begin{pmatrix} \text{orange} \\ \text{red} \\ \text{blue} \\ \text{green} \\ \text{cyan} \\ \text{magenta} \\ \text{gray} \\ \text{olive} \\ \text{pink} \end{pmatrix} \begin{pmatrix} 0 & 0 & 0 & 1 \\ t_1 & p_1 & z_1 & r_1 \\ t_2 & p_2 & z_2 & r_2 \\ t_3 & p_3 & z_3 & r_3 \\ t_4 & p_4 & z_4 & r_4 \\ t_5 & p_5 & z_5 & r_5 \\ t_6 & p_6 & z_6 & r_6 \\ t_7 & p_7 & z_7 & r_7 \\ t_8 & p_8 & z_8 & r_8 \end{pmatrix} = \begin{pmatrix} 0 & 0 & 0 & 1 \\ 0.1340386074 & 3.1415926536 & 0 & 1.9968460919 \\ 0.9885304723 & 8.935720912 & 52.4184492661 & 5.2372850077 \\ 0.1980244551 & 6.5141019897 & -104.9962792026 & 3.6442883187 \\ 2.8887794879 & 7.4839649434 & 19.9435158792 & 6.5161684705 \\ 0.0757583182 & 7.6139153147 & -104.3899059477 & 0.4752981793 \\ 2.8781560166 & 4.1485778781 & -52.8146593911 & 20.552226246 \\ 3.0594316634 & -2.2277784914 & -89.7093523819 & 0.1986672299 \\ 2.9356002757 & 3.083593533 & -83.933833203 & 0.1986672299 \end{pmatrix}$$

$$P = \begin{pmatrix} 0 & 1 & 1 & 1 & 1 & 1 & 1 & 1 & 1 \\ 1 & 0 & 1 & 1 & 1 & -1 & 1 & 1 & -1 \\ 1 & 1 & 0 & -1 & -1 & -1 & 1 & 1 & 1 \\ 1 & 1 & -1 & 0 & -1 & 1 & 1 & -1 & -1 \\ 1 & 1 & -1 & -1 & 0 & 1 & -1 & -1 & 1 \\ 1 & -1 & -1 & 1 & 1 & 0 & 1 & -1 & -1 \\ 1 & 1 & 1 & 1 & -1 & 1 & 0 & -1 & -1 \\ 1 & 1 & 1 & -1 & -1 & -1 & -1 & 0 & -1 \\ 1 & -1 & 1 & -1 & 1 & -1 & -1 & -1 & 0 \end{pmatrix} \qquad R = \begin{pmatrix} 0 & 2 & 2 & 10 & 2 & 2 & 5 & 5 & 2 \\ 6 & 0 & 9 & 6 & 9 & 6 & 6 & 6 & 6 \\ 0 & 0 & 0 & 0 & 0 & 0 & 0 & 0 & 0 \\ 0 & 0 & 0 & 0 & 0 & 0 & 0 & 0 & 0 \\ 2 & 2 & 10 & 5 & 0 & 2 & 2 & 2 & 5 \\ 1 & 1 & 1 & 1 & 6 & 0 & 6 & 1 & 1 \\ 0 & 0 & 0 & 0 & 0 & 0 & 0 & 0 & 0 \\ 4 & 4 & 7 & 4 & 5 & 8 & 11 & 0 & 5 \\ 7 & 5 & 4 & 11 & 5 & 8 & 4 & 4 & 0 \end{pmatrix}$$

$|P| = -52$

$I3(P,R) = 24.0193570498$

$I3(-P,R) = 21.6857084104$

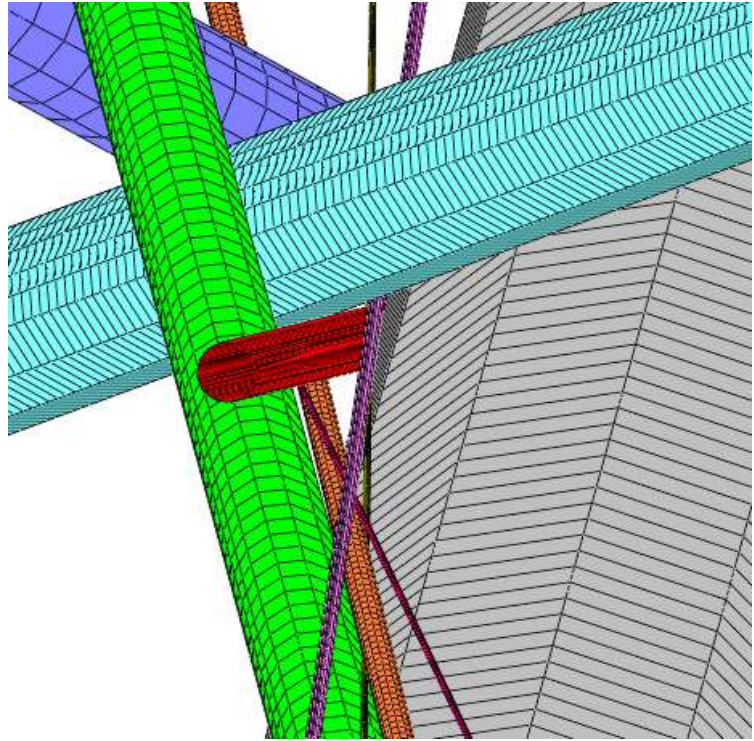

Appendix 2
9-cross
Configuration label a16a

$$\begin{pmatrix} \text{orange} \\ \text{red} \\ \text{blue} \\ \text{green} \\ \text{cyan} \\ \text{magenta} \\ \text{gray} \\ \text{olive} \\ \text{pink} \end{pmatrix} \begin{pmatrix} 0 & 0 & 0 & 1 \\ t_1 & p_1 & z_1 & r_1 \\ t_2 & p_2 & z_2 & r_2 \\ t_3 & p_3 & z_3 & r_3 \\ t_4 & p_4 & z_4 & r_4 \\ t_5 & p_5 & z_5 & r_5 \\ t_6 & p_6 & z_6 & r_6 \\ t_7 & p_7 & z_7 & r_7 \\ t_8 & p_8 & z_8 & r_8 \end{pmatrix} = \begin{pmatrix} 0 & 0 & 0 & 1 \\ 1.3624573075 & 3.1415926536 & 0 & 1.1253423662 \\ 2.0370199172 & 8.9782965227 & 4.3279705034 & 0.8522641276 \\ 0.3292355244 & 6.7599867497 & -11.9983435684 & 1.1253423662 \\ 1.3815920753 & 8.3589327697 & 1.8067539981 & 1.1253423662 \\ 1.2522006245 & 8.8692438739 & -9.3268299409 & 8.2294935906 \\ 4.0502510864 & 5.3820397842 & -7.9050531083 & 14.4127228981 \\ 2.8883760433 & -1.5014296653 & -14.346597525 & 0.3941363942 \\ 2.7819958764 & 3.0392103464 & -7.4368784066 & 0.1345213406 \end{pmatrix}$$

$$P = \begin{pmatrix} 0 & 1 & 1 & 1 & 1 & 1 & -1 & 1 & 1 \\ 1 & 0 & 1 & 1 & 1 & -1 & -1 & 1 & -1 \\ 1 & 1 & 0 & -1 & -1 & -1 & -1 & 1 & 1 \\ 1 & 1 & -1 & 0 & -1 & 1 & -1 & -1 & -1 \\ 1 & 1 & -1 & -1 & 0 & 1 & -1 & -1 & 1 \\ 1 & -1 & -1 & 1 & 1 & 0 & -1 & -1 & -1 \\ -1 & -1 & -1 & -1 & -1 & -1 & 0 & 1 & 1 \\ 1 & 1 & 1 & -1 & -1 & -1 & 1 & 0 & -1 \\ 1 & -1 & 1 & -1 & 1 & -1 & 1 & -1 & 0 \end{pmatrix} \qquad R = \begin{pmatrix} 0 & 3 & 3 & 9 & 6 & 3 & 9 & 6 & 3 \\ 5 & 0 & 10 & 5 & 9 & 7 & 6 & 7 & 5 \\ 0 & 0 & 0 & 0 & 0 & 0 & 0 & 0 & 0 \\ 0 & 0 & 0 & 0 & 0 & 0 & 0 & 0 & 0 \\ 2 & 2 & 10 & 5 & 0 & 2 & 2 & 2 & 5 \\ 0 & 0 & 0 & 0 & 0 & 0 & 0 & 0 & 0 \\ 0 & 0 & 0 & 0 & 0 & 0 & 0 & 0 & 0 \\ 3 & 3 & 6 & 3 & 3 & 9 & 9 & 0 & 6 \\ 6 & 6 & 3 & 9 & 3 & 9 & 3 & 3 & 0 \end{pmatrix}$$

$|P| = 16$

$I3(P, R) = 27.6275057092$

$I3(-P, R) = 29.7336970312$

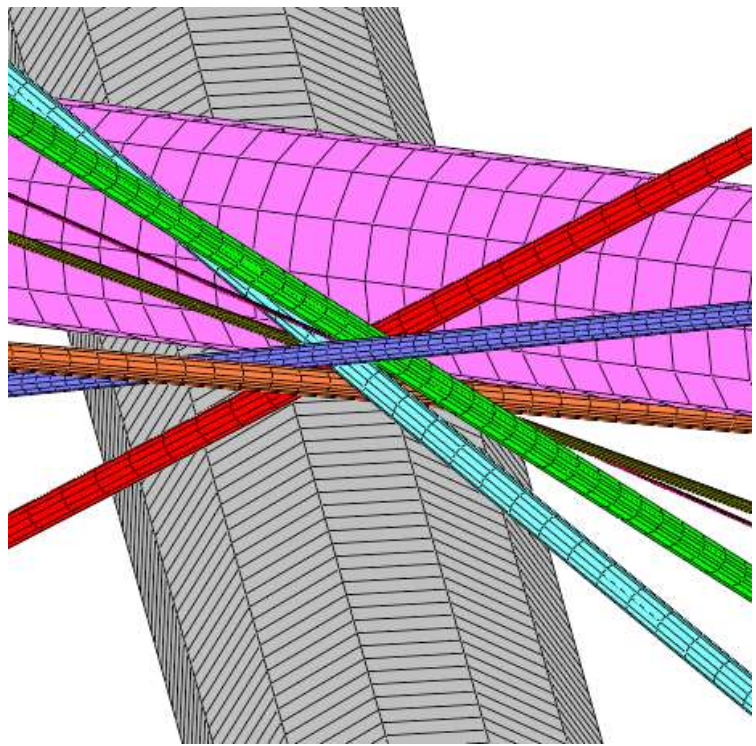

Appendix 2
9-cross
Configuration label a4b

$$\begin{pmatrix} \text{orange} \\ \text{red} \\ \text{blue} \\ \text{green} \\ \text{cyan} \\ \text{magenta} \\ \text{gray} \\ \text{olive} \\ \text{pink} \end{pmatrix} \begin{pmatrix} 0 & 0 & 0 & 1 \\ t_1 & p_1 & z_1 & r_1 \\ t_2 & p_2 & z_2 & r_2 \\ t_3 & p_3 & z_3 & r_3 \\ t_4 & p_4 & z_4 & r_4 \\ t_5 & p_5 & z_5 & r_5 \\ t_6 & p_6 & z_6 & r_6 \\ t_7 & p_7 & z_7 & r_7 \\ t_8 & p_8 & z_8 & r_8 \end{pmatrix} = \begin{pmatrix} 0 & 0 & 0 & 1 \\ -1.1018726105 & 3.1415926536 & 0 & 1.9405904123 \\ 1.0118250445 & 0.4254500737 & 1.1108627496 & 0.1202631601 \\ -2.271819041 & 4.06534099 & -4.2952056126 & 1.1996132639 \\ 1.2772092852 & 5.806599857 & 1.3988710502 & 0.929108154 \\ -0.6877061769 & 1.8653041045 & -0.1083027478 & 0.340256436 \\ 1.2305863236 & 0.6622592017 & 4.4561526786 & 1.1305019554 \\ -1.4950886842 & 0.3253920705 & 2.804864627 & 0.3976966419 \\ 1.4565105106 & 0.5152116172 & -0.3430133381 & 1.2351175625 \end{pmatrix}$$

$$P = \begin{pmatrix} 0 & -1 & 1 & -1 & 1 & -1 & 1 & -1 & 1 \\ -1 & 0 & -1 & 1 & 1 & -1 & -1 & -1 & 1 \\ 1 & -1 & 0 & -1 & 1 & -1 & -1 & -1 & 1 \\ -1 & 1 & -1 & 0 & 1 & 1 & -1 & -1 & -1 \\ 1 & 1 & 1 & 1 & 0 & 1 & -1 & 1 & 1 \\ -1 & -1 & -1 & 1 & 1 & 0 & -1 & 1 & 1 \\ 1 & -1 & -1 & -1 & -1 & -1 & 0 & 1 & -1 \\ -1 & -1 & -1 & -1 & 1 & 1 & 1 & 0 & -1 \\ 1 & 1 & 1 & -1 & 1 & 1 & -1 & -1 & 0 \end{pmatrix} \qquad R = \begin{pmatrix} 0 & 9 & 3 & 9 & 6 & 6 & 3 & 3 & 3 \\ 0 & 0 & 0 & 0 & 0 & 0 & 0 & 0 & 0 \\ 8 & 9 & 0 & 7 & 6 & 9 & 6 & 7 & 8 \\ 0 & 0 & 0 & 0 & 0 & 0 & 0 & 0 & 0 \\ 3 & 3 & 3 & 9 & 0 & 3 & 9 & 6 & 6 \\ 2 & 5 & 2 & 2 & 2 & 0 & 2 & 5 & 10 \\ 0 & 0 & 0 & 0 & 0 & 0 & 0 & 0 & 0 \\ 4 & 4 & 4 & 3 & 3 & 3 & 12 & 0 & 3 \\ 0 & 0 & 0 & 0 & 0 & 0 & 0 & 0 & 0 \end{pmatrix}$$

$|P| = 4$

$I3(P, R) = 27.0174639332$

$I3(-P, R) = 29.4031890661$

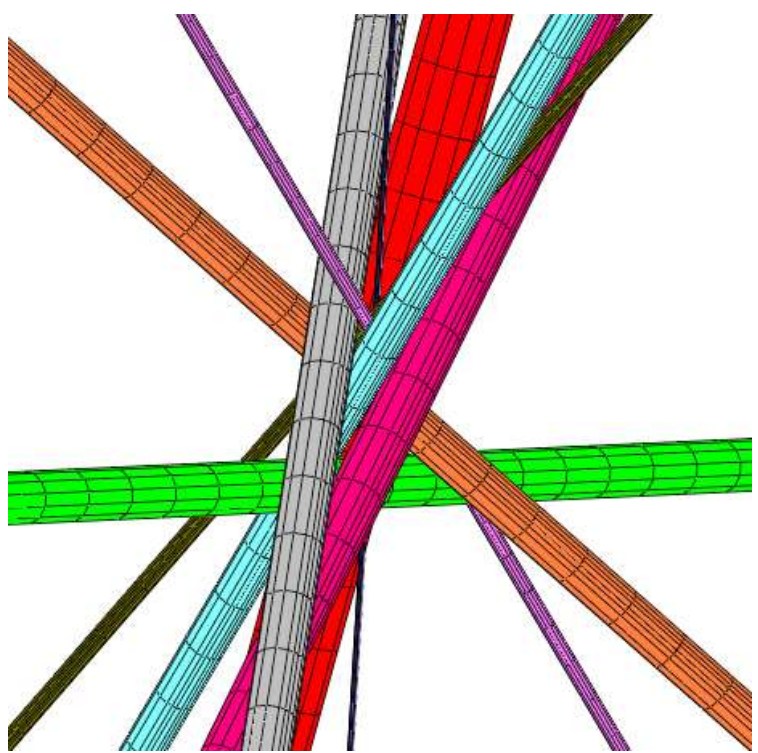

Appendix 2
9-cross
Configuration label a4a

$$\begin{pmatrix} \text{orange} \\ \text{red} \\ \text{blue} \\ \text{green} \\ \text{cyan} \\ \text{magenta} \\ \text{gray} \\ \text{olive} \\ \text{pink} \end{pmatrix} \begin{pmatrix} 0 & 0 & 0 & 1 \\ t_1 & p_1 & z_1 & r_1 \\ t_2 & p_2 & z_2 & r_2 \\ t_3 & p_3 & z_3 & r_3 \\ t_4 & p_4 & z_4 & r_4 \\ t_5 & p_5 & z_5 & r_5 \\ t_6 & p_6 & z_6 & r_6 \\ t_7 & p_7 & z_7 & r_7 \\ t_8 & p_8 & z_8 & r_8 \end{pmatrix} = \begin{pmatrix} 0 & 0 & 0 & 1 \\ 0.9071368704 & 3.1415926536 & 0 & 0.4424161588 \\ 1.6477467187 & 9.0151175871 & 1.9436660559 & 0.7245086968 \\ 1.2499869865 & 6.8746445895 & -2.7377466039 & 2.7041739474 \\ 1.7373910372 & 8.1593100218 & 0.6658298675 & 0.7245086968 \\ 1.014197371 & 8.8282882751 & -3.9600929517 & 1.7102533871 \\ 4.2979678016 & 5.3257731897 & 0.8738672527 & 3.0091345212 \\ 2.3670638427 & -1.9462015992 & -2.0960004929 & 0.7245086968 \\ 2.0886101107 & 2.6814005708 & -1.1042053298 & 0.1548918413 \end{pmatrix}$$

$$P = \begin{pmatrix} 0 & 1 & 1 & 1 & 1 & 1 & -1 & 1 & 1 \\ 1 & 0 & 1 & 1 & 1 & -1 & -1 & 1 & -1 \\ 1 & 1 & 0 & -1 & -1 & -1 & -1 & 1 & -1 \\ 1 & 1 & -1 & 0 & -1 & 1 & -1 & -1 & -1 \\ 1 & 1 & -1 & -1 & 0 & 1 & -1 & -1 & 1 \\ 1 & -1 & -1 & 1 & 1 & 0 & -1 & -1 & -1 \\ -1 & -1 & -1 & -1 & -1 & -1 & 0 & 1 & 1 \\ 1 & 1 & 1 & -1 & -1 & -1 & 1 & 0 & -1 \\ 1 & -1 & -1 & -1 & 1 & -1 & 1 & -1 & 0 \end{pmatrix} \qquad R = \begin{pmatrix} 0 & 3 & 3 & 9 & 6 & 3 & 9 & 6 & 3 \\ 5 & 0 & 10 & 5 & 9 & 7 & 6 & 7 & 5 \\ 0 & 0 & 0 & 0 & 0 & 0 & 0 & 0 & 0 \\ 0 & 0 & 0 & 0 & 0 & 0 & 0 & 0 & 0 \\ 2 & 2 & 10 & 5 & 0 & 2 & 2 & 2 & 5 \\ 0 & 0 & 0 & 0 & 0 & 0 & 0 & 0 & 0 \\ 0 & 0 & 0 & 0 & 0 & 0 & 0 & 0 & 0 \\ 3 & 3 & 6 & 3 & 3 & 9 & 9 & 0 & 6 \\ 7 & 5 & 5 & 11 & 4 & 8 & 4 & 4 & 0 \end{pmatrix}$$

$|P| = 4$

$I3(P, R) = 25.3527263569$

$I3(-P, R) = 31.380556605$

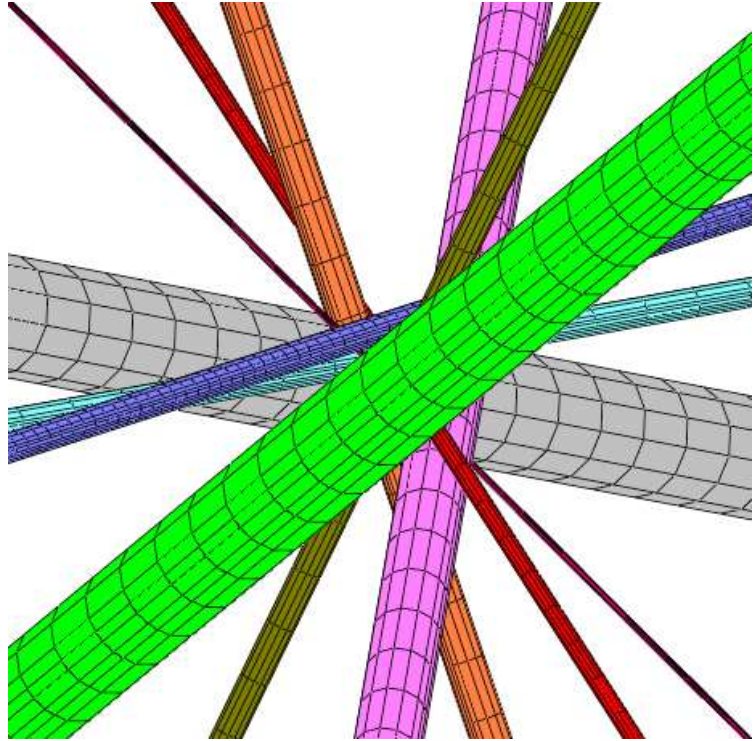

Appendix 2
9-cross
Configuration label a0c

$$\begin{pmatrix} \text{orange} \\ \text{red} \\ \text{blue} \\ \text{green} \\ \text{cyan} \\ \text{magenta} \\ \text{gray} \\ \text{olive} \\ \text{pink} \end{pmatrix} \begin{pmatrix} 0 & 0 & 0 & 1 \\ t1 & p1 & z1 & r1 \\ t2 & p2 & z2 & r2 \\ t3 & p3 & z3 & r3 \\ t4 & p4 & z4 & r4 \\ t5 & p5 & z5 & r5 \\ t6 & p6 & z6 & r6 \\ t7 & p7 & z7 & r7 \\ t8 & p8 & z8 & r8 \end{pmatrix} = \begin{pmatrix} 0 & 0 & 0 & 1 \\ 0.4532166313 & 3.1415926536 & 0 & 0.0595950528 \\ 0.5093229517 & 9.4474659246 & -1.4407510339 & 0.5213665534 \\ 2.4080159401 & 7.2426313562 & -1.6649518316 & 0.0325546844 \\ 0.54156115 & 9.2147404402 & 0.5053387975 & 0.8156030417 \\ 0.4990690202 & 9.2412839341 & -0.7621071696 & 0.0817082645 \\ 2.6060717325 & 6.0125681504 & 2.0433513998 & 0.6183797348 \\ 2.6499770645 & -0.4179085499 & 0.0947805221 & 0.4518842431 \\ 0.515649382 & 3.0830422326 & -0.3442752742 & 0.0194867152 \end{pmatrix}$$

$$P = \begin{pmatrix} 0 & 1 & 1 & 1 & 1 & 1 & 1 & 1 & 1 \\ 1 & 0 & 1 & 1 & 1 & -1 & 1 & 1 & -1 \\ 1 & 1 & 0 & 1 & 1 & 1 & -1 & -1 & 1 \\ 1 & 1 & 1 & 0 & -1 & 1 & 1 & -1 & -1 \\ 1 & 1 & 1 & -1 & 0 & 1 & -1 & -1 & 1 \\ 1 & -1 & 1 & 1 & 1 & 0 & 1 & -1 & -1 \\ 1 & 1 & -1 & 1 & -1 & 1 & 0 & -1 & -1 \\ 1 & 1 & -1 & -1 & -1 & -1 & -1 & 0 & -1 \\ 1 & -1 & 1 & -1 & 1 & -1 & -1 & -1 & 0 \end{pmatrix} \qquad R = \begin{pmatrix} 0 & 2 & 2 & 10 & 2 & 2 & 5 & 5 & 2 \\ 5 & 0 & 7 & 8 & 11 & 4 & 4 & 4 & 5 \\ 0 & 0 & 0 & 0 & 0 & 0 & 0 & 0 & 0 \\ 0 & 0 & 0 & 0 & 0 & 0 & 0 & 0 & 0 \\ 0 & 0 & 0 & 0 & 0 & 0 & 0 & 0 & 0 \\ 4 & 5 & 11 & 5 & 8 & 0 & 7 & 4 & 4 \\ 0 & 0 & 0 & 0 & 0 & 0 & 0 & 0 & 0 \\ 4 & 4 & 3 & 4 & 3 & 3 & 12 & 0 & 3 \\ 8 & 7 & 8 & 8 & 7 & 8 & 7 & 7 & 0 \end{pmatrix}$$

$|P| = 1.2490009027 \times 10^{-14}$

$I3(P,R) = 10.3290278722$

$I3(-P,R) = 46.4133242692$

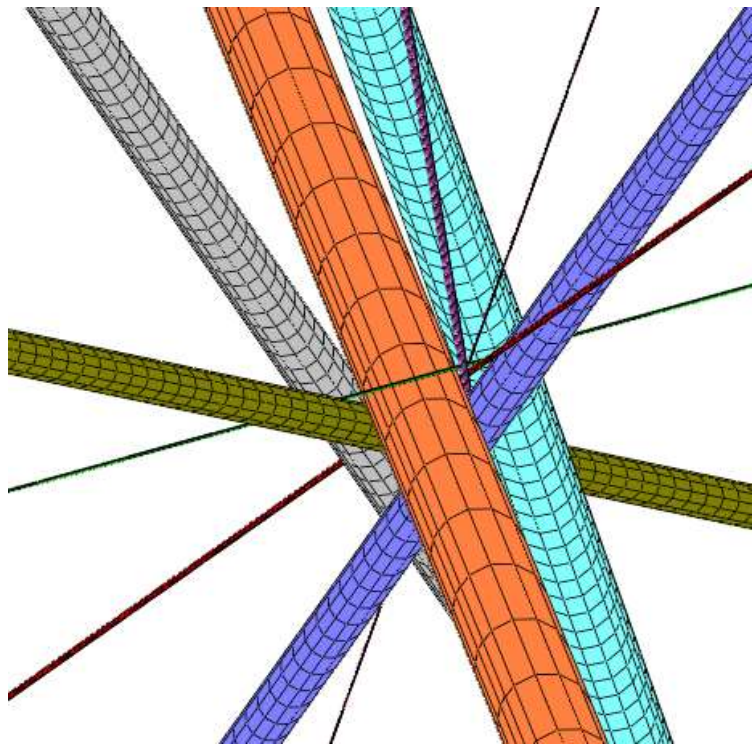

Appendix 2
9-cross
Configuration label a0b

$$\begin{pmatrix} \text{orange} \\ \text{red} \\ \text{blue} \\ \text{green} \\ \text{cyan} \\ \text{magenta} \\ \text{gray} \\ \text{olive} \\ \text{pink} \end{pmatrix} \begin{pmatrix} 0 & 0 & 0 & 1 \\ t_1 & p_1 & z_1 & r_1 \\ t_2 & p_2 & z_2 & r_2 \\ t_3 & p_3 & z_3 & r_3 \\ t_4 & p_4 & z_4 & r_4 \\ t_5 & p_5 & z_5 & r_5 \\ t_6 & p_6 & z_6 & r_6 \\ t_7 & p_7 & z_7 & r_7 \\ t_8 & p_8 & z_8 & r_8 \end{pmatrix} = \begin{pmatrix} 0 & 0 & 0 & 1 \\ 1.551872711 & 3.1415926536 & 0 & 1 \\ 2.4476796431 & 9.1311336 & 9.0148909954 & 1 \\ 0.2900183853 & 7.1697832107 & -5.0012529896 & 1.168728773 \\ 2.5861135333 & 8.6737152637 & 6.4096113692 & 0.9862923782 \\ 2.0499877011 & 9.0332257445 & -5.0636943264 & 1.5961769794 \\ 3.7181039642 & 5.6043646099 & -2.7210319717 & 9.4891794667 \\ 2.7897221676 & -1.2088148737 & -11.8027988457 & 0.4673591939 \\ 2.7841216303 & 2.8258085399 & -0.9015573909 & 0.0745888161 \end{pmatrix}$$

$$P = \begin{pmatrix} 0 & 1 & 1 & 1 & 1 & 1 & -1 & 1 & 1 \\ 1 & 0 & 1 & 1 & 1 & -1 & -1 & 1 & -1 \\ 1 & 1 & 0 & -1 & -1 & -1 & -1 & 1 & -1 \\ 1 & 1 & -1 & 0 & -1 & 1 & -1 & -1 & -1 \\ 1 & 1 & -1 & -1 & 0 & 1 & 1 & -1 & 1 \\ 1 & -1 & -1 & 1 & 1 & 0 & -1 & -1 & -1 \\ -1 & -1 & -1 & -1 & 1 & -1 & 0 & 1 & 1 \\ 1 & 1 & 1 & -1 & -1 & -1 & 1 & 0 & -1 \\ 1 & -1 & -1 & -1 & 1 & -1 & 1 & -1 & 0 \end{pmatrix} \qquad R = \begin{pmatrix} 0 & 3 & 3 & 9 & 6 & 3 & 9 & 6 & 3 \\ 5 & 0 & 10 & 5 & 9 & 7 & 6 & 7 & 5 \\ 0 & 0 & 0 & 0 & 0 & 0 & 0 & 0 & 0 \\ 0 & 0 & 0 & 0 & 0 & 0 & 0 & 0 & 0 \\ 2 & 2 & 10 & 5 & 0 & 2 & 2 & 2 & 5 \\ 0 & 0 & 0 & 0 & 0 & 0 & 0 & 0 & 0 \\ 0 & 0 & 0 & 0 & 0 & 0 & 0 & 0 & 0 \\ 3 & 3 & 6 & 3 & 3 & 9 & 9 & 0 & 6 \\ 8 & 4 & 4 & 11 & 4 & 7 & 5 & 5 & 0 \end{pmatrix}$$

$|P| = 8.7430063189 \times 10^{-15}$

$I3(P, R) = 31.2786803919$

$I3(-P, R) = 25.4237876721$

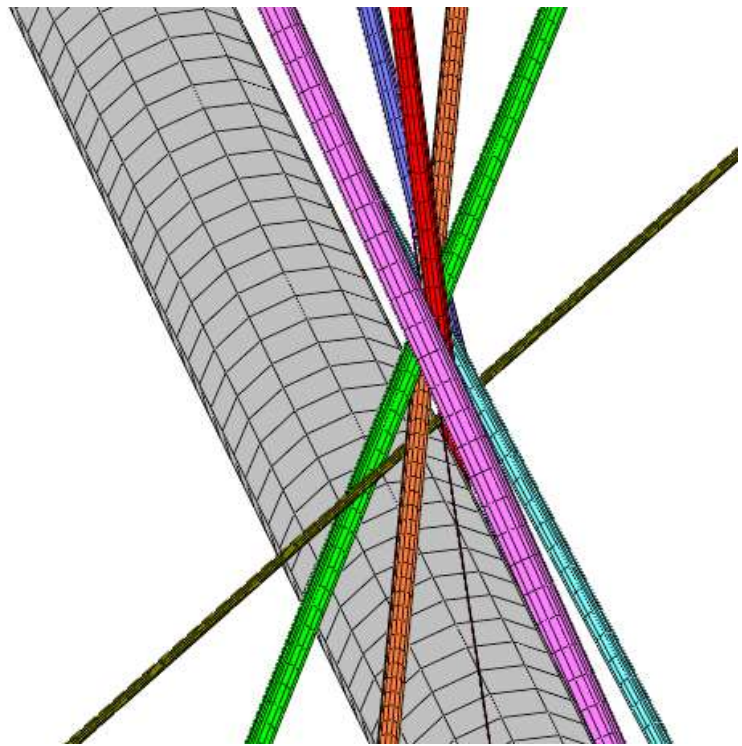

Appendix 2
9-cross
Configuration label a0a

$$\begin{pmatrix} \text{orange} \\ \text{red} \\ \text{blue} \\ \text{green} \\ \text{cyan} \\ \text{magenta} \\ \text{gray} \\ \text{olive} \\ \text{pink} \end{pmatrix} \begin{pmatrix} 0 & 0 & 0 & 1 \\ t_1 & p_1 & z_1 & r_1 \\ t_2 & p_2 & z_2 & r_2 \\ t_3 & p_3 & z_3 & r_3 \\ t_4 & p_4 & z_4 & r_4 \\ t_5 & p_5 & z_5 & r_5 \\ t_6 & p_6 & z_6 & r_6 \\ t_7 & p_7 & z_7 & r_7 \\ t_8 & p_8 & z_8 & r_8 \end{pmatrix} = \begin{pmatrix} 0 & 0 & 0 & 1 \\ 0.8296258347 & 3.1415926536 & 0 & 0.8987178684 \\ 1.7236362897 & 9.0894456276 & 6.0179979344 & 2.12645746 \\ 0.6686193455 & 6.8669620361 & -6.5792016629 & 1.3328198731 \\ 1.9157990476 & 8.086741106 & 2.3741681888 & 2.12645746 \\ 0.5396392432 & 8.248158423 & -5.4862825633 & 0.3132242642 \\ 2.2923604875 & 4.5916901697 & -3.718234883 & 2.0493522703 \\ 2.6789741053 & -1.8887006681 & -5.6307818502 & 0.1639896318 \\ 2.2608670154 & 2.7610732152 & -2.9725791599 & 0.3005267637 \end{pmatrix}$$

$$P = \begin{pmatrix} 0 & 1 & 1 & 1 & 1 & 1 & 1 & 1 & 1 \\ 1 & 0 & 1 & 1 & 1 & -1 & 1 & 1 & -1 \\ 1 & 1 & 0 & -1 & -1 & -1 & 1 & 1 & -1 \\ 1 & 1 & -1 & 0 & -1 & 1 & 1 & -1 & -1 \\ 1 & 1 & -1 & -1 & 0 & 1 & -1 & -1 & 1 \\ 1 & -1 & -1 & 1 & 1 & 0 & 1 & -1 & -1 \\ 1 & 1 & 1 & 1 & -1 & 1 & 0 & -1 & -1 \\ 1 & 1 & 1 & -1 & -1 & -1 & -1 & 0 & -1 \\ 1 & -1 & -1 & -1 & 1 & -1 & -1 & -1 & 0 \end{pmatrix} \qquad R = \begin{pmatrix} 0 & 2 & 2 & 10 & 2 & 2 & 5 & 5 & 2 \\ 6 & 0 & 9 & 6 & 9 & 6 & 6 & 6 & 6 \\ 0 & 0 & 0 & 0 & 0 & 0 & 0 & 0 & 0 \\ 0 & 0 & 0 & 0 & 0 & 0 & 0 & 0 & 0 \\ 2 & 2 & 10 & 5 & 0 & 2 & 2 & 2 & 5 \\ 1 & 1 & 1 & 1 & 6 & 0 & 6 & 1 & 1 \\ 0 & 0 & 0 & 0 & 0 & 0 & 0 & 0 & 0 \\ 4 & 4 & 7 & 4 & 5 & 8 & 11 & 0 & 5 \\ 8 & 4 & 4 & 11 & 4 & 7 & 5 & 5 & 0 \end{pmatrix}$$

$|P| = 8.7430063189 \times 10^{-15}$

$I3(P, R) = 20.8373932182$

$I3(-P, R) = 23.8438467564$

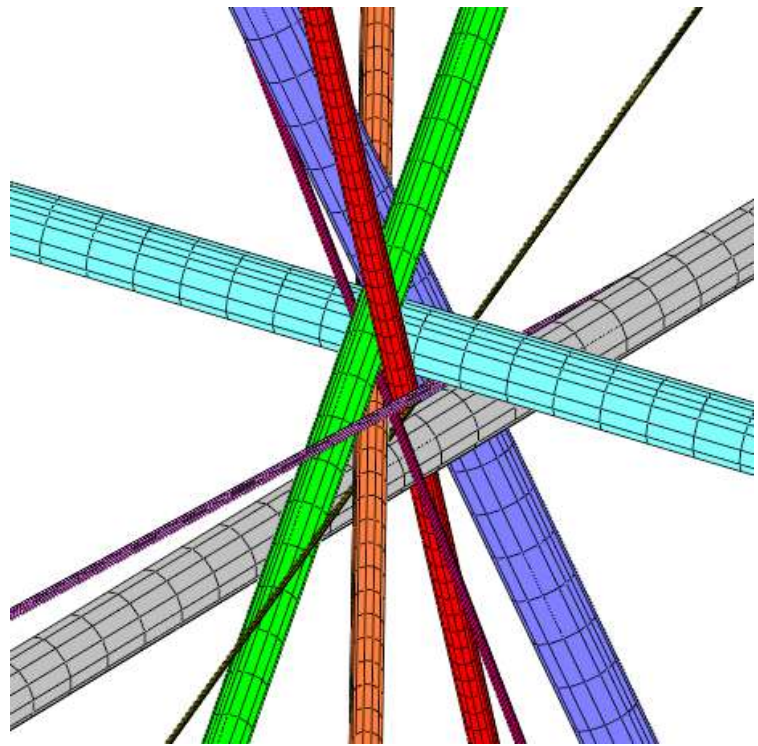

Appendix 3
Configurations with asymptotically equal radii.

7-cross
Configuration label **a59**

$$\begin{matrix} \text{orange} \\ \text{red} \\ \text{blue} \\ \text{green} \\ \text{cyan} \\ \text{magenta} \\ \text{gray} \end{matrix} \begin{pmatrix} 0 & 0 & 0 & 1 \\ t1 & p1 & z1 & r1 \\ t2 & p2 & z2 & r2 \\ t3 & p3 & z3 & r3 \\ t4 & p4 & z4 & r4 \\ t5 & p5 & z5 & r5 \\ t6 & p6 & z6 & r6 \end{pmatrix} = \begin{pmatrix} 0 & 0 & 0 & 1 \\ 1.5881497698 & 3.1415926536 & 0 & 0.9997507739 \\ 2.1985172094 & 8.6482186227 & 3.4610198657 & 0.9999832757 \\ 1.2029297351 & 7.4655177816 & 1.0648061394 & 1.0002747231 \\ 2.1983795887 & 8.6481369457 & -2.2451838529 & 0.9998822682 \\ 0.943086971 & 5.5064378732 & 0.6085524422 & 1.0001305012 \\ 1.9385960527 & -1.9590548291 & -2.7015773987 & 0.9999564662 \end{pmatrix}$$

$$P = \begin{pmatrix} 0 & 1 & 1 & 1 & 1 & 1 & 1 \\ 1 & 0 & 1 & 1 & -1 & 1 & 1 \\ 1 & 1 & 0 & -1 & -1 & 1 & 1 \\ 1 & 1 & -1 & 0 & 1 & 1 & -1 \\ 1 & -1 & -1 & 1 & 0 & 1 & -1 \\ 1 & 1 & 1 & 1 & 1 & 0 & -1 \\ 1 & 1 & 1 & -1 & -1 & -1 & 0 \end{pmatrix} \qquad R = \begin{pmatrix} 0 & 2 & 2 & 6 & 2 & 3 & 3 \\ 3 & 0 & 6 & 3 & 2 & 2 & 2 \\ 0 & 0 & 0 & 0 & 0 & 0 & 0 \\ 0 & 0 & 0 & 0 & 0 & 0 & 0 \\ 0 & 0 & 0 & 0 & 0 & 0 & 0 \\ 0 & 0 & 0 & 0 & 0 & 0 & 0 \\ 2 & 2 & 3 & 2 & 6 & 3 & 0 \end{pmatrix}$$

$I3(-P, R) = 20.4390243902$

$|P| = -2$

$I3(P, R) = 6.8780487805$

$-(x+2) \cdot \left(x^3 - x^2 - 9 \cdot x + 1\right)^2$

$$\text{eigenvals}(P) = \begin{pmatrix} -2 \\ 3.4939592074 \\ 3.4939592074 \\ 0.1099162642 \\ 0.1099162642 \\ -2.6038754716 \\ -2.6038754716 \end{pmatrix}$$

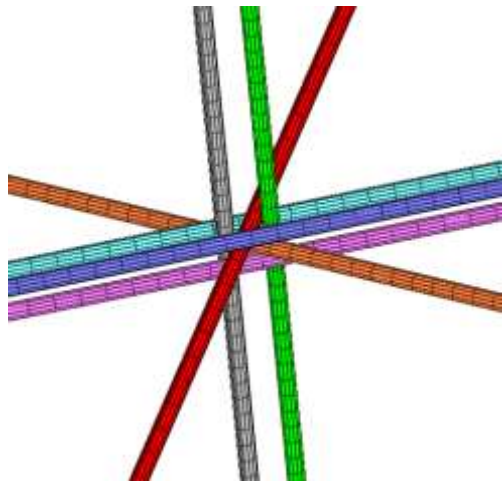

(xx, yy, zz), (xx1, yy1, zz1), (xx2, yy2, zz2), (xx3, yy3, zz3), (xx4, yy4, zz4), (xx5, yy5, zz5), (xx6, yy6, zz6)

Appendix 3

7-cross
Configuration label **m18c**

$$\begin{pmatrix} \text{orange} \\ \text{red} \\ \text{blue} \\ \text{green} \\ \text{cyan} \\ \text{magenta} \\ \text{gray} \end{pmatrix} \begin{pmatrix} 0 & 0 & 0 & 1 \\ t1 & p1 & z1 & r1 \\ t2 & p2 & z2 & r2 \\ t3 & p3 & z3 & r3 \\ t4 & p4 & z4 & r4 \\ t5 & p5 & z5 & r5 \\ t6 & p6 & z6 & r6 \end{pmatrix} = \begin{pmatrix} 0 & 0 & 0 & 1 \\ 1.1357762441 & 3.1415926536 & 0 & 0.9621027456 \\ 1.1309234184 & -0.952176621 & 0.8033402213 & 0.9746006523 \\ 1.4978712442 & 0.8343583044 & -0.0357113726 & 0.9606591048 \\ 1.6454323958 & -2.3085372364 & -2.0207026955 & 0.9819885055 \\ -1.115465217 & -0.017715365 & -2.0686754214 & 0.9823819141 \\ -1.1251391617 & 2.1686142806 & -2.8930318798 & 1.0377318223 \end{pmatrix}$$

$$P = \begin{pmatrix} 0 & 1 & 1 & 1 & 1 & -1 & -1 \\ 1 & 0 & 1 & 1 & 1 & 1 & 1 \\ 1 & 1 & 0 & 1 & -1 & -1 & -1 \\ 1 & 1 & 1 & 0 & -1 & 1 & -1 \\ 1 & 1 & -1 & -1 & 0 & 1 & 1 \\ -1 & 1 & -1 & 1 & 1 & 0 & -1 \\ -1 & 1 & -1 & -1 & 1 & -1 & 0 \end{pmatrix} \qquad R = \begin{pmatrix} 0 & 4 & 4 & 4 & 4 & 4 & 4 \\ 1 & 0 & 1 & 4 & 4 & 1 & 1 \\ 1 & 4 & 0 & 1 & 1 & 4 & 1 \\ 1 & 1 & 4 & 0 & 1 & 1 & 4 \\ 0 & 0 & 0 & 0 & 0 & 0 & 0 \\ 0 & 0 & 0 & 0 & 0 & 0 & 0 \\ 0 & 0 & 0 & 0 & 0 & 0 & 0 \end{pmatrix}$$

I3(P, R) = 18.2171052632

I3(−P, R) = 11.8421052632

$|P| = -18$

$$\text{eigenvals}(P) = \begin{pmatrix} -3.5950104757 \\ -2.6038754716 \\ -1.3567190173 \\ 3.4939592074 \\ 0.1099162642 \\ 2.4379957859 \\ 1.513733707 \end{pmatrix}$$

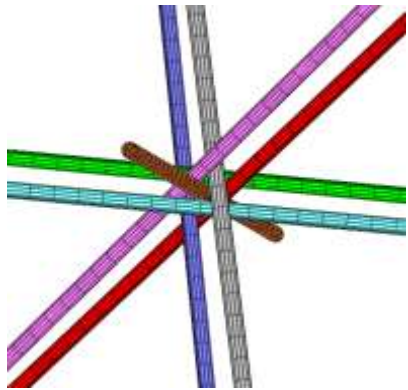

(xx,yy,zz),(xx1,yy1,zz1),(xx2,yy2,zz2),(xx3,yy3,zz3),(xx4,yy4,zz4),(xx5,yy5,zz5),(xx6,yy6,zz6)

# Appendix 3

7-cross
Configuration label **ma49**

$$\begin{pmatrix} \text{orange} \\ \text{red} \\ \text{blue} \\ \text{green} \\ \text{cyan} \\ \text{magenta} \\ \text{gray} \end{pmatrix} \begin{pmatrix} 0 & 0 & 0 & 1 \\ t_1 & p_1 & z_1 & r_1 \\ t_2 & p_2 & z_2 & r_2 \\ t_3 & p_3 & z_3 & r_3 \\ t_4 & p_4 & z_4 & r_4 \\ t_5 & p_5 & z_5 & r_5 \\ t_6 & p_6 & z_6 & r_6 \end{pmatrix} = \begin{pmatrix} 0 & 0 & 0 & 1 \\ 2.8828363798 & 3.1415926536 & 0 & 0.960000066 \\ 0.2855985908 & -0.3529409833 & 3.4513762605 & 0.9600000113 \\ 0.3001499523 & -0.1829409243 & -2.5317839264 & 1.0156698742 \\ 2.8976984933 & -0.2712958219 & -171.2549133418 & 1.0281841601 \\ 2.8413527737 & 2.9589099591 & -7.9436636703 & 1.04 \\ 0.2651785462 & 9.3323462987 & -155.0231975459 & 0.9600001115 \end{pmatrix}$$

$$P = \begin{pmatrix} 0 & 1 & 1 & 1 & 1 & 1 & 1 \\ 1 & 0 & -1 & -1 & 1 & -1 & -1 \\ 1 & -1 & 0 & 1 & 1 & -1 & -1 \\ 1 & -1 & 1 & 0 & -1 & 1 & 1 \\ 1 & 1 & 1 & -1 & 0 & -1 & 1 \\ 1 & -1 & -1 & 1 & -1 & 0 & 1 \\ 1 & -1 & -1 & 1 & 1 & 1 & 0 \end{pmatrix} \qquad R = \begin{pmatrix} 0 & 0 & 0 & 0 & 0 & 0 & 0 \\ 6 & 0 & 3 & 3 & 2 & 2 & 2 \\ 0 & 0 & 0 & 0 & 0 & 0 & 0 \\ 0 & 0 & 0 & 0 & 0 & 0 & 0 \\ 1 & 1 & 1 & 4 & 0 & 4 & 1 \\ 0 & 0 & 0 & 0 & 0 & 0 & 0 \\ 3 & 3 & 4 & 4 & 5 & 5 & 0 \end{pmatrix}$$

$I3(P, R) = -1$

$I3(-P, R) = 22.5714285714$

$|P| = -18$

$$\text{eigenvals}(P) = \begin{pmatrix} -3.5950104757 \\ -2.6038754716 \\ -1.3567190173 \\ 3.4939592074 \\ 0.1099162642 \\ 2.4379957859 \\ 1.513733707 \end{pmatrix}$$

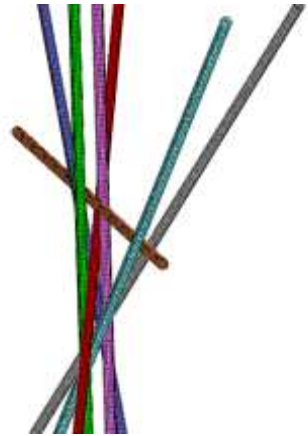

(xx, yy, zz), (xx1, yy1, zz1), (xx2, yy2, zz2), (xx3, yy3, zz3), (xx4, yy4, zz4), (xx5, yy5, zz5), (xx6, yy6, zz6)

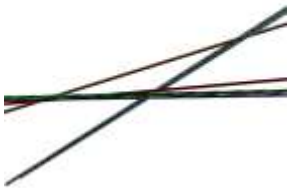

(xx, yy, zz), (xx1, yy1, zz1), (xx2, yy2, zz2), (xx3, yy3, zz3), (xx4, yy4, zz4), (xx5, yy5, zz5), (xx6, yy6, zz6)

# Appendix 3

7-cross
Configuration label **ma37**

$$\begin{matrix} \text{orange} \\ \text{red} \\ \text{blue} \\ \text{green} \\ \text{cyan} \\ \text{magenta} \\ \text{gray} \end{matrix} \begin{pmatrix} 0 & 0 & 0 & 1 \\ t1 & p1 & z1 & r1 \\ t2 & p2 & z2 & r2 \\ t3 & p3 & z3 & r3 \\ t4 & p4 & z4 & r4 \\ t5 & p5 & z5 & r5 \\ t6 & p6 & z6 & r6 \end{pmatrix} = \begin{pmatrix} 0 & 0 & 0 & 1 \\ 2.7512111711 & 3.1415926536 & 0 & 0.91 \\ 2.6173661055 & -0.2239255665 & -36.7105573466 & 1.09 \\ 0.4844036062 & -0.4366561104 & 0.5837311788 & 1.09 \\ 2.6554922801 & -0.6093305697 & -39.5047863813 & 0.9100000059 \\ 0.5229758968 & 9.1987535604 & -33.5304693544 & 1.09 \\ 0.4007840195 & 3.379773942 & -36.7508607087 & 0.9100000002 \end{pmatrix}$$

$$P = \begin{pmatrix} 0 & 1 & 1 & 1 & 1 & 1 & 1 \\ 1 & 0 & -1 & -1 & 1 & -1 & 1 \\ 1 & -1 & 0 & 1 & -1 & 1 & 1 \\ 1 & -1 & 1 & 0 & -1 & 1 & -1 \\ 1 & 1 & -1 & -1 & 0 & 1 & 1 \\ 1 & -1 & 1 & 1 & 1 & 0 & 1 \\ 1 & 1 & 1 & -1 & 1 & 1 & 0 \end{pmatrix} \quad R = \begin{pmatrix} 0 & 1 & 4 & 1 & 1 & 4 & 1 \\ 2 & 0 & 2 & 6 & 3 & 3 & 2 \\ 0 & 0 & 0 & 0 & 0 & 0 & 0 \\ 0 & 0 & 0 & 0 & 0 & 0 & 0 \\ 0 & 0 & 0 & 0 & 0 & 0 & 0 \\ 0 & 0 & 0 & 0 & 0 & 0 & 0 \\ 5 & 3 & 5 & 3 & 4 & 4 & 0 \end{pmatrix}$$

$I3(P, R) = 5.5714285714$

$I3(-P, R) = 16$

$|P| = -2$

$$\text{eigenvals}(P) = \begin{pmatrix} 3.4939592074 \\ 0.1099162642 \\ -2.6038754716 \\ -2 \\ 3.4939592074 \\ 0.1099162642 \\ -2.6038754716 \end{pmatrix}$$

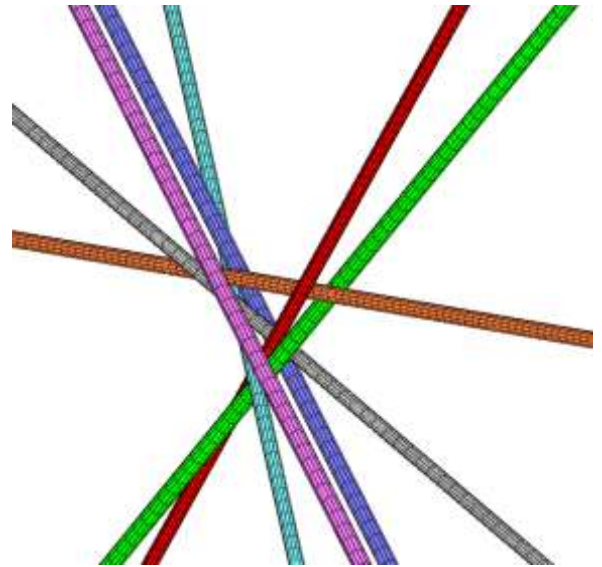

(xx,yy,zz),(xx1,yy1,zz1),(xx2,yy2,zz2),(xx3,yy3,zz3),(xx4,yy4,zz4),(xx5,yy5,zz5),(xx6,yy6,zz6)

Appendix 4

7-cross

Configuration label m162i

$$\begin{pmatrix} \text{orange} \\ \text{red} \\ \text{blue} \\ \text{green} \\ \text{cyan} \\ \text{magenta} \\ \text{gray} \end{pmatrix} \begin{pmatrix} 0 & 0 & 0 & 1 \\ t_1 & p_1 & z_1 & r_1 \\ t_2 & p_2 & z_2 & r_2 \\ t_3 & p_3 & z_3 & r_3 \\ t_4 & p_4 & z_4 & r_4 \\ t_5 & p_5 & z_5 & r_5 \\ t_6 & p_6 & z_6 & r_6 \end{pmatrix} = \begin{pmatrix} 0 & 0 & 0 & 1 \\ 0.3747323784 & 3.1415926536 & 0 & 0.7366370569 \\ 0.6854821951 & 2.4218241763 & 3.9974032529 & 5.3744730944 \\ 2.7102939927 & 0.2703059882 & -0.0135384855 & 1.0577819073 \\ 0.0745075176 & 0.3708281152 & -205.8554967521 & 0.666770171 \\ 0.441159476 & 3.2477176986 & 1.8977838095 & 0.1141904885 \\ 1.5169716323 & 0.4237811037 & -11.586361368 & 0.384131353 \end{pmatrix}$$

$$P = \begin{pmatrix} 0 & 1 & 1 & 1 & 1 & 1 & 1 \\ 1 & 0 & 1 & 1 & -1 & -1 & 1 \\ 1 & 1 & 0 & -1 & -1 & 1 & -1 \\ 1 & 1 & -1 & 0 & 1 & -1 & 1 \\ 1 & -1 & -1 & 1 & 0 & -1 & -1 \\ 1 & -1 & 1 & -1 & -1 & 0 & 1 \\ 1 & 1 & -1 & 1 & -1 & 1 & 0 \end{pmatrix} \qquad R = \begin{pmatrix} 0 & 0 & 0 & 0 & 0 & 0 & 0 \\ 4 & 0 & 1 & 1 & 4 & 1 & 1 \\ 0 & 0 & 0 & 0 & 0 & 0 & 0 \\ 0 & 0 & 0 & 0 & 0 & 0 & 0 \\ 0 & 0 & 0 & 0 & 0 & 0 & 0 \\ 3 & 4 & 5 & 3 & 5 & 0 & 4 \\ 4 & 3 & 5 & 5 & 4 & 3 & 0 \end{pmatrix}$$

$|P| = -162$

$I3(P, R) = 14.380952381$

$I3(-P, R) = 11.4285714286$

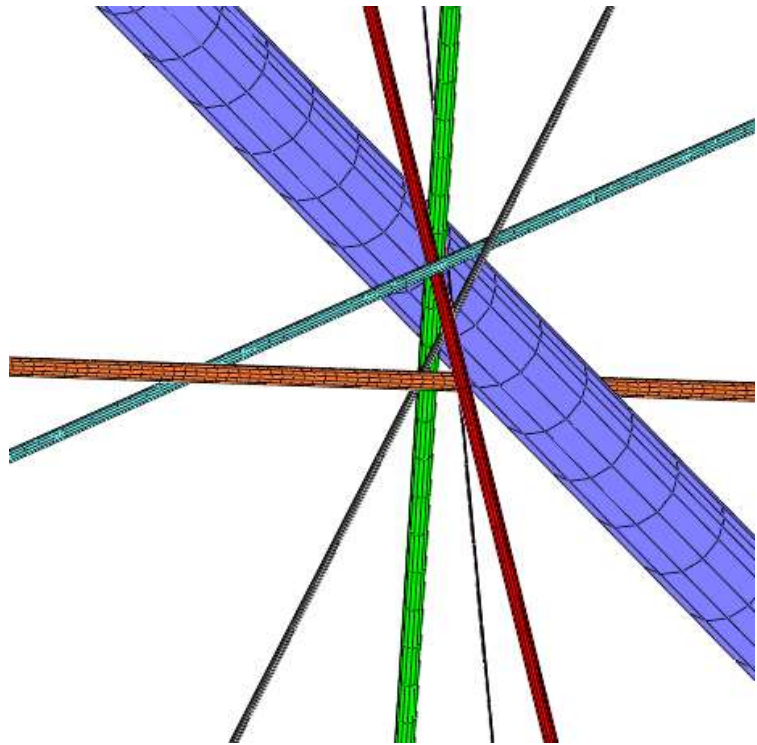

Appendix 4
7-cross
Configuration label m162b

$$\begin{pmatrix} \text{orange} \\ \text{red} \\ \text{blue} \\ \text{green} \\ \text{cyan} \\ \text{magenta} \\ \text{gray} \end{pmatrix} \begin{pmatrix} 0 & 0 & 0 & 1 \\ t_1 & p_1 & z_1 & r_1 \\ t_2 & p_2 & z_2 & r_2 \\ t_3 & p_3 & z_3 & r_3 \\ t_4 & p_4 & z_4 & r_4 \\ t_5 & p_5 & z_5 & r_5 \\ t_6 & p_6 & z_6 & r_6 \end{pmatrix} = \begin{pmatrix} 0 & 0 & 0 & 1 \\ 1.6434155029 & 3.1415926536 & 0 & 1.1577231217 \\ 1.0998072702 & -2.6044860781 & -3.941098406 & 1.4475439733 \\ 2.3800715444 & -2.6289055561 & -5.224973938 & 9.9999999997 \times 10^{-3} \\ 3.6816679207 & 2.0368644543 & -5.7315520185 & 0.0433739843 \\ -0.14761465 & 8.9144091378 & -2.6023755826 & 0.1696781504 \\ 1.6089944073 & 3.1242334187 & -4.8820649577 & 0.01 \end{pmatrix}$$

$$P = \begin{pmatrix} 0 & 1 & 1 & 1 & -1 & -1 & 1 \\ 1 & 0 & 1 & 1 & 1 & 1 & -1 \\ 1 & 1 & 0 & -1 & -1 & 1 & -1 \\ 1 & 1 & -1 & 0 & -1 & 1 & 1 \\ -1 & 1 & -1 & -1 & 0 & -1 & -1 \\ -1 & 1 & 1 & 1 & -1 & 0 & -1 \\ 1 & -1 & -1 & 1 & -1 & -1 & 0 \end{pmatrix} \qquad R = \begin{pmatrix} 0 & 0 & 0 & 0 & 0 & 0 & 0 \\ 0 & 0 & 0 & 0 & 0 & 0 & 0 \\ 1 & 4 & 0 & 1 & 1 & 1 & 4 \\ 1 & 1 & 4 & 0 & 1 & 4 & 1 \\ 1 & 4 & 1 & 4 & 0 & 1 & 1 \\ 3 & 2 & 2 & 2 & 6 & 0 & 3 \\ 4 & 1 & 1 & 4 & 1 & 1 & 0 \end{pmatrix}$$

$|P| = -162$

$I3(P, R) = 12.557319224$

$I3(-P, R) = 15.5291005291$

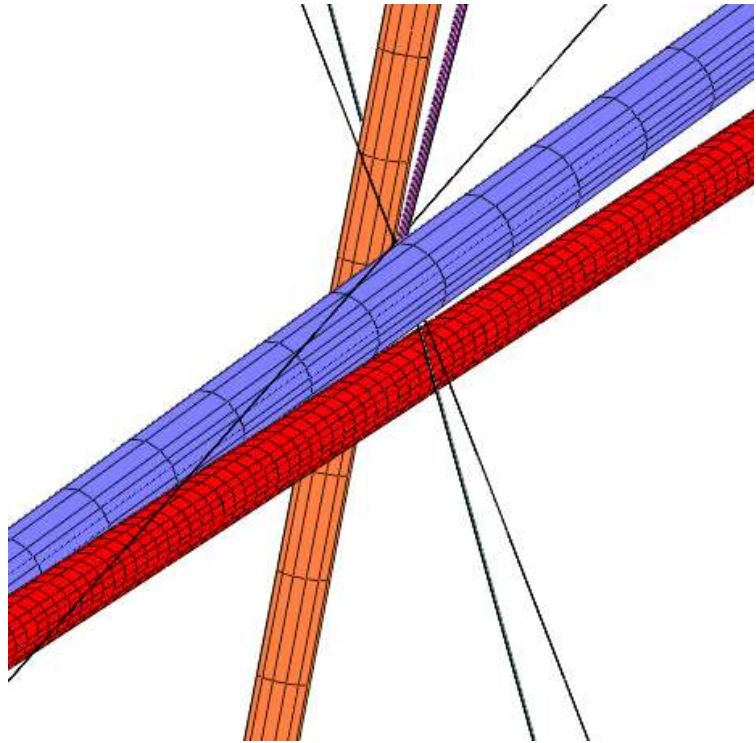

Appendix 4
7-cross
Configuration label m150q

$$\begin{pmatrix} \text{orange} \\ \text{red} \\ \text{blue} \\ \text{green} \\ \text{cyan} \\ \text{magenta} \\ \text{gray} \end{pmatrix} \begin{pmatrix} 0 & 0 & 0 & 1 \\ t1 & p1 & z1 & r1 \\ t2 & p2 & z2 & r2 \\ t3 & p3 & z3 & r3 \\ t4 & p4 & z4 & r4 \\ t5 & p5 & z5 & r5 \\ t6 & p6 & z6 & r6 \end{pmatrix} = \begin{pmatrix} 0 & 0 & 0 & 1 \\ 1.032025825 & 3.1415926536 & 0 & 0.0100556308 \\ 1.4184546063 & 4.2805231873 & -8.1242670732 & 14.7519227969 \\ 2.7611320201 & -3.1349756304 & -9.4227998077 & 0.0472444269 \\ -0.890264654 & 7.5670921297 & -0.8465107618 & 0.6501925381 \\ 0.1375236628 & 3.2300598166 & 2.2536732163 & 1.1203876658 \\ 1.7773708562 & 1.1008796434 & -5.5388309667 & 2.8271023019 \end{pmatrix}$$

$$P = \begin{pmatrix} 0 & 1 & 1 & 1 & -1 & 1 & 1 \\ 1 & 0 & 1 & 1 & 1 & -1 & -1 \\ 1 & 1 & 0 & -1 & 1 & 1 & -1 \\ 1 & 1 & -1 & 0 & -1 & -1 & -1 \\ -1 & 1 & 1 & -1 & 0 & -1 & 1 \\ 1 & -1 & 1 & -1 & -1 & 0 & 1 \\ 1 & -1 & -1 & -1 & 1 & 1 & 0 \end{pmatrix} \qquad R = \begin{pmatrix} 0 & 1 & 4 & 1 & 1 & 1 & 4 \\ 4 & 0 & 1 & 1 & 1 & 4 & 1 \\ 0 & 0 & 0 & 0 & 0 & 0 & 0 \\ 2 & 2 & 2 & 0 & 3 & 3 & 6 \\ 2 & 6 & 2 & 3 & 0 & 2 & 3 \\ 1 & 1 & 4 & 1 & 4 & 0 & 1 \\ 0 & 0 & 0 & 0 & 0 & 0 & 0 \end{pmatrix}$$

$|P| = -150$

$I3(P, R) = 11.5339412361$

$I3(-P, R) = 14.5997973658$

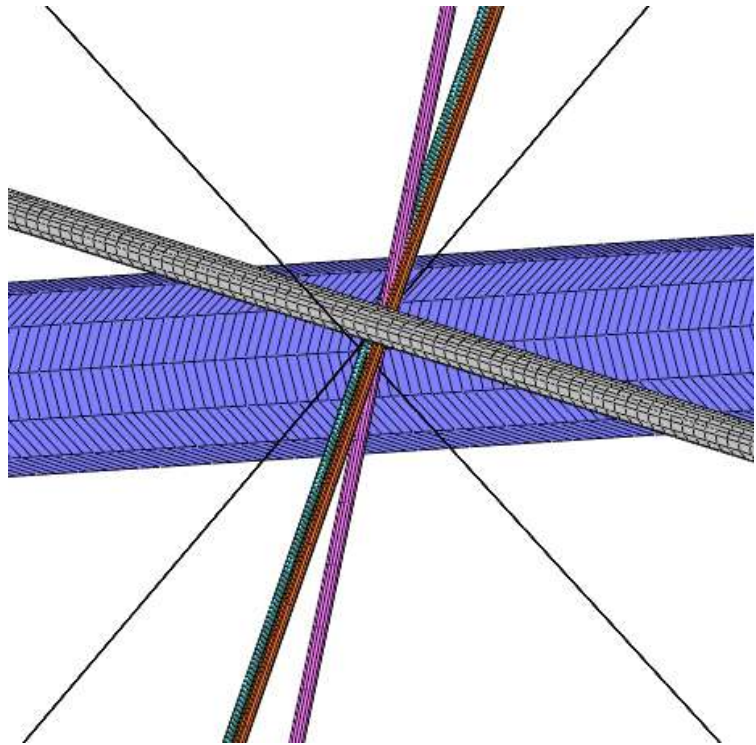

Appendix 4
7-cross
Configuration label a48m150

$$\begin{pmatrix} \text{orange} \\ \text{red} \\ \text{blue} \\ \text{green} \\ \text{cyan} \\ \text{magenta} \\ \text{gray} \end{pmatrix} \begin{pmatrix} 0 & 0 & 0 & 1 \\ t_1 & p_1 & z_1 & r_1 \\ t_2 & p_2 & z_2 & r_2 \\ t_3 & p_3 & z_3 & r_3 \\ t_4 & p_4 & z_4 & r_4 \\ t_5 & p_5 & z_5 & r_5 \\ t_6 & p_6 & z_6 & r_6 \end{pmatrix} = \begin{pmatrix} 0 & 0 & 0 & 1 \\ 0.7602445691 & 3.1415926536 & 0 & 1 \\ 2.0136844516 & 8.6074265279 & 4.055328572 & 1.6 \\ 1.0808245766 & 7.4891249587 & -0.6926447954 & 2 \\ 0.1913417177 & 7.7133032855 & -21.9787181008 & 0.531304411 \\ 2.4120873578 & 4.1256888405 & -5.0823953442 & 1 \\ 1.873884949 & 2.3624164662 & -4.7950350545 & 1 \end{pmatrix}$$

$$P = \begin{pmatrix} 0 & 1 & 1 & 1 & 1 & 1 & 1 \\ 1 & 0 & 1 & 1 & -1 & 1 & -1 \\ 1 & 1 & 0 & -1 & -1 & 1 & 1 \\ 1 & 1 & -1 & 0 & 1 & -1 & 1 \\ 1 & -1 & -1 & 1 & 0 & 1 & -1 \\ 1 & 1 & 1 & -1 & 1 & 0 & -1 \\ 1 & -1 & 1 & 1 & -1 & -1 & 0 \end{pmatrix} \qquad R = \begin{pmatrix} 0 & 0 & 0 & 0 & 0 & 0 & 0 \\ 5 & 0 & 4 & 5 & 3 & 3 & 4 \\ 0 & 0 & 0 & 0 & 0 & 0 & 0 \\ 0 & 0 & 0 & 0 & 0 & 0 & 0 \\ 1 & 1 & 1 & 4 & 0 & 4 & 1 \\ 1 & 1 & 4 & 1 & 1 & 0 & 4 \\ 4 & 1 & 1 & 1 & 4 & 1 & 0 \end{pmatrix}$$

$|P| = -150$

$I_3(P,R) = 19.2857142857$

$I_3(-P,R) = 11.6842105263$

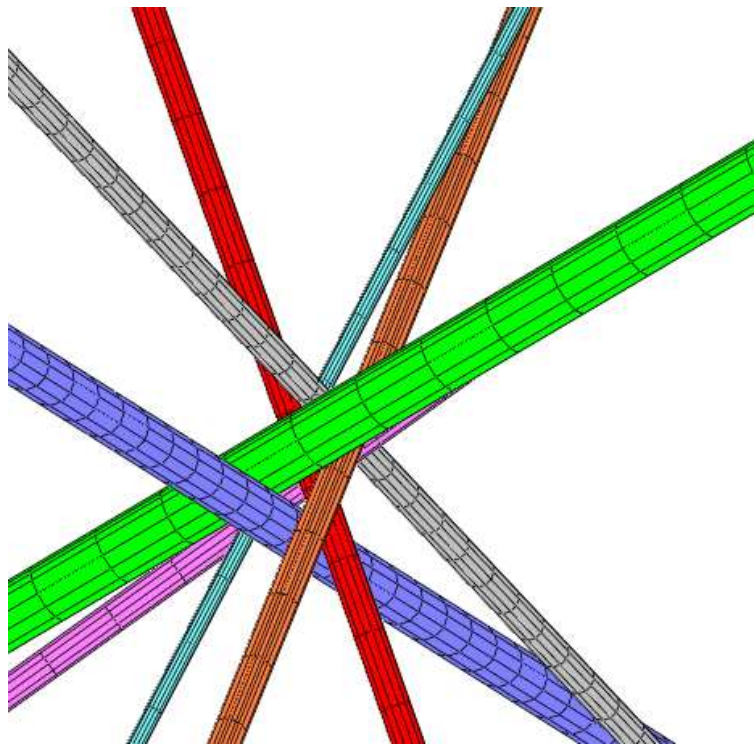

Appendix 4
7-cross
Configuration label b48

$$\begin{pmatrix} \text{orange} \\ \text{red} \\ \text{blue} \\ \text{green} \\ \text{cyan} \\ \text{magenta} \\ \text{gray} \end{pmatrix} \begin{pmatrix} 0 & 0 & 0 & 1 \\ t_1 & p_1 & z_1 & r_1 \\ t_2 & p_2 & z_2 & r_2 \\ t_3 & p_3 & z_3 & r_3 \\ t_4 & p_4 & z_4 & r_4 \\ t_5 & p_5 & z_5 & r_5 \\ t_6 & p_6 & z_6 & r_6 \end{pmatrix} = \begin{pmatrix} 0 & 0 & 0 & 1 \\ 0.1340386074 & 3.1415926536 & 0 & 1.9968460919 \\ 0.9885304727 & 2.6525356048 & 52.418449272 & 5.2372850077 \\ 2.8887794879 & 1.2007796362 & 19.9435158875 & 6.5161684705 \\ 0.0757583182 & 1.3307300075 & -104.3899059477 & 0.4752981793 \\ 2.8781560166 & 4.1485778781 & -52.8146593897 & 20.552226246 \\ 2.9356002757 & 3.083593533 & -83.933833201 & 0.1986672299 \end{pmatrix}$$

$$P = \begin{pmatrix} 0 & 1 & 1 & 1 & 1 & 1 & 1 \\ 1 & 0 & 1 & 1 & -1 & 1 & -1 \\ 1 & 1 & 0 & -1 & -1 & 1 & 1 \\ 1 & 1 & -1 & 0 & 1 & -1 & 1 \\ 1 & -1 & -1 & 1 & 0 & 1 & -1 \\ 1 & 1 & 1 & -1 & 1 & 0 & -1 \\ 1 & -1 & 1 & 1 & -1 & -1 & 0 \end{pmatrix} \qquad R = \begin{pmatrix} 0 & 0 & 0 & 0 & 0 & 0 & 0 \\ 4 & 0 & 4 & 4 & 4 & 4 & 4 \\ 0 & 0 & 0 & 0 & 0 & 0 & 0 \\ 1 & 1 & 4 & 0 & 1 & 1 & 4 \\ 1 & 1 & 1 & 4 & 0 & 4 & 1 \\ 0 & 0 & 0 & 0 & 0 & 0 & 0 \\ 4 & 1 & 1 & 1 & 4 & 1 & 0 \end{pmatrix}$$

$|P| = -150$

$I3(P, R) = 18.2763157895$

$I3(-P, R) = 11.7828947368$

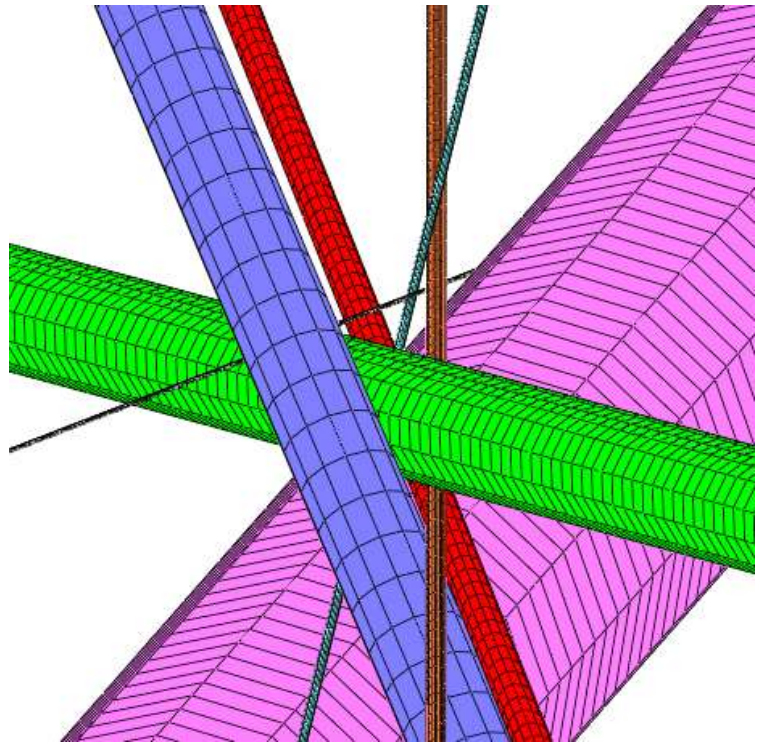

Appendix 4
7-cross
Configuration label mb28

$$\begin{pmatrix} \text{orange} \\ \text{red} \\ \text{blue} \\ \text{green} \\ \text{cyan} \\ \text{magenta} \\ \text{gray} \end{pmatrix} \begin{pmatrix} 0 & 0 & 0 & 1 \\ t_1 & p_1 & z_1 & r_1 \\ t_2 & p_2 & z_2 & r_2 \\ t_3 & p_3 & z_3 & r_3 \\ t_4 & p_4 & z_4 & r_4 \\ t_5 & p_5 & z_5 & r_5 \\ t_6 & p_6 & z_6 & r_6 \end{pmatrix} = \begin{pmatrix} 0 & 0 & 0 & 1 \\ 2.1530621814 & 3.1415926536 & 0 & 5.2372850077 \\ 2.9435681985 & 5.5632115758 & -157.4147284686 & 3.6442883187 \\ 0.2528131657 & 4.5933486221 & -32.4749333887 & 6.5161684705 \\ 3.0658343354 & 4.4633982508 & -156.8083552138 & 0.4752981793 \\ 0.263436637 & 1.6455503803 & -105.2331086578 & 20.552226246 \\ 0.2059923779 & 2.7105347254 & -136.3522824685 & 0.1986672299 \end{pmatrix}$$

$$P = \begin{pmatrix} 0 & 1 & 1 & 1 & 1 & 1 & 1 \\ 1 & 0 & 1 & 1 & 1 & -1 & -1 \\ 1 & 1 & 0 & 1 & -1 & -1 & 1 \\ 1 & 1 & 1 & 0 & -1 & 1 & -1 \\ 1 & 1 & -1 & -1 & 0 & -1 & 1 \\ 1 & -1 & -1 & 1 & -1 & 0 & 1 \\ 1 & -1 & 1 & -1 & 1 & 1 & 0 \end{pmatrix} \qquad R = \begin{pmatrix} 0 & 1 & 4 & 1 & 1 & 4 & 1 \\ 0 & 0 & 0 & 0 & 0 & 0 & 0 \\ 0 & 0 & 0 & 0 & 0 & 0 & 0 \\ 2 & 6 & 3 & 0 & 2 & 2 & 3 \\ 1 & 1 & 1 & 4 & 0 & 4 & 1 \\ 0 & 0 & 0 & 0 & 0 & 0 & 0 \\ 4 & 3 & 5 & 4 & 5 & 3 & 0 \end{pmatrix}$$

$|P| = -150$

$I3(P, R) = 25.9052132701$

$I3(-P, R) = 4.0568720379$

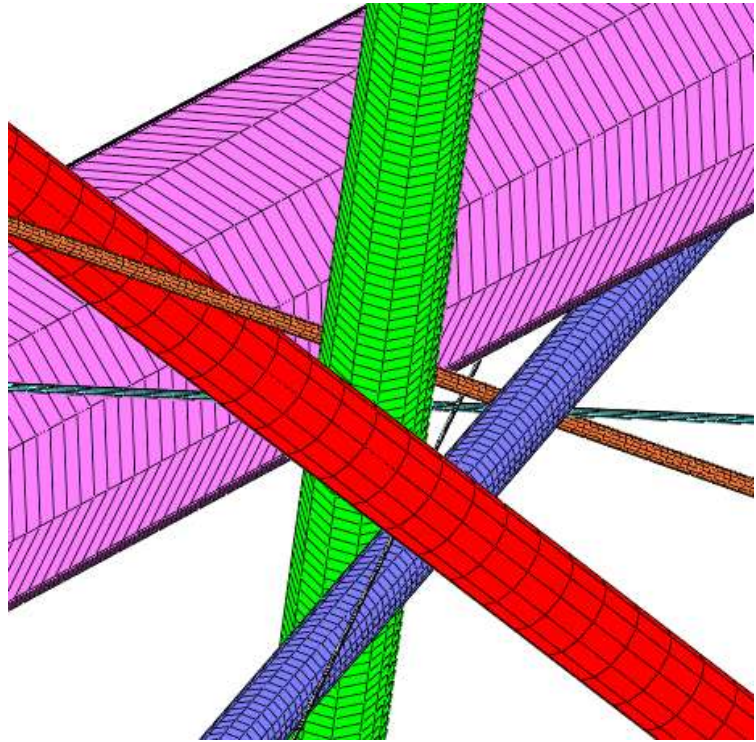

Appendix 4
7-cross
Configuration label b24

$$\begin{pmatrix} \text{orange} \\ \text{red} \\ \text{blue} \\ \text{green} \\ \text{cyan} \\ \text{magenta} \\ \text{gray} \end{pmatrix} \begin{pmatrix} 0 & 0 & 0 & 1 \\ t_1 & p_1 & z_1 & r_1 \\ t_2 & p_2 & z_2 & r_2 \\ t_3 & p_3 & z_3 & r_3 \\ t_4 & p_4 & z_4 & r_4 \\ t_5 & p_5 & z_5 & r_5 \\ t_6 & p_6 & z_6 & r_6 \end{pmatrix} = \begin{pmatrix} 0 & 0 & 0 & 1 \\ 0.9885304708 & 3.1415926536 & 0 & 5.2372850077 \\ 2.8887794877 & 1.6898366853 & -32.4749334028 & 6.5161684705 \\ 0.0757583182 & 1.8197870566 & -156.8083552138 & 0.4752981793 \\ 2.8781560166 & 4.6376349272 & -105.2331086586 & 20.552226246 \\ 3.0594316634 & 4.5444638649 & -142.1278016505 & 0.1986672299 \\ 2.9356002758 & 3.572650582 & -136.3522824705 & 0.1986672299 \end{pmatrix}$$

$$P = \begin{pmatrix} 0 & 1 & 1 & 1 & 1 & 1 & 1 \\ 1 & 0 & -1 & -1 & 1 & 1 & 1 \\ 1 & -1 & 0 & 1 & -1 & -1 & 1 \\ 1 & -1 & 1 & 0 & 1 & -1 & -1 \\ 1 & 1 & -1 & 1 & 0 & -1 & -1 \\ 1 & 1 & -1 & -1 & -1 & 0 & -1 \\ 1 & 1 & 1 & -1 & -1 & -1 & 0 \end{pmatrix} \qquad R = \begin{pmatrix} 0 & 0 & 0 & 0 & 0 & 0 & 0 \\ 0 & 0 & 0 & 0 & 0 & 0 & 0 \\ 1 & 4 & 0 & 1 & 1 & 1 & 4 \\ 1 & 1 & 4 & 0 & 4 & 1 & 1 \\ 0 & 0 & 0 & 0 & 0 & 0 & 0 \\ 4 & 5 & 3 & 4 & 5 & 0 & 3 \\ 4 & 1 & 1 & 4 & 1 & 1 & 0 \end{pmatrix}$$

$|P| = -150$

$I3(P, R) = 19.4962406015$

$I3(-P, R) = 11.4736842105$

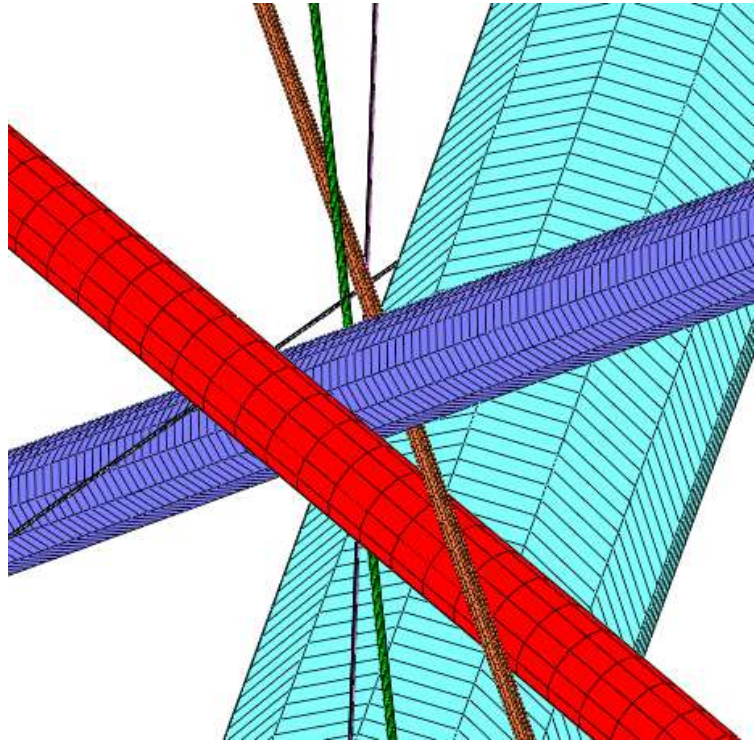

Appendix 4
7-cross
Configuration label m102n

$$\begin{pmatrix} \text{orange} \\ \text{red} \\ \text{blue} \\ \text{green} \\ \text{cyan} \\ \text{magenta} \\ \text{gray} \end{pmatrix} \begin{pmatrix} 0 & 0 & 0 & 1 \\ t1 & p1 & z1 & r1 \\ t2 & p2 & z2 & r2 \\ t3 & p3 & z3 & r3 \\ t4 & p4 & z4 & r4 \\ t5 & p5 & z5 & r5 \\ t6 & p6 & z6 & r6 \end{pmatrix} = \begin{pmatrix} 0 & 0 & 0 & 1 \\ 1.9268904145 & 3.1415926536 & 0 & 0.0291671339 \\ 2.6592769513 & 1.0716254986 & -8.5674790027 & 10.8656776836 \\ 0.2215826503 & 9.3975988661 & -10.3823963985 & 0.0102457996 \\ 2.3113253709 & 4.8684991633 & -6.7313879026 & 2.6150485477 \\ 3.4040128732 & -2.4359532899 & -3.9310417363 & 0.4604361542 \\ 2.5005016955 & 3.1380446509 & -5.4541921443 & 0.01 \end{pmatrix}$$

$$P = \begin{pmatrix} 0 & 1 & 1 & 1 & 1 & -1 & 1 \\ 1 & 0 & -1 & -1 & 1 & -1 & -1 \\ 1 & -1 & 0 & 1 & -1 & -1 & -1 \\ 1 & -1 & 1 & 0 & 1 & 1 & -1 \\ 1 & 1 & -1 & 1 & 0 & -1 & -1 \\ -1 & -1 & -1 & 1 & -1 & 0 & 1 \\ 1 & -1 & -1 & -1 & -1 & 1 & 0 \end{pmatrix} \quad R = \begin{pmatrix} 0 & 1 & 4 & 1 & 4 & 1 & 1 \\ 0 & 0 & 0 & 0 & 0 & 0 & 0 \\ 0 & 0 & 0 & 0 & 0 & 0 & 0 \\ 1 & 1 & 1 & 0 & 4 & 4 & 1 \\ 0 & 0 & 0 & 0 & 0 & 0 & 0 \\ 2 & 6 & 2 & 3 & 3 & 0 & 2 \\ 5 & 3 & 5 & 4 & 4 & 3 & 0 \end{pmatrix}$$

$|P| = -102$

$I3(P, R) = 16.2112676056$

$I3(-P, R) = 14$

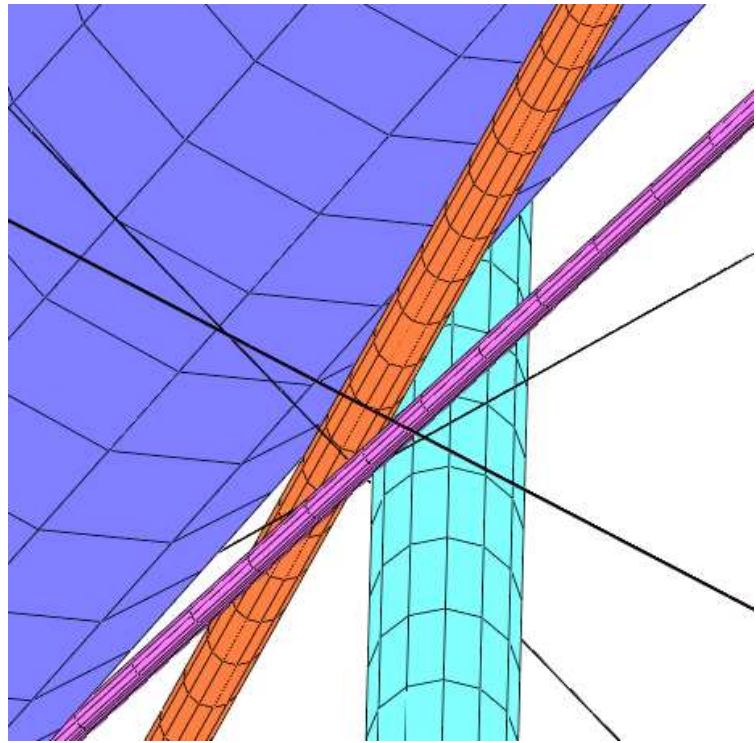

Appendix 4
7-cross
Configuration label m102c

$$\begin{pmatrix} \text{orange} \\ \text{red} \\ \text{blue} \\ \text{green} \\ \text{cyan} \\ \text{magenta} \\ \text{gray} \end{pmatrix} \begin{pmatrix} 0 & 0 & 0 & 1 \\ t_1 & p_1 & z_1 & r_1 \\ t_2 & p_2 & z_2 & r_2 \\ t_3 & p_3 & z_3 & r_3 \\ t_4 & p_4 & z_4 & r_4 \\ t_5 & p_5 & z_5 & r_5 \\ t_6 & p_6 & z_6 & r_6 \end{pmatrix} = \begin{pmatrix} 0 & 0 & 0 & 1 \\ 1.8432170544 & 3.1415926536 & 0 & 0.2318997748 \\ 1.1868288934 & 9.0822390103 & 2.0822980035 & 0.4861404601 \\ 1.3554894048 & 7.4773812744 & 0.5233414308 & 0.1750557856 \\ 0.0988620609 & 7.8069460197 & -18.5869081534 & 0.3509112444 \\ 2.741864816 & 3.7861897635 & -5.0378086436 & 2.3892037474 \\ -0.2651446037 & 1.8445521938 & 1.9460409559 & 0.2210707348 \end{pmatrix}$$

$$P = \begin{pmatrix} 0 & 1 & 1 & 1 & 1 & 1 & -1 \\ 1 & 0 & 1 & 1 & -1 & 1 & 1 \\ 1 & 1 & 0 & -1 & -1 & 1 & 1 \\ 1 & 1 & -1 & 0 & 1 & -1 & 1 \\ 1 & -1 & -1 & 1 & 0 & 1 & 1 \\ 1 & 1 & 1 & -1 & 1 & 0 & 1 \\ -1 & 1 & 1 & 1 & 1 & 1 & 0 \end{pmatrix} \qquad R = \begin{pmatrix} 0 & 0 & 0 & 0 & 0 & 0 & 0 \\ 4 & 0 & 4 & 4 & 4 & 4 & 4 \\ 0 & 0 & 0 & 0 & 0 & 0 & 0 \\ 1 & 4 & 4 & 0 & 1 & 1 & 1 \\ 1 & 1 & 1 & 4 & 0 & 4 & 1 \\ 0 & 0 & 0 & 0 & 0 & 0 & 0 \\ 6 & 2 & 2 & 3 & 3 & 2 & 0 \end{pmatrix}$$

$|P| = -102$

$I3(P, R) = 15.2205882353$

$I3(-P, R) = 15.9705882353$

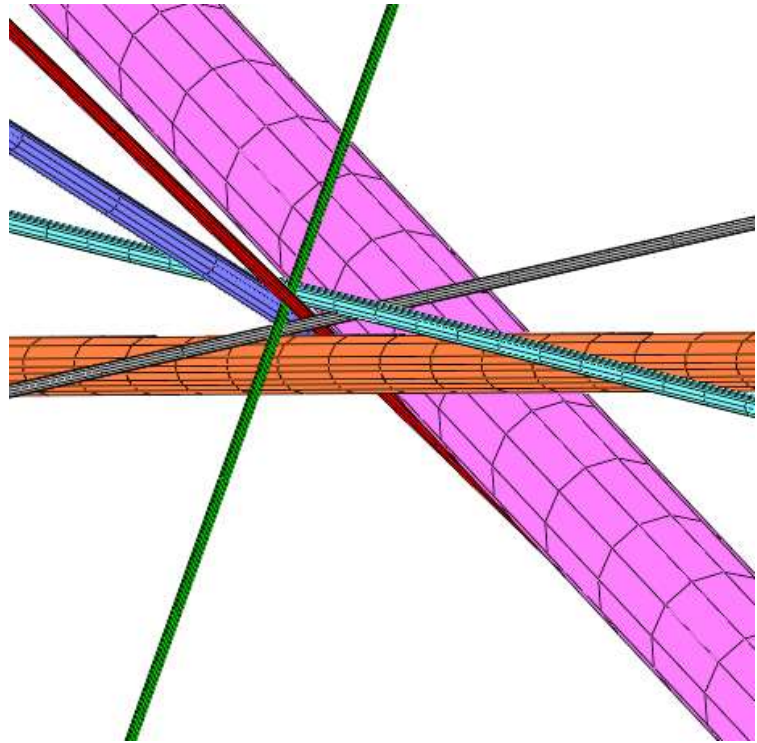

Appendix 4
7-cross
Configuration label m102b

$$\begin{pmatrix} \text{orange} \\ \text{red} \\ \text{blue} \\ \text{green} \\ \text{cyan} \\ \text{magenta} \\ \text{gray} \end{pmatrix} \begin{pmatrix} 0 & 0 & 0 & 1 \\ t_1 & p_1 & z_1 & r_1 \\ t_2 & p_2 & z_2 & r_2 \\ t_3 & p_3 & z_3 & r_3 \\ t_4 & p_4 & z_4 & r_4 \\ t_5 & p_5 & z_5 & r_5 \\ t_6 & p_6 & z_6 & r_6 \end{pmatrix} = \begin{pmatrix} 0 & 0 & 0 & 1 \\ 3.6746091143 & 3.1415926536 & 0 & 1 \\ 5.7669956198 & 15.5878482319 & -4.4892624669 & 0.076893307 \\ 3.6354801145 & 14.7491473689 & -5.2010477332 & 0.0100000003 \\ 0.4333839055 & 7.1797534377 & -7.0024858141 & 0.8501191328 \\ 5.9034511356 & 4.4842320311 & -7.647559772 & 0.3865582497 \\ 1.3905843778 & 2.1684064871 & -2.8835286668 & 0.1573134123 \end{pmatrix}$$

$$P = \begin{pmatrix} 0 & -1 & -1 & -1 & 1 & -1 & 1 \\ -1 & 0 & -1 & -1 & -1 & 1 & 1 \\ -1 & -1 & 0 & -1 & 1 & -1 & -1 \\ -1 & -1 & -1 & 0 & -1 & 1 & -1 \\ 1 & -1 & 1 & -1 & 0 & 1 & -1 \\ -1 & 1 & -1 & 1 & 1 & 0 & 1 \\ 1 & 1 & -1 & -1 & -1 & 1 & 0 \end{pmatrix} \quad R = \begin{pmatrix} 0 & 1 & 1 & 1 & 4 & 4 & 1 \\ 0 & 0 & 0 & 0 & 0 & 0 & 0 \\ 4 & 4 & 0 & 5 & 5 & 3 & 3 \\ 2 & 3 & 2 & 0 & 2 & 6 & 3 \\ 0 & 0 & 0 & 0 & 0 & 0 & 0 \\ 0 & 0 & 0 & 0 & 0 & 0 & 0 \\ 2 & 6 & 3 & 2 & 3 & 2 & 0 \end{pmatrix}$$

$|P| = -102$

$I3(P, R) = 15.8468899522$

$I3(-P, R) = 11.5885167464$

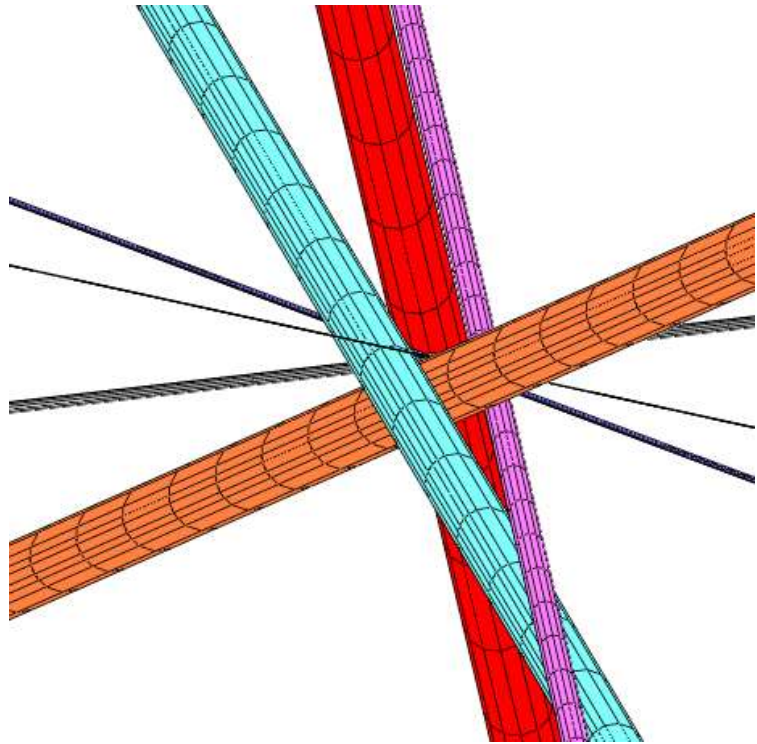

Appendix 4
7-cross
Configuration label m102a

$$\begin{pmatrix} \text{orange} \\ \text{red} \\ \text{blue} \\ \text{green} \\ \text{cyan} \\ \text{magenta} \\ \text{gray} \end{pmatrix} \begin{pmatrix} 0 & 0 & 0 & 1 \\ t_1 & p_1 & z_1 & r_1 \\ t_2 & p_2 & z_2 & r_2 \\ t_3 & p_3 & z_3 & r_3 \\ t_4 & p_4 & z_4 & r_4 \\ t_5 & p_5 & z_5 & r_5 \\ t_6 & p_6 & z_6 & r_6 \end{pmatrix} = \begin{pmatrix} 0 & 0 & 0 & 1 \\ 3.5440303969 & 3.1415926536 & 0 & 1 \\ 5.7471769155 & 15.6193882575 & -7.1037232467 & 0.076893307 \\ 3.4005236076 & 15.2554055886 & -9.814056102 & 0.0100000003 \\ 3.3054179744 & 11.0034036665 & 3.7627838767 & 1 \\ 5.6279326552 & 3.9803795596 & -8.1373582104 & 0.3865582497 \\ 2.9380161855 & 2.6580593753 & -4.2252770734 & 0.1573134123 \end{pmatrix}$$

$$P = \begin{pmatrix} 0 & -1 & -1 & -1 & -1 & -1 & 1 \\ -1 & 0 & -1 & -1 & 1 & 1 & 1 \\ -1 & -1 & 0 & -1 & -1 & -1 & -1 \\ -1 & -1 & -1 & 0 & 1 & 1 & -1 \\ -1 & 1 & -1 & 1 & 0 & -1 & 1 \\ -1 & 1 & -1 & 1 & -1 & 0 & 1 \\ 1 & 1 & -1 & -1 & 1 & 1 & 0 \end{pmatrix} \quad R = \begin{pmatrix} 0 & 0 & 0 & 0 & 0 & 0 & 0 \\ 0 & 0 & 0 & 0 & 0 & 0 & 0 \\ 4 & 1 & 0 & 4 & 1 & 1 & 1 \\ 3 & 4 & 3 & 0 & 4 & 5 & 5 \\ 0 & 0 & 0 & 0 & 0 & 0 & 0 \\ 2 & 2 & 3 & 2 & 6 & 0 & 3 \\ 1 & 4 & 4 & 1 & 1 & 1 & 0 \end{pmatrix}$$

$|P| = -102$

$I3(P, R) = 13.7023809524$

$I3(-P, R) = 17.4365079365$

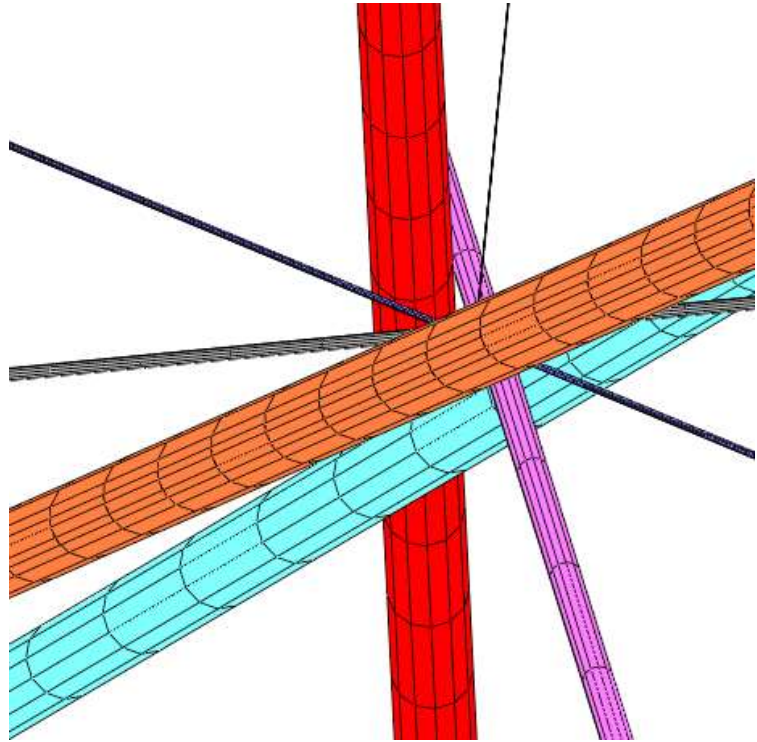

Appendix 4
7-cross
Configuration label b46

$$\begin{pmatrix} \text{orange} \\ \text{red} \\ \text{blue} \\ \text{green} \\ \text{cyan} \\ \text{magenta} \\ \text{gray} \end{pmatrix} \begin{pmatrix} 0 & 0 & 0 & 1 \\ t_1 & p_1 & z_1 & r_1 \\ t_2 & p_2 & z_2 & r_2 \\ t_3 & p_3 & z_3 & r_3 \\ t_4 & p_4 & z_4 & r_4 \\ t_5 & p_5 & z_5 & r_5 \\ t_6 & p_6 & z_6 & r_6 \end{pmatrix} = \begin{pmatrix} 0 & 0 & 0 & 1 \\ 0.1340386074 & 3.1415926536 & 0 & 1.9968460919 \\ 0.9885304709 & 2.6525356045 & 52.4184492388 & 5.2372850077 \\ 2.8887794876 & 1.2007796362 & 19.9435158464 & 6.5161684705 \\ 2.8781560166 & 4.1485778781 & -52.814659394 & 20.552226246 \\ 3.0594316634 & 4.0554068157 & -89.7093523806 & 0.1986672299 \\ 2.9356002757 & 3.083593533 & -83.933833203 & 0.1986672299 \end{pmatrix}$$

$$P = \begin{pmatrix} 0 & 1 & 1 & 1 & 1 & 1 & 1 \\ 1 & 0 & 1 & 1 & 1 & 1 & -1 \\ 1 & 1 & 0 & -1 & 1 & 1 & 1 \\ 1 & 1 & -1 & 0 & -1 & -1 & 1 \\ 1 & 1 & 1 & -1 & 0 & -1 & -1 \\ 1 & 1 & 1 & -1 & -1 & 0 & -1 \\ 1 & -1 & 1 & 1 & -1 & -1 & 0 \end{pmatrix} \qquad R = \begin{pmatrix} 0 & 0 & 0 & 0 & 0 & 0 & 0 \\ 4 & 0 & 4 & 4 & 4 & 4 & 4 \\ 0 & 0 & 0 & 0 & 0 & 0 & 0 \\ 1 & 1 & 4 & 0 & 1 & 1 & 4 \\ 0 & 0 & 0 & 0 & 0 & 0 & 0 \\ 2 & 2 & 2 & 3 & 6 & 0 & 3 \\ 0 & 0 & 0 & 0 & 0 & 0 & 0 \end{pmatrix}$$

$|P| = -102$

$I3(P, R) = 14.2941176471$

$I3(-P, R) = 10.0588235294$

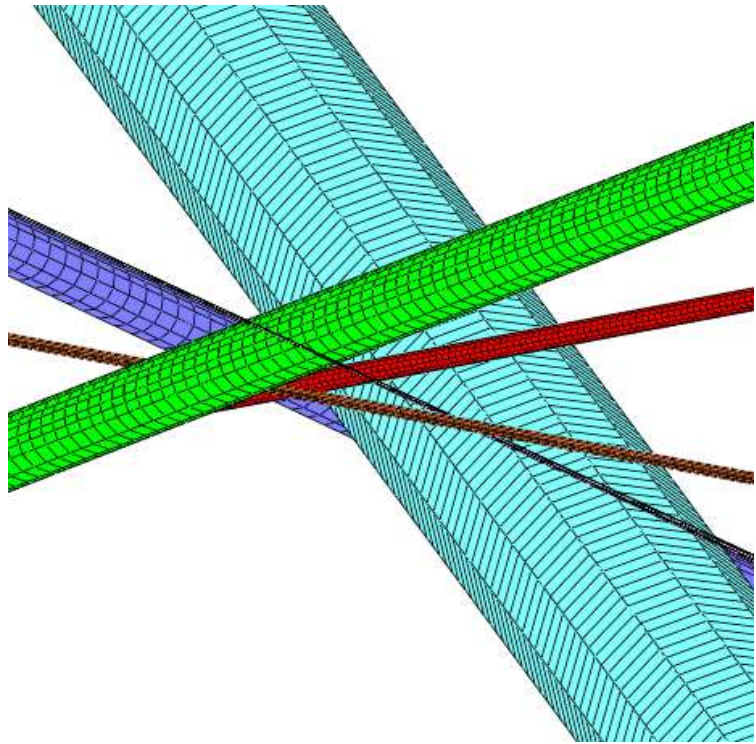

Appendix 4
7-cross
Configuration label b16

$$\begin{pmatrix} \text{orange} \\ \text{red} \\ \text{blue} \\ \text{green} \\ \text{cyan} \\ \text{magenta} \\ \text{gray} \end{pmatrix} \begin{pmatrix} 0 & 0 & 0 & 1 \\ t_1 & p_1 & z_1 & r_1 \\ t_2 & p_2 & z_2 & r_2 \\ t_3 & p_3 & z_3 & r_3 \\ t_4 & p_4 & z_4 & r_4 \\ t_5 & p_5 & z_5 & r_5 \\ t_6 & p_6 & z_6 & r_6 \end{pmatrix} = \begin{pmatrix} 0 & 0 & 0 & 1 \\ 0.8717111429 & 3.1415926536 & 0 & 2.622778505 \\ 0.329915725 & 0.677678256 & -60.7164431818 & 1.8250221354 \\ 2.9027654612 & 2.2778554462 & -34.3392956757 & 3.2632301993 \\ 2.7883138829 & 4.3693470432 & -34.9016969514 & 10.292343676 \\ 2.9462904054 & 4.0212968595 & -40.5582445842 & 0.0994905069 \\ -2.8016995043 & 0.5026702018 & -54.8335299675 & 0.0994905069 \end{pmatrix}$$

$$P = \begin{pmatrix} 0 & 1 & 1 & 1 & 1 & 1 & -1 \\ 1 & 0 & -1 & -1 & 1 & 1 & 1 \\ 1 & -1 & 0 & -1 & 1 & -1 & -1 \\ 1 & -1 & -1 & 0 & -1 & -1 & 1 \\ 1 & 1 & 1 & -1 & 0 & -1 & -1 \\ 1 & 1 & -1 & -1 & -1 & 0 & -1 \\ -1 & 1 & -1 & 1 & -1 & -1 & 0 \end{pmatrix} \qquad R = \begin{pmatrix} 0 & 4 & 4 & 4 & 4 & 4 & 4 \\ 0 & 0 & 0 & 0 & 0 & 0 & 0 \\ 0 & 0 & 0 & 0 & 0 & 0 & 0 \\ 2 & 6 & 3 & 0 & 2 & 2 & 3 \\ 0 & 0 & 0 & 0 & 0 & 0 & 0 \\ 2 & 2 & 2 & 3 & 6 & 0 & 3 \\ 3 & 2 & 6 & 3 & 2 & 2 & 0 \end{pmatrix}$$

$|P| = -102$

$I3(P,R) = 17.4590909091$

$I3(-P,R) = 12.2227272727$

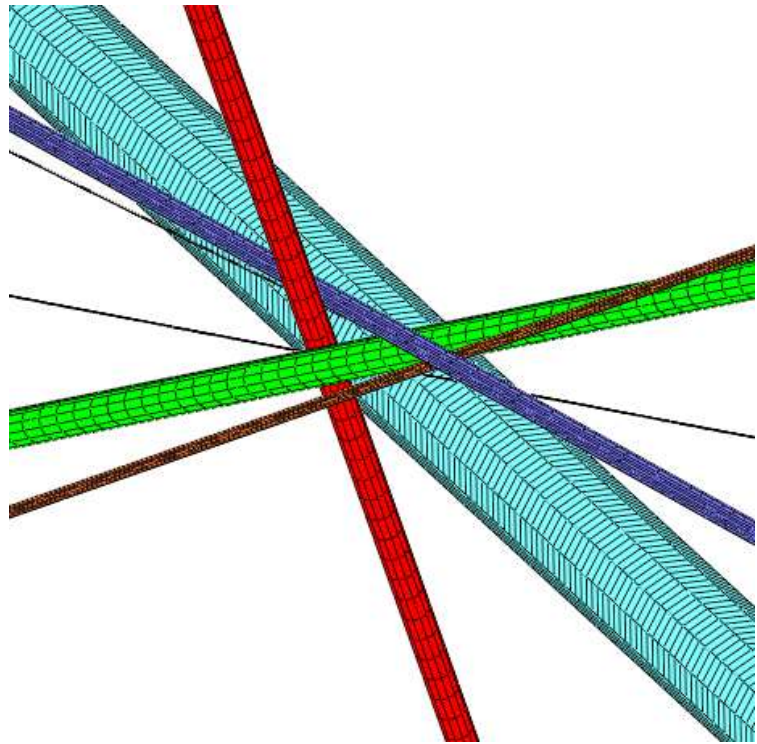

Appendix 4
7-cross
Configuration label b14

$$\begin{pmatrix} \text{orange} \\ \text{red} \\ \text{blue} \\ \text{green} \\ \text{cyan} \\ \text{magenta} \\ \text{gray} \end{pmatrix} \begin{pmatrix} 0 & 0 & 0 & 1 \\ t_1 & p_1 & z_1 & r_1 \\ t_2 & p_2 & z_2 & r_2 \\ t_3 & p_3 & z_3 & r_3 \\ t_4 & p_4 & z_4 & r_4 \\ t_5 & p_5 & z_5 & r_5 \\ t_6 & p_6 & z_6 & r_6 \end{pmatrix} = \begin{pmatrix} 0 & 0 & 0 & 1 \\ 0.8717111471 & 3.1415926536 & 0 & 2.622778505 \\ 2.9027654618 & 2.277855447 & -34.3392957179 & 3.2632301993 \\ -0.1688253538 & 4.1326812482 & -43.3563839515 & 0.2380244433 \\ 2.7883138827 & 4.3693470416 & -34.9016969488 & 10.292343676 \\ 2.9462904052 & 4.0212968579 & -40.5582445818 & 0.0994905069 \\ -2.8016995041 & 0.5026702006 & -54.8335299685 & 0.0994905069 \end{pmatrix}$$

$$P = \begin{pmatrix} 0 & 1 & 1 & -1 & 1 & 1 & -1 \\ 1 & 0 & -1 & -1 & 1 & 1 & 1 \\ 1 & -1 & 0 & 1 & -1 & -1 & 1 \\ -1 & -1 & 1 & 0 & 1 & -1 & -1 \\ 1 & 1 & -1 & 1 & 0 & -1 & -1 \\ 1 & 1 & -1 & -1 & -1 & 0 & -1 \\ -1 & 1 & 1 & -1 & -1 & -1 & 0 \end{pmatrix} \qquad R = \begin{pmatrix} 0 & 3 & 3 & 2 & 2 & 2 & 6 \\ 0 & 0 & 0 & 0 & 0 & 0 & 0 \\ 1 & 4 & 0 & 1 & 1 & 1 & 4 \\ 1 & 1 & 4 & 0 & 4 & 1 & 1 \\ 0 & 0 & 0 & 0 & 0 & 0 & 0 \\ 4 & 5 & 3 & 4 & 5 & 0 & 3 \\ 0 & 0 & 0 & 0 & 0 & 0 & 0 \end{pmatrix}$$

$|P| = -102$

$I3(P, R) = 14.7299270073$

$I3(-P, R) = 8.598540146$

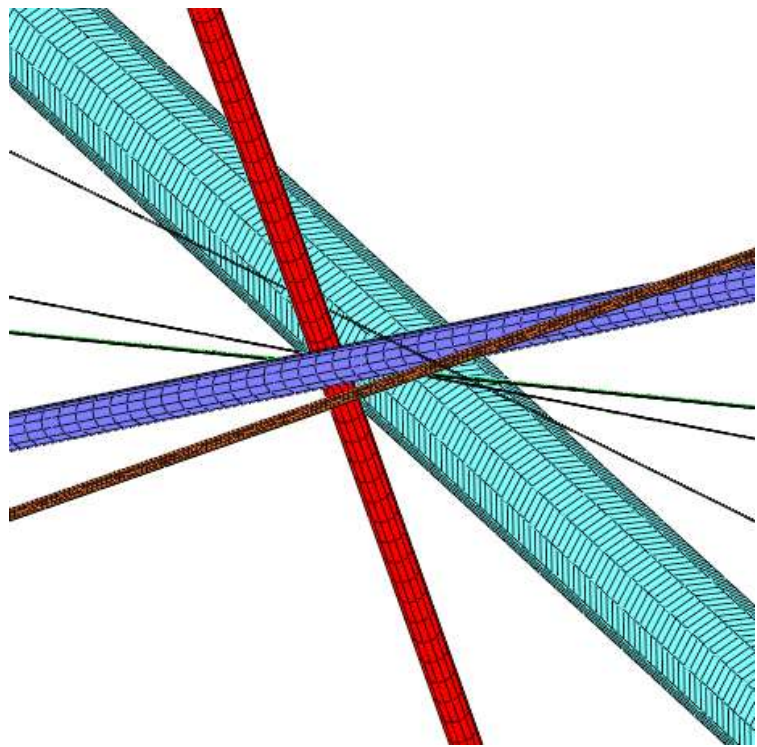

Appendix 4
7-cross
Configuration label x66e

$$\begin{pmatrix} \text{orange} \\ \text{red} \\ \text{blue} \\ \text{green} \\ \text{cyan} \\ \text{magenta} \\ \text{gray} \end{pmatrix} \begin{pmatrix} 0 & 0 & 0 & 1 \\ t1 & p1 & z1 & r1 \\ t2 & p2 & z2 & r2 \\ t3 & p3 & z3 & r3 \\ t4 & p4 & z4 & r4 \\ t5 & p5 & z5 & r5 \\ t6 & p6 & z6 & r6 \end{pmatrix} = \begin{pmatrix} 0 & 0 & 0 & 1 \\ 2.9924740175 & 3.1415926536 & 0 & 1.4786146441 \\ 2.8641954569 & 7.3554023648 & 8.3378482046 & 7.3883176609 \\ 1.8406812984 & 2.0063164249 & -10.4991508345 & 0.2212795487 \\ 1.2637301239 & 5.2421060448 & -4.3385055165 & 9.6314725076 \\ 3.387731002 & -2.6084354634 & -6.8881558748 & 0.1111150205 \\ 1.9513737611 & 3.8856365257 & -11.1174625263 & 0.222934851 \end{pmatrix}$$

$$P = \begin{pmatrix} 0 & 1 & 1 & 1 & 1 & -1 & 1 \\ 1 & 0 & -1 & -1 & -1 & -1 & 1 \\ 1 & -1 & 0 & 1 & -1 & -1 & -1 \\ 1 & -1 & 1 & 0 & 1 & -1 & -1 \\ 1 & -1 & -1 & 1 & 0 & -1 & -1 \\ -1 & -1 & -1 & -1 & -1 & 0 & -1 \\ 1 & 1 & -1 & -1 & -1 & -1 & 0 \end{pmatrix} \qquad R = \begin{pmatrix} 0 & 2 & 3 & 6 & 2 & 2 \\ 0 & 0 & 0 & 0 & 0 & 0 & 0 \\ 0 & 0 & 0 & 0 & 0 & 0 & 0 \\ 1 & 1 & 4 & 0 & 1 & 4 & 1 \\ 0 & 0 & 0 & 0 & 0 & 0 & 0 \\ 4 & 4 & 3 & 3 & 5 & 0 & 5 \\ 1 & 4 & 1 & 4 & 1 & 1 & 0 \end{pmatrix}$$

$|P| = 66$

$I3(P, R) = 15.3210526316$

$I3(-P, R) = 16.0052631579$

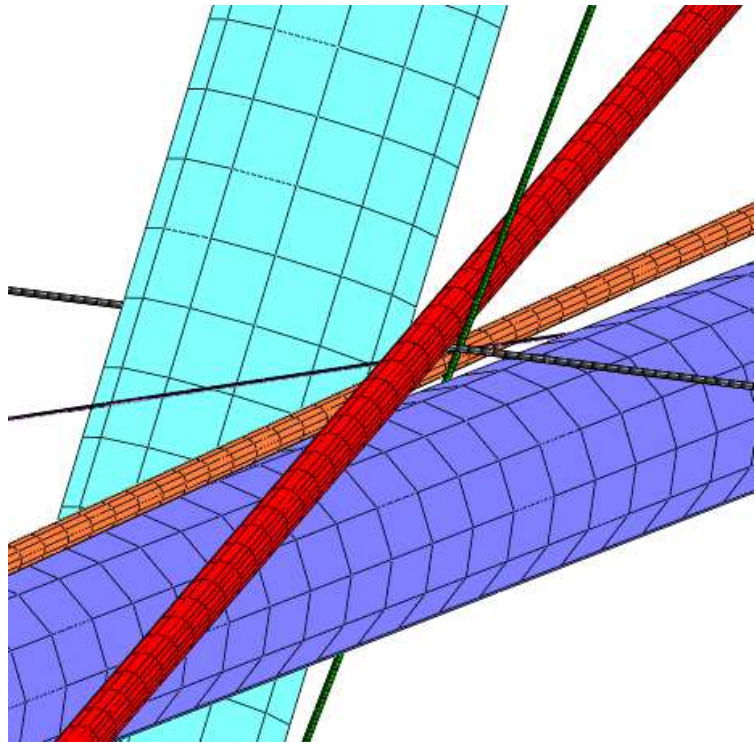

Appendix 4
7-cross
Configuration label x66a

$$\begin{pmatrix} \text{orange} \\ \text{red} \\ \text{blue} \\ \text{green} \\ \text{cyan} \\ \text{magenta} \\ \text{gray} \end{pmatrix} \begin{pmatrix} 0 & 0 & 0 & 1 \\ t_1 & p_1 & z_1 & r_1 \\ t_2 & p_2 & z_2 & r_2 \\ t_3 & p_3 & z_3 & r_3 \\ t_4 & p_4 & z_4 & r_4 \\ t_5 & p_5 & z_5 & r_5 \\ t_6 & p_6 & z_6 & r_6 \end{pmatrix} = \begin{pmatrix} 0 & 0 & 0 & 1 \\ -0.1843182352 & 3.1415926536 & 0 & 1.0189129916 \\ -0.0680640543 & 2.1012853513 & 4.4336779695 & 0.3463609234 \\ -0.1860966076 & -9.0412689873 & -8.6924684197 & 0.0100000003 \\ -2.5141106395 & -5.2220080312 & 2.2274924329 & 1.0075645913 \\ -2.7983081394 & 1.248594993 & -3.6430161967 & 0.3976220892 \\ -0.4456803977 & -2.8204972344 & -3.69006024 & 0.0461278404 \end{pmatrix}$$

$$P = \begin{pmatrix} 0 & -1 & -1 & -1 & -1 & -1 & -1 \\ -1 & 0 & -1 & 1 & -1 & -1 & 1 \\ -1 & -1 & 0 & -1 & -1 & -1 & -1 \\ -1 & 1 & -1 & 0 & -1 & -1 & 1 \\ -1 & -1 & -1 & -1 & 0 & 1 & 1 \\ -1 & -1 & -1 & -1 & 1 & 0 & 1 \\ -1 & 1 & -1 & 1 & 1 & 1 & 0 \end{pmatrix} \qquad R = \begin{pmatrix} 0 & 0 & 0 & 0 & 0 & 0 & 0 \\ 0 & 0 & 0 & 0 & 0 & 0 & 0 \\ 6 & 2 & 0 & 2 & 3 & 3 & 2 \\ 1 & 4 & 1 & 0 & 1 & 4 & 1 \\ 0 & 0 & 0 & 0 & 0 & 0 & 0 \\ 1 & 1 & 1 & 1 & 4 & 0 & 4 \\ 3 & 2 & 3 & 6 & 2 & 2 & 0 \end{pmatrix}$$

$|P| = 66$

$I3(P, R) = 11.2810810811$

$I3(-P, R) = 21.6324324324$

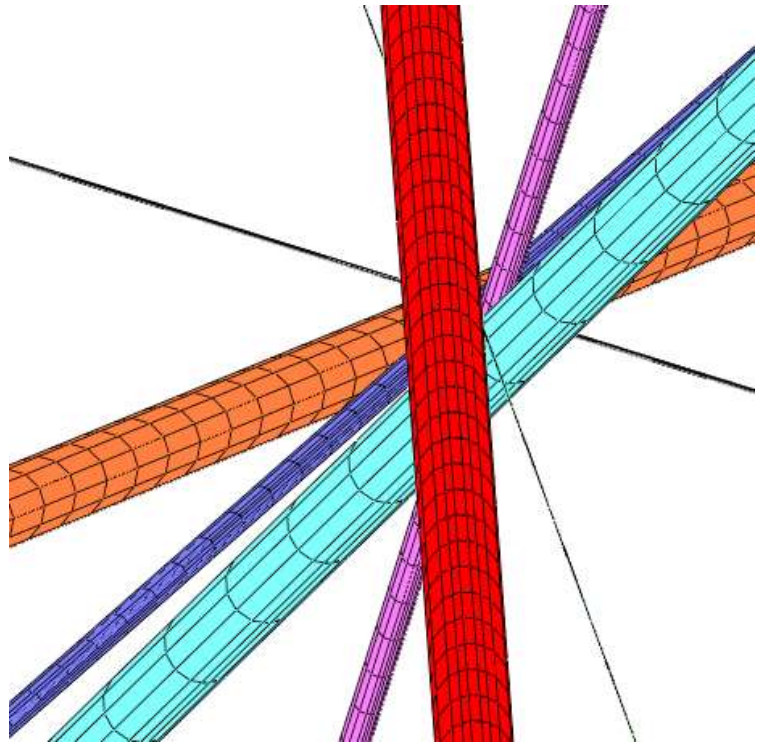

Appendix 4
7-cross
Configuration label a2466

$$\begin{pmatrix} \text{orange} \\ \text{red} \\ \text{blue} \\ \text{green} \\ \text{cyan} \\ \text{magenta} \\ \text{gray} \end{pmatrix} \begin{pmatrix} 0 & 0 & 0 & 1 \\ t_1 & p_1 & z_1 & r_1 \\ t_2 & p_2 & z_2 & r_2 \\ t_3 & p_3 & z_3 & r_3 \\ t_4 & p_4 & z_4 & r_4 \\ t_5 & p_5 & z_5 & r_5 \\ t_6 & p_6 & z_6 & r_6 \end{pmatrix} = \begin{pmatrix} 0 & 0 & 0 & 1 \\ -0.1854553046 & 3.1415926536 & 0 & 1 \\ -0.4244741528 & -9.159486875 & -4.1448223267 & 0.076893307 \\ -0.1859509103 & -9.0405410392 & -8.6489347965 & 0.0100000003 \\ -2.517183213 & -5.217042289 & 2.1813793581 & 1 \\ -2.8053153235 & 1.2526740191 & -3.7573092457 & 0.3865582497 \\ -0.440747322 & -2.833803516 & -3.7692828392 & 0.0547217375 \end{pmatrix}$$

$$P = \begin{pmatrix} 0 & -1 & -1 & -1 & -1 & -1 & -1 \\ -1 & 0 & 1 & 1 & -1 & -1 & 1 \\ -1 & 1 & 0 & 1 & 1 & 1 & -1 \\ -1 & 1 & 1 & 0 & -1 & -1 & 1 \\ -1 & -1 & 1 & -1 & 0 & 1 & 1 \\ -1 & -1 & 1 & -1 & 1 & 0 & 1 \\ -1 & 1 & -1 & 1 & 1 & 1 & 0 \end{pmatrix} \qquad R = \begin{pmatrix} 0 & 0 & 0 & 0 & 0 & 0 & 0 \\ 0 & 0 & 0 & 0 & 0 & 0 & 0 \\ 4 & 1 & 0 & 4 & 1 & 1 & 1 \\ 1 & 4 & 1 & 0 & 1 & 4 & 1 \\ 0 & 0 & 0 & 0 & 0 & 0 & 0 \\ 2 & 2 & 3 & 2 & 6 & 0 & 3 \\ 4 & 4 & 5 & 5 & 3 & 3 & 0 \end{pmatrix}$$

$|P| = 66$

$I_3(P, R) = 19.05907173$

$I_3(-P, R) = 10.4810126582$

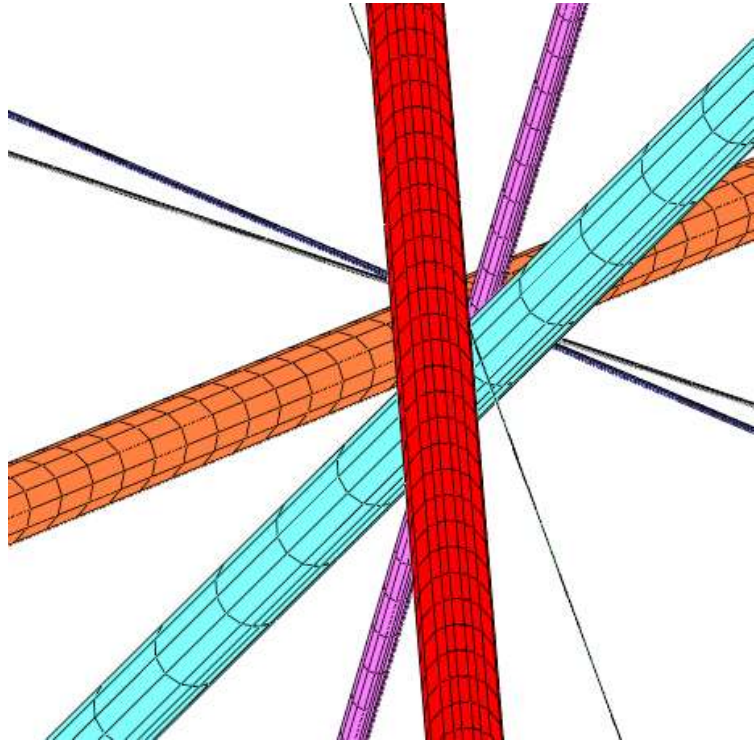

Appendix 4
7-cross
Configuration label m54f

$$\begin{pmatrix} \text{orange} \\ \text{red} \\ \text{blue} \\ \text{green} \\ \text{cyan} \\ \text{magenta} \\ \text{gray} \end{pmatrix} \begin{pmatrix} 0 & 0 & 0 & 1 \\ t_1 & p_1 & z_1 & r_1 \\ t_2 & p_2 & z_2 & r_2 \\ t_3 & p_3 & z_3 & r_3 \\ t_4 & p_4 & z_4 & r_4 \\ t_5 & p_5 & z_5 & r_5 \\ t_6 & p_6 & z_6 & r_6 \end{pmatrix} = \begin{pmatrix} 0 & 0 & 0 & 1 \\ 0.6472961325 & 3.1415926536 & 0 & 2.9301926741 \\ 2.1491323697 & 5.5928790695 & -1.5959129161 & 0.0100000015 \\ 2.3225368105 & 6.4596268383 & -2.6234027808 & 0.1901390518 \\ 2.261690526 & 5.9173285848 & -0.2860402847 & 0.8720483684 \\ -2.2777091438 & 9.1288097702 & -2.6653696393 & 0.0240152272 \\ 0.8289710387 & 2.8341073902 & -2.9409455479 & 0.0100000003 \end{pmatrix}$$

$$P = \begin{pmatrix} 0 & 1 & 1 & 1 & 1 & -1 & 1 \\ 1 & 0 & 1 & 1 & 1 & 1 & -1 \\ 1 & 1 & 0 & 1 & -1 & 1 & 1 \\ 1 & 1 & 1 & 0 & 1 & 1 & 1 \\ 1 & 1 & -1 & 1 & 0 & 1 & 1 \\ -1 & 1 & 1 & 1 & 1 & 0 & 1 \\ 1 & -1 & 1 & 1 & 1 & 1 & 0 \end{pmatrix} \qquad R = \begin{pmatrix} 0 & 2 & 3 & 6 & 3 & 2 & 2 \\ 0 & 0 & 0 & 0 & 0 & 0 & 0 \\ 2 & 2 & 0 & 3 & 6 & 3 & 2 \\ 0 & 0 & 0 & 0 & 0 & 0 & 0 \\ 0 & 0 & 0 & 0 & 0 & 0 & 0 \\ 4 & 4 & 3 & 5 & 3 & 0 & 5 \\ 1 & 4 & 4 & 1 & 1 & 1 & 0 \end{pmatrix}$$

$|P| = -54$

$I3(P, R) = 3.5344827586$

$I3(-P, R) = 25.8218390805$

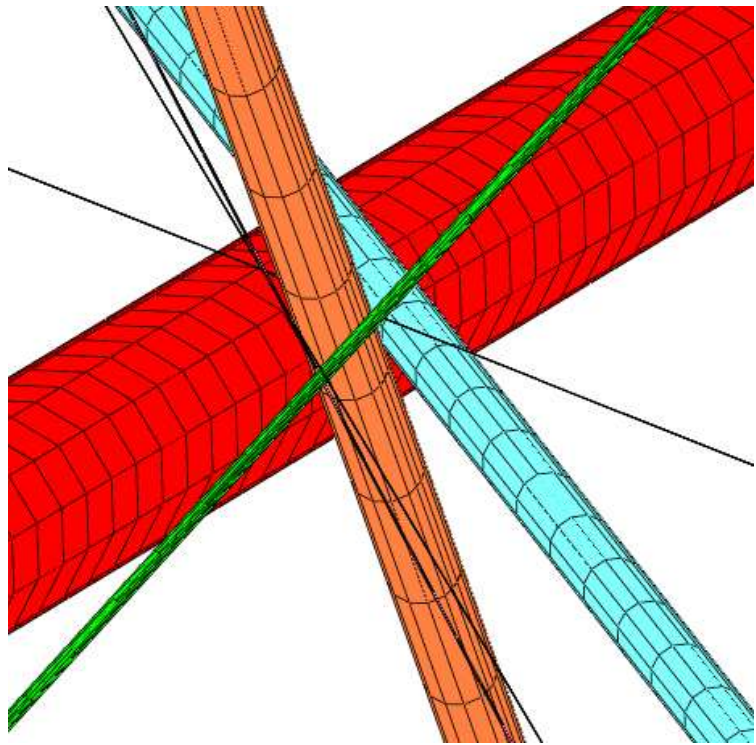

Appendix 4
7-cross
Configuration label x42ae

$$\begin{pmatrix} \text{orange} \\ \text{red} \\ \text{blue} \\ \text{green} \\ \text{cyan} \\ \text{magenta} \\ \text{gray} \end{pmatrix} \begin{pmatrix} 0 & 0 & 0 & 1 \\ t_1 & p_1 & z_1 & r_1 \\ t_2 & p_2 & z_2 & r_2 \\ t_3 & p_3 & z_3 & r_3 \\ t_4 & p_4 & z_4 & r_4 \\ t_5 & p_5 & z_5 & r_5 \\ t_6 & p_6 & z_6 & r_6 \end{pmatrix} = \begin{pmatrix} 0 & 0 & 0 & 1 \\ 2.4142965242 & 3.1415926536 & 0 & 1.0390205538 \\ 0.8306242378 & 4.7367811146 & -1.1792131003 & 0.0248835651 \\ 0.627524572 & 2.2448831596 & -7.1329308345 & 2.0587981649 \\ -0.9963787786 & 7.5828197129 & 1.3630038368 & 0.2705394048 \\ 2.3583918484 & 2.9174220888 & -2.2925037608 & 0.5352455258 \\ 0.7829673211 & 4.7272731801 & -0.9926218892 & 0.0100000002 \end{pmatrix}$$

$$P = \begin{pmatrix} 0 & 1 & 1 & 1 & -1 & 1 & 1 \\ 1 & 0 & 1 & -1 & 1 & -1 & 1 \\ 1 & 1 & 0 & 1 & 1 & -1 & 1 \\ 1 & -1 & 1 & 0 & -1 & -1 & 1 \\ -1 & 1 & 1 & -1 & 0 & -1 & -1 \\ 1 & -1 & -1 & -1 & -1 & 0 & -1 \\ 1 & 1 & 1 & 1 & -1 & -1 & 0 \end{pmatrix} \qquad R = \begin{pmatrix} 0 & 2 & 3 & 2 & 6 & 2 & 3 \\ 0 & 0 & 0 & 0 & 0 & 0 & 0 \\ 1 & 4 & 0 & 1 & 1 & 4 & 1 \\ 0 & 0 & 0 & 0 & 0 & 0 & 0 \\ 0 & 0 & 0 & 0 & 0 & 0 & 0 \\ 4 & 1 & 1 & 4 & 1 & 0 & 1 \\ 3 & 5 & 5 & 4 & 3 & 4 & 0 \end{pmatrix}$$

$|P| = 42$

$I3(P, R) = 11.375$

$I3(-P, R) = 17.9017857143$

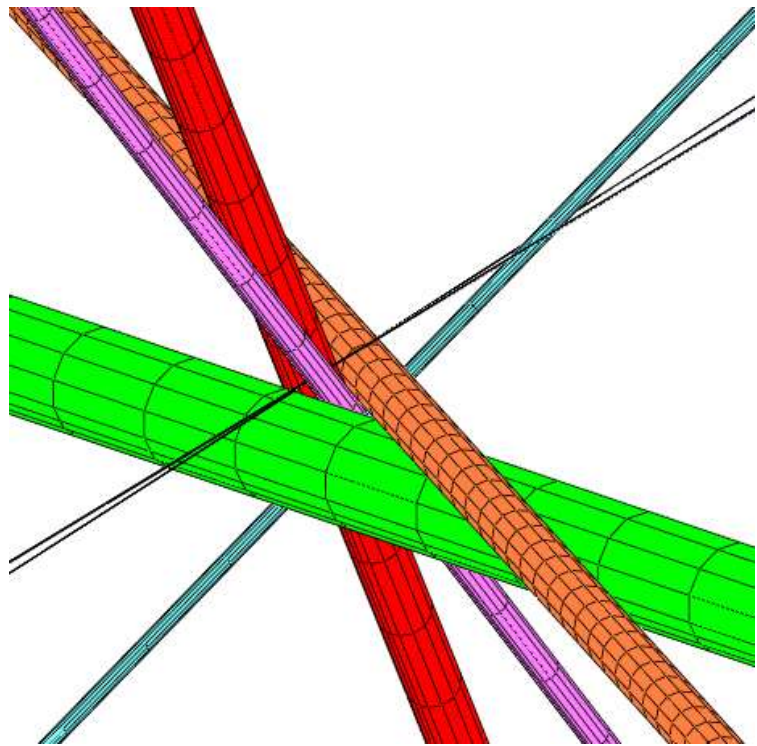

Appendix 4
7-cross
Configuration label mc78

$$\begin{pmatrix} \text{orange} \\ \text{red} \\ \text{blue} \\ \text{green} \\ \text{cyan} \\ \text{magenta} \\ \text{gray} \end{pmatrix} \begin{pmatrix} 0 & 0 & 0 & 1 \\ t_1 & p_1 & z_1 & r_1 \\ t_2 & p_2 & z_2 & r_2 \\ t_3 & p_3 & z_3 & r_3 \\ t_4 & p_4 & z_4 & r_4 \\ t_5 & p_5 & z_5 & r_5 \\ t_6 & p_6 & z_6 & r_6 \end{pmatrix} = \begin{pmatrix} 0 & 0 & 0 & 1 \\ 1.7791353451 & 3.1415926536 & 0 & 1.1253423662 \\ 1.1045727356 & 3.5880740918 & 4.3279705041 & 0.8522641276 \\ 2.8123571293 & 5.8063838647 & -11.9983435684 & 1.1253423662 \\ 1.760000576 & 4.2074378441 & 1.8067540017 & 1.1253423662 \\ 1.889392027 & 3.6971267403 & -9.3268299421 & 8.2294935906 \\ 0.3595967769 & 3.2439749608 & -7.4368784089 & 0.1345213406 \end{pmatrix}$$

$$P = \begin{pmatrix} 0 & 1 & 1 & 1 & 1 & 1 & 1 \\ 1 & 0 & -1 & -1 & -1 & 1 & 1 \\ 1 & -1 & 0 & 1 & 1 & 1 & -1 \\ 1 & -1 & 1 & 0 & 1 & -1 & 1 \\ 1 & -1 & 1 & 1 & 0 & -1 & -1 \\ 1 & 1 & 1 & -1 & -1 & 0 & 1 \\ 1 & 1 & -1 & 1 & -1 & 1 & 0 \end{pmatrix} \qquad R = \begin{pmatrix} 0 & 0 & 0 & 0 & 0 & 0 & 0 \\ 2 & 0 & 3 & 2 & 3 & 6 & 2 \\ 0 & 0 & 0 & 0 & 0 & 0 & 0 \\ 0 & 0 & 0 & 0 & 0 & 0 & 0 \\ 2 & 2 & 6 & 3 & 0 & 2 & 3 \\ 0 & 0 & 0 & 0 & 0 & 0 & 0 \\ 4 & 4 & 3 & 5 & 3 & 5 & 0 \end{pmatrix}$$

$|P| = 42$

$I3(P, R) = 9.2857142857$

$I3(-P, R) = 21.6857142857$

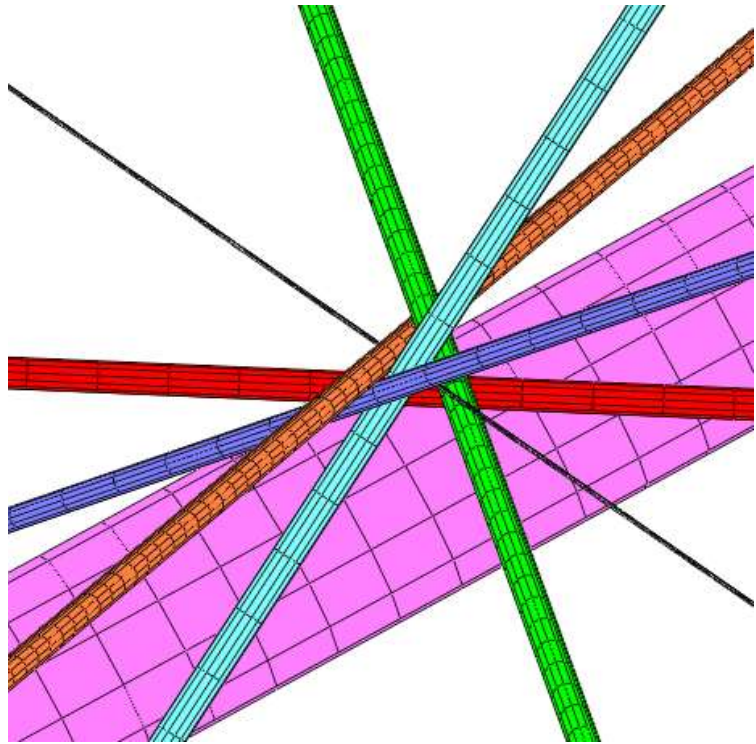

Appendix 4
7-cross
Configuration label mc67

$$\begin{pmatrix} \text{orange} \\ \text{red} \\ \text{blue} \\ \text{green} \\ \text{cyan} \\ \text{magenta} \\ \text{gray} \end{pmatrix} \begin{pmatrix} 0 & 0 & 0 & 1 \\ t_1 & p_1 & z_1 & r_1 \\ t_2 & p_2 & z_2 & r_2 \\ t_3 & p_3 & z_3 & r_3 \\ t_4 & p_4 & z_4 & r_4 \\ t_5 & p_5 & z_5 & r_5 \\ t_6 & p_6 & z_6 & r_6 \end{pmatrix} = \begin{pmatrix} 0 & 0 & 0 & 1 \\ 1.7791353468 & 3.1415926536 & 0 & 1.1253423662 \\ 1.104572737 & 3.5880740917 & 4.3279705027 & 0.8522641276 \\ 2.8123571291 & 5.8063838647 & -11.9983435697 & 1.1253423662 \\ 1.7600005788 & 4.2074378446 & 1.8067539969 & 1.1253423662 \\ 0.2532166103 & 1.5014296653 & -14.346597525 & 0.3941363942 \\ 0.3595967772 & 3.2439749608 & -7.4368784066 & 0.1345213406 \end{pmatrix}$$

$$P = \begin{pmatrix} 0 & 1 & 1 & 1 & 1 & 1 & 1 \\ 1 & 0 & -1 & -1 & -1 & -1 & 1 \\ 1 & -1 & 0 & 1 & 1 & -1 & -1 \\ 1 & -1 & 1 & 0 & 1 & 1 & 1 \\ 1 & -1 & 1 & 1 & 0 & 1 & -1 \\ 1 & -1 & -1 & 1 & 1 & 0 & 1 \\ 1 & 1 & -1 & 1 & -1 & 1 & 0 \end{pmatrix} \qquad R = \begin{pmatrix} 0 & 1 & 1 & 4 & 1 & 4 & 1 \\ 2 & 0 & 3 & 2 & 3 & 6 & 2 \\ 0 & 0 & 0 & 0 & 0 & 0 & 0 \\ 0 & 0 & 0 & 0 & 0 & 0 & 0 \\ 2 & 2 & 6 & 3 & 0 & 2 & 3 \\ 0 & 0 & 0 & 0 & 0 & 0 & 0 \\ 1 & 4 & 1 & 4 & 1 & 1 & 0 \end{pmatrix}$$

$|P| = 42$

$I3(P, R) = 10.188976378$

$I3(-P, R) = 17.8897637795$

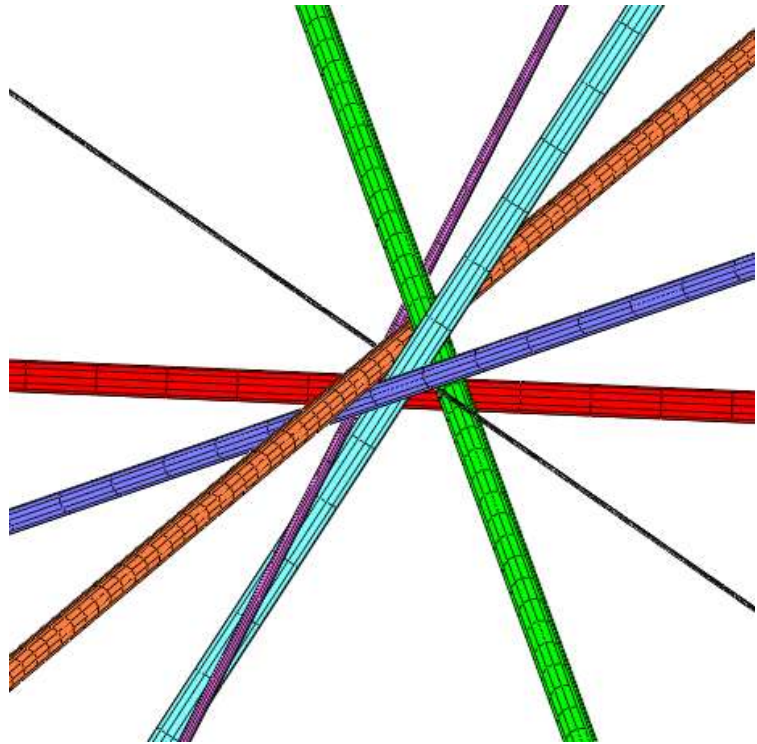

Appendix 4
7-cross
Configuration label md26

$$\begin{pmatrix} \text{orange} \\ \text{red} \\ \text{blue} \\ \text{green} \\ \text{cyan} \\ \text{magenta} \\ \text{gray} \end{pmatrix} \begin{pmatrix} 0 & 0 & 0 & 1 \\ t1 & p1 & z1 & r1 \\ t2 & p2 & z2 & r2 \\ t3 & p3 & z3 & r3 \\ t4 & p4 & z4 & r4 \\ t5 & p5 & z5 & r5 \\ t6 & p6 & z6 & r6 \end{pmatrix} = \begin{pmatrix} 0 & 0 & 0 & 1 \\ 1.4938459348 & 3.1415926536 & 0 & 0.7245086968 \\ 1.891605667 & 5.2820656512 & -4.6814126592 & 2.7041739475 \\ 1.4042016163 & 3.9974002189 & -1.2778361881 & 0.7245086968 \\ 1.9852175057 & 0.5477517437 & -1.0697988026 & 3.0091345213 \\ 0.7745288109 & 1.5365412255 & -4.0396665488 & 0.7245086969 \\ 1.0529825428 & 3.1921243627 & -3.0478713856 & 0.1548918413 \end{pmatrix}$$

$$P = \begin{pmatrix} 0 & 1 & 1 & 1 & 1 & 1 & 1 \\ 1 & 0 & 1 & 1 & -1 & -1 & 1 \\ 1 & 1 & 0 & 1 & -1 & 1 & 1 \\ 1 & 1 & 1 & 0 & -1 & 1 & -1 \\ 1 & -1 & -1 & -1 & 0 & 1 & 1 \\ 1 & -1 & 1 & 1 & 1 & 0 & 1 \\ 1 & 1 & 1 & -1 & 1 & 1 & 0 \end{pmatrix} \qquad R = \begin{pmatrix} 0 & 3 & 5 & 4 & 5 & 4 & 3 \\ 0 & 0 & 0 & 0 & 0 & 0 & 0 \\ 0 & 0 & 0 & 0 & 0 & 0 & 0 \\ 2 & 6 & 3 & 0 & 2 & 2 & 3 \\ 0 & 0 & 0 & 0 & 0 & 0 & 0 \\ 1 & 1 & 1 & 1 & 4 & 0 & 4 \\ 1 & 4 & 4 & 1 & 1 & 1 & 0 \end{pmatrix}$$

$|P| = 42$

$I3(P, R) = 0.598540146$

$I3(-P, R) = 22.7299270073$

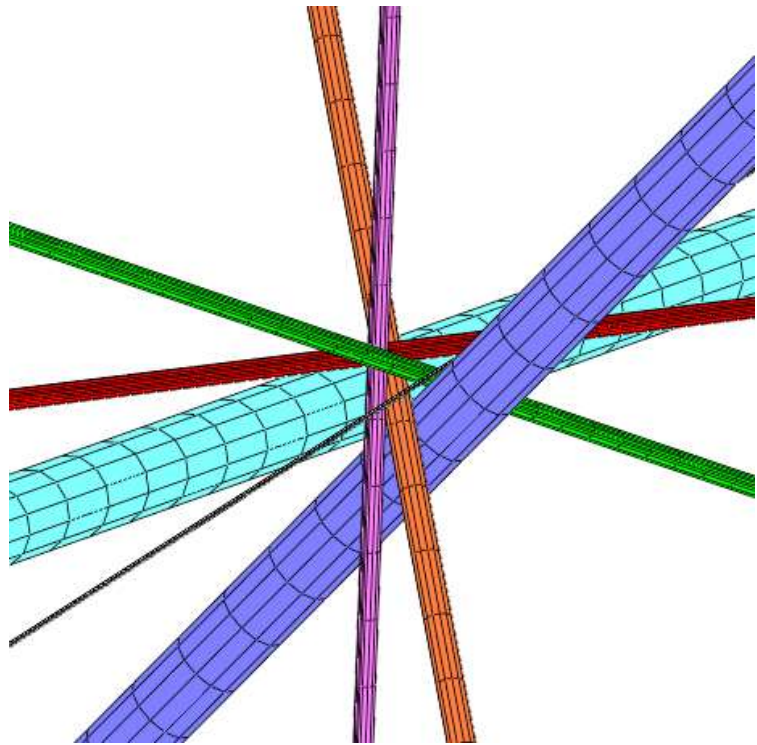

Appendix 4
7-cross
Configuration label d46

$$\begin{pmatrix} \text{orange} \\ \text{red} \\ \text{blue} \\ \text{green} \\ \text{cyan} \\ \text{magenta} \\ \text{gray} \end{pmatrix} \begin{pmatrix} 0 & 0 & 0 & 1 \\ t_1 & p_1 & z_1 & r_1 \\ t_2 & p_2 & z_2 & r_2 \\ t_3 & p_3 & z_3 & r_3 \\ t_4 & p_4 & z_4 & r_4 \\ t_5 & p_5 & z_5 & r_5 \\ t_6 & p_6 & z_6 & r_6 \end{pmatrix} = \begin{pmatrix} 0 & 0 & 0 & 1 \\ 0.9071368704 & 3.1415926536 & 0 & 0.4424161588 \\ 1.6477467186 & 2.73193228 & 1.9436660558 & 0.7245086968 \\ 1.737391037 & 1.8761247147 & 0.6658298673 & 0.7245086968 \\ -1.9852175055 & 5.3257731899 & 0.8738672531 & 3.0091345213 \\ 2.3670638426 & 4.3369837082 & -2.0960004924 & 0.7245086969 \\ 2.0886101106 & 2.6814005709 & -1.1042053296 & 0.1548918413 \end{pmatrix}$$

$$P = \begin{pmatrix} 0 & 1 & 1 & 1 & -1 & 1 & 1 \\ 1 & 0 & 1 & 1 & -1 & 1 & -1 \\ 1 & 1 & 0 & -1 & -1 & 1 & -1 \\ 1 & 1 & -1 & 0 & -1 & -1 & 1 \\ -1 & -1 & -1 & -1 & 0 & 1 & 1 \\ 1 & 1 & 1 & -1 & 1 & 0 & -1 \\ 1 & -1 & -1 & 1 & 1 & -1 & 0 \end{pmatrix} \qquad R = \begin{pmatrix} 0 & 1 & 1 & 4 & 4 & 1 & 1 \\ 3 & 0 & 5 & 4 & 4 & 5 & 3 \\ 0 & 0 & 0 & 0 & 0 & 0 & 0 \\ 1 & 1 & 4 & 0 & 1 & 1 & 4 \\ 0 & 0 & 0 & 0 & 0 & 0 & 0 \\ 1 & 1 & 1 & 1 & 4 & 0 & 4 \\ 0 & 0 & 0 & 0 & 0 & 0 & 0 \end{pmatrix}$$

$|P| = 42$

$I3(P, R) = -5.0526315789$

$I3(-P, R) = 26.6842105263$

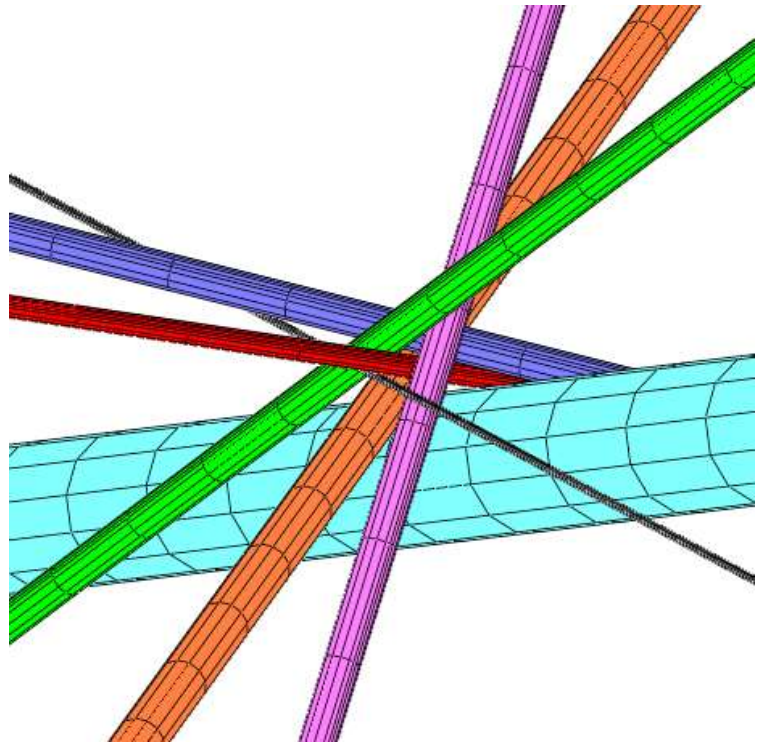

Appendix 4
7-cross
Configuration label d14

$$\begin{pmatrix} \text{orange} \\ \text{red} \\ \text{blue} \\ \text{green} \\ \text{cyan} \\ \text{magenta} \\ \text{gray} \end{pmatrix} \begin{pmatrix} 0 & 0 & 0 & 1 \\ t_1 & p_1 & z_1 & r_1 \\ t_2 & p_2 & z_2 & r_2 \\ t_3 & p_3 & z_3 & r_3 \\ t_4 & p_4 & z_4 & r_4 \\ t_5 & p_5 & z_5 & r_5 \\ t_6 & p_6 & z_6 & r_6 \end{pmatrix} = \begin{pmatrix} 0 & 0 & 0 & 1 \\ 0.8324454359 & 3.1415926536 & 0 & 1.6376180716 \\ 1.4390421887 & 2.4591034075 & -6.3139087372 & 1.6376180716 \\ -0.4975082062 & 5.3138676133 & -16.5868109729 & 3.8657118484 \\ -1.402980087 & 5.9879993996 & 3.0917189391 & 6.8015927119 \\ 1.8114013396 & 4.4425860304 & -3.3945233218 & 1.6376180718 \\ -1.2573334374 & 0.1491143695 & -7.4210428119 & 0.3501043943 \end{pmatrix}$$

$$P = \begin{pmatrix} 0 & 1 & 1 & -1 & -1 & 1 & -1 \\ 1 & 0 & -1 & -1 & -1 & 1 & -1 \\ 1 & -1 & 0 & 1 & -1 & -1 & 1 \\ -1 & -1 & 1 & 0 & -1 & -1 & -1 \\ -1 & -1 & -1 & -1 & 0 & 1 & 1 \\ 1 & 1 & -1 & -1 & 1 & 0 & -1 \\ -1 & -1 & 1 & -1 & 1 & -1 & 0 \end{pmatrix} \quad R = \begin{pmatrix} 0 & 4 & 5 & 3 & 4 & 3 & 5 \\ 0 & 0 & 0 & 0 & 0 & 0 & 0 \\ 1 & 4 & 0 & 1 & 1 & 1 & 4 \\ 0 & 0 & 0 & 0 & 0 & 0 & 0 \\ 0 & 0 & 0 & 0 & 0 & 0 & 0 \\ 3 & 4 & 3 & 5 & 5 & 0 & 4 \\ 0 & 0 & 0 & 0 & 0 & 0 & 0 \end{pmatrix}$$

$|P| = 42$

$I3(P, R) = -6.6$

$I3(-P, R) = 30.5$

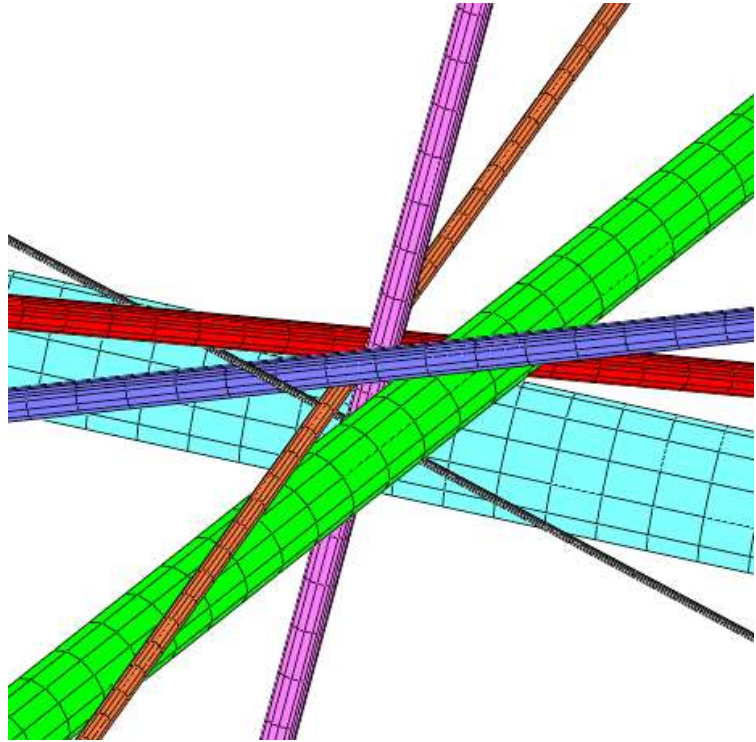

Appendix 4
7-cross
Configuration label d12

$$\begin{pmatrix} \text{orange} \\ \text{red} \\ \text{blue} \\ \text{green} \\ \text{cyan} \\ \text{magenta} \\ \text{gray} \end{pmatrix} \begin{pmatrix} 0 & 0 & 0 & 1 \\ t_1 & p_1 & z_1 & r_1 \\ t_2 & p_2 & z_2 & r_2 \\ t_3 & p_3 & z_3 & r_3 \\ t_4 & p_4 & z_4 & r_4 \\ t_5 & p_5 & z_5 & r_5 \\ t_6 & p_6 & z_6 & r_6 \end{pmatrix} = \begin{pmatrix} 0 & 0 & 0 & 1 \\ -2.1347933998 & 3.1415926536 & 0 & 3.7324244132 \\ -0.8534774438 & 3.6271451736 & 7.0593818229 & 1 \\ -0.6580278153 & 2.1631195781 & 17.9868279513 & 2.3605698518 \\ -0.6265189059 & 4.0945979308 & 18.2234874798 & 4.1533449282 \\ 1.5397269943 & 2.6767209096 & 10.8905648361 & 1.0000000001 \\ -0.4434477891 & 4.9414166212 & -0.1209552217 & 0.2137887951 \end{pmatrix}$$

$$P = \begin{pmatrix} 0 & -1 & -1 & -1 & -1 & 1 & -1 \\ -1 & 0 & -1 & 1 & -1 & -1 & -1 \\ -1 & -1 & 0 & 1 & -1 & -1 & 1 \\ -1 & 1 & 1 & 0 & -1 & -1 & -1 \\ -1 & -1 & -1 & -1 & 0 & 1 & 1 \\ 1 & -1 & -1 & -1 & 1 & 0 & -1 \\ -1 & -1 & 1 & -1 & 1 & -1 & 0 \end{pmatrix} \qquad R = \begin{pmatrix} 0 & 0 & 0 & 0 & 0 & 0 & 0 \\ 0 & 0 & 0 & 0 & 0 & 0 & 0 \\ 6 & 3 & 0 & 2 & 2 & 2 & 3 \\ 0 & 0 & 0 & 0 & 0 & 0 & 0 \\ 0 & 0 & 0 & 0 & 0 & 0 & 0 \\ 4 & 3 & 3 & 5 & 5 & 0 & 4 \\ 3 & 6 & 2 & 3 & 2 & 2 & 0 \end{pmatrix}$$

$|P| = 42$

$I3(P, R) = 5.6981132075$

$I3(-P, R) = 22.3018867925$

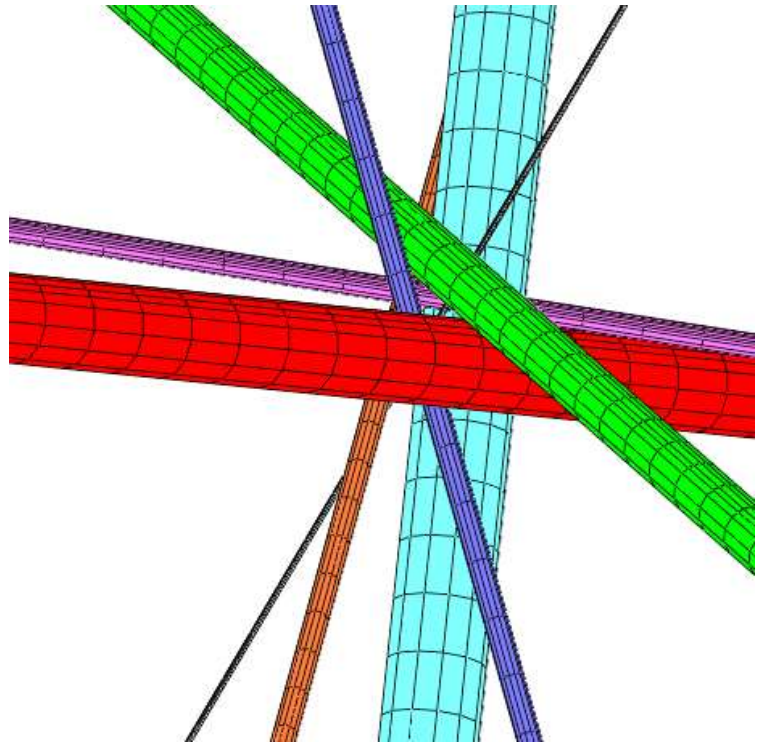

Appendix 4

7-cross
Configuration label mb57

$$\begin{pmatrix} \text{orange} \\ \text{red} \\ \text{blue} \\ \text{green} \\ \text{cyan} \\ \text{magenta} \\ \text{gray} \end{pmatrix} \begin{pmatrix} 0 & 0 & 0 & 1 \\ t1 & p1 & z1 & r1 \\ t2 & p2 & z2 & r2 \\ t3 & p3 & z3 & r3 \\ t4 & p4 & z4 & r4 \\ t5 & p5 & z5 & r5 \\ t6 & p6 & z6 & r6 \end{pmatrix} = \begin{pmatrix} 0 & 0 & 0 & 1 \\ 3.0075540462 & 3.1415926536 & 0 & 1.9968460919 \\ 2.1530621814 & 3.6306497025 & 52.4184492569 & 5.2372850077 \\ 2.9435681985 & 6.0522686246 & -104.9962792013 & 3.6442883187 \\ 3.0658343354 & 4.9524552997 & -104.3899059477 & 0.4752981793 \\ 0.0821609902 & 2.2277784914 & -89.7093523835 & 0.1986672299 \\ 0.2059923778 & 3.1995917742 & -83.9338332048 & 0.1986672299 \end{pmatrix}$$

$$P = \begin{pmatrix} 0 & 1 & 1 & 1 & 1 & 1 & 1 \\ 1 & 0 & -1 & -1 & 1 & -1 & 1 \\ 1 & -1 & 0 & 1 & 1 & -1 & -1 \\ 1 & -1 & 1 & 0 & -1 & 1 & 1 \\ 1 & 1 & 1 & -1 & 0 & 1 & 1 \\ 1 & -1 & -1 & 1 & 1 & 0 & 1 \\ 1 & 1 & -1 & 1 & 1 & 1 & 0 \end{pmatrix} \qquad R = \begin{pmatrix} 0 & 1 & 1 & 4 & 1 & 4 & 1 \\ 2 & 0 & 6 & 2 & 3 & 3 & 2 \\ 0 & 0 & 0 & 0 & 0 & 0 & 0 \\ 0 & 0 & 0 & 0 & 0 & 0 & 0 \\ 0 & 0 & 0 & 0 & 0 & 0 & 0 \\ 1 & 1 & 4 & 1 & 4 & 0 & 1 \\ 4 & 4 & 3 & 5 & 5 & 3 & 0 \end{pmatrix}$$

$|P| = -18$

$I3(P, R) = 7.401459854$

$I3(-P, R) = 15.9270072993$

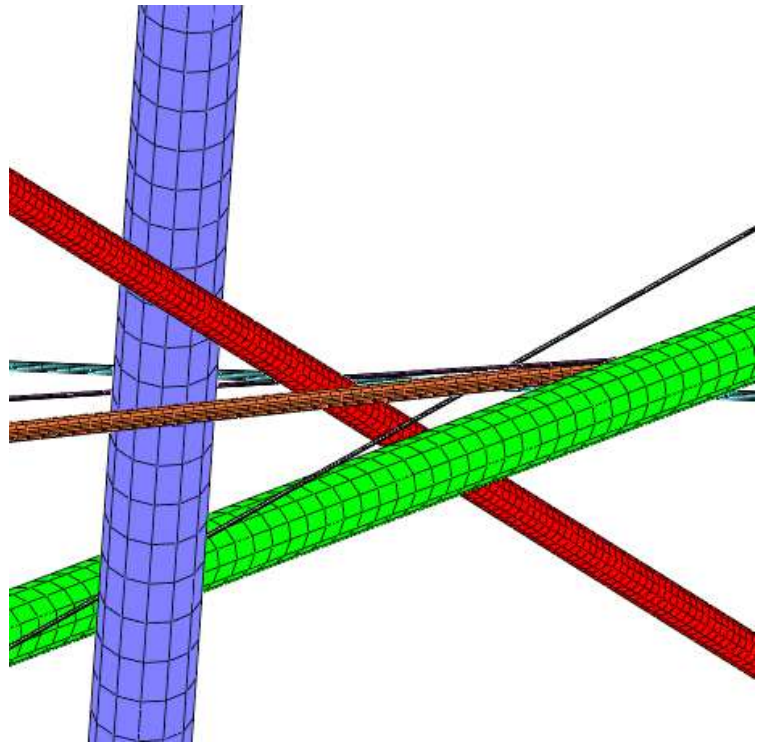

Appendix 4
7-cross
Configuration label b27

$$\begin{pmatrix} \text{orange} \\ \text{red} \\ \text{blue} \\ \text{green} \\ \text{cyan} \\ \text{magenta} \\ \text{gray} \end{pmatrix} \begin{pmatrix} 0 & 0 & 0 & 1 \\ t_1 & p_1 & z_1 & r_1 \\ t_2 & p_2 & z_2 & r_2 \\ t_3 & p_3 & z_3 & r_3 \\ t_4 & p_4 & z_4 & r_4 \\ t_5 & p_5 & z_5 & r_5 \\ t_6 & p_6 & z_6 & r_6 \end{pmatrix} = \begin{pmatrix} 0 & 0 & 0 & 1 \\ 0.9885304618 & 3.1415926536 & 0 & 5.2372850077 \\ 0.1980244566 & 0.7199737283 & -157.4147284274 & 3.6442883187 \\ 2.8887794891 & 1.6898366789 & -32.4749334353 & 6.5161684705 \\ 0.0757583181 & 1.8197870505 & -156.8083552138 & 0.4752981793 \\ 3.0594316632 & 4.5444638579 & -142.1278014044 & 0.1986672299 \\ 2.9356002737 & 3.5726505793 & -136.3522823262 & 0.1986672299 \end{pmatrix}$$

$$P = \begin{pmatrix} 0 & 1 & 1 & 1 & 1 & 1 & 1 \\ 1 & 0 & -1 & -1 & -1 & 1 & 1 \\ 1 & -1 & 0 & -1 & 1 & -1 & -1 \\ 1 & -1 & -1 & 0 & 1 & -1 & 1 \\ 1 & -1 & 1 & 1 & 0 & -1 & -1 \\ 1 & 1 & -1 & -1 & -1 & 0 & -1 \\ 1 & 1 & -1 & 1 & -1 & -1 & 0 \end{pmatrix} \qquad R = \begin{pmatrix} 0 & 1 & 4 & 1 & 1 & 4 & 1 \\ 0 & 0 & 0 & 0 & 0 & 0 & 0 \\ 0 & 0 & 0 & 0 & 0 & 0 & 0 \\ 2 & 6 & 3 & 0 & 2 & 2 & 3 \\ 0 & 0 & 0 & 0 & 0 & 0 & 0 \\ 1 & 4 & 1 & 1 & 4 & 0 & 1 \\ 4 & 3 & 5 & 4 & 5 & 3 & 0 \end{pmatrix}$$

$|P| = -18$

$I3(P, R) = 17.5035460993$

$I3(-P, R) = 7.6595744681$

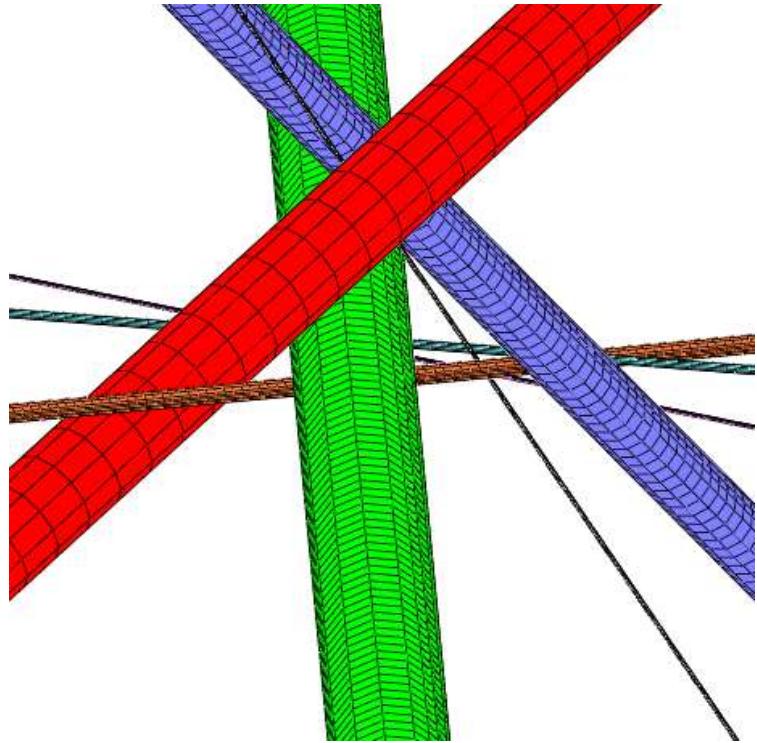

Appendix 4
7-cross
Configuration label ma49

$$\begin{pmatrix} \text{orange} \\ \text{red} \\ \text{blue} \\ \text{green} \\ \text{cyan} \\ \text{magenta} \\ \text{gray} \end{pmatrix} \begin{pmatrix} 0 & 0 & 0 & 1 \\ t_1 & p_1 & z_1 & r_1 \\ t_2 & p_2 & z_2 & r_2 \\ t_3 & p_3 & z_3 & r_3 \\ t_4 & p_4 & z_4 & r_4 \\ t_5 & p_5 & z_5 & r_5 \\ t_6 & p_6 & z_6 & r_6 \end{pmatrix} = \begin{pmatrix} 0 & 0 & 0 & 1 \\ 2.8818123122 & 3.1415926536 & 0 & 1 \\ 0.2833324327 & -0.3533072984 & 4.1314378435 & 1 \\ 0.2996794362 & -0.1842493099 & -2.1085601822 & 1 \\ 2.8945201973 & -0.2691324276 & -171.7820458028 & 1 \\ 2.8408722586 & 2.9601745121 & -8.1355881072 & 1 \\ 0.2703593444 & 9.329780127 & -155.2252341719 & 0.93 \end{pmatrix}$$

$$P = \begin{pmatrix} 0 & 1 & 1 & 1 & 1 & 1 & 1 \\ 1 & 0 & -1 & -1 & 1 & -1 & -1 \\ 1 & -1 & 0 & 1 & 1 & -1 & -1 \\ 1 & -1 & 1 & 0 & -1 & 1 & 1 \\ 1 & 1 & 1 & -1 & 0 & -1 & 1 \\ 1 & -1 & -1 & 1 & -1 & 0 & 1 \\ 1 & -1 & -1 & 1 & 1 & 1 & 0 \end{pmatrix} \quad R = \begin{pmatrix} 0 & 0 & 0 & 0 & 0 & 0 & 0 \\ 6 & 0 & 3 & 3 & 2 & 2 & 2 \\ 0 & 0 & 0 & 0 & 0 & 0 & 0 \\ 0 & 0 & 0 & 0 & 0 & 0 & 0 \\ 1 & 1 & 1 & 4 & 0 & 4 & 1 \\ 0 & 0 & 0 & 0 & 0 & 0 & 0 \\ 3 & 3 & 4 & 4 & 5 & 5 & 0 \end{pmatrix}$$

$|P| = -18$

$I3(P, R) = -1$

$I3(-P, R) = 22.5714285714$

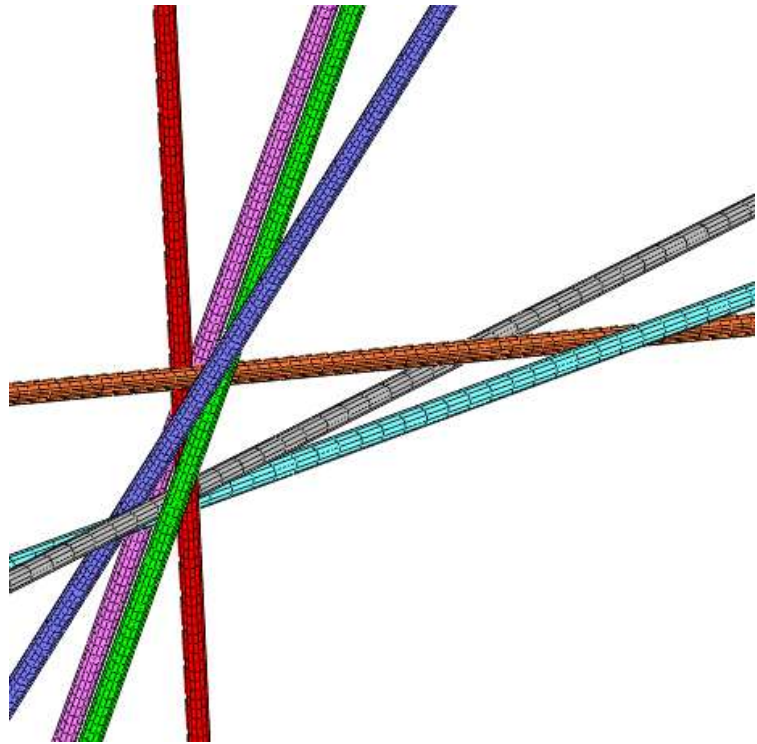

Appendix 4
7-cross
Configuration label ma29

$$\begin{pmatrix} \text{orange} \\ \text{red} \\ \text{blue} \\ \text{green} \\ \text{cyan} \\ \text{magenta} \\ \text{gray} \end{pmatrix} \begin{pmatrix} 0 & 0 & 0 & 1 \\ t1 & p1 & z1 & r1 \\ t2 & p2 & z2 & r2 \\ t3 & p3 & z3 & r3 \\ t4 & p4 & z4 & r4 \\ t5 & p5 & z5 & r5 \\ t6 & p6 & z6 & r6 \end{pmatrix} = \begin{pmatrix} 0 & 0 & 0 & 1 \\ 0.1751285285 & 3.1415926536 & 0 & 1 \\ 3.1045005059 & 4.939161367 & -107.4772787932 & 1 \\ 0.2442917266 & 3.4790879912 & -10.6870107377 & 1 \\ 2.8283262735 & 3.3177362857 & -98.9508955197 & 1 \\ 0.2694291198 & 6.7681766486 & -113.6989713016 & 1 \\ 0.0627211797 & 13.1751833334 & -156.2307520591 & 0.22 \end{pmatrix}$$

$$P = \begin{pmatrix} 0 & 1 & 1 & 1 & 1 & 1 & 1 \\ 1 & 0 & 1 & 1 & 1 & -1 & -1 \\ 1 & 1 & 0 & 1 & -1 & -1 & 1 \\ 1 & 1 & 1 & 0 & -1 & 1 & 1 \\ 1 & 1 & -1 & -1 & 0 & -1 & 1 \\ 1 & -1 & -1 & 1 & -1 & 0 & 1 \\ 1 & -1 & 1 & 1 & 1 & 1 & 0 \end{pmatrix} \qquad R = \begin{pmatrix} 0 & 2 & 6 & 2 & 2 & 3 & 3 \\ 0 & 0 & 0 & 0 & 0 & 0 & 0 \\ 0 & 0 & 0 & 0 & 0 & 0 & 0 \\ 1 & 4 & 4 & 0 & 1 & 1 & 1 \\ 1 & 1 & 1 & 4 & 0 & 4 & 1 \\ 0 & 0 & 0 & 0 & 0 & 0 & 0 \\ 3 & 4 & 3 & 4 & 5 & 5 & 0 \end{pmatrix}$$

$|P| = -18$

$I3(P, R) = 16.3780487805$

$I3(-P, R) = 7.487804878$

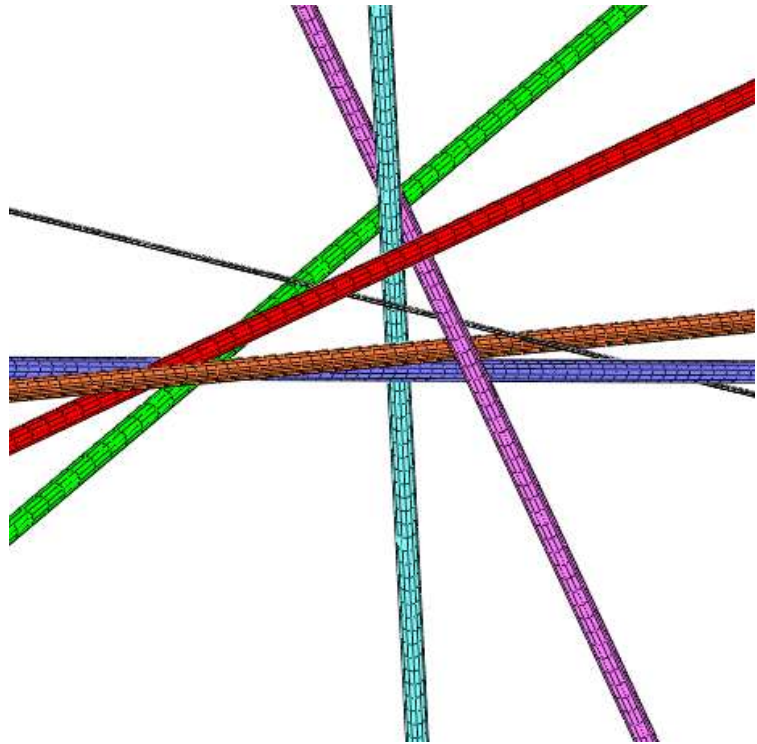

Appendix 4
7-cross
Configuration label ma28

$$\begin{pmatrix} \text{orange} \\ \text{red} \\ \text{blue} \\ \text{green} \\ \text{cyan} \\ \text{magenta} \\ \text{gray} \end{pmatrix} \begin{pmatrix} 0 & 0 & 0 & 1 \\ t_1 & p_1 & z_1 & r_1 \\ t_2 & p_2 & z_2 & r_2 \\ t_3 & p_3 & z_3 & r_3 \\ t_4 & p_4 & z_4 & r_4 \\ t_5 & p_5 & z_5 & r_5 \\ t_6 & p_6 & z_6 & r_6 \end{pmatrix} = \begin{pmatrix} 0 & 0 & 0 & 1 \\ 1.0023676966 & 3.1415926536 & 0 & 1 \\ 3.0945408013 & 4.5981544412 & -79.7066858018 & 1 \\ 1.0158405491 & 3.2257787125 & -2.2995748006 & 1 \\ 2.4506215342 & 3.18554618 & -84.1428805766 & 0.9 \\ 0.6874149396 & 6.4104635276 & -87.5452912918 & 1 \\ 0.5219298206 & 9.4399202347 & -34.9304686969 & 0.2404687146 \end{pmatrix}$$

$$P = \begin{pmatrix} 0 & 1 & 1 & 1 & 1 & 1 & 1 \\ 1 & 0 & 1 & 1 & 1 & -1 & 1 \\ 1 & 1 & 0 & 1 & -1 & -1 & 1 \\ 1 & 1 & 1 & 0 & -1 & 1 & -1 \\ 1 & 1 & -1 & -1 & 0 & -1 & 1 \\ 1 & -1 & -1 & 1 & -1 & 0 & 1 \\ 1 & 1 & 1 & -1 & 1 & 1 & 0 \end{pmatrix} \qquad R = \begin{pmatrix} 0 & 1 & 4 & 1 & 1 & 4 & 1 \\ 0 & 0 & 0 & 0 & 0 & 0 & 0 \\ 0 & 0 & 0 & 0 & 0 & 0 & 0 \\ 2 & 6 & 3 & 0 & 2 & 2 & 3 \\ 1 & 1 & 1 & 4 & 0 & 4 & 1 \\ 0 & 0 & 0 & 0 & 0 & 0 & 0 \\ 5 & 3 & 5 & 3 & 4 & 4 & 0 \end{pmatrix}$$

$|P| = -18$

$I3(P, R) = 20.9738219895$

$I3(-P, R) = 7.1518324607$

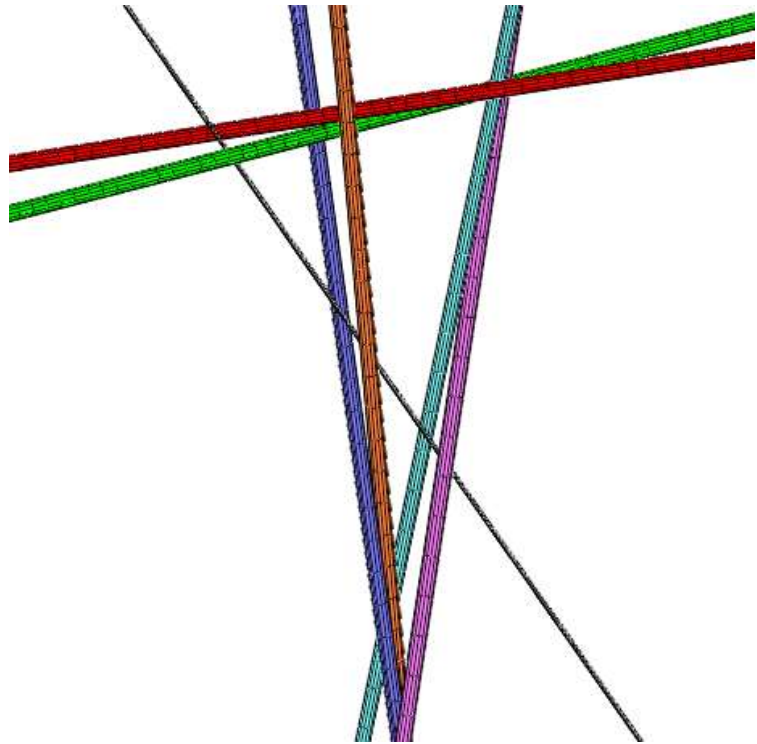

Appendix 4
7-cross
Configuration label a89

$$\begin{pmatrix} \text{orange} \\ \text{red} \\ \text{blue} \\ \text{green} \\ \text{cyan} \\ \text{magenta} \\ \text{gray} \end{pmatrix} \begin{pmatrix} 0 & 0 & 0 & 1 \\ t_1 & p_1 & z_1 & r_1 \\ t_2 & p_2 & z_2 & r_2 \\ t_3 & p_3 & z_3 & r_3 \\ t_4 & p_4 & z_4 & r_4 \\ t_5 & p_5 & z_5 & r_5 \\ t_6 & p_6 & z_6 & r_6 \end{pmatrix} = \begin{pmatrix} 0 & 0 & 0 & 1 \\ 0.3255445253 & 3.1415926536 & 0 & 1 \\ 0.3800224809 & 9.3515809945 & 13.4293807485 & 1 \\ 2.4244469645 & 6.4336249988 & -23.5948605815 & 1 \\ 0.3316040913 & 9.1928411285 & 4.8892717374 & 1 \\ 0.2843728713 & 9.2680129687 & -9.9091055237 & 1 \\ 2.8406082265 & 6.0133479798 & -10.6464648447 & 1 \end{pmatrix}$$

$$P = \begin{pmatrix} 0 & 1 & 1 & 1 & 1 & 1 & 1 \\ 1 & 0 & 1 & 1 & 1 & -1 & 1 \\ 1 & 1 & 0 & -1 & -1 & -1 & 1 \\ 1 & 1 & -1 & 0 & -1 & 1 & 1 \\ 1 & 1 & -1 & -1 & 0 & 1 & -1 \\ 1 & -1 & -1 & 1 & 1 & 0 & 1 \\ 1 & 1 & 1 & 1 & -1 & 1 & 0 \end{pmatrix} \qquad R = \begin{pmatrix} 0 & 1 & 1 & 4 & 1 & 1 & 4 \\ 4 & 0 & 4 & 4 & 4 & 4 & 4 \\ 0 & 0 & 0 & 0 & 0 & 0 & 0 \\ 0 & 0 & 0 & 0 & 0 & 0 & 0 \\ 1 & 1 & 4 & 4 & 0 & 1 & 1 \\ 1 & 1 & 1 & 1 & 4 & 0 & 4 \\ 0 & 0 & 0 & 0 & 0 & 0 & 0 \end{pmatrix}$$

$|P| = -18$

$I3(P, R) = 14.0731707317$

$I3(-P, R) = 7.8048780488$

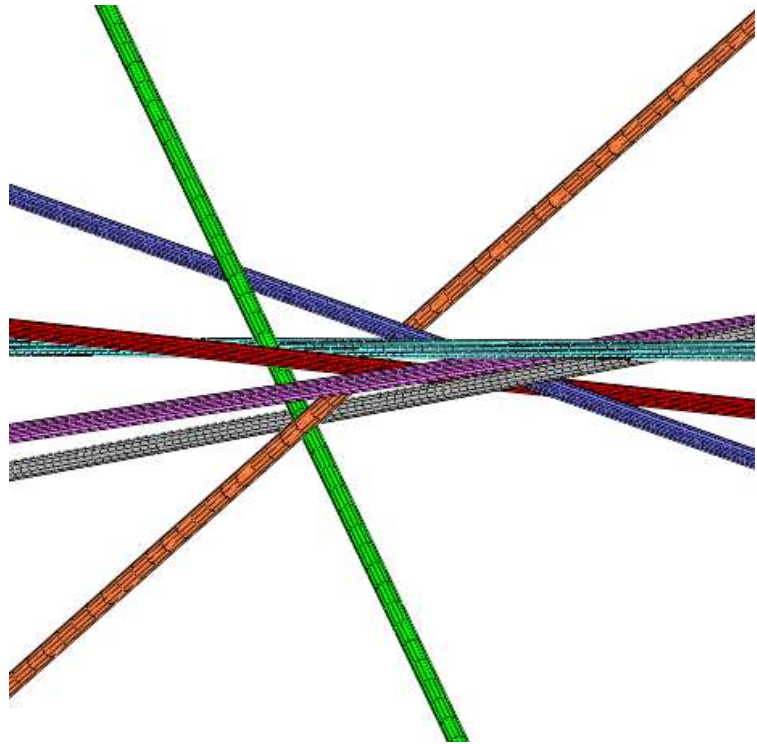

Appendix 4
7-cross
Configuration label a48

$$\begin{pmatrix} \text{orange} \\ \text{red} \\ \text{blue} \\ \text{green} \\ \text{cyan} \\ \text{magenta} \\ \text{gray} \end{pmatrix} \begin{pmatrix} 0 & 0 & 0 & 1 \\ t_1 & p_1 & z_1 & r_1 \\ t_2 & p_2 & z_2 & r_2 \\ t_3 & p_3 & z_3 & r_3 \\ t_4 & p_4 & z_4 & r_4 \\ t_5 & p_5 & z_5 & r_5 \\ t_6 & p_6 & z_6 & r_6 \end{pmatrix} = \begin{pmatrix} 0 & 0 & 0 & 1 \\ 1.4170009063 & 3.1415926536 & 0 & 1 \\ 2.2046039322 & 8.8661719545 & 4.1778592744 & 1 \\ 0.7272528342 & 7.604397779 & -0.1528808167 & 1 \\ 0.0187774206 & 8.0985087711 & -193.1372673589 & 1 \\ 2.4256456939 & 4.4531366885 & -4.7163990297 & 1 \\ 2.2099340599 & 2.580794305 & -3.3522304668 & 0.95 \end{pmatrix}$$

$$P = \begin{pmatrix} 0 & 1 & 1 & 1 & 1 & 1 & 1 \\ 1 & 0 & 1 & 1 & -1 & 1 & -1 \\ 1 & 1 & 0 & -1 & -1 & 1 & -1 \\ 1 & 1 & -1 & 0 & 1 & -1 & 1 \\ 1 & -1 & -1 & 1 & 0 & 1 & -1 \\ 1 & 1 & 1 & -1 & 1 & 0 & -1 \\ 1 & -1 & -1 & 1 & -1 & -1 & 0 \end{pmatrix} \qquad R = \begin{pmatrix} 0 & 0 & 0 & 0 & 0 & 0 & 0 \\ 4 & 0 & 4 & 4 & 4 & 4 & 4 \\ 0 & 0 & 0 & 0 & 0 & 0 & 0 \\ 1 & 1 & 4 & 0 & 1 & 1 & 4 \\ 1 & 1 & 1 & 4 & 0 & 4 & 1 \\ 0 & 0 & 0 & 0 & 0 & 0 & 0 \\ 4 & 1 & 1 & 1 & 4 & 1 & 0 \end{pmatrix}$$

$|P| = -18$

$I3(P, R) = 15.0592105263$

$I3(-P, R) = 15$

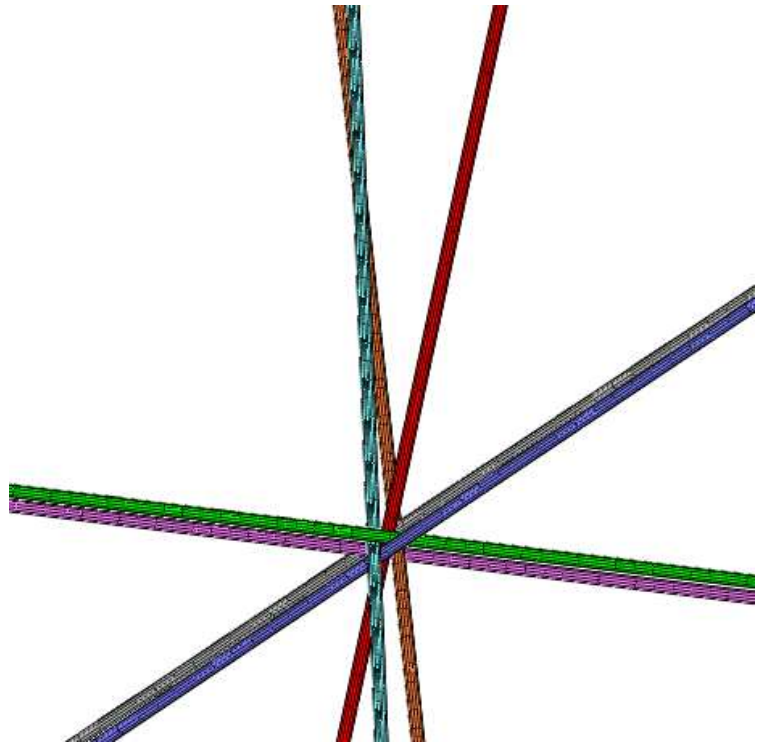

Appendix 4
7-cross
Configuration label  a24

$$\begin{pmatrix} \text{orange} \\ \text{red} \\ \text{blue} \\ \text{green} \\ \text{cyan} \\ \text{magenta} \\ \text{gray} \end{pmatrix} \begin{pmatrix} 0 & 0 & 0 & 1 \\ t_1 & p_1 & z_1 & r_1 \\ t_2 & p_2 & z_2 & r_2 \\ t_3 & p_3 & z_3 & r_3 \\ t_4 & p_4 & z_4 & r_4 \\ t_5 & p_5 & z_5 & r_5 \\ t_6 & p_6 & z_6 & r_6 \end{pmatrix} = \begin{pmatrix} 0 & 0 & 0 & 1 \\ 1.5801077274 & 3.1415926536 & 0 & 1.4 \\ 1.7285598483 & 2.1133775241 & -3.1117809296 & 1.8493691875 \\ 0.358943564 & 2.145262318 & -10.6342239331 & 0.0426287245 \\ 2.5806677038 & -1.6226476996 & -7.9931295822 & 0.8896820927 \\ 2.7473102984 & -8.0487328444 & -9.3755515913 & 0.1207057652 \\ 2.1099298228 & -3.1935065022 & -7.1447434841 & 0.3123436692 \end{pmatrix}$$

$$P = \begin{pmatrix} 0 & 1 & 1 & 1 & 1 & 1 & 1 \\ 1 & 0 & -1 & -1 & 1 & 1 & -1 \\ 1 & -1 & 0 & 1 & -1 & -1 & 1 \\ 1 & -1 & 1 & 0 & 1 & -1 & -1 \\ 1 & 1 & -1 & 1 & 0 & -1 & -1 \\ 1 & 1 & -1 & -1 & -1 & 0 & -1 \\ 1 & -1 & 1 & -1 & -1 & -1 & 0 \end{pmatrix} \qquad R = \begin{pmatrix} 0 & 0 & 0 & 0 & 0 & 0 & 0 \\ 0 & 0 & 0 & 0 & 0 & 0 & 0 \\ 1 & 4 & 0 & 1 & 1 & 1 & 4 \\ 1 & 1 & 4 & 0 & 4 & 1 & 1 \\ 0 & 0 & 0 & 0 & 0 & 0 & 0 \\ 4 & 5 & 3 & 4 & 5 & 0 & 3 \\ 4 & 1 & 1 & 4 & 1 & 1 & 0 \end{pmatrix}$$

$|P| = -18$

$I3(P, R) = 16.7368421053$

$I3(-P, R) = 14.2330827068$

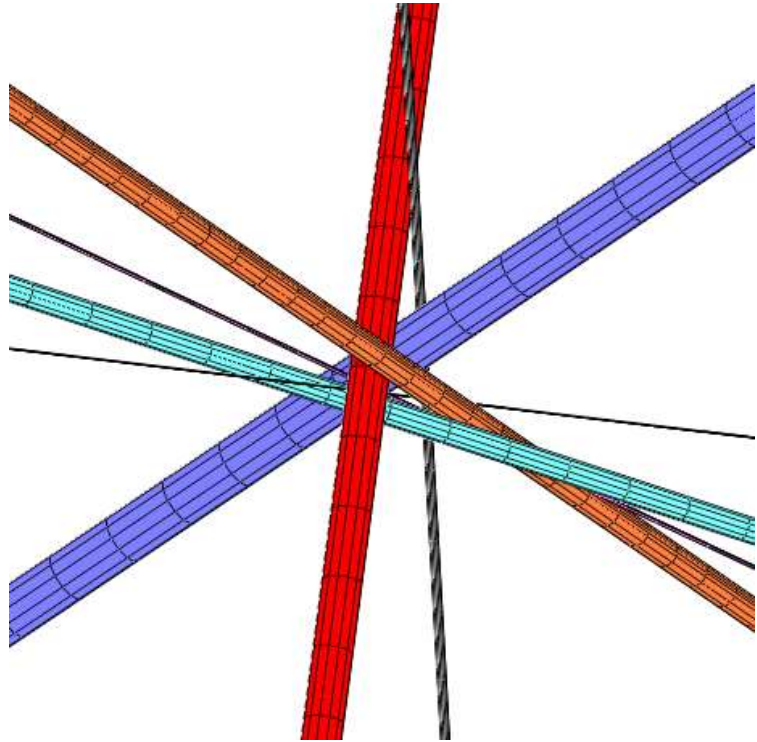

Appendix 4
7-cross
Configuration label me69

$$\begin{pmatrix} \text{orange} \\ \text{red} \\ \text{blue} \\ \text{green} \\ \text{cyan} \\ \text{magenta} \\ \text{gray} \end{pmatrix} \begin{pmatrix} 0 & 0 & 0 & 1 \\ t_1 & p_1 & z_1 & r_1 \\ t_2 & p_2 & z_2 & r_2 \\ t_3 & p_3 & z_3 & r_3 \\ t_4 & p_4 & z_4 & r_4 \\ t_5 & p_5 & z_5 & r_5 \\ t_6 & p_6 & z_6 & r_6 \end{pmatrix} = \begin{pmatrix} 0 & 0 & 0 & 1 \\ 0.8434337352 & 3.1415926536 & 0 & 0.6625731317 \\ 0.4275092599 & 3.3139417528 & 11.1325761992 & 1.5860208535 \\ 2.852323249 & 5.5226770137 & -4.1006057873 & 1.0061146618 \\ 0.3535237438 & 4.1474629414 & 8.4210644889 & 1.1847711232 \\ 2.7630717221 & 0.6702247257 & -2.3991298051 & 4.6496979387 \\ 0.3164930051 & 1.599688017 & -9.6898972492 & 0.4028951014 \end{pmatrix}$$

$$P = \begin{pmatrix} 0 & 1 & 1 & 1 & 1 & 1 & 1 \\ 1 & 0 & -1 & -1 & -1 & -1 & -1 \\ 1 & -1 & 0 & 1 & 1 & -1 & -1 \\ 1 & -1 & 1 & 0 & 1 & -1 & 1 \\ 1 & -1 & 1 & 1 & 0 & 1 & 1 \\ 1 & -1 & -1 & -1 & 1 & 0 & 1 \\ 1 & -1 & -1 & 1 & 1 & 1 & 0 \end{pmatrix} \qquad R = \begin{pmatrix} 0 & 3 & 3 & 5 & 4 & 5 & 4 \\ 3 & 0 & 5 & 3 & 4 & 4 & 5 \\ 0 & 0 & 0 & 0 & 0 & 0 & 0 \\ 0 & 0 & 0 & 0 & 0 & 0 & 0 \\ 1 & 1 & 4 & 4 & 0 & 1 & 1 \\ 0 & 0 & 0 & 0 & 0 & 0 & 0 \\ 0 & 0 & 0 & 0 & 0 & 0 & 0 \end{pmatrix}$$

$|P| = -10$

$I3(P, R) = 7.9$

$I3(-P, R) = 16$

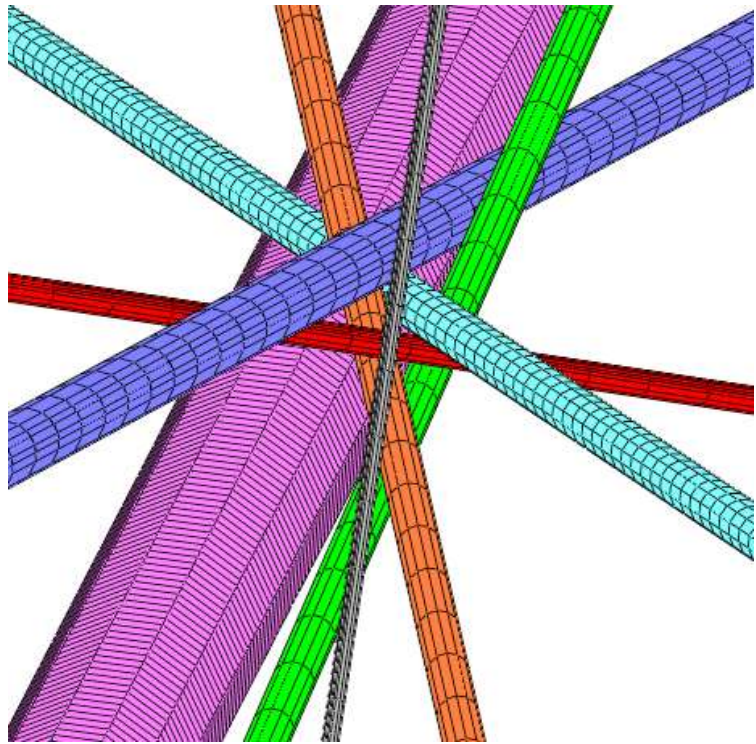

Appendix 4
7-cross
Configuration label  e46

$$\begin{pmatrix} \text{orange} \\ \text{red} \\ \text{blue} \\ \text{green} \\ \text{cyan} \\ \text{magenta} \\ \text{gray} \end{pmatrix} \begin{pmatrix} 0 & 0 & 0 & 1 \\ t_1 & p_1 & z_1 & r_1 \\ t_2 & p_2 & z_2 & r_2 \\ t_3 & p_3 & z_3 & r_3 \\ t_4 & p_4 & z_4 & r_4 \\ t_5 & p_5 & z_5 & r_5 \\ t_6 & p_6 & z_6 & r_6 \end{pmatrix} = \begin{pmatrix} 0 & 0 & 0 & 1 \\ 2.2981589182 & 3.1415926536 & 0 & 0.6625731317 \\ 2.7140833937 & 2.9692435543 & 11.1325762004 & 1.5860208535 \\ 2.7880689098 & 2.1357223657 & 8.4210644902 & 1.1847711232 \\ 0.3785209314 & 5.6129605814 & -2.3991298043 & 4.6496979387 \\ 2.8250996486 & 4.68349729 & -9.6898972492 & 0.4028951014 \\ 2.8783916438 & 2.9483195774 & -2.2945677397 & 0.0926035906 \end{pmatrix}$$

$$P = \begin{pmatrix} 0 & 1 & 1 & 1 & 1 & 1 & 1 \\ 1 & 0 & 1 & 1 & 1 & 1 & -1 \\ 1 & 1 & 0 & -1 & 1 & 1 & -1 \\ 1 & 1 & -1 & 0 & -1 & -1 & 1 \\ 1 & 1 & 1 & -1 & 0 & -1 & -1 \\ 1 & 1 & 1 & -1 & -1 & 0 & -1 \\ 1 & -1 & -1 & 1 & -1 & -1 & 0 \end{pmatrix} \qquad R = \begin{pmatrix} 0 & 1 & 1 & 4 & 4 & 1 & 1 \\ 3 & 0 & 5 & 4 & 4 & 5 & 3 \\ 0 & 0 & 0 & 0 & 0 & 0 & 0 \\ 1 & 1 & 4 & 0 & 1 & 1 & 4 \\ 0 & 0 & 0 & 0 & 0 & 0 & 0 \\ 1 & 1 & 1 & 1 & 4 & 0 & 4 \\ 0 & 0 & 0 & 0 & 0 & 0 & 0 \end{pmatrix}$$

$|P| = -10$

$I3(P, R) = 3.8947368421$

$I3(-P, R) = 17.7368421053$

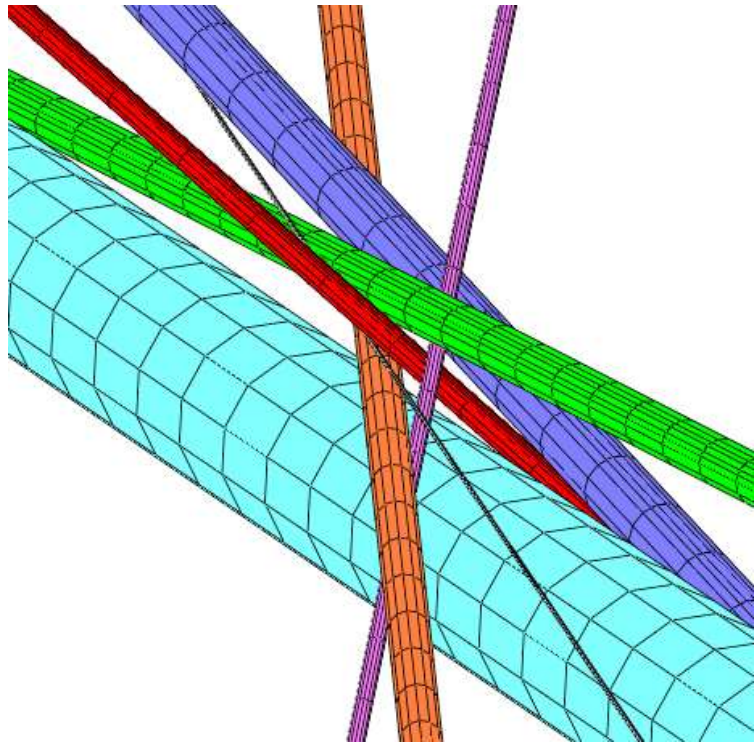

Appendix 4
7-cross
Configuration label  me39

$$\begin{pmatrix} \text{orange} \\ \text{red} \\ \text{blue} \\ \text{green} \\ \text{cyan} \\ \text{magenta} \\ \text{gray} \end{pmatrix} \begin{pmatrix} 0 & 0 & 0 & 1 \\ t1 & p1 & z1 & r1 \\ t2 & p2 & z2 & r2 \\ t3 & p3 & z3 & r3 \\ t4 & p4 & z4 & r4 \\ t5 & p5 & z5 & r5 \\ t6 & p6 & z6 & r6 \end{pmatrix} = \begin{pmatrix} 0 & 0 & 0 & 1 \\ 0.8434337352 & 3.1415926536 & 0 & 0.6625731317 \\ 2.852323249 & 5.5226770137 & -4.1006057872 & 1.0061146618 \\ 0.3535237438 & 4.1474629413 & 8.4210644883 & 1.1847711232 \\ 0.6471527664 & 3.4925662464 & -4.8510541922 & 1.6441357842 \\ 2.7630717221 & 0.6702247257 & -2.3991298052 & 4.6496979387 \\ 0.3164930051 & 1.5996880169 & -9.6898972492 & 0.4028951014 \end{pmatrix}$$

$$P = \begin{pmatrix} 0 & 1 & 1 & 1 & 1 & 1 & 1 \\ 1 & 0 & -1 & -1 & 1 & -1 & -1 \\ 1 & -1 & 0 & 1 & -1 & -1 & 1 \\ 1 & -1 & 1 & 0 & -1 & 1 & 1 \\ 1 & 1 & -1 & -1 & 0 & -1 & 1 \\ 1 & -1 & -1 & 1 & -1 & 0 & 1 \\ 1 & -1 & 1 & 1 & 1 & 1 & 0 \end{pmatrix} \qquad R = \begin{pmatrix} 0 & 3 & 5 & 4 & 3 & 5 & 4 \\ 2 & 0 & 2 & 6 & 3 & 2 & 3 \\ 0 & 0 & 0 & 0 & 0 & 0 & 0 \\ 0 & 0 & 0 & 0 & 0 & 0 & 0 \\ 0 & 0 & 0 & 0 & 0 & 0 & 0 \\ 0 & 0 & 0 & 0 & 0 & 0 & 0 \\ 1 & 1 & 1 & 1 & 4 & 4 & 0 \end{pmatrix}$$

$|P| = -10$

$I3(P, R) = -0.3448275862$

$I3(-P, R) = 22.2068965517$

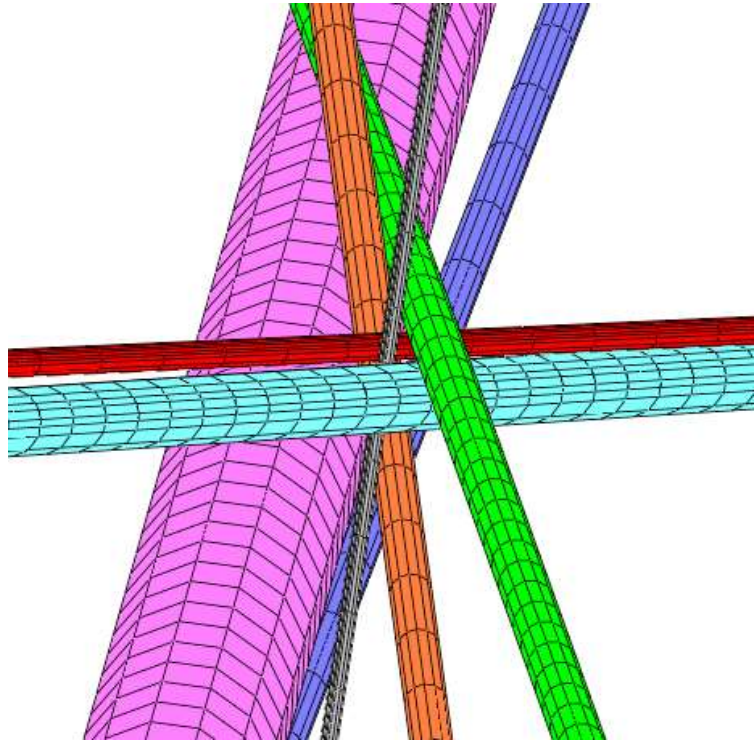

Appendix 4
7-cross
Configuration label  me36

$$\begin{pmatrix} \text{orange} \\ \text{red} \\ \text{blue} \\ \text{green} \\ \text{cyan} \\ \text{magenta} \\ \text{gray} \end{pmatrix} \begin{pmatrix} 0 & 0 & 0 & 1 \\ t1 & p1 & z1 & r1 \\ t2 & p2 & z2 & r2 \\ t3 & p3 & z3 & r3 \\ t4 & p4 & z4 & r4 \\ t5 & p5 & z5 & r5 \\ t6 & p6 & z6 & r6 \end{pmatrix} = \begin{pmatrix} 0 & 0 & 0 & 1 \\ 0.8434337352 & 3.1415926536 & 0 & 0.6625731317 \\ 2.852323249 & 5.5226770137 & -4.1006057872 & 1.0061146618 \\ 0.3535237438 & 4.1474629413 & 8.421064488 & 1.1847711232 \\ 2.7630717221 & 0.6702247257 & -2.3991298051 & 4.6496979387 \\ 0.3164930051 & 1.5996880169 & -9.6898972492 & 0.4028951014 \\ 0.2632010099 & 3.3348657297 & -2.2945677391 & 0.0926035906 \end{pmatrix}$$

$$P = \begin{pmatrix} 0 & 1 & 1 & 1 & 1 & 1 & 1 \\ 1 & 0 & -1 & -1 & -1 & -1 & 1 \\ 1 & -1 & 0 & 1 & -1 & 1 & 1 \\ 1 & -1 & 1 & 0 & 1 & 1 & -1 \\ 1 & -1 & -1 & 1 & 0 & 1 & 1 \\ 1 & -1 & 1 & 1 & 1 & 0 & 1 \\ 1 & 1 & 1 & -1 & 1 & 1 & 0 \end{pmatrix} \qquad R = \begin{pmatrix} 0 & 3 & 5 & 4 & 5 & 4 & 3 \\ 1 & 0 & 1 & 4 & 1 & 4 & 1 \\ 0 & 0 & 0 & 0 & 0 & 0 & 0 \\ 0 & 0 & 0 & 0 & 0 & 0 & 0 \\ 0 & 0 & 0 & 0 & 0 & 0 & 0 \\ 1 & 1 & 1 & 1 & 4 & 0 & 4 \\ 2 & 3 & 6 & 3 & 2 & 2 & 0 \end{pmatrix}$$

$|P| = -10$

$I3(P, R) = 10.8601398601$

$I3(-P, R) = 17.8041958042$

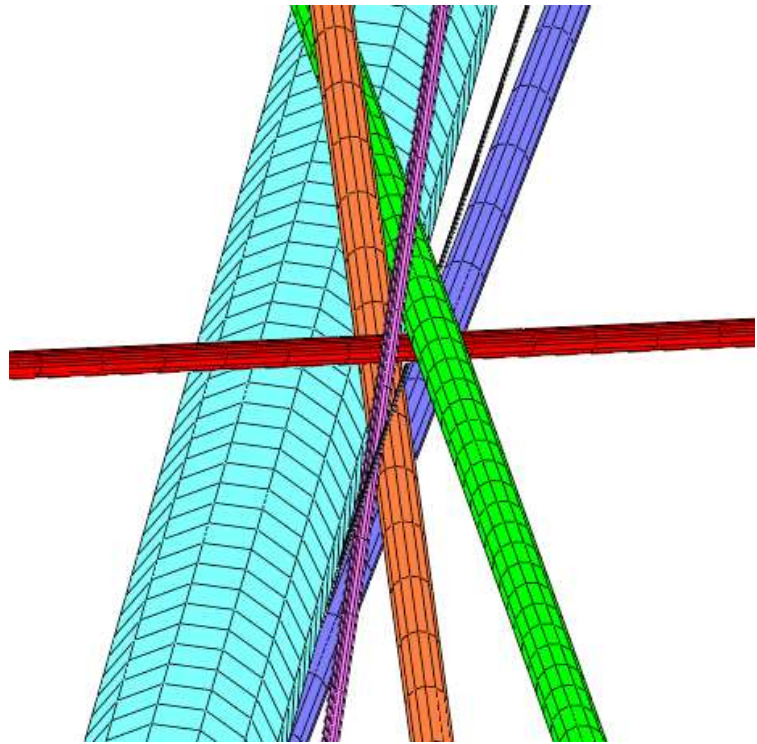

Appendix 4
7-cross
Configuration label  md67

$$\begin{pmatrix} \text{orange} \\ \text{red} \\ \text{blue} \\ \text{green} \\ \text{cyan} \\ \text{magenta} \\ \text{gray} \end{pmatrix} \begin{pmatrix} 0 & 0 & 0 & 1 \\ t_1 & p_1 & z_1 & r_1 \\ t_2 & p_2 & z_2 & r_2 \\ t_3 & p_3 & z_3 & r_3 \\ t_4 & p_4 & z_4 & r_4 \\ t_5 & p_5 & z_5 & r_5 \\ t_6 & p_6 & z_6 & r_6 \end{pmatrix} = \begin{pmatrix} 0 & 0 & 0 & 1 \\ 2.2344557833 & 3.1415926536 & 0 & 0.4424161588 \\ 1.493845935 & 3.5512530273 & 1.9436660557 & 0.7245086968 \\ 1.891605667 & 5.6917260248 & -2.7377466038 & 2.7041739475 \\ 1.4042016164 & 4.4070605926 & 0.6658298674 & 0.7245086968 \\ 0.7745288109 & 1.9462015992 & -2.0960004926 & 0.7245086969 \\ 1.0529825429 & 3.6017847364 & -1.1042053297 & 0.1548918413 \end{pmatrix}$$

$$P = \begin{pmatrix} 0 & 1 & 1 & 1 & 1 & 1 & 1 \\ 1 & 0 & -1 & -1 & -1 & -1 & 1 \\ 1 & -1 & 0 & 1 & 1 & -1 & 1 \\ 1 & -1 & 1 & 0 & 1 & 1 & 1 \\ 1 & -1 & 1 & 1 & 0 & 1 & -1 \\ 1 & -1 & -1 & 1 & 1 & 0 & 1 \\ 1 & 1 & 1 & 1 & -1 & 1 & 0 \end{pmatrix} \qquad R = \begin{pmatrix} 0 & 1 & 1 & 4 & 1 & 4 & 1 \\ 2 & 0 & 3 & 2 & 3 & 6 & 2 \\ 0 & 0 & 0 & 0 & 0 & 0 & 0 \\ 0 & 0 & 0 & 0 & 0 & 0 & 0 \\ 2 & 2 & 6 & 3 & 0 & 2 & 3 \\ 0 & 0 & 0 & 0 & 0 & 0 & 0 \\ 2 & 3 & 3 & 6 & 2 & 2 & 0 \end{pmatrix}$$

$|P| = -10$

$I3(P, R) = 10.1384615385$

$I3(-P, R) = 15.8461538462$

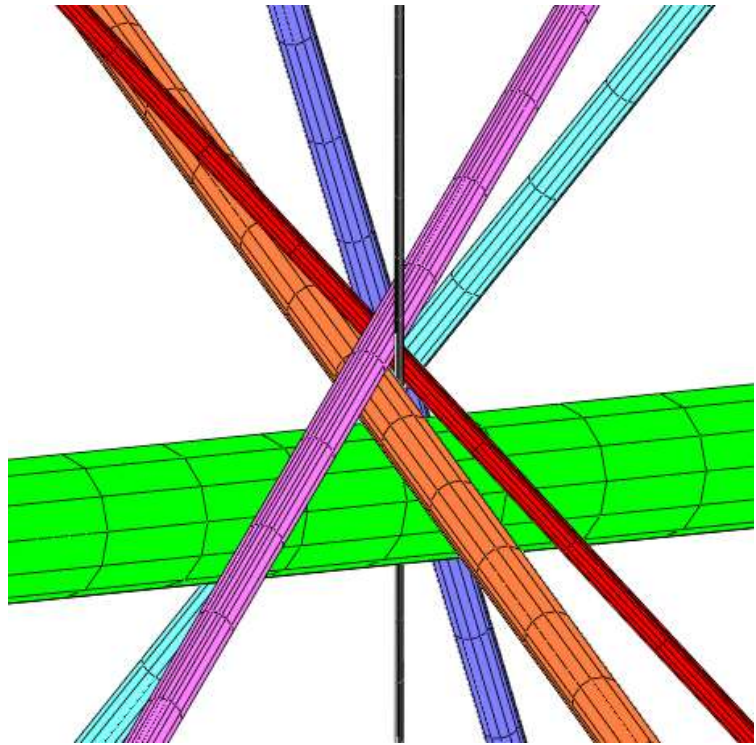

Appendix 4
7-cross
Configuration label  ma69

$$\begin{pmatrix} \text{orange} \\ \text{red} \\ \text{blue} \\ \text{green} \\ \text{cyan} \\ \text{magenta} \\ \text{gray} \end{pmatrix} \begin{pmatrix} 0 & 0 & 0 & 1 \\ t1 & p1 & z1 & r1 \\ t2 & p2 & z2 & r2 \\ t3 & p3 & z3 & r3 \\ t4 & p4 & z4 & r4 \\ t5 & p5 & z5 & r5 \\ t6 & p6 & z6 & r6 \end{pmatrix} = \begin{pmatrix} 0 & 0 & 0 & 1 \\ 2.4350672766 & 3.1415926536 & 0 & 1 \\ 2.3306212912 & -3.0803765257 & 7.6740161121 & 1 \\ 1.2837798184 & -0.1197866094 & -14.3813154944 & 1 \\ 2.4618520543 & -2.9203250131 & 2.6628906163 & 1 \\ 0.6223654767 & 0.2513273739 & -7.2349214808 & 1 \\ 0.4102096255 & 6.6075059922 & -17.6852161963 & 0.42 \end{pmatrix}$$

$$P = \begin{pmatrix} 0 & 1 & 1 & 1 & 1 & 1 & 1 \\ 1 & 0 & -1 & -1 & -1 & -1 & -1 \\ 1 & -1 & 0 & 1 & 1 & -1 & -1 \\ 1 & -1 & 1 & 0 & 1 & -1 & 1 \\ 1 & -1 & 1 & 1 & 0 & 1 & 1 \\ 1 & -1 & -1 & -1 & 1 & 0 & 1 \\ 1 & -1 & -1 & 1 & 1 & 1 & 0 \end{pmatrix} \qquad R = \begin{pmatrix} 0 & 2 & 2 & 6 & 2 & 3 & 3 \\ 4 & 0 & 4 & 4 & 4 & 4 & 4 \\ 0 & 0 & 0 & 0 & 0 & 0 & 0 \\ 0 & 0 & 0 & 0 & 0 & 0 & 0 \\ 1 & 1 & 4 & 4 & 0 & 1 & 1 \\ 0 & 0 & 0 & 0 & 0 & 0 & 0 \\ 1 & 1 & 1 & 1 & 4 & 4 & 0 \end{pmatrix}$$

$|P| = -10$

$I3(P, R) = 11.2151898734$

$I3(-P, R) = 12.1772151899$

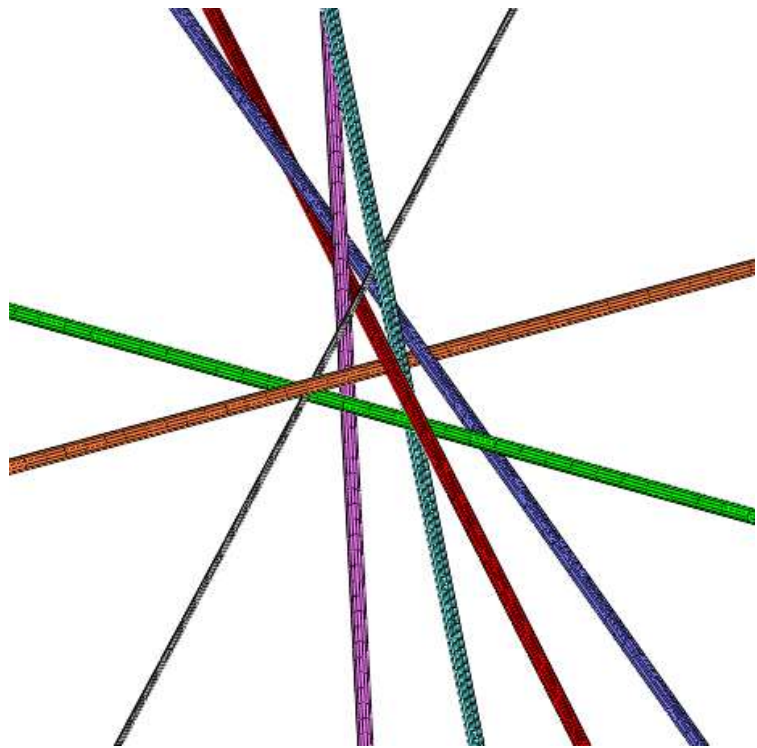

Appendix 4
7-cross
Configuration label  ma67

$$\begin{pmatrix} \text{orange} \\ \text{red} \\ \text{blue} \\ \text{green} \\ \text{cyan} \\ \text{magenta} \\ \text{gray} \end{pmatrix} \begin{pmatrix} 0 & 0 & 0 & 1 \\ t1 & p1 & z1 & r1 \\ t2 & p2 & z2 & r2 \\ t3 & p3 & z3 & r3 \\ t4 & p4 & z4 & r4 \\ t5 & p5 & z5 & r5 \\ t6 & p6 & z6 & r6 \end{pmatrix} = \begin{pmatrix} 0 & 0 & 0 & 1 \\ 2.9680871265 & 3.1415926536 & 0 & 1 \\ 2.4656809729 & -2.6433270855 & 16.0186869199 & 1 \\ 3.0215413671 & -1.0230355358 & -21.0490886869 & 1 \\ 0.7595430869 & -2.0193885871 & 8.1348382419 & 1 \\ 0.1044514939 & 7.9156563623 & -30.1645938871 & 1 \\ 0.2122535362 & 3.7432271956 & -8.8710230251 & 0.22 \end{pmatrix}$$

$$P = \begin{pmatrix} 0 & 1 & 1 & 1 & 1 & 1 & 1 \\ 1 & 0 & -1 & -1 & -1 & -1 & 1 \\ 1 & -1 & 0 & 1 & 1 & -1 & 1 \\ 1 & -1 & 1 & 0 & 1 & 1 & 1 \\ 1 & -1 & 1 & 1 & 0 & 1 & -1 \\ 1 & -1 & -1 & 1 & 1 & 0 & 1 \\ 1 & 1 & 1 & 1 & -1 & 1 & 0 \end{pmatrix} \qquad R = \begin{pmatrix} 0 & 1 & 1 & 4 & 1 & 4 & 1 \\ 2 & 0 & 3 & 2 & 3 & 6 & 2 \\ 0 & 0 & 0 & 0 & 0 & 0 & 0 \\ 0 & 0 & 0 & 0 & 0 & 0 & 0 \\ 2 & 2 & 6 & 3 & 0 & 2 & 3 \\ 0 & 0 & 0 & 0 & 0 & 0 & 0 \\ 3 & 2 & 2 & 6 & 2 & 3 & 0 \end{pmatrix}$$

$|P| = -10$

$I3(P, R) = 8.3770491803$

$I3(-P, R) = 15.9672131148$

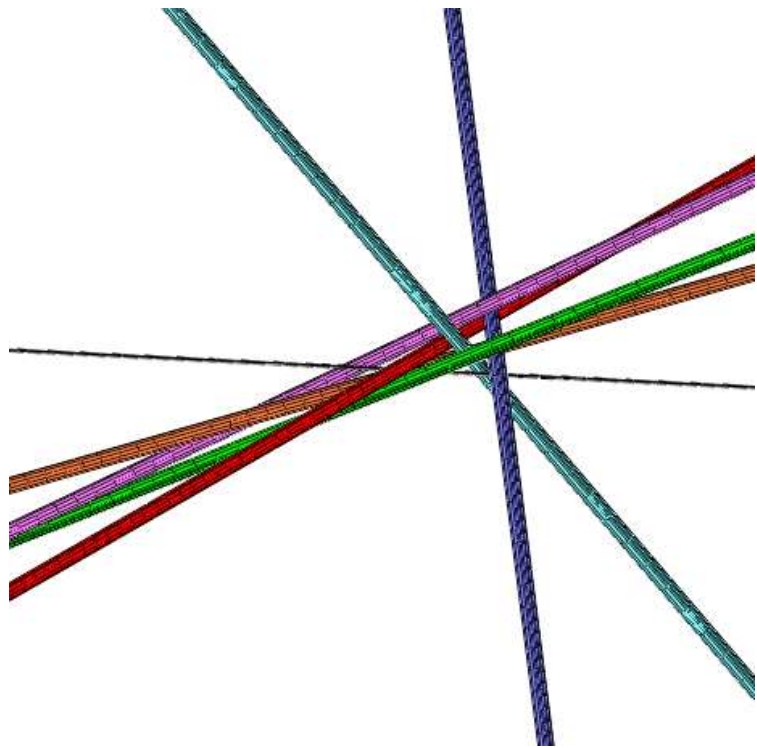

Appendix 4
7-cross
Configuration label  ma39

$$\begin{pmatrix} \text{orange} \\ \text{red} \\ \text{blue} \\ \text{green} \\ \text{cyan} \\ \text{magenta} \\ \text{gray} \end{pmatrix} \begin{pmatrix} 0 & 0 & 0 & 1 \\ t_1 & p_1 & z_1 & r_1 \\ t_2 & p_2 & z_2 & r_2 \\ t_3 & p_3 & z_3 & r_3 \\ t_4 & p_4 & z_4 & r_4 \\ t_5 & p_5 & z_5 & r_5 \\ t_6 & p_6 & z_6 & r_6 \end{pmatrix} = \begin{pmatrix} 0 & 0 & 0 & 1 \\ 2.7374809407 & 3.1415926536 & 0 & 1 \\ 1.0463047815 & -0.2698607451 & -12.3689131406 & 1 \\ 2.703525161 & -2.7910472102 & 3.5146536054 & 1 \\ 2.7596167402 & -2.8672770393 & -6.7151801635 & 1 \\ 0.4190435343 & 0.3754305768 & -4.1179263214 & 1 \\ 0.2758141584 & 6.7826908518 & -14.5848370812 & 0.5015750769 \end{pmatrix}$$

$$P = \begin{pmatrix} 0 & 1 & 1 & 1 & 1 & 1 & 1 \\ 1 & 0 & -1 & -1 & 1 & -1 & -1 \\ 1 & -1 & 0 & 1 & -1 & -1 & 1 \\ 1 & -1 & 1 & 0 & -1 & 1 & 1 \\ 1 & 1 & -1 & -1 & 0 & -1 & 1 \\ 1 & -1 & -1 & 1 & -1 & 0 & 1 \\ 1 & -1 & 1 & 1 & 1 & 1 & 0 \end{pmatrix} \qquad R = \begin{pmatrix} 0 & 2 & 6 & 2 & 2 & 3 & 3 \\ 3 & 0 & 3 & 6 & 2 & 2 & 2 \\ 0 & 0 & 0 & 0 & 0 & 0 & 0 \\ 0 & 0 & 0 & 0 & 0 & 0 & 0 \\ 1 & 1 & 1 & 4 & 0 & 4 & 1 \\ 0 & 0 & 0 & 0 & 0 & 0 & 0 \\ 2 & 2 & 2 & 3 & 3 & 6 & 0 \end{pmatrix}$$

$|P| = -10$

$I3(P, R) = 2.5573770492$

$I3(-P, R) = 21.7868852459$

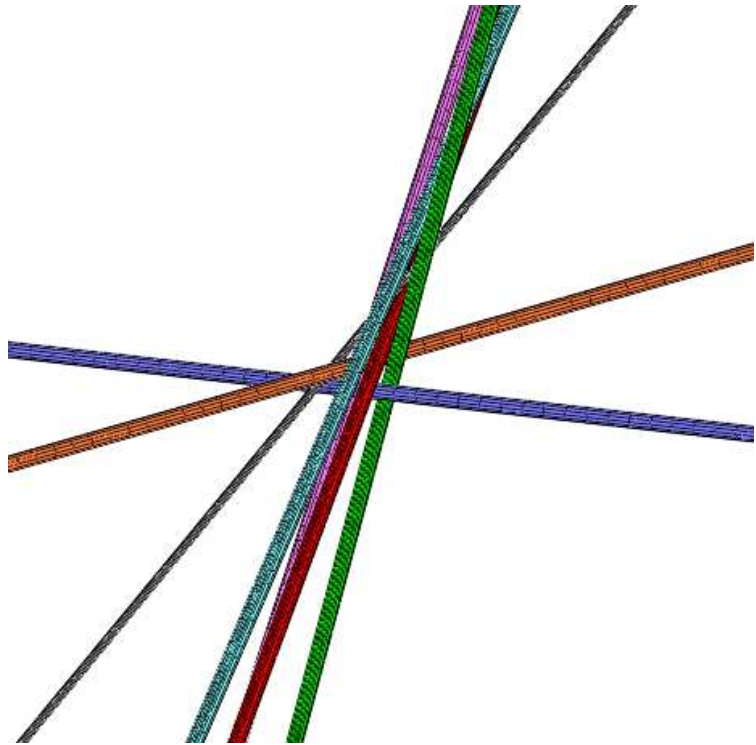

Appendix 4
7-cross
Configuration label  ma36

$$\begin{pmatrix} \text{orange} \\ \text{red} \\ \text{blue} \\ \text{green} \\ \text{cyan} \\ \text{magenta} \\ \text{gray} \end{pmatrix} \begin{pmatrix} 0 & 0 & 0 & 1 \\ t_1 & p_1 & z_1 & r_1 \\ t_2 & p_2 & z_2 & r_2 \\ t_3 & p_3 & z_3 & r_3 \\ t_4 & p_4 & z_4 & r_4 \\ t_5 & p_5 & z_5 & r_5 \\ t_6 & p_6 & z_6 & r_6 \end{pmatrix} = \begin{pmatrix} 0 & 0 & 0 & 1 \\ 1.815649312 & 3.1415926536 & 0 & 1 \\ 3.1010695435 & -0.5013348391 & -91.0269771184 & 1 \\ 1.4528433228 & -1.9434430426 & 1.9306092387 & 1 \\ 0.0563311292 & 1.674572818 & -67.2844777746 & 1 \\ 0.03185402 & 8.1135418711 & -92.7747858406 & 0.1561241855 \\ 1.4075891494 & 4.3135291794 & -2.1104319155 & 0.9 \end{pmatrix}$$

$$P = \begin{pmatrix} 0 & 1 & 1 & 1 & 1 & 1 & 1 \\ 1 & 0 & -1 & -1 & -1 & -1 & 1 \\ 1 & -1 & 0 & 1 & -1 & 1 & 1 \\ 1 & -1 & 1 & 0 & 1 & 1 & -1 \\ 1 & -1 & -1 & 1 & 0 & 1 & 1 \\ 1 & -1 & 1 & 1 & 1 & 0 & 1 \\ 1 & 1 & 1 & -1 & 1 & 1 & 0 \end{pmatrix} \qquad R = \begin{pmatrix} 0 & 2 & 6 & 2 & 3 & 3 & 2 \\ 2 & 0 & 2 & 6 & 3 & 3 & 2 \\ 0 & 0 & 0 & 0 & 0 & 0 & 0 \\ 0 & 0 & 0 & 0 & 0 & 0 & 0 \\ 0 & 0 & 0 & 0 & 0 & 0 & 0 \\ 2 & 2 & 2 & 3 & 6 & 0 & 3 \\ 2 & 3 & 6 & 3 & 2 & 2 & 0 \end{pmatrix}$$

$|P| = -10$

$I3(P, R) = 10.7094017094$

$I3(-P, R) = 16.3418803419$

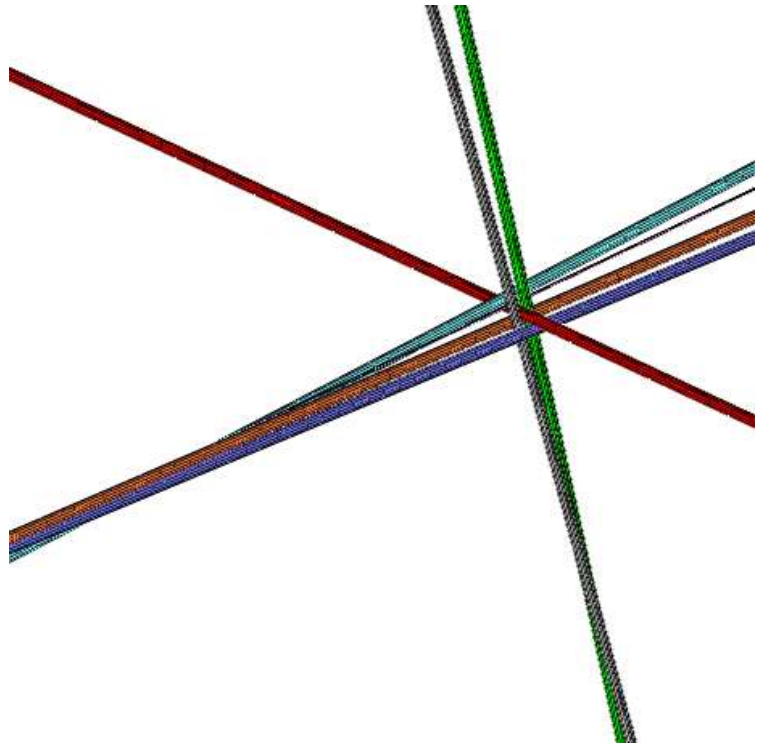

Appendix 4
7-cross
Configuration label  ma27

$$\begin{pmatrix} \text{orange} \\ \text{red} \\ \text{blue} \\ \text{green} \\ \text{cyan} \\ \text{magenta} \\ \text{gray} \end{pmatrix} \begin{pmatrix} 0 & 0 & 0 & 1 \\ t1 & p1 & z1 & r1 \\ t2 & p2 & z2 & r2 \\ t3 & p3 & z3 & r3 \\ t4 & p4 & z4 & r4 \\ t5 & p5 & z5 & r5 \\ t6 & p6 & z6 & r6 \end{pmatrix} = \begin{pmatrix} 0 & 0 & 0 & 1 \\ 0.8670725547 & 3.1415926536 & 0 & 1 \\ 2.9471707151 & 5.0396314482 & -22.2887966512 & 1 \\ 2.3055759111 & 4.3795424785 & -3.7010800277 & 1 \\ 2.9620830068 & 3.9095881359 & -21.5844134055 & 1 \\ 0.1930653908 & 14.456658727 & -15.5915792155 & 0.2018724285 \\ 0.4265892635 & 9.5112644416 & -10.8742649147 & 0.35 \end{pmatrix}$$

$$P = \begin{pmatrix} 0 & 1 & 1 & 1 & 1 & 1 & 1 \\ 1 & 0 & 1 & 1 & 1 & -1 & 1 \\ 1 & 1 & 0 & 1 & -1 & 1 & 1 \\ 1 & 1 & 1 & 0 & -1 & 1 & -1 \\ 1 & 1 & -1 & -1 & 0 & 1 & 1 \\ 1 & -1 & 1 & 1 & 1 & 0 & 1 \\ 1 & 1 & 1 & -1 & 1 & 1 & 0 \end{pmatrix} \qquad R = \begin{pmatrix} 0 & 1 & 4 & 1 & 1 & 4 & 1 \\ 0 & 0 & 0 & 0 & 0 & 0 & 0 \\ 0 & 0 & 0 & 0 & 0 & 0 & 0 \\ 2 & 6 & 3 & 0 & 2 & 2 & 3 \\ 0 & 0 & 0 & 0 & 0 & 0 & 0 \\ 1 & 4 & 1 & 1 & 4 & 0 & 1 \\ 5 & 3 & 5 & 3 & 4 & 4 & 0 \end{pmatrix}$$

$|P| = -10$

$I3(P, R) = -0.7073170732$

$I3(-P, R) = 24.5731707317$

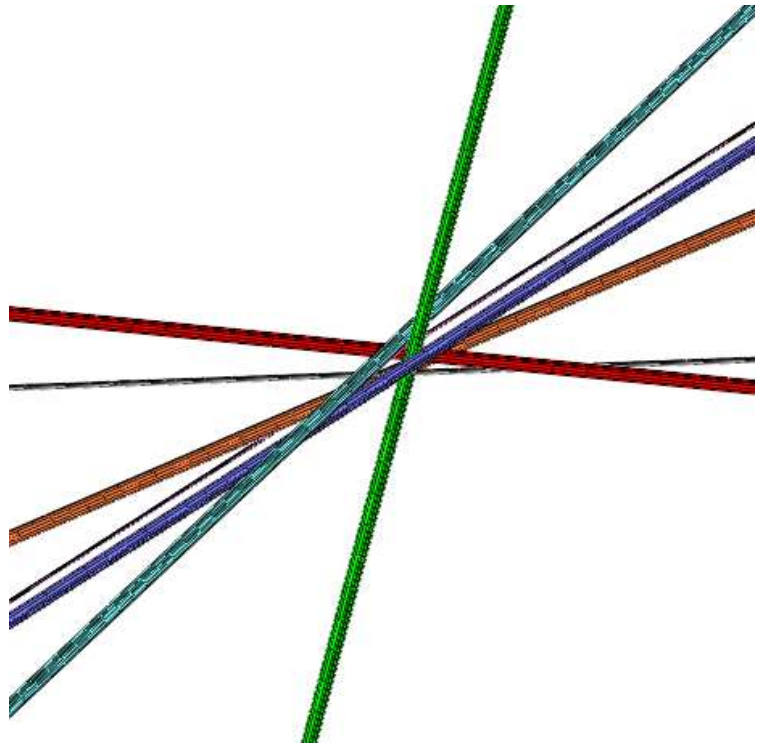

Appendix 4
7-cross
Configuration label ma26

$$\begin{pmatrix} \text{orange} \\ \text{red} \\ \text{blue} \\ \text{green} \\ \text{cyan} \\ \text{magenta} \\ \text{gray} \end{pmatrix} \begin{pmatrix} 0 & 0 & 0 & 1 \\ t_1 & p_1 & z_1 & r_1 \\ t_2 & p_2 & z_2 & r_2 \\ t_3 & p_3 & z_3 & r_3 \\ t_4 & p_4 & z_4 & r_4 \\ t_5 & p_5 & z_5 & r_5 \\ t_6 & p_6 & z_6 & r_6 \end{pmatrix} = \begin{pmatrix} 0 & 0 & 0 & 1 \\ 2.9992963375 & 3.1415926536 & 0 & 1 \\ 3.0697117671 & 5.3858609763 & -135.4537848705 & 1 \\ 2.3792005688 & 3.3643725978 & -55.8140165278 & 1 \\ 0.0507457798 & 7.6599867482 & -70.8575344312 & 1 \\ 0.0312215801 & 14.0694186276 & -94.2720843883 & 0.1255909636 \\ 2.3827817956 & 9.5926170298 & -58.3001663582 & 0.8 \end{pmatrix}$$

$$P = \begin{pmatrix} 0 & 1 & 1 & 1 & 1 & 1 & 1 \\ 1 & 0 & 1 & 1 & -1 & -1 & 1 \\ 1 & 1 & 0 & 1 & -1 & 1 & 1 \\ 1 & 1 & 1 & 0 & 1 & 1 & -1 \\ 1 & -1 & -1 & 1 & 0 & 1 & 1 \\ 1 & -1 & 1 & 1 & 1 & 0 & 1 \\ 1 & 1 & 1 & -1 & 1 & 1 & 0 \end{pmatrix} \qquad R = \begin{pmatrix} 0 & 2 & 6 & 2 & 3 & 3 & 2 \\ 0 & 0 & 0 & 0 & 0 & 0 & 0 \\ 0 & 0 & 0 & 0 & 0 & 0 & 0 \\ 2 & 6 & 3 & 0 & 2 & 2 & 3 \\ 0 & 0 & 0 & 0 & 0 & 0 & 0 \\ 2 & 2 & 2 & 3 & 6 & 0 & 3 \\ 2 & 3 & 6 & 3 & 2 & 2 & 0 \end{pmatrix}$$

$|P| = -10$

$I3(P, R) = 12.6612903226$

$I3(-P, R) = 15.2903225806$

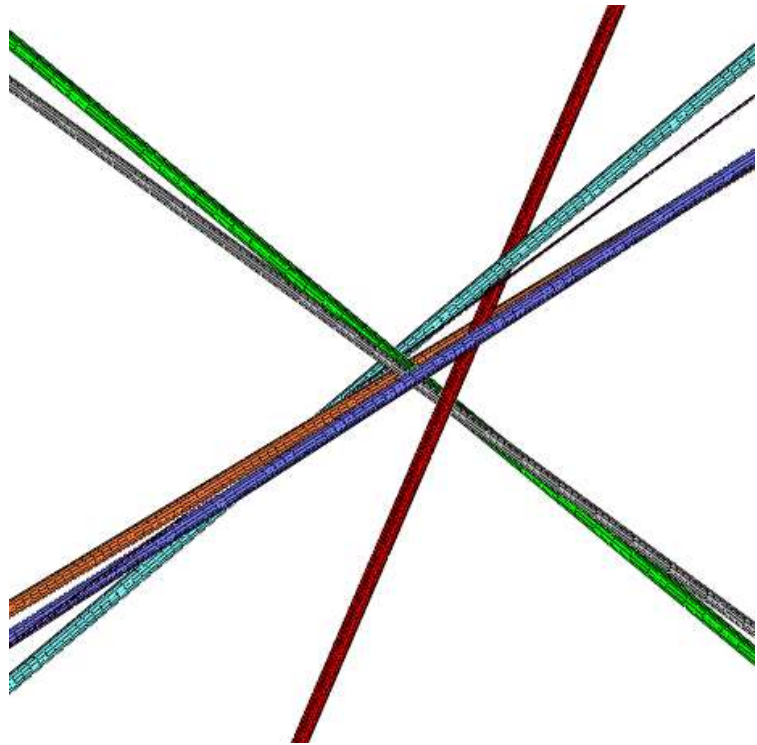

Appendix 4
7-cross
Configuration label  a58

$$\begin{pmatrix} \text{orange} \\ \text{red} \\ \text{blue} \\ \text{green} \\ \text{cyan} \\ \text{magenta} \\ \text{gray} \end{pmatrix} \begin{pmatrix} 0 & 0 & 0 & 1 \\ t_1 & p_1 & z_1 & r_1 \\ t_2 & p_2 & z_2 & r_2 \\ t_3 & p_3 & z_3 & r_3 \\ t_4 & p_4 & z_4 & r_4 \\ t_5 & p_5 & z_5 & r_5 \\ t_6 & p_6 & z_6 & r_6 \end{pmatrix} = \begin{pmatrix} 0 & 0 & 0 & 1 \\ 0.5096303691 & 3.1415926536 & 0 & 1 \\ 0.3766959277 & 8.5250650362 & 1.2958990871 & 1 \\ 0.4150556985 & 6.6831650665 & -19.5123392101 & 1 \\ 0.2379693866 & 8.2983149218 & -15.9506530822 & 1 \\ 2.8293243125 & 4.9213200305 & -10.1489077064 & 1 \\ 0.4796952391 & 2.1870970378 & -6.3368384738 & 0.32 \end{pmatrix}$$

$$P = \begin{pmatrix} 0 & 1 & 1 & 1 & 1 & 1 & 1 \\ 1 & 0 & 1 & 1 & -1 & 1 & -1 \\ 1 & 1 & 0 & -1 & -1 & 1 & -1 \\ 1 & 1 & -1 & 0 & 1 & 1 & -1 \\ 1 & -1 & -1 & 1 & 0 & 1 & -1 \\ 1 & 1 & 1 & 1 & 1 & 0 & -1 \\ 1 & -1 & -1 & -1 & -1 & -1 & 0 \end{pmatrix} \qquad R = \begin{pmatrix} 0 & 1 & 1 & 4 & 1 & 4 & 1 \\ 2 & 0 & 6 & 2 & 3 & 3 & 2 \\ 0 & 0 & 0 & 0 & 0 & 0 & 0 \\ 0 & 0 & 0 & 0 & 0 & 0 & 0 \\ 0 & 0 & 0 & 0 & 0 & 0 & 0 \\ 0 & 0 & 0 & 0 & 0 & 0 & 0 \\ 5 & 3 & 3 & 5 & 4 & 4 & 0 \end{pmatrix}$$

$|P| = -10$

$I3(P, R) = 1.7142857143$

$I3(-P, R) = 19.8571428571$

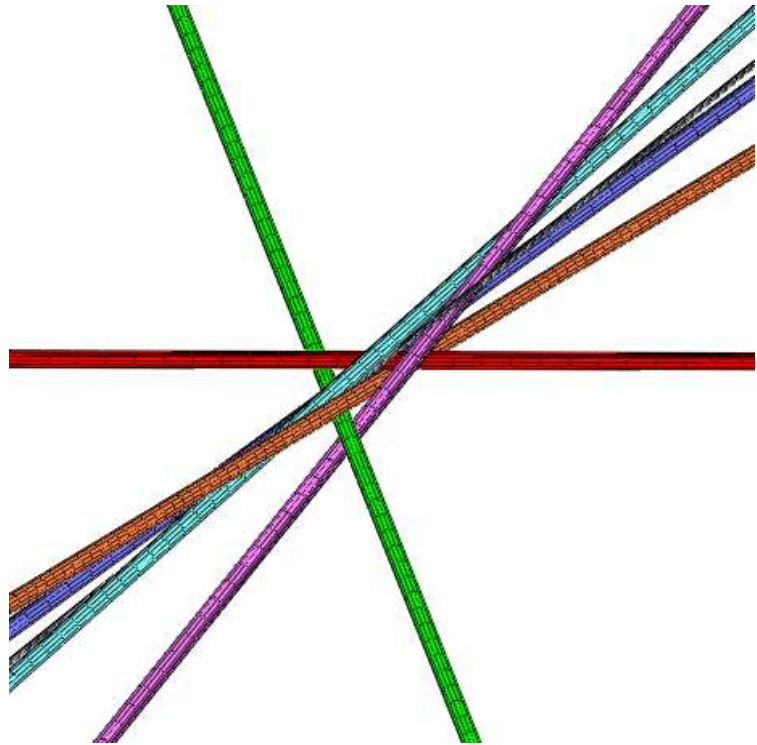

Appendix 4
7-cross
Configuration label  a56

$$\begin{pmatrix} \text{orange} \\ \text{red} \\ \text{blue} \\ \text{green} \\ \text{cyan} \\ \text{magenta} \\ \text{gray} \end{pmatrix} \begin{pmatrix} 0 & 0 & 0 & 1 \\ t_1 & p_1 & z_1 & r_1 \\ t_2 & p_2 & z_2 & r_2 \\ t_3 & p_3 & z_3 & r_3 \\ t_4 & p_4 & z_4 & r_4 \\ t_5 & p_5 & z_5 & r_5 \\ t_6 & p_6 & z_6 & r_6 \end{pmatrix} = \begin{pmatrix} 0 & 0 & 0 & 1 \\ 1.1109456268 & 3.1415926536 & 0 & 1 \\ 0.5277928906 & 8.4500165895 & 1.7854883041 & 1 \\ 0.2927269363 & 6.6606983808 & -19.3012389549 & 1 \\ 2.8485296782 & 4.3985915219 & -11.7210902945 & 1 \\ 2.955444361 & -2.0452911715 & -16.0333171763 & 0.18 \\ 0.840433757 & 1.8181873884 & -4.910955917 & 1 \end{pmatrix}$$

$$P = \begin{pmatrix} 0 & 1 & 1 & 1 & 1 & 1 & 1 \\ 1 & 0 & 1 & 1 & 1 & 1 & -1 \\ 1 & 1 & 0 & -1 & 1 & 1 & -1 \\ 1 & 1 & -1 & 0 & 1 & -1 & -1 \\ 1 & 1 & 1 & 1 & 0 & -1 & -1 \\ 1 & 1 & 1 & -1 & -1 & 0 & -1 \\ 1 & -1 & -1 & -1 & -1 & -1 & 0 \end{pmatrix} \qquad R = \begin{pmatrix} 0 & 2 & 2 & 6 & 3 & 3 & 2 \\ 2 & 0 & 6 & 2 & 3 & 3 & 2 \\ 0 & 0 & 0 & 0 & 0 & 0 & 0 \\ 0 & 0 & 0 & 0 & 0 & 0 & 0 \\ 0 & 0 & 0 & 0 & 0 & 0 & 0 \\ 1 & 1 & 1 & 1 & 4 & 0 & 4 \\ 2 & 3 & 3 & 6 & 2 & 2 & 0 \end{pmatrix}$$

$|P| = -10$

$I3(P, R) = 8.4242424242$

$I3(-P, R) = 19.4025974026$

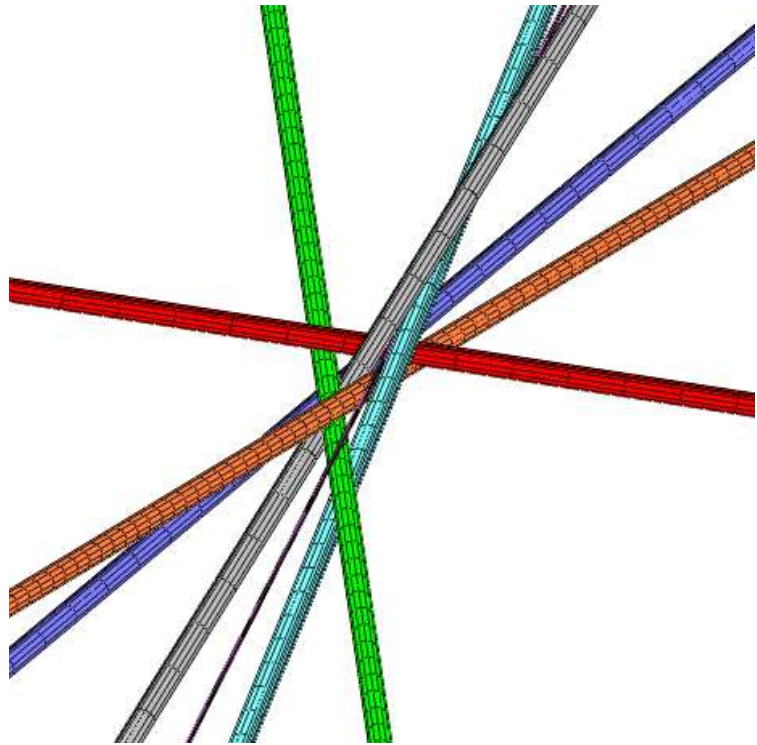

Appendix 4
7-cross
Configuration label e13

$$\begin{pmatrix} \text{orange} \\ \text{red} \\ \text{blue} \\ \text{green} \\ \text{cyan} \\ \text{magenta} \\ \text{gray} \end{pmatrix} \begin{pmatrix} 0 & 0 & 0 & 1 \\ t_1 & p_1 & z_1 & r_1 \\ t_2 & p_2 & z_2 & r_2 \\ t_3 & p_3 & z_3 & r_3 \\ t_4 & p_4 & z_4 & r_4 \\ t_5 & p_5 & z_5 & r_5 \\ t_6 & p_6 & z_6 & r_6 \end{pmatrix} = \begin{pmatrix} 0 & 0 & 0 & 1 \\ 2.4842823378 & 3.1415926536 & 0 & 1.5184960175 \\ 0.7040583282 & 5.4868383704 & 4.6701097909 & 1.788136383 \\ -0.3066658752 & 2.057068686 & -6.7132406797 & 2.4814404713 \\ 2.5573886738 & 2.3846938085 & 20.0511391079 & 7.0176373237 \\ 0.8781268932 & 0.0885453641 & 10.0384682835 & 0.6080764253 \\ -0.5867999455 & 2.723580373 & 4.7750878577 & 0.1397635765 \end{pmatrix}$$

$$P = \begin{pmatrix} 0 & 1 & 1 & -1 & 1 & 1 & -1 \\ 1 & 0 & -1 & 1 & 1 & -1 & -1 \\ 1 & -1 & 0 & 1 & -1 & -1 & 1 \\ -1 & 1 & 1 & 0 & 1 & -1 & -1 \\ 1 & 1 & -1 & 1 & 0 & -1 & -1 \\ 1 & -1 & -1 & -1 & -1 & 0 & -1 \\ -1 & -1 & 1 & -1 & -1 & -1 & 0 \end{pmatrix} \qquad R = \begin{pmatrix} 0 & 2 & 6 & 3 & 2 & 3 & 2 \\ 0 & 0 & 0 & 0 & 0 & 0 & 0 \\ 0 & 0 & 0 & 0 & 0 & 0 & 0 \\ 0 & 0 & 0 & 0 & 0 & 0 & 0 \\ 0 & 0 & 0 & 0 & 0 & 0 & 0 \\ 2 & 2 & 2 & 3 & 6 & 0 & 3 \\ 2 & 6 & 2 & 2 & 3 & 3 & 0 \end{pmatrix}$$

$|P| = -2$

$I3(P, R) = 15.9166666667$

$I3(-P, R) = 14.4166666667$

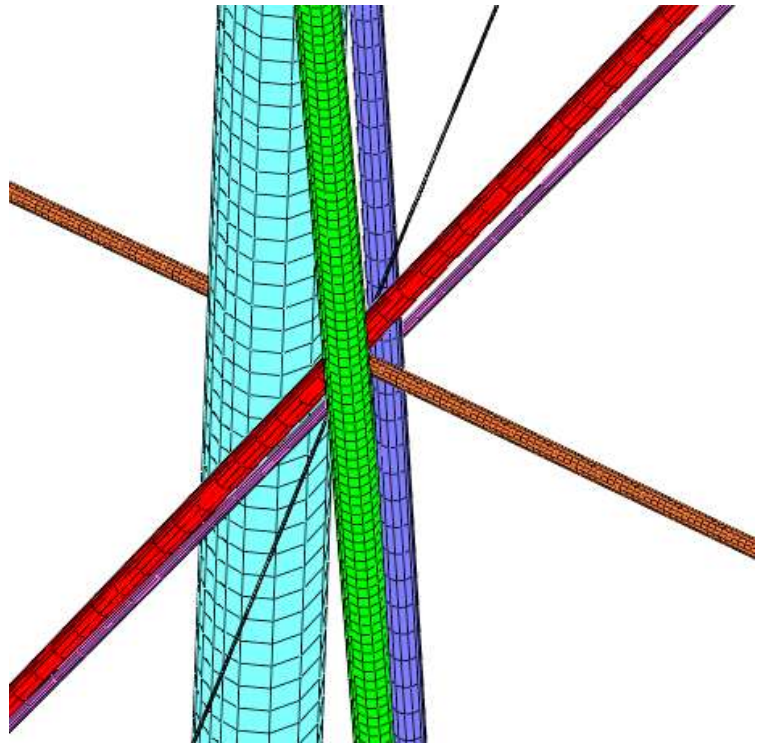

Appendix 4
7-cross
Configuration label c69

$$\begin{pmatrix} \text{orange} \\ \text{red} \\ \text{blue} \\ \text{green} \\ \text{cyan} \\ \text{magenta} \\ \text{gray} \end{pmatrix} \begin{pmatrix} 0 & 0 & 0 & 1 \\ t_1 & p_1 & z_1 & r_1 \\ t_2 & p_2 & z_2 & r_2 \\ t_3 & p_3 & z_3 & r_3 \\ t_4 & p_4 & z_4 & r_4 \\ t_5 & p_5 & z_5 & r_5 \\ t_6 & p_6 & z_6 & r_6 \end{pmatrix} = \begin{pmatrix} 0 & 0 & 0 & 1 \\ 1.3624573075 & 3.1415926536 & 0 & 1.1253423662 \\ 2.0370199171 & 2.6951112155 & 4.3279705032 & 0.8522641276 \\ 0.3292355244 & 0.4768014425 & -11.9983435682 & 1.1253423662 \\ 1.3815920754 & 2.0757474625 & 1.8067539983 & 1.1253423662 \\ -2.2329342208 & 5.3820397842 & -7.9050531083 & 14.4127228981 \\ 2.8883760432 & 4.7817556419 & -14.346597525 & 0.3941363942 \end{pmatrix}$$

$$P = \begin{pmatrix} 0 & 1 & 1 & 1 & 1 & -1 & 1 \\ 1 & 0 & 1 & 1 & 1 & -1 & 1 \\ 1 & 1 & 0 & -1 & -1 & -1 & 1 \\ 1 & 1 & -1 & 0 & -1 & -1 & -1 \\ 1 & 1 & -1 & -1 & 0 & -1 & -1 \\ -1 & -1 & -1 & -1 & -1 & 0 & 1 \\ 1 & 1 & 1 & -1 & -1 & 1 & 0 \end{pmatrix} \qquad R = \begin{pmatrix} 0 & 3 & 3 & 5 & 4 & 5 & 4 \\ 3 & 0 & 5 & 3 & 4 & 4 & 5 \\ 0 & 0 & 0 & 0 & 0 & 0 & 0 \\ 0 & 0 & 0 & 0 & 0 & 0 & 0 \\ 1 & 1 & 4 & 4 & 0 & 1 & 1 \\ 0 & 0 & 0 & 0 & 0 & 0 & 0 \\ 0 & 0 & 0 & 0 & 0 & 0 & 0 \end{pmatrix}$$

$|P| = -2$

$I3(P, R) = 7.8$

$I3(-P, R) = 16.1$

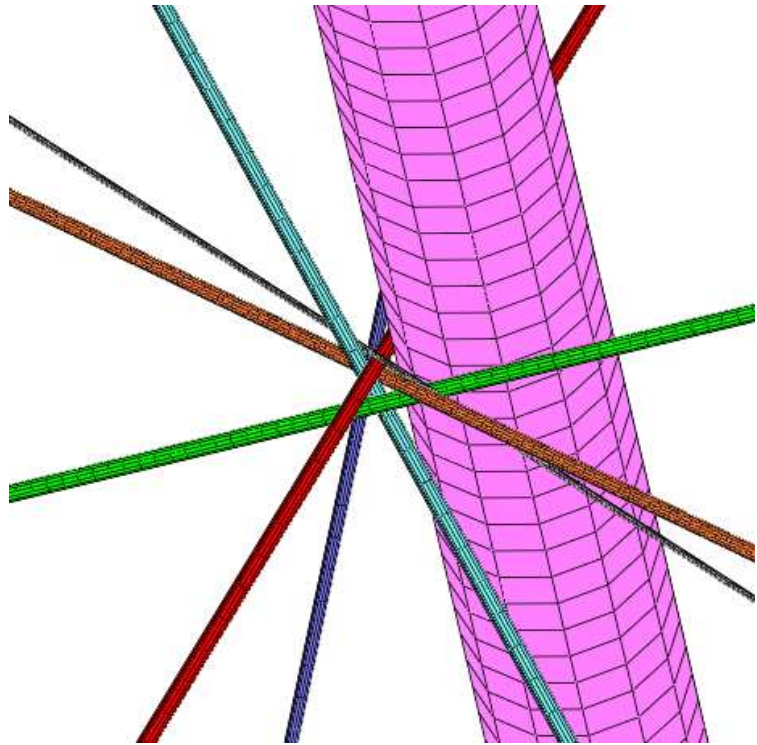

Appendix 4
7-cross
Configuration label  ma47

$$\begin{pmatrix} \text{orange} \\ \text{red} \\ \text{blue} \\ \text{green} \\ \text{cyan} \\ \text{magenta} \\ \text{gray} \end{pmatrix} \begin{pmatrix} 0 & 0 & 0 & 1 \\ t_1 & p_1 & z_1 & r_1 \\ t_2 & p_2 & z_2 & r_2 \\ t_3 & p_3 & z_3 & r_3 \\ t_4 & p_4 & z_4 & r_4 \\ t_5 & p_5 & z_5 & r_5 \\ t_6 & p_6 & z_6 & r_6 \end{pmatrix} = \begin{pmatrix} 0 & 0 & 0 & 1 \\ 3.0452079427 & 3.1415926536 & 0 & 1 \\ 0.1723459538 & -0.8583318575 & 0.8586515965 & 1 \\ 0.1667920833 & -0.5523132023 & -7.5962213713 & 1 \\ 2.4571207529 & -0.7156053079 & -96.3740086147 & 1 \\ 2.6455827091 & 9.0771697776 & -46.417569502 & 1 \\ 0.1382831705 & 5.4714409861 & -24.9579225639 & 0.41 \end{pmatrix}$$

$$P = \begin{pmatrix} 0 & 1 & 1 & 1 & 1 & 1 & 1 \\ 1 & 0 & -1 & -1 & 1 & -1 & 1 \\ 1 & -1 & 0 & 1 & 1 & -1 & 1 \\ 1 & -1 & 1 & 0 & -1 & 1 & -1 \\ 1 & 1 & 1 & -1 & 0 & 1 & 1 \\ 1 & -1 & -1 & 1 & 1 & 0 & 1 \\ 1 & 1 & 1 & -1 & 1 & 1 & 0 \end{pmatrix} \qquad R = \begin{pmatrix} 0 & 0 & 0 & 0 & 0 & 0 & 0 \\ 4 & 0 & 4 & 4 & 4 & 4 & 4 \\ 0 & 0 & 0 & 0 & 0 & 0 & 0 \\ 1 & 1 & 4 & 0 & 1 & 1 & 4 \\ 0 & 0 & 0 & 0 & 0 & 0 & 0 \\ 1 & 1 & 4 & 1 & 4 & 0 & 1 \\ 4 & 1 & 1 & 1 & 4 & 1 & 0 \end{pmatrix}$$

$|P| = -2$

$I3(P, R) = 16.5365853659$

$I3(-P, R) = 5.3414634146$

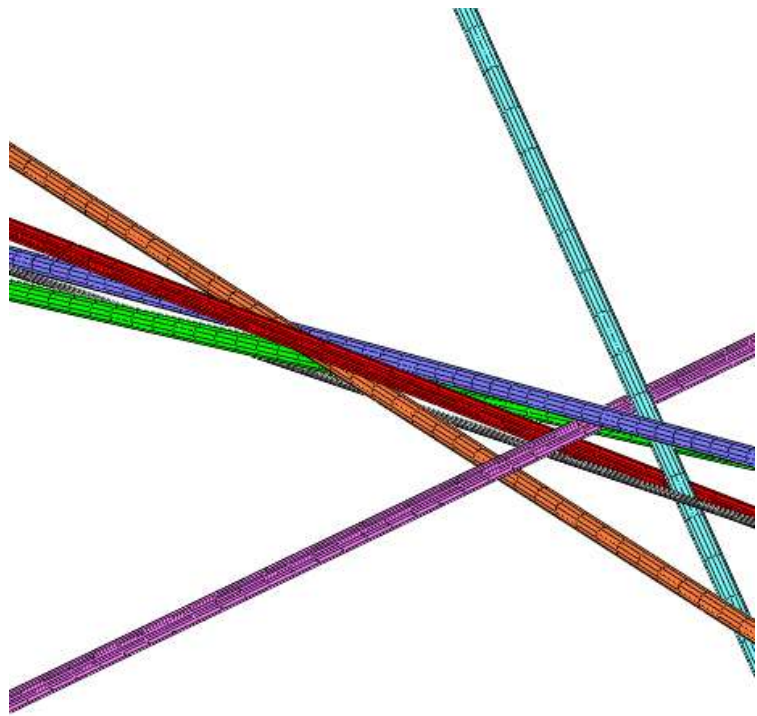

Appendix 4
7-cross
Configuration label ma45

$$\begin{pmatrix} \text{orange} \\ \text{red} \\ \text{blue} \\ \text{green} \\ \text{cyan} \\ \text{magenta} \\ \text{gray} \end{pmatrix} \begin{pmatrix} 0 & 0 & 0 & 1 \\ t_1 & p_1 & z_1 & r_1 \\ t_2 & p_2 & z_2 & r_2 \\ t_3 & p_3 & z_3 & r_3 \\ t_4 & p_4 & z_4 & r_4 \\ t_5 & p_5 & z_5 & r_5 \\ t_6 & p_6 & z_6 & r_6 \end{pmatrix} = \begin{pmatrix} 0 & 0 & 0 & 1 \\ 2.4711453573 & 3.1415926536 & 0 & 1 \\ 2.2139336809 & -2.9787958435 & 6.5758804565 & 1 \\ 2.6258628753 & -2.801182907 & -5.7553128549 & 1 \\ 0.6091778375 & 0.5647215202 & -4.687691704 & 1 \\ 0.404281182 & 6.903913655 & -8.8838241141 & 0.13 \\ 2.2129197866 & 3.3044798459 & -2.2323959752 & 0.25 \end{pmatrix}$$

$$P = \begin{pmatrix} 0 & 1 & 1 & 1 & 1 & 1 & 1 \\ 1 & 0 & -1 & 1 & -1 & -1 & 1 \\ 1 & -1 & 0 & 1 & -1 & -1 & 1 \\ 1 & 1 & 1 & 0 & -1 & 1 & 1 \\ 1 & -1 & -1 & -1 & 0 & 1 & 1 \\ 1 & -1 & -1 & 1 & 1 & 0 & 1 \\ 1 & 1 & 1 & 1 & 1 & 1 & 0 \end{pmatrix} \qquad R = \begin{pmatrix} 0 & 0 & 0 & 0 & 0 & 0 & 0 \\ 3 & 0 & 6 & 2 & 2 & 2 & 3 \\ 0 & 0 & 0 & 0 & 0 & 0 & 0 \\ 0 & 0 & 0 & 0 & 0 & 0 & 0 \\ 0 & 0 & 0 & 0 & 0 & 0 & 0 \\ 3 & 3 & 4 & 5 & 5 & 0 & 4 \\ 4 & 1 & 1 & 4 & 1 & 1 & 0 \end{pmatrix}$$

$|P| = -2$

$I3(P, R) = 18.8965517241$

$I3(-P, R) = 2.9655172414$

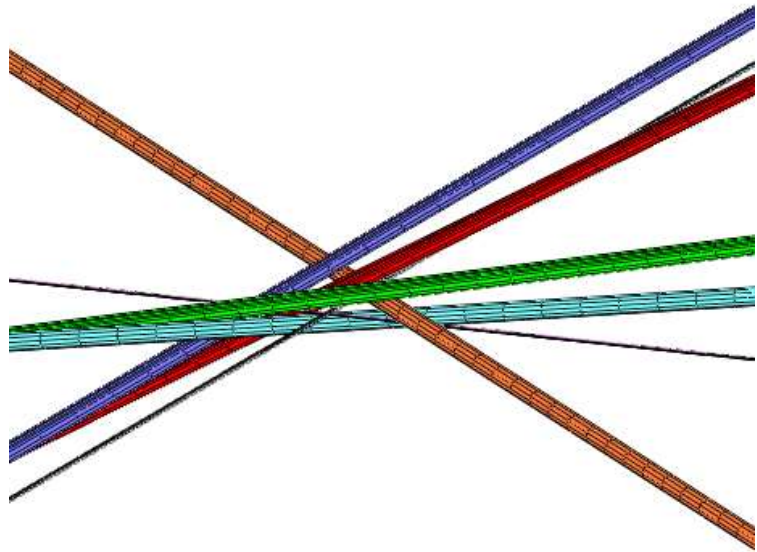

Appendix 4
7-cross
Configuration label  a15

$$\begin{pmatrix} \text{orange} \\ \text{red} \\ \text{blue} \\ \text{green} \\ \text{cyan} \\ \text{magenta} \\ \text{gray} \end{pmatrix} \begin{pmatrix} 0 & 0 & 0 & 1 \\ t_1 & p_1 & z_1 & r_1 \\ t_2 & p_2 & z_2 & r_2 \\ t_3 & p_3 & z_3 & r_3 \\ t_4 & p_4 & z_4 & r_4 \\ t_5 & p_5 & z_5 & r_5 \\ t_6 & p_6 & z_6 & r_6 \end{pmatrix} = \begin{pmatrix} 0 & 0 & 0 & 1 \\ 0.3670539738 & 3.1415926536 & 0 & 1 \\ 0.9702020678 & 3.7534162051 & 9.5819794215 & 1 \\ 2.6862801939 & 0.3531540227 & -0.2307994636 & 1 \\ 2.6065227644 & 0.2251887379 & -6.7421034487 & 1 \\ 2.4654065903 & 6.8101965651 & -3.4888901462 & 0.25 \\ 2.7496307564 & 6.0829608965 & -0.6114696657 & 0.1233979826 \end{pmatrix}$$

$$P = \begin{pmatrix} 0 & 1 & 1 & 1 & 1 & 1 & 1 \\ 1 & 0 & -1 & 1 & 1 & 1 & 1 \\ 1 & -1 & 0 & -1 & 1 & -1 & 1 \\ 1 & 1 & -1 & 0 & -1 & 1 & -1 \\ 1 & 1 & 1 & -1 & 0 & -1 & 1 \\ 1 & 1 & -1 & 1 & -1 & 0 & 1 \\ 1 & 1 & 1 & -1 & 1 & 1 & 0 \end{pmatrix} \quad R = \begin{pmatrix} 0 & 6 & 3 & 2 & 2 & 2 & 3 \\ 0 & 0 & 0 & 0 & 0 & 0 & 0 \\ 0 & 0 & 0 & 0 & 0 & 0 & 0 \\ 0 & 0 & 0 & 0 & 0 & 0 & 0 \\ 0 & 0 & 0 & 0 & 0 & 0 & 0 \\ 3 & 4 & 3 & 5 & 5 & 0 & 4 \\ 2 & 2 & 6 & 2 & 3 & 3 & 0 \end{pmatrix}$$

$|P| = -2$

$I3(P, R) = 7.4242424242$

$I3(-P, R) = 23.9393939394$

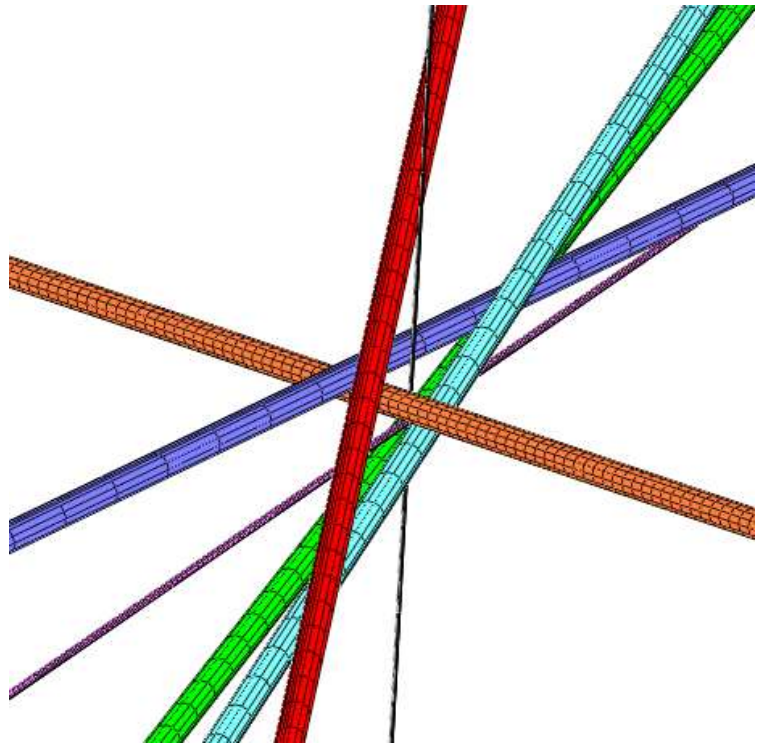

Appendix 4
7-cross
Configuration label  a12

$$\begin{pmatrix} \text{orange} \\ \text{red} \\ \text{blue} \\ \text{green} \\ \text{cyan} \\ \text{magenta} \\ \text{gray} \end{pmatrix} \begin{pmatrix} 0 & 0 & 0 & 1 \\ t_1 & p_1 & z_1 & r_1 \\ t_2 & p_2 & z_2 & r_2 \\ t_3 & p_3 & z_3 & r_3 \\ t_4 & p_4 & z_4 & r_4 \\ t_5 & p_5 & z_5 & r_5 \\ t_6 & p_6 & z_6 & r_6 \end{pmatrix} = \begin{pmatrix} 0 & 0 & 0 & 1 \\ 0.9229270329 & 3.1415926536 & 0 & 1 \\ 2.7035309819 & 3.6787289084 & 10.4129688582 & 1 \\ 2.0500812188 & -3.2822441705 & 18.2448210165 & 1 \\ 1.9623658867 & 3.2573260042 & 20.3670756694 & 1 \\ 1.3305510155 & 9.348020398 & 13.3869154539 & 0.2 \\ 2.9211508015 & 5.8644979027 & -17.2995624019 & 0.049076264 \end{pmatrix}$$

$$P = \begin{pmatrix} 0 & 1 & 1 & 1 & 1 & 1 & 1 \\ 1 & 0 & -1 & 1 & -1 & 1 & -1 \\ 1 & -1 & 0 & 1 & 1 & 1 & 1 \\ 1 & 1 & 1 & 0 & -1 & 1 & -1 \\ 1 & -1 & 1 & -1 & 0 & -1 & 1 \\ 1 & 1 & 1 & 1 & -1 & 0 & 1 \\ 1 & -1 & 1 & -1 & 1 & 1 & 0 \end{pmatrix} \qquad R = \begin{pmatrix} 0 & 0 & 0 & 0 & 0 & 0 & 0 \\ 0 & 0 & 0 & 0 & 0 & 0 & 0 \\ 6 & 3 & 0 & 2 & 2 & 2 & 3 \\ 1 & 1 & 4 & 0 & 4 & 1 & 1 \\ 0 & 0 & 0 & 0 & 0 & 0 & 0 \\ 5 & 4 & 3 & 4 & 5 & 0 & 3 \\ 2 & 6 & 2 & 2 & 3 & 3 & 0 \end{pmatrix}$$

$|P| = -2$

$I3(P, R) = 3.1278688525$

$I3(-P, R) = 24.6032786885$

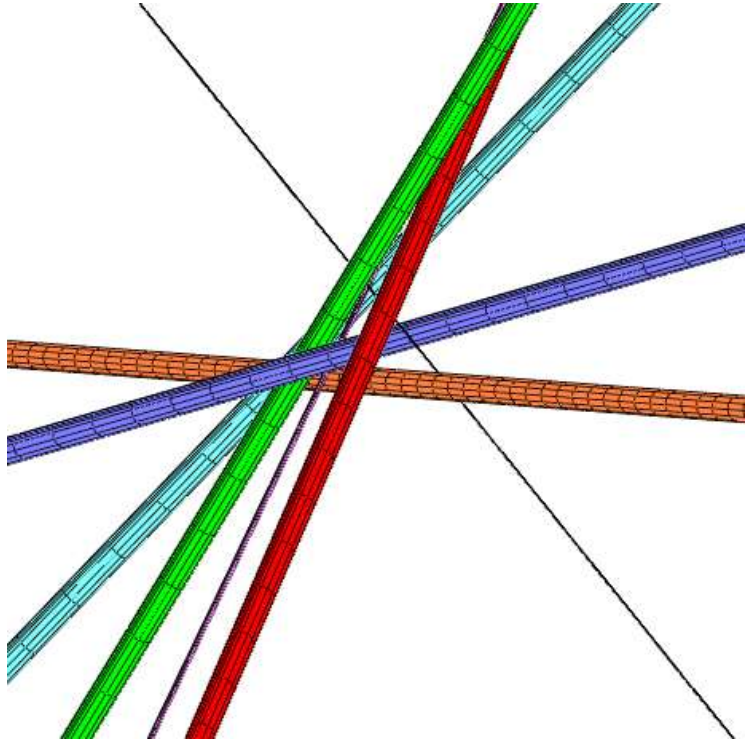

Appendix 5

$$\text{md}(x,h) := \mod\left[\frac{\left(x - \mod\left(x, 2^h\right)\right)}{2^h}, 2\right]$$

$$\text{tonum}(x,n) := \begin{array}{|l} \text{for } i \in 0..n-1 \\ \quad a_{0,(n-1-i)} \leftarrow \text{md}(x,i) \\ a \end{array}$$

$$\text{sum}(x,n) := \begin{array}{|l} a \leftarrow 0 \\ \text{for } i \in 0..n-1 \\ \quad a \leftarrow a + \text{md}(x,i) \\ a \end{array}$$

$$\text{mat}(x,n) := \begin{array}{|l} p \leftarrow 0 \\ \text{for } i \in 0..n-1 \\ \quad \text{for } j \in 0..n-1 \\ \quad\quad a_{i,j} \leftarrow \text{if}\left[(i = j) + (i > j), 0, \text{sign}\left[0.5 - \text{tonum}\left[x, \frac{n \cdot (n-1)}{2}\right]_{0,p}\right]\right] \\ \quad\quad p \leftarrow \text{if}(j > i, p+1, p) \\ a \leftarrow a + a^T \\ a \end{array}$$

$$\text{mat1}(A, x, n1, n2) := \begin{array}{|l} p \leftarrow 0 \\ n \leftarrow n2 - n1 \\ \text{for } i \in 0..n2-1 \\ \quad \text{for } j \in 0..n2-1 \\ \quad\quad nn \leftarrow \frac{n2 \cdot (n2-1)}{2} - \frac{n1 \cdot (n1-1)}{2} \\ \quad\quad a_{i,j} \leftarrow \text{if}\left[(i < j), \text{if}\left(i < n, \text{sign}(0.5 - \text{tonum}(x,nn)_{0,p}), A_{i-n, j-n}\right), 0\right] \\ \quad\quad p \leftarrow \text{if}[(j > i) \cdot (i < n), p+1, p] \\ a \leftarrow a + a^T \\ a \end{array}$$

$$\text{outmat}(A, V, k) := \begin{array}{|l} \text{for } i \in 0..k-1 \\ \quad \text{for } j \in 0..k-1 \\ \quad\quad a_{i,j} \leftarrow A_{V_i, V_j} \\ a \end{array}$$

$$\text{tonum2}(n,k) := \begin{array}{|l} p \leftarrow 0 \\ \text{for } i \in 1..2^n \\ \quad \text{if sum}(i,n) = k \\ \quad \quad \begin{array}{|l} p \leftarrow p + 1 \\ f \leftarrow 0 \\ w \leftarrow \text{tonum}(i,n) \\ \text{for } j \in 0..n-1 \\ \quad \text{if } w_{0,j} \neq 0 \\ \quad \quad \begin{array}{|l} v_f \leftarrow j \\ f \leftarrow f + 1 \end{array} \\ a_p \leftarrow v \end{array} \\ a_0 \leftarrow p \\ a \end{array} \qquad \text{checkforKx}(A, TN, M, x) := \begin{array}{|l} a \leftarrow 0 \\ \text{for } i \in 1..TN_0 \\ \quad \begin{array}{|l} O \leftarrow \text{outmat}(A, TN_i, x) \\ \text{for } j \in 0..2^x - 1 \\ \quad \text{if } O = M_j \\ \quad \quad \begin{array}{|l} a \leftarrow 1 \\ \text{break} \end{array} \\ \text{break if } a \end{array} \\ a \end{array}$$

$$\text{matgenK5}(n) := \begin{array}{|l} \text{for } i \in 0..n-1 \\ \quad \text{for } j \in 0..n-1 \\ \quad \quad \begin{array}{|l} P1_{i,j} \leftarrow \text{if}(i = j, 0, 1) \\ P2_{i,j} \leftarrow \text{if}(i = j, 0, -1) \\ U_{i,j} \leftarrow \text{if}(i = j, 1, 0) \end{array} \\ \text{for } i \in 0..2^4 - 1 \\ \quad \begin{array}{|l} \text{for } k \in 1..n-1 \\ \quad U_{k,k} \leftarrow \text{sign}(0.5 - \text{md}(i, k-1)) \\ N_i \leftarrow U \cdot P1 \cdot U \end{array} \\ \text{for } i \in 0..2^4 - 1 \\ \quad \begin{array}{|l} \text{for } k \in 1..n-1 \\ \quad U_{k,k} \leftarrow \text{sign}(0.5 - \text{md}(i, k-1)) \\ N_{i+2^4} \leftarrow U \cdot P2 \cdot U \end{array} \\ N \end{array}$$

$MK5 := \text{matgenK5}(5)$

$TN := \text{tonum2}(9, 5)$

$$\text{matgen}(n) := \begin{vmatrix} N_{0,0} \leftarrow 0 \\ \text{for } i \in 0..n-1 \\ \quad K_i \leftarrow 0 \\ N_{1,1} \leftarrow K \\ \text{for } i \in 0..2^{\frac{n \cdot (n-1)}{2}} - 1 \\ \quad \begin{vmatrix} P \leftarrow \text{mat}(i,n) \\ E \leftarrow \text{eigenvals}(P) \\ p \leftarrow 1 \\ \text{for } j \in 1..N_{0,0} \\ \quad \text{if } (N_{1,j} - E) \cdot (N_{1,j} - E) < 0.0001 \\ \quad \quad \begin{vmatrix} p \leftarrow 0 \\ \text{break} \end{vmatrix} \\ \text{if } p \\ \quad \begin{vmatrix} N_{0,0} \leftarrow N_{0,0} + 1 \\ N_{0,N_{0,0}} \leftarrow P \\ N_{1,N_{0,0}} \leftarrow E \end{vmatrix} \end{vmatrix} \\ N \end{vmatrix}$$

m6 := matgen(6)

| m6 = | 0 | 1 | 2 | 3 | 4 | 5 | 6 |
|---|---|---|---|---|---|---|---|
| 0 | 16 | [6, 6] | [6, 6] | [6, 6] | [6, 6] | [6, 6] | [6, 6] |
| 1 | 0 | [6, 1] | [6, 1] | [6, 1] | [6, 1] | [6, 1] | ... |

$$\text{MatrixGenerator}(A, n) := \begin{array}{|l} N_{0,0} \leftarrow 0 \\ \text{for } i \in 0..n \\ \quad K_i \leftarrow 0 \\ N_{1,1} \leftarrow K \\ \text{for } t \in 1..A_{0,0} \\ \quad \begin{array}{|l} B \leftarrow A_{0,t} \\ \text{for } i \in 0..2^n - 1 \\ \quad \begin{array}{|l} P \leftarrow \text{mat1}(B, i, n, n+1) \\ E \leftarrow \text{eigenvals}(P) \\ p \leftarrow 1 \\ \text{for } j \in 1..N_{0,0} \\ \quad \text{if } (N_{1,j} - E) \cdot (N_{1,j} - E) < 0.0001 \\ \quad \quad \begin{array}{|l} p \leftarrow 0 \\ \text{break} \end{array} \\ \text{if } p \\ \quad \begin{array}{|l} N_{0,0} \leftarrow N_{0,0} + 1 \\ N_{0, N_{0,0}} \leftarrow P \\ N_{1, N_{0,0}} \leftarrow E \end{array} \end{array} \end{array} \\ N \end{array}$$

m7 := MatrixGenerator(m6, 6)

|   |   | 0  | 1      | 2      | 3      | 4      | 5      | 6      |
|---|---|----|--------|--------|--------|--------|--------|--------|
| m7 = | 0 | 54 | [7, 7] | [7, 7] | [7, 7] | [7, 7] | [7, 7] | [7, 7] |
|   | 1 | 0  | [7, 1] | [7, 1] | [7, 1] | [7, 1] | [7, 1] | ...    |

$\text{MatrixGeneratorNoK5}(A, n, s) :=$
$\quad$ MK5 ← matgenK5(5)
$\quad$ TN ← tonum2(n + 1, 5)
$\quad$ if s
$\quad\quad$ $T_0 \leftarrow 0$
$\quad\quad$ for $i \in 1..TN_0$
$\quad\quad\quad$ if $(TN_i)_0 = 0$
$\quad\quad\quad\quad$ $T_0 \leftarrow T_0 + 1$
$\quad\quad\quad\quad$ $T_{(T_0)} \leftarrow TN_i$
$\quad\quad$ TN ← T
$\quad$ $N_{0,0} \leftarrow 0$
$\quad$ for $i \in 0..n$
$\quad\quad$ $K_i \leftarrow 0$
$\quad$ $N_{1,1} \leftarrow K$
$\quad$ for $t \in 1..A_{0,0}$
$\quad\quad$ $B \leftarrow A_{0,t}$
$\quad\quad$ for $i \in 0..2^n - 1$
$\quad\quad\quad$ $P \leftarrow \text{mat1}(B, i, n, n+1)$
$\quad\quad\quad$ $E \leftarrow \text{eigenvals}(P)$
$\quad\quad\quad$ $p \leftarrow 1$
$\quad\quad\quad$ for $j \in 1..N_{0,0}$
$\quad\quad\quad\quad$ if $(N_{1,j} - E) \cdot (N_{1,j} - E) < 0.0001$
$\quad\quad\quad\quad\quad$ $p \leftarrow 0$
$\quad\quad\quad\quad\quad$ break
$\quad\quad\quad$ $p \leftarrow 0$ if checkforKx(P, TN, MK5, 5) if p
$\quad\quad\quad$ if p
$\quad\quad\quad\quad$ $N_{0,0} \leftarrow N_{0,0} + 1$
$\quad\quad\quad\quad$ $N_{0, N_{0,0}} \leftarrow P$
$\quad\quad\quad\quad$ $N_{1, N_{0,0}} \leftarrow E$
$\quad$ N

m7noK5 := MatrixGeneratorNoK5(m6, 6, 0)

| | 0 | 1 | 2 | 3 | 4 | 5 | 6 |
|---|---|---|---|---|---|---|---|
| m7noK5 = 0 | 22 | [7, 7] | [7, 7] | [7, 7] | [7, 7] | [7, 7] | [7, 7] |
| 1 | 0 | [7, 1] | [7, 1] | [7, 1] | [7, 1] | [7, 1] | ... |

m8noK5 := MatrixGeneratorNoK5(m7noK5, 7, 1)

| | 0 | 1 | 2 | 3 | 4 | 5 | 6 |
|---|---|---|---|---|---|---|---|
| m8noK5 = 0 | 51 | [8, 8] | [8, 8] | [8, 8] | [8, 8] | [8, 8] | [8, 8] |
| 1 | 0 | [8, 1] | [8, 1] | [8, 1] | [8, 1] | [8, 1] | ... |

m9noK5 := MatrixGeneratorNoK5(m8noK5, 8, 1)

| | 0 | 1 | 2 | 3 | 4 | 5 | 6 |
|---|---|---|---|---|---|---|---|
| m9noK5 = 0 | 105 | [9, 9] | [9, 9] | [9, 9] | [9, 9] | [9, 9] | [9, 9] |
| 1 | 0 | [9, 1] | [9, 1] | [9, 1] | [9, 1] | [9, 1] | ... |

m10noK5 := MatrixGeneratorNoK5(m9noK5, 9, 1)

| | 0 | 1 | 2 | 3 | 4 |
|---|---|---|---|---|---|
| m10noK5 = 0 | 172 | [10, 10] | [10, 10] | [10, 10] | [10, 10] |
| 1 | 0 | [10, 1] | [10, 1] | [10, 1] | ... |

m11noK5 := MatrixGeneratorNoK5(m10noK5, 10, 1)

| | 0 | 1 | 2 | 3 | 4 |
|---|---|---|---|---|---|
| m11noK5 = 0 | 142 | [11, 11] | [11, 11] | [11, 11] | [11, 11] |
| 1 | 0 | [11, 1] | [11, 1] | [11, 1] | ... |

m12noK5 := MatrixGeneratorNoK5(m11noK5, 11, 1)

| | 0 | 1 | 2 | 3 | 4 |
|---|---|---|---|---|---|
| m12noK5 = 0 | 61 | [12, 12] | [12, 12] | [12, 12] | [12, 12] |
| 1 | 0 | [12, 1] | [12, 1] | [12, 1] | ... |

m13noK5 := MatrixGeneratorNoK5(m12noK5, 12, 1)

$$m13noK5 = \begin{pmatrix} 8 & \{13,13\} & \{13,13\} & \{13,13\} & \{13,13\} & \{13,13\} & \{13,13\} & \{13,13\} & \{13,13\} \\ 0 & \{13,1\} & \{13,1\} & \{13,1\} & \{13,1\} & \{13,1\} & \{13,1\} & \{13,1\} & \{13,1\} \end{pmatrix}$$

m14noK5 := MatrixGeneratorNoK5(m13noK5, 13, 1)

$$m14noK5 = \begin{pmatrix} 5 & \{14,14\} & \{14,14\} & \{14,14\} & \{14,14\} & \{14,14\} \\ 0 & \{14,1\} & \{14,1\} & \{14,1\} & \{14,1\} & \{14,1\} \end{pmatrix}$$

m15noK5 := MatrixGeneratorNoK5(m14noK5, 14, 1)

$$m15noK5 = \begin{pmatrix} 2 & \{15,15\} & \{15,15\} \\ 0 & \{15,1\} & \{15,1\} \end{pmatrix}$$

m16noK5 := MatrixGeneratorNoK5(m15noK5, 15, 1)

$$\text{m16noK5} = \begin{pmatrix} 1 & \{16,16\} \\ 0 & \{16,1\} \end{pmatrix}$$

m17noK5 := MatrixGeneratorNoK5(m16noK5, 16, 1)

$$\text{m17noK5} = \begin{pmatrix} 1 & \{17,17\} \\ 0 & \{17,1\} \end{pmatrix}$$

m18noK5 := MatrixGeneratorNoK5(m17noK5, 17, 1)

$$\text{m18noK5} = \begin{pmatrix} 1 & \{18,18\} \\ 0 & \{18,1\} \end{pmatrix}$$

m19noK5 := MatrixGeneratorNoK5(m18noK5, 18, 1)

$$\text{m19noK5} = \begin{pmatrix} 0 & 0 \\ 0 & \{19,1\} \end{pmatrix}$$

$\text{m18noK5}_{0,1} =$

|    | 3  | 4  | 5  | 6  | 7  | 8  | 9  | 10 | 11 | 12 |
|----|----|----|----|----|----|----|----|----|----|----|
| 1  | 1  | -1 | -1 | 1  | -1 | 1  | 1  | -1 | 1  | 1  |
| 2  | 1  | 1  | -1 | -1 | 1  | -1 | -1 | 1  | -1 | 1  |
| 3  | 0  | 1  | 1  | 1  | -1 | 1  | -1 | -1 | -1 | 1  |
| 4  | 1  | 0  | 1  | 1  | 1  | -1 | 1  | 1  | -1 | 1  |
| 5  | 1  | 1  | 0  | 1  | 1  | 1  | -1 | -1 | 1  | -1 |
| 6  | 1  | 1  | 1  | 0  | 1  | 1  | 1  | -1 | -1 | 1  |
| 7  | -1 | 1  | 1  | 1  | 0  | 1  | 1  | 1  | -1 | -1 |
| 8  | 1  | -1 | 1  | 1  | 1  | 0  | 1  | 1  | 1  | -1 |
| 9  | -1 | 1  | -1 | 1  | 1  | 1  | 0  | 1  | 1  | 1  |
| 10 | -1 | 1  | -1 | -1 | 1  | 1  | 1  | 0  | 1  | 1  |
| 11 | -1 | -1 | 1  | -1 | -1 | 1  | 1  | 1  | 0  | 1  |
| 12 | 1  | 1  | -1 | 1  | -1 | -1 | 1  | 1  | 1  | 0  |
| 13 | 1  | -1 | 1  | -1 | 1  | 1  | -1 | 1  | 1  | 1  |
| 14 | -1 | -1 | -1 | 1  | 1  | -1 | 1  | -1 | -1 | ...|